\documentclass[12pt]{article}
\usepackage[a4paper]{geometry}
\usepackage{amssymb, latexsym, amsmath, exscale}

\newsavebox{\ttt}
\sbox{\ttt}{}
\newcommand{\startsection}[1]
    {\section[#1]{#1}
    \sbox{\ttt}{\thesection\ \ \textsc{#1}}
    \thispagestyle{plain}
}

\newenvironment{evlist}[2]{
\begin{list}{}{
\setlength{\topsep}{0.5ex plus0.2ex minus0.1ex} 
\setlength{\leftmargin}{#1}
\setlength{\itemsep}{#2 plus0.2ex}
\setlength{\parsep}{0ex plus0.2ex} }}
{\end{list}}
\newcommand{\Self}[1]{\mathrm{T}_{#1}}
\newcommand{\vonn}{\sigma}
\newcommand{\Selfb}[1]{\mathrm{S}_{#1}}
\newcommand{\ordass}{\varrho}

\newcommand{\Nat}{\mathbb{N}}
\newcommand{\Int}{\mathbb{Z}}
\newcommand{\Bool}{\mathbb{B}}

\newcommand{\id}{\mathrm{id}}
\newcommand{\twogroup}{\mathrm{F}_2}
\newcommand{\definition}[1]{\textit{#1}}
\newcommand{\proof}{{\textit{Proof}\enspace}}
\newcommand{\ssuc}[1]{\textsf{s}_{\,#1}}
\newcommand{\esuc}[1]{\textsf{e}_{\,#1}}
\newcommand{\standard}{$\Nat$-like}

\newcommand{\eop}{\ \vbox{\hrule
                       \hbox{\vrule
                             \hskip 6pt
                             \vrule height 6pt width 0pt
                             \vrule}%
                       \hrule}%
                     \vspace{\medskipamount}
                }

\newtheorem{lemma}{Lemma}[section]
\newtheorem{proposition}{Proposition}[section]
\newtheorem{theorem}{Theorem}[section]
\newcommand{\proper}[1]{#1^p}
\newcommand{\prop}{\mathsf{P}}

\newcommand{\lind}[1]{\mathcal{I}_0(#1)}
\newcommand{\finite}[1]{\mathsf{Fin}(#1)}

\newcommand{\fin}{\mathsf{Fin}}

\newcommand{\class}[1]{\mathsf{#1}}

\newcommand{\onept}{\varepsilon}

\newcommand{\val}{\omega}

\newlength{\extextwidth}
\newenvironment{exframe}{\begin{minipage}{\extextwidth}{}}
{\end{minipage}}

\begin{document}

\begin{titlepage}

\begin{center}
\phantom{q}
\vspace{60pt}

{\huge Some Notes on Finite Sets}

\bigskip
{\LARGE Chris Preston}

\bigskip
    {\large April 2018}
\end{center}

\bigskip
\bigskip
\bigskip
\bigskip
\begin{quote}

\end{quote}

\end{titlepage}

\thispagestyle{empty}

\addtocontents{toc}{\vskip 20pt}
\addtolength{\parskip}{-5pt}
\tableofcontents
\addtolength{\parskip}{5pt}

\startsection{Introduction}

\label{intro}

For sets $E$ and $F$ we write $E \approx F$ if there exists a bijective mapping $f : E \to F$.
The standard definition of a set $A$ being finite is that $A \approx [n]$ for some $n \in \Nat$, where 
$\Nat = \{0,1,\ldots\,\}$ is the set of natural numbers, $[0] = \varnothing$, and $[n] = \{0,1,\ldots,n-1\}$ for each $n \in \Nat \setminus \{0\}$.

This seems straightforward enough until it is perhaps recalled how much effort is required to give a rigorous definition of the 
sets $[m],\,m \in \Nat$. Moreover, it is often not that easy to apply. For example, try giving, without too much thought, a  proof 
of the fact that any injective mapping of a finite set into itself is bijective.

A more fundamental objection to  the standard definition is its use of the infinite set $\Nat$, which
might even give the impression that finite sets can only be defined with the help of such a set. This is certainly not the case and
in  what follows we are going to work with one of several possible definitions not involving the natural numbers. 
The definition introduced below is usually called \definition{Kuratowski-finiteness} \cite{kura} and it is  
essentially that employed by Whitehead  and Russell in \textit{Principia Mathematica} \cite{wr}. 
We will also make use of a characterisation of finite sets due to Tarski \cite{tarski}.
There are similar approaches which also appeared in the first decades of the previous century and excellent 
treatments of this topic can be found, for example, in Levy \cite{levy} and Suppes \cite{suppes}.

These notes aim to give a gentle account of one approach to the theory of finite sets without making use of the natural numbers or any other infinite set.
They were written to be used as the basis for a student seminar.
There are no real prerequisites except for a certain familiarity with the kind of mathematics seen in the first couple of years of a university
mathematics course and the language used to describe sets in such a course will be employed here.
It is assumed that the reader is acquainted with such basic concepts as equivalence relations and partial and totally ordered sets as well as basic structures
such as monoids and groups.
It is  not assumed that the reader has taken a course on axiomatic set theory. However, since the subject of these notes is finite sets, we need to briefly discuss 
one of the axioms
of set theory, namely the axiom of infinity. This is needed to ensure that an infinite set 
exists. There is a version of set theory in which the axiom of infinity is replaced by its negation and in this version of set theory
all sets are finite and so the natural numbers do not exist (at least as a set). 
In general we do not assume that an infinite set exists. However, the symbol $\Nat$ always refers to the first element of a triple
$(\Nat,\mathsf{s},0)$ satisfying the Peano axioms, where $\Nat$ is the set of natural numbers, $\mathsf{s} : \Nat\to \Nat$ is the successor 
mapping with $\mathsf{s}(n) = n+1$ for all $n \in \Nat$ and $0$ is the initial element of $\Nat$. 
This means that the mapping $\mathsf{s}$ is injective, $\mathsf{s}(n) \ne 0$ for all $n \in \Nat$ and that $N = \Nat$ whenever $N$ is a subset of $\Nat$
containing $0$ and such that $\mathsf{s}(n) \in N$ for all $n \in N$.
Whenever $(\Nat,\mathsf{s},0)$  occurs it is assumed that it exists and hence that the negation of the axiom of infinity is not in force.

We assume that the reader has heard of the axiom of choice. This is not involved
when dealing only with finite sets.
 Its use will be pointed out on the couple of occasions when a statement (involving infinite sets) depends on this axiom.

There is one important point which should be mentioned. We will be dealing with mappings defined on the collection of all finite sets and this collection
is too large to be considered a set, meaning that treating it as a set might possibly lead to various paradoxes. Such a large collection is called a 
\definition{proper class} and something which is either a proper class or a set is referred to as a \definition{class}. The proper class of all finite sets will be denoted by $\fin$. However, as far as what is to be found in these notes, 
there is no problem in treating proper classes as if they were just sets. 

It could be objected that, after having rejected the infinite set $\Nat$ as a means of introducing finite sets, we now resort to objects which are so large that
they cannot even be considered to be sets.
But mappings defined on proper classes can make sense without involving an infinite set. An important example is the mapping
$\vonn : \fin \to \fin$ defined by $\vonn(A) = A \cup \{A\}$ for each finite set $A$. This mapping $\vonn$ will form the basis for defining the finite
ordinals.

The power set of a set $E$, i.e., the set of all its subsets, will be denoted by $\mathcal{P}(E)$ and the set of non-empty subsets by
$\mathcal{P}_0(E)$.
If $E$ is a set and $\mathcal{S}$ is a subset of $\mathcal{P}(E)$ then $\proper{\mathcal{S}}$ will be used to denote the set
of subsets in $\mathcal{S}$ which are proper subsets of $E$.

Finite sets will be defined here in terms of what is known as an inductive system, where
a subset $\mathcal{S}$ of $\mathcal{P}(E)$ is called an \definition{inductive $E$-system} if 
$\varnothing \in \mathcal{S}$ and $F \cup \{e\} \in \mathcal{S}$ for all $F \in \proper{\mathcal{S}}$, $e \in E \setminus F$. 
In particular, $\mathcal{P}(E)$ is itself an inductive $E$-system.

The definition of being finite which will be used here is the following; \cite{kura}, \cite{wr}:

\medskip

\fbox{\begin{exframe}
\begin{center}A set $E$ is defined to be \definition{finite} if $\mathcal{P}(E)$ is the only inductive $E$-system.\end{center} 
\end{exframe}}
\medskip

In Theorem~\ref{theorem_intro_2} we show that the above definition of being finite is equivalent to the standard
definition given in terms of the natural numbers.

\begin{lemma}\label{lemma_intro_1}
The empty set $\varnothing$ is finite. Moreover, for each finite set  $A$ and each element $a \notin A$ the set $A \cup \{a\}$
is finite.
\end{lemma}

\proof 
The empty set $\varnothing$ is finite since $\mathcal{P}(\varnothing) = \{\varnothing\}$ is the only
subset of $\mathcal{P}(\varnothing)$ containing $\varnothing$. Now consider a finite set $A$ and $a \notin A$. Put $B = A \cup \{a\}$ and 
let $\mathcal{R}$ be an inductive $B$-system. 
Then $\mathcal{S} = \mathcal{R} \cap \mathcal{P}(A)$ is clearly an inductive $A$-system and 
thus $\mathcal{S} = \mathcal{P}(A)$ (since $A$ is finite), i.e., $\mathcal{P}(A) \subset \mathcal{R}$. Moreover, 
$A' \cup \{a\} \in \mathcal{R}$ for all $A' \in \mathcal{P}(A)$, since $\mathcal{R}$ is an inductive $B$-system and 
$\mathcal{P}(A) \subset \mathcal{R}$. This implies that $\mathcal{R} = \mathcal{P}(B)$ and hence that
$B = A \cup \{a\}$ is finite. \eop

Note that an arbitrary intersection of inductive $E$-systems is again an inductive $E$-system  
and so there is a least inductive $E$-system (namely the intersection of all such $E$-systems).
Thus if the least inductive $E$-system is denoted by $\lind{E}$
then a set $E$ is finite if and only if $\lind{E} = \mathcal{P}(E)$.

The set of finite subsets of a set $E$ will be denoted by $\finite{E}$.

\begin{lemma}\label{lemma_intro_2}
For each set $E$ the least inductive $E$-system is exactly the set of finite subsets of $E$, i.e., $\lind{E} = \finite{E}$.
In particular, a set $E$ is finite if and only if every inductive $E$-system contains $E$.
\end{lemma}

\proof By Lemma~\ref{lemma_intro_1} $\finite{E}$ is an inductive $E$-system and hence $\lind{E} \subset \finite{E}$. 
Conversely, if $A \in \finite{E}$ then $A \in \mathcal{P}(A) = \lind{A} \subset \lind{E}$ and therefore also
$\finite{E} \subset \lind{E}$. 
In particular, it follows that if every inductive $E$-system contains $E$ then $E \in  \lind{E} = \finite{E}$, and hence $E$ is finite. Clearly 
if $E$ is finite then 
$E \in \mathcal{P}(E) = \lind{E}$ and so every inductive $E$-system contains $E$. 
\eop

\begin{proposition}\label{prop_intro_1}
Every subset of a finite set $A$ is finite.
\end{proposition}

\proof 
By Lemma~\ref{lemma_intro_2}  $\finite{A} = \mathcal{P}(A)$, and hence every subset of $A$ is finite. 
\eop

We next show that the definition of being finite employed here is equivalent to the standard
definition. 
To help distinguish between these two definitions 
let us call sets which are finite according to the standard definition \definition{$\Nat$-finite}. Thus a set
$A$ is $\Nat$-finite if and only if there exists a bijective mapping $h : [n] \to A$ for some $n \in \Nat$.
(When working with this definition we assume the reader is familiar with the properties of the sets $[n]$, $n \in \Nat$.)

\begin{theorem}\label{theorem_intro_2}
A set is finite if and only if it is $\Nat$-finite.
\end{theorem}

\proof
Let $A$ be a finite set and let $\mathcal{S} = \{ B \in \mathcal{P}(A) : \mbox{$B$ is $\Nat$-finite} \}$. Clearly $\varnothing \in \mathcal{S}$, so consider 
$B \in \proper{\mathcal{S}}$, let $a \in A \setminus B$ and put $B' = B \cup \{a\}$. By assumption there exists $n \in \Nat$  and a bijective mapping $h : [n] \to B$
and the mapping $h$ can be extended to a bijective mapping $h' : [n+1] \to B'$ by putting $h'(n) = a$; hence $B' \in \mathcal{S}$. It follows that $\mathcal{S}$ 
is an inductive $A$-system and hence $\mathcal{S} = \mathcal{P}(A)$, since $A$ is finite. In particular $A \in \mathcal{S}$, i.e., $A$ is $\Nat$-finite. 
This shows that each finite set is $\Nat$-finite.

Now let $A$ be $\Nat$-finite, so there exists $n \in \Nat$  and a bijective mapping $h : [n] \to A$. Let $\mathcal{S}$ be an inductive $A$-system.
For each $k \in [n + 1] = [n] \cup \{n\}$ put $A_k = h([k])$.
Then $A_0 = h(\varnothing) = \varnothing$, $A_n = h([n]) = A$ and for each $k \in [n]$ 
\[A_{k+1} = h([k+1]) = h([k]) \cup h(\{k\}) = A_k \cup \{a_k\}\,,\]
where $a_k = h(k)$. Thus $A_0 = \varnothing \in \mathcal{S}$ and if $A_k \in \mathcal{S}$ for some $k \in [n]$ 
then also $A_{k+1} = A_k \cup \{a_k\} \in \mathcal{S}$, since $\mathcal{S}$ is an inductive $A$-system. This means that if we 
put $J = \{ k \in [n+1] : A_k \in \mathcal{S}\}$ then $0 \in J$ and $k + 1 \in J$ whenever $k \in J$ for some $k \in [n]$.
It follows that $J = [n+1]$ (insert your own proof of this fact here) and in particular $n \in J$, i.e., $A = A_n \in \mathcal{S}$.
This shows that every inductive $A$-system contains $A$ and therefore by Lemma~\ref{lemma_intro_2} $A$ is finite. \eop

There is a further characterisation of finite sets due to Tarski \cite{tarski} which will be very useful for establishing properties of such sets.
If $\mathcal{C}$ is a non-empty subset of $\mathcal{P}(E)$ then $C \in \mathcal{C}$ is said to be 
\definition{minimal} if $D \notin \mathcal{C}$ for each proper subset $D$ of $C$.

\begin{proposition}\label{prop_intro_2}
A set $E$ is finite if and only if each non-empty subset of $\mathcal{P}(E)$ contains a minimal element.
\end{proposition}

\proof 
Let $A$ be a finite set and let $\mathcal{S}$ be the set consisting of those elements $B \in \mathcal{P}(A)$ such that each non-empty subset
of $\mathcal{P}(B)$ contains a minimal element. Then $\varnothing \in \mathcal{S}$, since the only non-empty 
subset of $\mathcal{P}(\varnothing)$ is $\{\varnothing\}$ and then $\varnothing$ is the required minimal element.
Let $B \in \proper{\mathcal{S}}$ and $a \in A \setminus B$, and  $\mathcal{C}$ be a non-empty subset of 
$\mathcal{P}(B \cup \{a\})$. Put $\mathcal{D} = \mathcal{C} \cap \mathcal{P}(B)$; there are two cases:

($\alpha$)\enskip 
$\mathcal{D} \ne \varnothing$. Here $\mathcal{D}$ is a non-empty subset of $\mathcal{P}(B)$ and thus contains 
a minimal element $D$ which is then a minimal element of $\mathcal{C}$, since each set in 
$\mathcal{C} \setminus \mathcal{D}$ contains $a$ and so is not a proper subset of $D$.
 
($\beta$)\enskip 
$\mathcal{D} = \varnothing$ (and so each set in $\mathcal{C}$ contains $a$). Put 
$\mathcal{F} = \{ C \subset B : C \cup \{a\} \in \mathcal{C} \}$; then $\mathcal{F}$ is a non-empty subset 
of $\mathcal{P}(B)$ and thus contains a minimal element $F$. It follows that $F' = F \cup \{a\}$ is a minimal 
element of $\mathcal{C}$: A proper subset of $F'$ has either the form $C$ with $C \subset F$, in which case 
$C \notin \mathcal{C}$ (since each set in $\mathcal{C}$ contains $a$), or has the form $C \cup \{a\}$ with $C$ a 
proper subset of $F$ and here $C \cup \{a\} \notin \mathcal{C}$, since $C \notin \mathcal{F}$.

This shows that $B \cup \{a\} \in \mathcal{S}$ and thus that $\mathcal{S}$ is an inductive $A$-system. Hence 
$\mathcal{S} = \mathcal{P}(A)$ and in particular $A \in \mathcal{S}$, i.e., each non-empty subset of $\mathcal{P}(E)$ 
contains a minimal element.

Conversely, suppose $E$ is not finite and let 
$\mathcal{C} = \{ C \in \mathcal{P}(E) : \mbox{$C$ is not finite} \}$; then $\mathcal{C}$ is non-empty since it 
contains $E$. However $\mathcal{C}$ cannot contain a minimal element: If $D$ were a minimal element of 
$\mathcal{C}$ then $D \ne \varnothing$, since $\varnothing$ is finite. Choose $d \in D$; then $D \setminus \{d\}$
is a proper subset of $D$ and thus $D \setminus \{d\} \notin \mathcal{C}$, i.e., $D \setminus \{d\}$ is finite. 
But then by Lemma~\ref{lemma_intro_1} $D = (D \setminus \{d\}) \cup \{d\}$ would be finite. 
\eop

If $\mathcal{C}$ is a non-empty subset of $\mathcal{P}(E)$ then $C \in \mathcal{C}$ is said to be 
\definition{maximal} if $C' = C$ whenever $C' \in \mathcal{C}$ with $C \subset C' \subset E$.

\begin{proposition}\label{prop_intro_3}
If $A$ is finite then  every non-empty subset $\mathcal{C}$ of $\mathcal{P}(A)$ contains a maximal element.
\end{proposition}

\proof
The set $\mathcal{D} = \{ A \setminus C :  C \in \mathcal{C} \}$ is also a non-empty subset of $\mathcal{P}(A)$ and therefore by  
Proposition~\ref{prop_intro_2} it contains a minimal element which has the form $A \setminus C$ with $C \in \mathcal{C}$. Hence
$\{ D \in \mathcal{C} : C \subset D \subset A \} = \{C\}$ and so $C$ is maximal. \eop

\bigskip
\ref{posets}, 
We end the Introduction by outlining some of the main results to be found in these notes.
\medskip 

The notes are divided into three parts. The first part consists of Sections~\ref{fsets}, \ref{perm}, \ref{binom}, \ref{posets} and \ref{enums} and deals
with results solely involving finite sets.

In Section~\ref{fsets} we establish the basic properties of finite sets.
Most of these simply confirm that finite sets are closed under the usual set-theoretic operations.
More precisely, if $A$ and $B$ are finite sets then their union $A \cup B$, their product $A \times B$ and 
$B^A$ (the set of all mappings from $A$ to $B$) are all finite sets.
Moreover, the power set $\mathcal{P}(A)$ is finite and (Proposition~\ref{prop_intro_1}) any subset of a finite set is finite.

Apart from these closure properties there are two
properties which depend crucially on the set involved being finite. The first is  given in 
Theorem~\ref{theorem_fs_1} which states that if $A$ is a  finite set and $f : A \to A$ is a mapping then  $f$ is injective if and only if it is surjective (and 
thus if and 
only if it is bijective). 

Theorem~\ref{theorem_fs_1} implies the set $\Nat$ of natural numbers is infinite (i.e., it is not finite),
since the successor mapping $\mathsf{s} : \Nat \to \Nat$ with $\mathsf{s}(n) = n + 1$ for all $n \in \Nat$ is injective but not surjective. 

If $E, \,F$ are any sets then we write $E \approx F$ if there exists a bijective mapping $f : E \to F$. A direct corollary
of Theorem~\ref{theorem_fs_1} (Theorem~\ref{theorem_fs_2}) is
that if $B$ is a subset of a finite set $A$ with $B \approx A$ then $B = A$.

The second important property involving finite sets is
Proposition~\ref{prop_fs_2}, which states that if $A$ and $E$ are sets with $A$ finite
and if there exists either  an injective mapping $f : E \to A$ 
or a surjective mapping $f : A \to E$ then $E$ is finite.
\bigskip

For a given finite set $A$ Section~\ref{perm} looks at the group $\Selfb{A}$ of  bijective mappings $\sigma : A \to A$ (with functional composition $\circ$ as group 
operation and $\id_A$ as identity element.  The elements of $S_A$ are called \definition{permutations}.

An element $\tau$ of $\Selfb{A}$ is a \definition{transposition} 
if there exist $b,\,c \in A$ with $b \ne c$ such that
\[ 
   \tau(x) = \left\{ \begin{array}{cl}
                  c &\ \mbox{if}\ x = b\;,\\
                  b &\ \mbox{if}\ x = c\;,\\
                  x   &\  \mbox{otherwise}\;.\\
\end{array} \right. 
\]

Denote by $\twogroup$ the multiplicative group $\{+,-\}$ with ${+} \cdot {+} = {-} \cdot {-} = {+}$ and ${-} \cdot {+} = {+} \cdot {-} = {-}$.
For each element $s \in \twogroup$ the other element will be denoted by $-s$.
A  mapping $\sigma : \Selfb{A} \to \twogroup$ is a \definition{signature} if $\sigma(\id_A) = {+}\,$ and
$\sigma(\tau \circ f) = -\sigma(f)$ for each $f \in \Selfb{A}$ and for each transposition $\tau$. 

Theorem~\ref{theorem_perm_1} gives a proof of the fundamental fact that there is a unique signature $\sigma : \Selfb{A} \to \twogroup$ and that
$\sigma$ is then a group homomorphism.
\bigskip
\medskip

For each finite set $A$ and each $B \subset A$ denote the set $ \{ C \in \mathcal{P}(A) : C\approx B \}$  by  $A \,\Delta\, B$.   
The set $A \,\Delta\, B$ plays the role of a \definition{binomial coefficient}: If $|A| = n$ (with $|A|$ the usual cardinality of the set $A$)  
and $|B|= k$ then $|A\,\Delta\,B | = {n\choose k}$. In Section~\ref{binom} we establish results which correspond to some of the usual identities for binomial
coefficients.

If $A$, $B$ and $C$ are finite sets then we write $C \approx A \amalg B$ if there exist disjoint sets $A'$ and $B'$ with $A \approx A'$, $B \approx B'$ and
$C\approx A' \cup B'$. 
Theorem~\ref{theorem_binom_1} corresponds to the identity
${{n+1}\choose{k+1}} {n\choose{k+1}} + {n\choose k}$
used to generate Pascal's triangle. It states that
 if $A$ is a finite set, $B$ is a proper subset of $A$, $a\notin A$  and $b \in A \setminus B$ then 
\[(A \cup \{a\})\,\Delta\,(B \cup \{a\})\approx (A\,\Delta\, B) \amalg (A \,\Delta\, (B\cup \{b\})).\] 

If $C$ is a finite set then, as before, $\Selfb{C}$ denotes the group of bijections $h :C \to C$. If $B$ is a subset of a finite set $A$
then $I_{B,A}$ denotes the set of injective mappings $k : B \to A$.
Theorem~\ref{theorem_binom_2} states that 
if $B$ is a subset of a finite set $A$ then
$I_{B,A} \approx (A\,\Delta\, B) \times \Selfb{B}$.

Theorem~\ref{theorem_binom_5} states that if
$B$ is a subset of a finite set $A$ then $\Selfb{A} \approx I_{B,A} \times \Selfb{A \setminus B}$.

Theorem~\ref{theorem_binom_6} corresponds to the usual expression for binomial coefficients:
\[{n\choose m}  = \frac{n!}{m!\cdot (n-m)!}\] 
and states that if
$B$ is a subset of a finite set $A$ then $\Selfb{A} \times (A\, \Delta \,B) \approx \Selfb{B} \times \Selfb{A \setminus B}$.

Theorem~\ref{theorem_binom_7} corresponds to the following identity for binomial coefficients:

\[{n\choose m} {m\choose k} = {n\choose k} {{n- k}\choose {m - k}}\;.\]

It states that
if $A$, $B$ and $C$ are finite sets with $C \subset B \subset A$ then
\[(A \,\Delta\, B) \times( B\, \Delta\, C) \approx (A\,\Delta\,C) \times ((A \setminus C)\,\Delta\,(B \setminus C))\;.  \] 

\bigskip

In Section~\ref{posets} we prove Dilworth's decomposition theorem~\cite{dilworth} by modifying a proof due to
Galvin \cite{galvin} to work with the present treatment of finite sets.
This theorem states that if
$\le$ is a partial order on a finite set $A$  
then there exists a chain-partition $\mathcal{C}$ of $A$ and an antichain $D$
such that $D \approx \mathcal{C}$.

We apply Dilworth's theorem to give a proof of Hall's theorem on the existence of a system of distinct representatives \cite{hall} (known as the
Marriage Theorem).

\bigskip

In Section~\ref{enums} we give a further characterisation of a set being finite. This
can be seen as having something to do with enumerating the elements in the set. 
Let $E$ be a set.
A subset $\mathcal{U}$ of $\mathcal{P}(E)$ is called an \definition{$E$-selector} if $\varnothing \in \mathcal{U}$ and  for each 
$U \in \proper{\mathcal{U}}$ there exists  a unique element $e \in E \setminus U$ such that $U \cup \{e\} \in \mathcal{U}$. 
If $\mathcal{U}$ is an $E$-selector then a subset $\mathcal{V}$ of $\mathcal{U}$ is said to be \definition{invariant} if
$\varnothing \in \mathcal{V}$ and if $V \cup \{e\} \in \mathcal{V}$ for all $V \in \proper{\mathcal{V}}$, where 
$e$ is the unique element of $E \setminus V$ with $V \cup \{e\} \in \mathcal{U}$. 
In other words, a subset $\mathcal{V}$ of $\mathcal{U}$ is invariant if and only if it is itself an $E$-selector.
An $E$-selector  $\mathcal{U}$ is said to be \definition{minimal} if the only invariant subset of $\,\mathcal{U}$ is $\,\mathcal{U}$ itself.
A minimal $E$-selector containing $E$ will be called an \definition{$E$-enumerator}.

Theorem~\ref{theorem_enums_21} states that
an $E$-enumerator exists  if and only if $E$ is finite
and Theorem~\ref{theorem_enums_1} then
shows that if $A$ is a finite set then
an $A$-selector is minimal if and only if it is totally ordered (with respect to inclusion) and thus it is an $A$-enumerator if and only if it is totally ordered.
This implies every $A$-enumerator is totally ordered.

For each $A$-selector $\mathcal{U}$ let $\esuc{\mathcal{U}} : \proper{\mathcal{U}} \to A$ 
and $\ssuc{\mathcal{U}} : \proper{\mathcal{U}}  \to \mathcal{U} \setminus \{\varnothing\}$ be the mappings with
$\esuc{\mathcal{U}}(U) = e$ and $\ssuc{\mathcal{U}}(U) = U \cup \{e\}$, where $e$ is the unique 
element in $A \setminus U$ such that $U \cup \{e\} \in \mathcal{U}$. 
Proposition~\ref{prop_enums_1} states that
if $\,\mathcal{U}$ is an $A$-enumerator then the mappings
$\ssuc{\mathcal{U}} : \proper{\mathcal{U}} \to \mathcal{U} \setminus \{\varnothing\}$ 
and $\esuc{\mathcal{U}} : \proper{\mathcal{U}}  \to A$ are  both bijections.
In particular, if $\,\mathcal{U}$ and $\mathcal{V}$ are  $A$-enumerators then $\mathcal{U} \approx \mathcal{V}$. 
(This means, somewhat imprecisely, that any $A$-enumerator contains one more element than $A$.)

If $\mathcal{U}$ is an $A$-enumerator then for each $U \in \mathcal{U}$ the set $\mathcal{U} \cap \mathcal{P}(U)$
will be denoted by $\mathcal{U}_U$. 
Note that, as far as the definition of $\mathcal{U}_U^p$ is concerned,
$\mathcal{U}_U$ is considered here to be a subset of $\mathcal{P}(U)$ and so
$\mathcal{U}_U^p  =\{ U' \in \mathcal{U} : \mbox{$U'$ is a proper subset of $U$} \}$. 
$\mathcal{U}_U$ 
is in fact  a $U$-enumerator.

Theorem~\ref{theorem_enums_2} states that
if $\,\mathcal{U}$ and $\,\mathcal{V}$ are any two $A$-enumerators then there exists a unique
mapping $\pi : \mathcal{U} \to \mathcal{V}$ with
$\pi(\varnothing) = \varnothing$ such that $\pi(\proper{\mathcal{U}}) \subset \proper{\mathcal{V}}$ and that 
$\pi(\ssuc{\mathcal{U}}(U)) = \ssuc{\mathcal{V}}(\pi(U))$ for all $U \in \proper{\mathcal{U}}$. 
Moreover, the mapping $\pi$ is a bijection and $\pi(A) = A$.

\bigskip

Part II of the notes consists of Sections~\ref{iterators}, \ref{ordinals}, \ref{finiteiter}, \ref{am}, \ref{mam} and \ref{assoc} and studies various aspects of
what we call an iterator.
 
Recall that $\fin$ denotes the proper class of all finite sets (which is itself too large to be considered a set).
In Section~\ref{iterators} we introduce what will be called an assignment of finite sets in a triple $\,\mathbf{I} = (X, f, x_0)$, where 
$X$ is some class of
objects, $f : X \to X$ is a mapping of the class $X$ into itself and $x_0$ is an object of $X$. Such a triple will be called
an \definition{iterator}. The results in this section will
will be applied in Section~\ref{ordinals} to define the finite ordinals. In this case $(X,f,x_0) = (\fin,\vonn, \varnothing)$, where 
$\vonn : \fin \to \fin$ is the mapping given by $\vonn(A) = A \cup\{A\}$  for each finite set $A$.

The archetypal example of an iterator whose first component is a set is $(\Nat,\textsf{s},0)$.
However, we will also be dealing with examples in which the first component is a class or a finite set.
 Let $\,\mathbf{I} = (X,f, x_0)$ be an iterator.
A mapping $\val : \fin \to X$ will be called  an \definition{assignment of finite sets in $\,\mathbf{I}$}  
if $\val(\varnothing) = x_0$ and 
$\val(A \cup \{a\}) = f(\val(A))$ 
for each finite set $A$ and each element $a \notin A$. 

For each finite set $A$ denote the cardinality of $A$ by $|A|$ (with $|A|$ defined as usual in terms of  $\Nat$). Then it is clear that the mapping
$|\cdot| : \fin \to \Nat$  defines an assignment of finite sets in $(\Nat,\textsf{s},0)$ and that this is the 
unique  such assignment.

Theorem~\ref{theorem_iterators_1} states that  for each iterator $\,\mathbf{I}$ there exists a unique assignment $\val$ of finite sets in $\,\mathbf{I}$ and that  
if $A$ and $B$ are finite sets with $A \approx B$ then $\val(A) = \val(B)$. 

Let $\,X_0 = \{ x \in X : \mbox{$x = \val(A)$ for some finite set $A$} \}$. 
A subclass $Y$ of $X$ is said to be \definition{$f$-invariant} if $f(y) \in Y$ for all $y \in Y$. 
Lemma~\ref{lemma_iterators_3} shows that $X_0$ is the least $f$-invariant subclass of $X$ containing $x_0$.

The iterator $\,\mathbf{I}$ is said to be \definition{minimal} if the only 
$f$-invariant subclass of $X$ containing $x_0$ is $X$ itself, thus $\,\mathbf{I}$ is minimal if and only if
$X_0 = X$. In particular, it is easy to see that the Principle of Mathematical Induction is exactly the requirement that the iterator 
$(\Nat,\textsf{s},0)$ be minimal. If $\,\mathbf{I}$ is minimal then Lemma~\ref{lemma_iterators_3} implies that
for each $x \in X$ there exists a finite set $A$ with $x = \val(A)$.

The iterator $\,\mathbf{I}$ will be called \definition{regular} if $B_1 \approx B_2$ whenever $B_1$ and $B_2$ are finite sets with 
$\val(B_1) = \val(B_2)$.
If $\,\mathbf{I}$ is regular then by Theorem~\ref{theorem_iterators_1} $B_1 \approx B_2$ holds if and only if $\val(B_1) = \val(B_2)$.

$\,\mathbf{I}$ will be called a \definition{Peano iterator} if it is minimal
and \standard, where \definition{\standard}{} means that the mapping $f$ is injective and $x_0 \notin f(X)$.
The Peano axioms thus require $(\Nat,\textsf{s},0)$ to be a Peano iterator.
If $\,\mathbf{I} = (X,f,x_0)$ is a Peano iterator then by Theorem~\ref{theorem_fs_1} $X$ cannot b a finite set.
Theorem~\ref{theorem_iterators_2} states that a minimal iterator is regular if and only if it is a Peano iterator. A corollary is
Theorem~\ref{theorem_iterators_3} (the recursion theorem for Peano iterators). It states that
if $\,\mathbf{I}$ is a Peano iterator then  for each iterator $\,\mathbf{J} = (Y,g,y_0)$ there exists a unique mapping
$ \pi: X \to Y$ with $\pi(x_0) = y_0$ such that $\pi \circ f = g \circ \pi$.
The Peano axioms require $(\Nat,\mathsf{s},0)$ to be a Peano iterator and hence  for each iterator $\,\mathbf{J} = (Y,g,y_0)$ there exists 
a unique mapping
$ \pi: \Nat \to Y$ with $\pi(0) = y_0$ such that $\pi \circ \mathsf{s} = g \circ \pi$.
Theorem~\ref{theorem_iterators_4} shows that
if $\,\mathbf{I} = (X,f,x_0)$ is minimal
 then the class $X$ is  a finite set if and only if $\,\mathbf{I}$ is not regular.
Theorems~\ref{theorem_iterators_2} and \ref{theorem_iterators_4} imply that for a minimal iterator $\,\mathbf{I} = (X,f,x_0)$ there are two mutually exclusive possibilities: Either $\mathbf{I}$ is a Peano iterator or $X$ is a finite set.

Let $\,\mathbf{I} =(X,f,x_0)$ and 
$\,\mathbf{J}= (Y,g,y_0)$ be iterators; a mapping $\mu : X \to Y$ is called a \definition{morphism} from 
$\,\mathbf{I}$ to $\,\mathbf{J}$ if $\mu(x_0) = y_0$ and $g\circ \mu = \mu \circ f$.
 This will also be expressed by 
saying that $\mu : \mathbf{I} \to \mathbf{J}$ is a morphism. 
The iterators $\,\mathbf{I}$ and $\,\mathbf{J}$ are said to be 
\definition{isomorphic} if there exists a  morphism $\mu : \mathbf{I} \to \mathbf{J} $ and a morphism $\nu : \mathbf{J} \to \mathbf{I}$ such that 
$\nu \circ \mu = \id_{X}$ and  $\mu \circ \nu = \id_{Y}$. In particular, the mappings $\mu$ and $\nu$ are then both bijections.   
The iterator $\,\mathbf{I}$ is said to be \definition{initial} if for each iterator $\,\mathbf{J}$ there 
is a unique morphism from $\,\mathbf{I}$ to $\,\mathbf{J}$. Theorem~\ref{theorem_iterators_3} (the recursion theorem) thus states that a Peano iterator is initial.
Let  $\,\mathbf{I}$ be initial  and $\pi : \mathbf{I} \to \mathbf{J}$ be the unique morphism.
If $\,\mathbf{J}$ is initial then $\pi$ is an isomorphism and so $\,\mathbf{I}$ and $\,\mathbf{J}$ are isomorphic.
Conversely, if $\pi$ is an isomorphism then $\,\mathbf{J}$ is initial.
This shows that, up to isomorphism, there is a unique initial iterator. Of course, this is only
true if an initial iterator exists, and $(\Nat,\mathsf{s},0)$ is an initial iterator. In fact, 
in Proposition~\ref{prop_iterators_222} we introduce a Peano (and thus initial) iterator $\,\mathbf{H} = (H,h,\varnothing)$
which is defined only in terms of finite sets. The elements occurring in the sets of $H$ are the hereditarily finite sets.
Moreover, in Section~\ref{ordinals} we will exhibit another Peano iterator
$\mathbf{O} = (O,\vonn,\varnothing)$,which again is defined only in terms of finite sets.
The elements of $O$ are the finite ordinals. 

Theorem~\ref{theorem_iterators_6} is a result of Lawvere \cite{lawvere} which  shows the converse of the recursion theorem holds. 
That is, an initial iterator $\,\mathbf{I} = (X,f,x_0)$ is a Peano iterator.

Let $\,\mathbf{I} = (X,f,x_0)$ be an iterator. A total order $\le$ on $X$ will be called \definition{compatible with $\,\mathbf{I}$} if
$x \le f(x)$ for all $x \in X$.

\medskip

Theorem~\ref{theorem_iterators_199} states that if
$\,\mathbf{I} = (X,f,x_0)$ is a Peano iterator
then there exists a unique total order $\le$ on $X$ compatible with $\,\mathbf{I}$ and the following hold:

(1)\enskip $x_0 \le x$ and $x <f(x)$ for all $x \in X$, where as usual we write $y < x$ if $y \le x$ but $y \ne x$.

(2)\enskip If $y \le x$ then $f(y) \le f(x)$.

(3)\enskip Let $x,\,y \in X$. Then $y \le f(x)$ if and only if $y \le x$ or $y = f(x)$. Thus $y < f(x)$ if and only if $y \le x$.

(4)\enskip For each $x \in X$ let $L_x = \{ y \in X : y < x \}$. Then $L_{x_0} = \varnothing$ and
$L_{f(x)}$ is the disjoint union of $\{x_0\}$ and $f(L_x)$ for each $x \in X$.

(5) \enskip The subclass  $L_x$ is a finite set for each $x \in X$.

(6) \enskip Each non-empty subclass $Y$ of $X$ contains a minimum element, i.e., an element $x$ with $x \le y$ for all $y \in Y$.

Let $\,\mathbf{I} = (X,f,x_0)$ be a Peano iterator with $\le$ the unique total order
on $X$ compatible with $\,\mathbf{I}$. Define an iterator $\,\mathbf{I}_\le = (X_\le,f_\le,\varnothing)$
by letting $X_\le = \{L_x : x \in X\}$ and with $f_\le : X_\le \to X_\le$ given by $f_\le(L_x) = L_{f(x)}$ for all $x \in X$. Also define $\pi_\le : X \to X_\le$
by $\pi_\le(x) = L_x$ for all $x \in X$. 
Then Theorem~\ref{theorem_iterators_99} states that:

(1)\enskip $\pi_\le : \mathbf{I} \to \mathbf{I}_\le$ is an isomorphism.

(2)\enskip $\,\mathbf{I}_\le$ is a Peano iterator.

(3)\enskip $f_\le(L_x)$ is the disjoint union of $\{x_0\}$ and $f(L_x)$ for each $x \in X$.

(4)\enskip The sets in $X_\le$ are totally ordered by inclusion and inclusion is the unique total order on $X_\le$ compatible with $\,\mathbf{I}_\le$.

(5)\enskip Let $\val_\le : \fin \to X_\le$ be the assignment of finite sets in $\,\mathbf{I}_\le$.
Then $\val_\le(A) \approx A$ for all $A \in \fin$.
 
We call $\,\mathbf{I}_\le$ the finite segment iterator associated with $\,\mathbf{I}$.

\medskip

Theorem~\ref{theorem_iterators_11} is a result which guarantees the existence of mappings 'defined by recursion'. It states that
if $\,\mathbf{I} = (X,f,x_0)$ is a Peano iterator, $Y$ and $Z$ are classes and $\beta : X \times Y \to Z$
and $\alpha : X \times Y \times Z \to Z$ are  mappings then there is a unique mapping
$\pi : X \times Y \to Z$ with $\pi(x_0,y) = \beta(y)$ for all $y \in Y$ such that
$\pi(f(x),y) = \alpha(x,y,\pi(x,y))$ for all $x \in X$, $y \in Y$.

In Section~\ref{ordinals}  we study finite ordinals using the standard approach introduced by von Neumann \cite{vonn}. 
Let $\vonn : \fin \to \fin$ be the mapping given by $\vonn(A) = A \cup \{A\}$ for 
each finite set $A$. 
If we iterate the operation $\vonn$ starting with the empty set and label the resulting sets using the natural numbers then we obtain the following:

$0 =\varnothing$,\\
$1 = \vonn(0) = 0 \cup\{0\} = \varnothing \cup \{0\} = \{0\}$,\\ 
$2 = \vonn(1) = 1 \cup \{1\} = \{0\} \cup \{1\} = \{0,1\}$,\\
$3 = \vonn(2) = 2 \cup\{2\} = \{0,1\} \cup \{2\} = \{0,1,2\}$,\\
$4 = \vonn(3) = 3 \cup\{3\} = \{0,1,2\} \cup \{3\} = \{0,1,2,3\}$,\\
$5 = \vonn(4) = 4 \cup \{4\} = \{0,1,2,3\} \cup \{4\} = \{0,1,2,3,4\}$,
 
$n+1 =  \vonn(n) = n \cup \{n\} = \{0,1,2,\ldots,n-1\} \cup\{n\} = \{0,1,2,\ldots,n\}$ .
\bigskip

Denote by $\,\mathbf{O}'$ the iterator $(\fin,\vonn,\varnothing)$.
Then by Theorem~\ref{theorem_iterators_1} there exists a unique assignment $\ordass$
of finite sets in $\,\mathbf{O}'$. Thus $\ordass : \fin \to \fin$  is the unique mapping with 
$\ordass(\varnothing) = \varnothing$ and such that
$\ordass(A \cup \{a\}) = \vonn(\ordass(A))$
for each finite set $A$ and each element $a \notin A$.
Moreover, if $A$ and $B$ are finite sets with $A \approx B$ then $\ordass(A) = \ordass(B)$. 
Theorem~\ref{theorem_ordinal_1} states that $\ordass(A) \approx A$ for each finite set and thus  $\ordass(A) = \ordass(B)$ if and only if $A \approx B$.
Let $\,O = \{ B \in \fin : \mbox{$B = \ordass(A)$ for some finite set $A$} \}$. The elements of $O$ will be called 
\definition{finite ordinals}. Thus for  each finite set $A$  there exists  a unique  $o \in O$  with $o = \ordass(A)$.
By Lemma~\ref{lemma_iterators_3} $O$ is the least $\vonn$-invariant subclass of $\fin$ containing $\varnothing$.
We denote the restriction of $\vonn$ to $O \to O$ again by $\vonn$. Then $\,\mathbf{O} = (O,\vonn,\varnothing)$ is a minimal 
iterator.

In fact,Theorem~\ref{theorem_ordinal_2}
shows that$\,\mathbf{O}$ is a Peano iterator.

Proposition~\ref{prop_ordinal_1} states that for each finite set $A$
\[\ordass(A) = \{ o \in O : o = \ordass(A') \mbox{ for some proper subset $A'$ of $A$} \}.\]

It follows that if $A$ is a finite set and $B\subset A$; then $\ordass(B) \subset \ordass(A)$.
It also follows that for each $o \in O$ 
\[o = \{ o' \in O : \mbox{$o'$ is a proper subset of $o$} \},\]
\[\vonn(o) = \{ o' \in O :\mbox{ $o'$ is a subset of $o$} \}.\]

Proposition~\ref{prop_ordinal_11} states that if
$o,\, o'\in O$ with $o \ne o'$.  then either $o$ is a proper subset of $o'$ or $o'$ is a proper subset of
$o$.

\bigskip

We also  consider a situation which is somewhat more general than that occurring with the iterator $\,\mathbf{O}$.
Let $\,\mathbf{J} = (T,h,\varnothing)$ be a minimal iterator with $T$ a subclass of $\fin$. We call $\,\mathbf{J}$ an \definition{ordinal iterator}
if for each $B \in T$ there exists an element $b \notin B$ such that $h(B) = B \cup \{b\}$. Note that if it is not assumed that $\,\mathbf{J}$ is minimal
then the associated minimal iterator $\,\mathbf{J}_0$ will be an ordinal iterator. The archetypal example of an ordinal iterator is of course 
$\,\mathbf{O}$. 
Moreover, if $\,\mathbf{I}$ is a Peano iterator and $\,\mathbf{I}_\le$ is the initial segment iterator associated with $\,\mathbf{I}$
then Theorem~\ref{theorem_iterators_199} (4) shows that  $\,\mathbf{I}_\le$ 
will be an ordinal iterator.

Let $\,\mathbf{J} = (T,h,\varnothing)$ be an ordinal iterator and let $\tau : \fin \to T$ be the evaluation of finite sets in $\,\mathbf{J}$. 

Theorem~\ref{theorem_ordinal_3} states that

(1)\enskip $\tau(A) \approx A$ for all $A \in \fin$ and therefore $\tau(A) = \tau(A')$ if and only if $A \approx A'$. In particular
$\tau(\tau(A)) = \tau(A)$ for all $A \in \fin$ and $\tau(B) = B$ for all $B \in T$ (since $\tau$ is surjective).

(2)\enskip $\,\mathbf{J}$ is  a Peano iterator.

(3)\enskip If $A,\,A' \in \fin$ with $A \subset A'$ then $\tau(A) \subset \tau(A')$.

(4)\enskip If $A,\,A' \in \fin$ with $A \preceq A'$ then $\tau(A) \subset \tau(A')$.

(5)\enskip For all $B,\,B' \in T$ either $B \subset B'$ or $B' \subset B$. Thus $T$ is totally ordered by inclusion. Moreover, inclusion is the unique total
order on $T$ compatible with $\,\mathbf{J}$.

Again let $\,\mathbf{J} = (T,h,\varnothing)$ be an ordinal iterator.
Let $T^*$ be the class consisting of all elements $c$ for which there exists $B \in T$ such that $h(B) = B \cup \{c\}$.
Define $\gamma: T \to T^*$ by letting $\gamma(B) = c$, where $h(B) = B \cup \{c\}$. Thus $h(B) = B \cup \{\gamma(B)\}$  for all $B \in T$.

Theorem~\ref{theorem_ordinal_4} states that:

(1)\enskip
The mapping $\gamma : T \to T^*$ is a bijection.

(2)\enskip Define  a mapping $h^*: T^* \to T^*$ by $h^* = \gamma \circ h \circ \gamma^{-1}$ and an iterator by $\,\mathbf{J}^* = (T^*,h^*,t^*_0)$, where $t^*_0 = \gamma(\varnothing)$ 
and so $t^*_0$ is the single element in $h(\varnothing)$.  Then $\gamma : \,\mathbf{J} \to \mathbf{J}^*$ is an isomorphism.

(3)\enskip $\,\mathbf{J}^*$ is a Peano iterator.

(4)\enskip $h(\gamma^{-1}(t)) = \gamma^{-1}(t) \cup \{t\}$ for all $t \in T^*$.

(5)\enskip  Define $\le$ on $T^*$ by letting $t' \le t$ if and only if
$\gamma^{-1}(t') \subset \gamma^{-1}(t)$. 
Then $\le$ is the unique a total order on $T^*$      compatible with $\,\mathbf{J}^*$.  

(6)\enskip $L_t =  \gamma^{-1}(t)$ for all $t \in T^*$.

(7)\enskip $h(L_t) = L_t \cup \{t\}$ for all $t \in T^*$.

We call the iterator $\,\mathbf{J}^*$ the \definition{dual iterator} to the ordinal iterator $\,\mathbf{J}$.

The canonical ordinal iterator $\mathbf{O}$ is its own dual.

Let $\,\mathbf{I} = (X,f,x_0)$ be a Peano iterator and let $\,\mathbf{I}_\le$ be the finite segment iterator associated with $\,\mathbf{I}$.
Also let $\,\mathbf{I}_\le^*$ be the dual  iterator  to the ordinal iterator $\,\mathbf{I}_\le$. 
Proposition~\ref{prop_ordinal_22} shows that $\,\mathbf{I}_\le^* = \,\mathbf{I}$.

Let $\,\mathbf{J} = (T,h,\varnothing)$ be an ordinal iterator and let $\,\mathbf{J}^* = (T^*,h^*,t^*_0)$ be the dual iterator to $\,\mathbf{J}$.
Also let $\,\mathbf{J}^*_\le$ be the finite segment iterator associated with $\,\mathbf{J}^*$.
Then Proposition~\ref{prop_ordinal_23} states that 
$\,\mathbf{J} = \mathbf{J}^*_\le$.

\medskip

Theorems~\ref{theorem_iterators_2} and \ref{theorem_iterators_4} imply that for a minimal iterator $\,\mathbf{I} = (X,f,x_0)$ there are 
two mutually 
exclusive possibilities: Either $\,\mathbf{I}$ is  a Peano iterator or $X$ is a finite set. In Section~\ref{finiteiter} we deal with the case when 
$X$ is a finite set.

Let $\,\mathbf{I} = (X,f,x_0)$ be a minimal iterator with $X$ a finite set. 
For each $x \in X$ let $X_x$ be the least $f$-invariant subset of $X$ containing $x$ and let $f_x : X_x \to X_x$ be the restriction of $f$ to $X_x$. Thus
$\,\mathbf{I}_x = (X_x,f_x,x)$ is a minimal iterator. 
An element $x \in X$ is said to be \definition{periodic} if $x \in f_x(X_x)$ and so by Proposition~\ref{prop_iterators_1} and Theorem~\ref{theorem_fs_1} $x$ is
 periodic if and only if $f_x$ is a bijection.  Theorem~\ref{theorem_finiteiter_1} states that:

(1) \enskip Let $X_P = \{ x \in X : \mbox{ $x$ is periodic} \}$.
Then $X_P$ is non-empty and  $X_x = X_y$ for all $x,\,y \in X_P$.  Thus $f$ maps $X_P$ bijectively onto itself.

(2) \enskip Let $X_N =\{ x \in X : \mbox{ $x$ is not periodic} \}$ and suppose $X_N \ne \varnothing$.
Then $f$ is injective on $X_N$ and there exists a unique element $u \in X_N$ such that $f(u)$ is periodic. Moreover, there exists a unique element 
$v \in X_P$ such that $f(v) = f(u)$ and $u$ and $v$ are the unique elements of $X$ with $u \ne v$ such that $f(u) = f(v)$.

Statement (1) corresponds to the elementary fact that a mapping $f : X \to X$ with $X$ a finite set is eventually periodic.
\bigskip

In Section~\ref{am} we show how an addition and a multiplication can be defined for any minimal iterator. 
These operations are associative and commutative and can be specified by the rules (a0), (a1), (m0) and (m1) 
below, which are usually employed when defining the operations on $\Nat$ via the Peano axioms.

Let $\,\mathbf{I} = (X,f,x_0)$ be a minimal iterator  with $\val$ the assignment of finite sets in $\,\mathbf{I}$.
Theorem~\ref{theorem_am_1}
states that there exists a unique binary operation $\oplus$ on $X$ such that
\[ 
\val(A) \oplus \val(B) = \val(A \cup B) 
\]
whenever $A$ and $B$ are disjoint finite sets. This operation $\oplus$ is both associative and commutative, 
$x \oplus x_0 = x$ for all $x \in X$ and for all $x_1,\,x_2 \in X$ there is a $x \in X$ such that either 
$x_1 = x_2 \oplus x$ or $x_2 = x_1 \oplus x$. Moreover, $\oplus$ is the unique binary operation $\oplus$ on $X$ 
such that
\begin{evlist}{15pt}{6pt}
\item[(a0)]  $x \oplus x_0 = x$ for all $x \in X$. 

\item[(a1)]  $x \oplus f(x') = f(x \oplus x')$ for all $x,\,x' \in X$.
\end{evlist}

\medskip

If $f$ is injective then the cancellation law holds for $\oplus$ (meaning  that $x_1 = x_2$ whenever 
$x_1 \oplus x = x_2 \oplus x$ for some $x \in X$).

If $x_0 \in f(X)$ then by Theorem~\ref{theorem_iterators_4} and Proposition~\ref{prop_iterators_2} $X$ is finite and  
$f$ is bijective and here $(X,\oplus,x_0)$ is the cyclic group.generated by the 
element $f(x_0)$. 

If $\,\mathbf{I}$ is a Peano iterator and $\le$ is the unique total order on $X$ compatible with $\,\mathbf{I}$ then
$y \le x$ if and only if $x = y \oplus z$ for some $z\in X$.
\bigskip

Theorem~\ref{theorem_am_2} states that 
there exists a unique binary operation $\otimes$ on $X$ such that
\[ 
\val(A) \otimes \val(B) = \val(A\times B) 
\]
for all finite sets $A$ and $B$. This operation $\otimes$ is both associative and commutative, 
$x \otimes x_0 = x_0$ for all $x\in X$ and $x \otimes f(x_0) = x$ for all $x \in X$ with $ x \ne x_0$ (and so $f(x_0)$ is a 
multiplicative identity) and the distributive law holds for $\oplus$ and $\otimes$: 
\[x \otimes (x_1 \oplus x_2) = (x \otimes x_1) \oplus (x \otimes x_2)\]
for all $x,\,x_1,\,x_2 \in X$.
Moreover, $\otimes$ is the unique binary operation on $X$ such that 
\begin{evlist}{15pt}{6pt}
\item[(m0)]  $x \otimes x_0 = x_0$ for all $x \in X$. 

\item[(m1)]  $x \otimes f(x') = x \oplus (x \otimes x')$ for all $x,\,x' \in X$.
\end{evlist}

We also look at the operation of exponentiation. Here we have to be more careful: For example, 
$2 \cdot 2 \cdot 2 = 2$ in $\Int_3$ and so $2^3$ is not well-defined if the exponent $3$ is considered as an 
element of $\Int_3$ (since we would also have to have $2^0 = 1$). However, $2^3$ does make sense if $2$ is 
considered as an element of $\Int_3$ and the exponent $3$ as an element of $\Nat$.

In general it is the case that if $\,\mathbf{J} = (Y,g,y_0)$ is a Peano iterator then we can define an element of $X$ which is 
`$x$ to the power of $y$' for each $x \in X$ and each $y \in Y$ and this operation has the properties which might 
be expected. 

Let $\,\mathbf{J} = (Y,g,y_0)$ be a Peano iterator with $\val'$ the assignment of finite sets in $\,\mathbf{J}$. (As 
before $ \,\mathbf{I} = (X,f,,x_0)$ is assumed to be minimal with $\val$ the assignment of finite sets in $\,\mathbf{I}$.)
Also let $\oplus$ and $\otimes$ be the operations given in Theorems~\ref{theorem_am_1} and \ref{theorem_am_2} for the iterator
$\,\mathbf{I}$.

Theorem~\ref{theorem_am_3} states that there exists a unique operation ${\uparrow} : X \times Y \to X$ such that  
\[ 
\val(A) \uparrow \val'(B) = \val(A^B) 
\]
for all finite sets $A$ and $B$. This operation ${\uparrow}$ satisfies
\[     
x \uparrow (y_1 \oplus y_2) = (x \uparrow y_1) \otimes (x \uparrow y_2) 
\]
for all $x \in X$ and all $y_1,\,y_2 \in Y$ and
\[     
(x_1 \otimes x_2) \uparrow y = (x_1 \uparrow y) \otimes (x_2 \uparrow y) 
\]
for all $x_1,\,x_2 \in X$ and $y \in Y$. Moreover, ${\uparrow}$ is the unique operation such that
\begin{evlist}{15pt}{6pt}
\item[(e0)]  $x \uparrow y_0 = f(x_0)$ for all $x \in X$. 

\item[(e1)]  $x \uparrow g(y) = x \otimes (x \uparrow y)$ for all $x \in X$, $y\in Y$.
\end{evlist}

\bigskip

In Section~\ref{mam} we give alternative proofs for Theorem~\ref{theorem_am_1} and Theorem~\ref{theorem_am_2}. 

\bigskip

Section~\ref{assoc} really belongs in Part 1 of the notes, but the results in Section~\ref{finiteiter} are needed here.

Let $\bullet$ be a binary operation on a set $X$, written using infix notation, so $x_1 \bullet x_2$ denotes the product of $x_1$ and $x_2$. The large majority of 
such operations occurring in mathematics are \definition{associative}, meaning that
$(x_1 \bullet x_2) \bullet x_3 = x_1 \bullet ( x_2 \bullet x_3)$ for all $x_1,\,x_2,\,x_3 \in X$. If $\bullet$ is associative and 
$x_1,\,x_2,\,\ldots,\,x_n \in X$ then the product $x_1 \bullet x_2 \bullet \cdots \bullet x_n$ is well-defined, meaning its value does not depend
on the order in which the operations are carried out.

This result will be  established in Section~\ref{assoc}. We first define a particular order of carrying out the operations. This is the order in which at each
stage the product of the current first and second components are taken. For example, the product of the 6 components $x_1,\,x_2,\,x_3,\, x_4,\,x_5, \, x_6$
 evaluated using this order results in the value
$\bullet(x_1,\ldots,x_6) = (((((x_1 \bullet x_2) \bullet x_3) \bullet x_4)  \bullet x_5)\bullet x_6)$.
In general, the corresponding product of $n$ terms will be denoted by $\bullet(x_1,\ldots,x_n)$.
Theorem~\ref{theorem_assoc_1} states that if $\bullet$ is associative then
$\bullet(x_1,\ldots,x_m,x_{m+1},\ldots,x_n) = \alpha \bullet \beta$,
where $\alpha = \bullet(x_1,\ldots,x_m)$ and $ \beta = \bullet(x_{m+1},\ldots,x_n)$.
This is a weak form of the generalised associative law, although it is one which is often all that is needed.

Theorem~\ref{theorem_assoc_2} gives the general form of the generalised associative law and states that if $\bullet$ is associative then
$\bullet(x_1,\ldots,x_n) = \bullet_{\mathbf{R}}(x_1,\ldots,x_n)$ for each $\,\mathbf{R}$ from the set of prescriptions describing how the operations are
carried out. The main task is to give a rigorous definition of this set. We do this using partitions of intervals of the form
$\{ k \in \Int : m \le  k \le n \}$ in which each element in the partition is also an interval of this form.

By a \definition{partition} of a set $S$ we mean  a subset $\mathcal{Q}$ of $\mathcal{P}_0(S)$ such that for each $s \in S$ there exists a unique 
$Q \in \mathcal{Q}$ such that $s \in Q$. Thus, different elements in a partition of $S$ are disjoint and their union is $S$. 

Consider the product $((x_1 \bullet x_2) \bullet ((x_3 \bullet x_4) \bullet x_5))$.
The order of operations involved here can described with the help of the following sequence of partitions of the set $\{1,2,3,4,5\}$:

$\{\{1\},\{2\},\{3\},\{4\},\{5\}\}$\\
$\{\{1\},\{2\},\{3,4\},\{5\}\}$\\
$\{\{1,2\},\{3,4\},\{5\}\}$\\
$\{\{1,2\},\{3,4,5\}\}$\\
$\{\{1,2,3,4,5\}\}$

For each of these partitions (except the last one) the next partition  is obtained by amalgamating two adjacent partitions. Corresponding to
these partitions there  is a sequence of partial evaluations:

$\{\{x_1\},\{x_2\},\{x_3\},\{x_4\},\{x_5\}\}$\\
$\{\{x_1\},\{x_2\},\{(x_3\bullet x_4)\},\{x_5\}\}$\\
$\{\{(x_1\bullet x_2)\},\{(x_3\bullet x_4)\},\{x_5\}\}$\\
$\{\{(x_1 \bullet x_2)\},\{((x_3 \bullet x_4) \bullet x_5)\}\}$\\
$\{\{((x_1 \bullet x_2 )\bullet ((x_3 \bullet x_4) \bullet x_5))\}\}$

and the final expression is essentially the product we started with.

\bigskip

Part III of the notes consists of Section \ref{lists} and presents an approach to dealing with (finite) lists taking their values in some fixed class $E$. 
We first formulate things using the natural numbers but
then employ a  general Peano iterator in its place. For each $n \in \Nat$ lists of length $n$ are mappings from $L_n =\{0,1.\ldots,n-1\}$. to $E$.
A list $t : L_n \to E$ with $t(k) = e_k$ for all $k \in L_n$ will be represented in the form $[e_0,e_1,\ldots,e_{n-1}]$. In particular, $[]$
represents the empty list, i.e., the list with no elements.

When using the natural numbers the sets $L_n = \{0,1,\ldots,n-1\}$ play an important role and these correspond to the sets in the finite segment iterator associated
with a general Peano iterator. The natural numbers will only be used when giving examples. Let us fix a Peano iterator
$\,\mathbf{I} = (N,s,n_0)$  with $\le$ the unique total order on $N$ compatible with $\,\mathbf{I}$
and let $\,\mathbf{I}_\le = (N_\le,s_\le,\varnothing)$ be the finite segment iterator associated with $\,\mathbf{I}$.
Thus $N_\le = \{L_n : n \in N \}$ with $L_n = \{ m \in N : m < n \}$ and $s_\le(L_n) = L_{s(n)}$ for all $n \in N$.

For each $n\in N$ denote by $E^*_n$ the class of all mappings from $L_n$ to $E$. In particular, 
$L_{n_0} = \varnothing$ and so $E^*_{n_0}$ consists of the unique mapping from $\varnothing$ to $E$, which we denote by $\varepsilon$. 

We call $E^*_n$  \definition{the class of lists based on $L_n$ with
values in $E$}. The single element $\varepsilon \in E^*_{n_0}$ will be referred to as the empty list.
Let $m,\,n \in N$ with $m \ne n$; then $L_m \ne L_n$ and hence $E^*_m$ and $E^*_n$ are disjoint.
Put $E^* = \bigcup_{n\in N} E^*_n$; this is the \definition{class of all lists with values in $E$}. Also let
$E^*_+ = E^* \setminus \{\varepsilon\}$.

We define mappings $\triangleleft : E \times E^* \to E^*_+$ 
and $\triangleright : E^*_+ \to E \times E^*$. 
The list $\triangleleft(e,r)$ is obtained by adding the element $e$ to the beginning of the list $r$ and so
$\triangleleft(e,[e_0,e_1,\ldots,e_{n-1}]) = [e,e_0,e_1\ldots,e_{n-1}]$. In particular $\triangleleft(e,[\;]) = [e]$ is a non-empty list.
If $r$ is a non-empty list then the first component of $\triangleright(r)$ is the first element of $r$  (the head of the list) and the second component of $\triangleright(r)$ is the
rest of the list (its tail). Thus $\triangleright([e_0,e_1,\ldots,e_{n-1}]) = (e_0,[e_1,\ldots,e_{n-1}])$.
Proposition~\ref{prop_lists_1} states that
he mappings $\triangleleft : E \times E^* \to E^*_+$ 
and $\triangleright : E^*_+ \to E \times E^*$ are inverse to each other and in particular they are both bijections.

A triple $(X,f,x_0)$ with $X$ a class, $f : E \times X \to X$ a mapping and $x_0$ an element of
$X$ will be called a \definition{list algebra}. Thus $(E^*,\triangleleft,\varepsilon)$ is a list algebra. If
$(X,f,x_0)$ is a list algebra then for each $e \in E$ let $f_e : X \to X$ be the mapping with
$f_e(x) = f(e,x)$ for all $x \in X$. There is then the iterator $(X,f_e,x_0)$.
A subclass $X_0$ of $X$ is said to be \definition{$f$-invariant} if $f_e(X_0) \subset X_0$ for all $e \in E$ and
$(X,f,x_0)$ is said to be \definition{minimal} if $X$ itself is the only $f$-invariant subclass of $X$ containing $x_0$.

Lemma~\ref{lemma_lists_1} states that
the list algebra $(E^*,\triangleleft,\varepsilon)$ is minimal.

Let $(X,f,x_0)$ and $(Y,g,y_0)$ be list algebras.
A  morphism $\pi : (X,f,x_0)) \to (Y,g,y_0)$ is then a mapping $\pi : X \to Y$ with $\pi(x_0) = y_0$ such that
$  g_e \circ\ \pi = \pi \circ f_e$ for all $e \in E$.  Thus $\pi : (X,f,x_0) \to (Y,g,y_0)$ being a morphism means exactly that $\pi (X,f_e,x_0) \to (Y,g_e,y_0)$
is a morphism of iterators for each $e \in E$.
A list algebra $(X,f,x_0)$ is said to be \definition{initial} if for each list algebra $(Y,g,y_0)$ there exists a unique morphism
$\pi : (X,f,x_0) \to (Y,g,y_0)$.

Theorem~\ref{theorem_lists_1} states that the list algebra $(E^*,\triangleleft,\varepsilon)$ is initial.
An important application of  this result is the following:
Let $r \in E^*$. Then there is a unique morphism $\psi_r :(E^*,\triangleleft,\varnothing) \to (E^*,\triangleleft,r)$. Thus $\psi_r$ is the unique mapping 
$\psi_r : E^* \to E^*$ with
$\psi_r(\varepsilon) = r$ and such that $\psi_r(\triangleleft(e,t)) = \triangleleft(e,\psi_r(t))$ for all $(e,t) \in E \times E^*$.
 
The mapping $\psi_r$ appends the list $r$ to its argument. This is verified by giving an explicit expression for the mapping $\psi_r$.

For $r,\,t \in E^*$ we write $t \bowtie r$ instead of $\psi_r(t)$. 
Thus $\psi_r(t) = t \bowtie r$ for all $t, \,r \in E^*$. We consider $\bowtie$ as an infix operation on $E^*$. By the uniqueness of $\psi_r$ it follows that
$\bowtie$ is the unique operation on $E^*$ with   $\varepsilon \bowtie r = r$ for all $r \in E^*$ 
such that $\triangleleft(e,t) \bowtie r = \triangleleft(e, (t \bowtie r))$ for all $(e,t) \in E \times E^*$ and all $r \in E^*$.

Now since  $t \bowtie r$ appends the list $r$ to the list $t$,  it is to be expected that $\bowtie$ is associative, i.e., that
$t\bowtie (u \bowtie v) = (t \bowtie u) \bowtie v$ for all $t,\,u,\,v \in E^*$. This is shown to be the case.

A list algebra $(X,f,x_0)$ will be called \definition{unambiguous} if the mapping $f_e$ is 
injective for each $e \in E$ and the classes $f_e(X)$, $e \in E$, are disjoint and 
$x_0 \notin \bigcup_{s\in S} f_e(X)$. Being unambiguous corresponds to being $\Nat$-like for an iterator.

By Lemmas~\ref{lemma_lists_1} and \ref{lemma_lists_2} the list algebra $(E^*,\triangleleft,\varepsilon)$ is both minimal and unambiguous
and by Theorem~\ref{theorem_lists_1} $(E^*,\triangleleft,\varepsilon)$ is initial.

Theorem~\ref{theorem_lists_2} states that a list algebra is initial if and only if it is both minimal and unambiguous. This corresponds to
Theorem~\ref{theorem_iterators_2} for iterators.

\medskip

We end the Introduction with the following important fact:

\begin{lemma}\label{lemma_intro_3}
For each set $A$ there exists an element $a$ not in $A$. If $A$ is finite then by Lemma~\ref{lemma_intro_1} 
$A \cup \{a\}$ will also be finite. 
\end{lemma}

\proof
In fact there must exist an element in $\mathcal{P}(A) \setminus A$. If this were  not the case then 
$\mathcal{P}(A) \subset A$, and we could define a surjective mapping $f : A \to \mathcal{P}(A)$ by letting 
$f(x) = x$ if $x \in \mathcal{P}(A)$ and $f(x) = \varnothing$ otherwise. But by Cantor's diagonal argument (which 
states that a mapping $f : X \to \mathcal{P}(X)$ cannot be surjective) this is not possible.
\eop


\startsection{Finite sets}

\label{fsets}

Recall that a set $A$ is defined to be \definition{finite} if $\mathcal{P}(A)$ is the only inductive $A$-system, 
where a subset $\mathcal{S}$ of the power set $\mathcal{P}(A)$ is an inductive $A$-system if 
$\varnothing \in \mathcal{S}$ and $B \cup \{a\} \in \mathcal{S}$ for all $B \in \proper{\mathcal{S}}$, $a \in A \setminus B$ (and where $\proper{\mathcal{S}}$
denotes the set of subsets in $\mathcal{S}$ which are proper subsets of $A$).

In this section we establish the basic properties of finite sets.
Most of these simply confirm that finite sets are closed under the usual set-theoretic operations so let us state what these operations are. For
arbitrary sets $X$ and $Y$ there is a set $X \cup Y$ (their union), a set $X \times Y$ (their product), a set $Y ^X$ (the set of all mappings from $X$ to $Y$)
and the power set $\mathcal{P}(X)$ of $X$ ( the set of all subsets of $X$). The union $X \cup Y$ is the set of all elements which are members of $X$ or $Y$.
The cartesian product  $X \times Y$ is the set of all ordered pairs $(x,y)$ with $x \in X$ and $y \in Y$ and where for all elements $u,\,v$ there exists
an element $(u,v)$ (the ordered pair) such that $(u,v) = (u',v')$ if and only if $u = u'$ and $v = v'$.

Mappings will be defined in terms of their graphs.
Let $f : X \to Y$ be a mapping. Then the 
\definition{graph of $f$} is the subset $\Gamma_f = \{ (x,y) \in X \times Y : \mbox{ $y = f(x)$ for some $x \in X$} \}$ of $\mathcal{P}(X \times Y)$. The set
$\Gamma_f$ has the property that for each $x \in X$ there exists a unique $y \in Y$ with $y = f(x)$ and a set with this property will be called an 
\definition{$X\times Y$-graph}.
Thus if $f : X \to Y$ is a mapping then $\Gamma_f$ is an $X \times Y$-graph. Now this fact can be used as the definition of a mapping by stipulating that
each $X \times Y$-graph $G$ defines a mapping $f : X \to Y$ with of course $\Gamma_f = G$.

We will see that if $A$ and $B$ are finite sets then their union $A \cup B$, their product $A \times B$ and 
$B^A$ (the set of all mappings from $A$ to $B$) are all finite sets.
Moreover, the power set $\mathcal{P}(A)$ is finite and (Proposition~\ref{prop_intro_1}) any subset of a finite set is finite.

Now, although mappings are identified with their graphs, the present definition of being finite allows us to be much more restrictive about
defining mappings between finite sets.  We assume that the following statements are  valid for the mappings to be considered here: 
(The sets occurring below are all finite.)

(1)\enskip If $f : A \to B$ is a mapping then to  each $a \in A$ there is associated a unique element $f(a)$ of $B$ (the value of $f$ at $a$).
In particular, this implies that the set $\Gamma_f = \{ (a,b) \in A \times B : \mbox{$b = f(a)$ for some $a \in A$} \}$ is an $A \times B$-graph.
It also implies that if $A$ is non-empty then there can be no mapping $f : A \to \varnothing$.

(2)\enskip Mappings are determined by their values, meaning that if $f,\,g : A \to B$ are mappings with $f(a) = g(a)$ for all $a \in A$ then $f = g$.
Equivalently, if $\Gamma_f = \Gamma_g$ then $f = g$. 

We assume that if $f : A \to B$ and $g : A \to C$ with $f(a) = g(a)$ for all  $a \in A$ then $f = g$, i.e., we do not insist that the codomains have to be
equal for the mappings to be equal.

(3)\enskip 
For each set $A$ there is the identity mapping $\id_A: A\to A$ with $\id_A(a) = a$ for all $a \in A$; these mappings are bijections.
Note that $\id_\varnothing : \varnothing \to \varnothing$ is the unique mapping $f : \varnothing \to \varnothing$ (since mappings are determined by their
values).
 
(4)\enskip If $f : A \to B$ and $g : B \to C$ are two mappings then there is a mapping $g \circ f : A \to C$ (their composition) satisfying $(g\circ f)(a) = g(f(a))$ 
for all $a \in A$.

(5)\enskip Mappings can be defined by explicitly giving their values. For example:

(5.1)\enskip If $a$ is an element and $B$ a non-empty finite set then for each $b \in B$ there exists a constant mapping $h : \{a\} \to B$ with $h(a) = b$.

(5.2)\enskip If $A_1$ and $A_2$ are disjoint and $f_1 : A_1 \to B$ and $f_2 : A_2 \to B$ are mappings then there is a mapping $f : A_1 \cup A_2 \to B$ satisfying 
$f(a_1) = f_1(a_1)$ if $a_1 \in A_1$ and $f(a_2) =  f_2(a_2)$ if $a_2 \in A_2$.
In particular, if $f : A \to B$ is a mapping and $a \notin A$ then $f$ can be extended to a mapping $f' : A \cup \{a\}\to B$
with $f'(a)$ chosen to be any element in $B$.

(5.3)\enskip If $f : A \to B$ is a mapping and $C \subset A$ then there is the restriction mapping $f_{|C} : C \to B$ satisfying $f_{|C}(c) = f(c)$ for all $c \in C$.

(5.4)\enskip If $f : A \to B$ is a mapping and $C$ is a finite set with $f(A) \subset C$ then there is the  extension mapping $f^{|C}: A \to C$ satisfying 
$f^{|C}(a) = f(a)$ 
for all $a \in A$.  (If $f : A \to B$ and $A' \subset A$ then  as usual the set $\{ f(c) : c \in A' \}$ is denoted by $f(A')$.)
In particular, for each finite set $B$ there is a unique mapping $f_\varnothing^B : \varnothing \to B$. This mapping is unique since mappings are determined by their
values.

(5.5) If $f_1 : A_1 \to B_1$ and $f_2 : A_2 \to B_2$ are mappings then there is a mapping $g :A_1 \times B_1 \to A_2 \times B_2$ satisfying
$g((a,b) = (f_1(a),f_2(a))$ for all $(a,b) \in A_1 \times B_1$.

(5.6)\enskip A mapping can be defined by a formula involving a finite number of cases (where finite means an explicit number such as two or three).
For example, if $b,\,c \in E$ with $b \ne c$ then a transposition $\tau : E \to E$ can be defined by 
\[ 
   \tau(x) = \left\{ \begin{array}{cl}
                  c &\ \mbox{if}\ x = b\;,\\
                  b &\ \mbox{if}\ x = c\;,\\
                  x   &\  \mbox{otherwise}\;.\\
\end{array} \right. 
\]    
(6)\enskip A mapping can be defined by modifying a previously defined  mapping in finitely many places.

For example, if $a \notin A$, $A'$ is a non-empty subset of $A$ and 
$f : A \cup \{a\} \to B$ has been defined previously then a new mapping $g : A \to B$ can be defined by
\[ g(c) = \left\{ \begin{array}{cl}
                  f(c) &\ \mbox{if}\ c \in A \setminus A'\;,\\
                  f(a)   &\  \mbox{if}\ c \in A'\;.\\
\end{array} \right. \]

\medskip

We will see later in Proposition~\ref{prop_fs_333} that if $G$ is an $A \times B$-graph then there is a mapping $f : A \to B$ obtained using only the above 
statements such that $\Gamma_f = G$. 
However, we only make use of Proposition~\ref{prop_fs_333} in the proof of Proposition~\ref{prop_fs_77}, and even then it is not really necessary. In all other cases it is much easier to apply the above 
statements directly.

\bigskip

We start by looking at a fundamental property which depends crucially on the set involved being finite. One reason for presenting
the result at this point is to convince the reader that
the definition of being finite employed here leads to rather straightforward proofs.

\begin{theorem}\label{theorem_fs_1}
Let $A$ be finite and $f : A \to A$ be a mapping. Then  $f$ is injective if and only if it is surjective (and thus if and only if it is bijective). Therefore
the mapping $f : A \to A$ is either a bijection or it is neither injective nor surjective.
\end{theorem}

\proof 
We first show that an injective mapping is bijective. Let $\mathcal{S}$ be the set  consisting of those $B \in \mathcal{P}(A)$ having 
the property that every injective mapping $p : B \to B$ is bijective. Then $\varnothing \in \mathcal{S}$, since the 
only mapping $p : \varnothing \to \varnothing$ is bijective. Let $B \in \proper{\mathcal{S}}$ and $a \in A \setminus B$; 
consider an injective mapping $p : B \cup \{a\} \to B \cup \{a\}$. There are two cases:

($\alpha$)\enskip 
$p(B) \subset B$. Then the restriction $p_{|B} : B \to B$ of $p$ to $B$ is injective
and hence bijective, since $B \in \mathcal{S}$. If $p(a) \in B$ then $p(b) = p_{|B}(b) = p(a)$ for 
some $b \in B$, since $p_{|B}$ is surjective, which contradicts the fact that $p$ is injective. Thus $p(a) = a$, and 
it follows that $p$ is bijective.

($\beta$)\enskip 
$p(B) \not\subset B$. In this case there exists $b \in B$ with $p(b) = a$ and, since $p$ is injective, we must have 
$p(c) \in B$ for all $c \in B \setminus \{b\}$ and $p(a) \in B$. This means there is an injective mapping 
$q : B \to B$ defined by letting
\[ q(c) = \left\{ \begin{array}{cl}
                  p(c) &\ \mbox{if}\ c \in B \setminus \{b\}\;,\\
                  p(a)   &\  \mbox{if}\ c = b\;\\
\end{array} \right. \]
and then $q$ is bijective, since $B \in \mathcal{S}$. Therefore $p$ is again bijective.
 
This shows that $B \cup \{a\} \in \mathcal{S}$ and thus that $\mathcal{S}$ is an inductive $A$-system. Hence 
$\mathcal{S} = \mathcal{P}(A)$, since $A$ is finite. In particular, $A \in \mathcal{S}$ and so every injective mapping 
$f : A \to A$ is bijective.

We now show  that a surjective mapping is bijective, and here let $\mathcal{S}$ be the set consisting of those elements $B \in \mathcal{P}(A)$ 
having the property that every surjective mapping $p : B \to B$ is bijective. Then $\varnothing \in \mathcal{S}$, 
again since the only mapping $p : \varnothing \to \varnothing$ is bijective. Let $B \in \proper{\mathcal{S}}$ and 
$a \in A \setminus B$; consider a surjective mapping $p : B \cup \{a\} \to B \cup \{a\}$. Let 
$D = \{ b \in B : p(b) = a \}$; there are three cases:

($\alpha$)\enskip 
$D = \varnothing$. Then $p(a) = a$, since $p$ is surjective, thus the restriction $p_{|B} : B \to B$ of $p$ to 
$B$ is surjective and hence bijective (since $B \in \mathcal{S}$), and this means $p$ is bijective. 

($\beta$)\enskip 
$D \ne \varnothing$ and $p(a) \in B$. Here we can define a surjective mapping $q : B \to B$ by letting
\[ 
   q(c) = \left\{ \begin{array}{cl}
                  p(c) &\ \mbox{if}\ c \in B \setminus D\;,\\
                  p(a)   &\  \mbox{if}\ c \in D\;.\\
\end{array} \right. 
\] 
Thus $q$ is bijective (since $B \in \mathcal{S}$), which implies that $D = \{b\}$ for some $b \in C$ and in 
particular $p$ is also injective. 

($\gamma$)\enskip 
$D \ne \varnothing$ and $p(a) = a$. This is not possible since then $p(B \setminus D) = B$ and so, choosing any 
$b \in D$, the mapping $h : B \to B$ with
\[ 
   q(c) = \left\{ \begin{array}{cl}
                  p(c) &\ \mbox{if}\ c \in B \setminus D\;,\\
                  b   &\  \mbox{if}\ c \in D\\
\end{array} \right. 
\] 
would be surjective but not injective (since there also exists $c \in B \setminus D$ with $p(c) = b$).

This shows that $B \cup \{a\} \in \mathcal{S}$ and thus that $\mathcal{S}$ is an inductive $A$-system. Hence 
$\mathcal{S} = \mathcal{P}(A)$, since $A$ is finite. In particular, $A \in \mathcal{S}$ and so every surjective mapping $f : A \to A$ is bijective.
\eop

Note that Theorem~\ref{theorem_fs_1} implies the set $\Nat$ of natural numbers is infinite (i.e., it is not finite),
since the successor mapping $\mathsf{s} : \Nat \to \Nat$ with $\mathsf{s}(n) = n + 1$ for all $n \in \Nat$ is injective but not surjective.

If $E, \,F$ are any sets then we write $E \approx F$ if there exists a bijective mapping $f : E \to F$. 
The following result is a direct corollary
of Theorem~\ref{theorem_fs_1}:

\begin{theorem}\label{theorem_fs_2}
If $B$ is a subset of a finite set $A$ with $B \approx A$ then $B = A$.
\end{theorem}

\proof 
There exists a bijective mapping $f : A \to B$ and the restriction $f_{|B} : B \to B$ of $f$ to $B$ is then 
injective; thus by Theorem~\ref{theorem_fs_1} $f_{|B}$ is bijective. But this is only possible if $B = A$, 
since if $a \in A \setminus B$ then $f(a) \notin f_{|B}(B)$.
\eop

The form of the proof of Theorem~\ref{theorem_fs_1} is repeated in practically every proof which follows:
In general we will start with some statement $\prop$ about finite sets, meaning  for each finite set $A$ we have a
statement $\prop(A)$. (For example, $\prop(A)$ could be the statement that any injective mapping $f : A \to A$
is surjective.) The aim is then to establish that $\prop$ is a property of finite sets, i.e., to establish that $\prop(A)$ holds 
for every finite set $A$. To accomplish this we fix a finite set $A$ and consider the set 
$\mathcal{S} = \{ B \in \mathcal{P}(A) : \mbox{$\prop(B)$ holds} \}$
(recalling from Proposition~\ref{prop_intro_1} that each subset of $A$ is finite). 
We then show that $\mathcal{S}$ is an inductive
$A$-system (i.e., show that $\prop(\varnothing)$ holds and  $\prop(B \cup \{a\})$ holds whenever $B \in \proper{\mathcal{S}}$ and
$a \in A \setminus B$) to conclude that $\mathcal{S} = \mathcal{P}(A)$, since $A$ is finite. In particular  $A \in \mathcal{S}$, i.e.,
$\prop(A)$ holds.

This template for proving facts about finite sets can be regarded as a `local' version of the following \definition{induction principle 
for finite sets} which first appeared in a 1909 paper of Zermelo \cite{zermelo}:

\begin{theorem}\label{theorem_fs_3}
Let $\,\prop$ be a statement about finite sets. Suppose
$\,\prop(\varnothing)$ holds and that $\,\prop(A \cup \{a\})$ holds for each element $a \notin A$ whenever $\,\prop(A)$ holds for a finite set $A$.
Then $\,\prop$ is a property of finite sets, i.e., $\prop(A)$ holds for every finite set $A$.
\end{theorem}

\proof
Let $A$ be a finite set and recall from Proposition~\ref{prop_intro_1} that each subset of $A$ is finite. 
Put $\mathcal{S} = \{ B \in \mathcal{P}(A) : \mbox{$\prop(B)$ holds} \}$; then $\mathcal{S}$ is an inductive $A$-system
and hence $\mathcal{S} = \mathcal{P}(A)$. In particular $A \in \mathcal{S}$, i.e., $\prop(A)$ holds. \eop

The above proof of Theorem~\ref{theorem_fs_1} and nearly all the proofs which follow can easily
be converted into proofs based on Theorem~\ref{theorem_fs_3}. 
As an example, we give two proofs of Proposition~\ref{prop_fs_1} below.
Proofs based on Theorem~\ref{theorem_fs_3} seem to be more elegant (although this really a matter of taste).
However, we prefer to continue with the style used in the proof of Theorem~\ref{theorem_fs_1}, since such proofs are internal
to the finite set being considered, and thus appear to be more concrete.
These proofs almost always end with a mantra of the form:

\begin{evlist}{10pt}{6pt}
\item[] \textit{It follows that  $\mathcal{S}$ is an inductive $A$-system. Thus $\mathcal{S} = \mathcal{P}(A)$, since $A$ is finite. In particular, 
$A \in \mathcal{S}$ and so the statement about $A$ holds.}
\end{evlist}

and this will be shortened to the following:

\begin{evlist}{10pt}{6pt}
\item[] \textit{It follows that $\mathcal{S}$ is an inductive $A$-system. Thus $A \in \mathcal{S}$ and so the statement about $A$ holds.}
\end{evlist}

We now establish the usual properties of finite sets mentioned above. The proofs are mostly very straightforward and, since they all
follow the same pattern, they tend to become somewhat monotonous.

\begin{proposition}\label{prop_fs_1}
If $A$ and $B$ are finite sets then so is $A \cup B$.
\end{proposition}

\proof
Consider the set $\mathcal{S} = \{ C \in \mathcal{P}(A) : \mbox{$C \cup B$ is finite} \}$. Then 
$\varnothing \in \mathcal{S}$, since by assumption $\varnothing \cup B = B$ is finite and if $C \in \proper{\mathcal{S}}$
(i.e., $C \cup B$ is finite) and $a \in A \setminus C$ then by Lemma~\ref{lemma_intro_1} 
$(C \cup \{a\}) \cup B = (C \cup B) \cup \{a\} \in \mathcal{S}$. If follows that $\mathcal{S}$ is an inductive $A$-system. 
Thus $A \in \mathcal{S}$ and so $A \cup B$ is finite. 

Here is a proof based on Theorem~\ref{theorem_fs_3}:
Consider the finite set $B$ to be fixed and for each finite set $A$ let $\prop(A)$ be the statement that $A \cup B$ is finite.
Then $\prop(\varnothing)$ holds, since by assumption $\varnothing \cup B = B$ is finite. Moreover, if $\prop(A)$ holds
(i.e., $A \cup B$ is finite) and $a \notin A$ then by Lemma~\ref{lemma_intro_1} 
$(A \cup \{a\}) \cup B = (A \cup B) \cup \{a\}$ is finite, i.e., $\prop(A \cup \{a\})$ holds.
Thus by Theorem~\ref{theorem_fs_3} $A \cup B$ is finite for every finite set $A$.
\eop

\begin{proposition}\label{prop_fs_2}
Let $A$ and $E$ be sets with $A$ finite.

(1)\enskip 
If there exists an injective mapping $f : E \to A$ then $E$ is also finite.

(2)\enskip 
If there exists a surjective mapping $f : A \to E$ then $E$ is again finite.

(3) \enskip If $E$ and $F$ are any sets with $E \approx F$ then $E$ is finite if and only if $F$ is. 
\end{proposition}

\proof
(1)\enskip
Let $\mathcal{S}$ be the set consisting of those elements $C \in \mathcal{P}(A)$ such that if $D$ is any set for which there exists an 
injective mapping $p : D \to C$ then $D$ is finite. Then $\varnothing \in \mathcal{S}$, since there can only 
exist a mapping $p : D \to \varnothing$ if $D = \varnothing$ and the empty set $\varnothing$ is finite. Let 
$C \in \proper{\mathcal{S}}$ and $a \in A \setminus C$. Consider a set $D$ for which there exists an injective mapping 
$p : D \to C \cup \{a\}$. There are two cases:

($\alpha$)\enskip $p(d) \in C$ for all $d \in D$. Here we can consider $p$ as a mapping from $D$ to $C$ and as 
such it is still injective. Thus $D$ is finite since  $C \in \mathcal{S}$.

($\beta$)\enskip 
There exists an element $b \in D$ with $p(b) = a$. Put $D' = D \setminus \{b\}$. Now since $p$ is injective it 
follows that $p(d) \ne a$ for all $d  \in D'$, and thus we can define a mapping $q : D' \to C$ by letting 
$q(d) = f(d)$ for all $d \in D'$. Then $g : D' \to C$ is also injective (since $p : D \to C \cup \{a\}$ is) and 
therefore $D'$ is finite since $C \in \mathcal{S}$. Hence by Lemma~\ref{lemma_intro_1} $D = D' \cup \{b\}$ 
is finite.

This shows that $C \cup \{a\} \in \mathcal{S}$ and therefore $\mathcal{S}$ is an inductive $A$-system. 
Thus $A \in \mathcal{S}$, which means that if there exists an injective mapping $f : E \to A$ then $E$ is also finite.

(2)\enskip
Let $\mathcal{S}$ be the set consisting of those elements $C \in \mathcal{P}(A)$ such that if $D$ is any set for which there exists a 
surjective mapping $p : C \to D$ then $D$ is finite. Then $\varnothing \in \mathcal{S}$, since there can only 
exist a surjective mapping $p : \varnothing \to D$ if $D = \varnothing$ and the empty set $\varnothing$ is finite.
Let $C \in \proper{\mathcal{S}}$ and $a \in A \setminus C$. Consider a set $D$ for which there exists a surjective mapping 
$p : C \cup \{a\} \to D$. There are again two cases:

($\alpha$)\enskip The restriction $p_{|C} : C \to D$ of $p$ to $C$ is still surjective. Then $D$ is finite since 
$C \in \mathcal{S}$.

($\beta$)\enskip The restriction $p_{|C}$ is not surjective. Put $b = p(a)$ and $D' = D \setminus \{b\}$. Then 
$p(c) \ne b$ for all $c \in C$ (since $f_{|C}$ is not surjective) and therefore we can define a mapping 
$q : C \to D'$ by letting $q(c) = p(c)$ for all $c \in C$. But $p : C \cup \{a\} \to D$ is surjective and hence 
$q : C \to D'$ is also surjective. Thus $D'$ is finite since $C \in \mathcal{S}$ holds, and so by 
Lemma~\ref{lemma_intro_1} $D = D' \cup \{b\}$ is finite.

This shows that $C \cup \{a\} \in \mathcal{S}$ and it follows that $\mathcal{S}$ is an inductive $A$-system. 
Thus $A \in \mathcal{S}$, which means that if there exists a surjective  mapping $f : A \to E$ then $E$ is also finite. 

(3)\enskip This is now clear.
\eop

Note that Proposition~\ref{prop_fs_2} (2) is still valid if the set $E$ is replaced by a class $X$. That is, if $A$ is a finite set and $X$ a class and there
exists a surjective mapping $f : A \to X$ then $X$ is a finite set. (Just replace the set $D$ in the above proof by a class $Y$ and the statement $D$ is finite
by the statement $Y$ is a finite set.)

\begin{proposition}\label{prop_fs_3}
If $A$ is a finite set then so is the power set $\mathcal{P}(A)$.
\end{proposition}

\proof 
Let $\mathcal{S}$ be the set consisting of those elements $B \in \mathcal{P}(A)$ for which the power set $\mathcal{P}(B)$ is finite.
Then by Lemma~\ref{lemma_intro_1} $\varnothing \in \mathcal{S}$, since 
$\mathcal{P}(\varnothing) = \{\varnothing\} = \varnothing \cup \{\varnothing\}$. 
Thus consider $B \in \proper{\mathcal{S}}$ and $a \in A \setminus B$.
Then $\mathcal{P}(B \cup \{a\}) = \mathcal{P}(B) \cup \mathcal{P}_a(B)$,
where $\mathcal{P}_a(B) = \{ C \cup \{a\} : C \in \mathcal{P}(B) \}$,
and the mapping $C \mapsto C \cup \{a\}$ from
$\mathcal{P}(B)$ to $\mathcal{P}_a(B)$ is surjective. It follows from Proposition~\ref{prop_fs_2}~(2) that 
$\mathcal{P}_a(B)$ is finite and so by Proposition~\ref{prop_fs_1} $\mathcal{P}(B \cup \{a\})$ is finite, 
i.e., $B \cup \{a\} \in \mathcal{S}$. This shows that $\mathcal{S}$ is an inductive $A$-system. 
Hence $A \in \mathcal{S}$ and so the power set $\mathcal{P}(A)$ is finite. 
\eop

\begin{proposition}\label{prop_fs_4}

If $A$ and $B$ are finite sets then so is their product $A \times B$.
\end{proposition}

\proof
Consider $\mathcal{S} = \{ C \in \mathcal{P}(A) : \mbox{$C \times B$ is finite}\}$. Then $\varnothing \in \mathcal{S}$, 
since $\varnothing \times B = \varnothing$. Let $C \in \proper{\mathcal{S}}$ and $a \in A \setminus C$. Then 
$(C \cup \{a\}) \times B = (C \times B) \cup (\{a\} \times B)$ and by Proposition~\ref{prop_fs_2}~(2) 
$\{a\} \times B$ is finite since the mapping $f : B \to \{a\} \times B$ with $f(b) = (a,b)$ for all $b \in B$ is 
surjective. Thus by Proposition~\ref{prop_fs_1} $(C \cup \{a\}) \times B$ is finite, i.e., 
$C \cup \{a\} \in \mathcal{S}$. This shows that $\mathcal{S}$ is an inductive $A$-system. Hence 
$A \in \mathcal{S}$ and so $A \times B$ is finite. 
\eop

\begin{proposition}\label{prop_fs_5}
If $A$ and $B$ are finite sets then so is $B^A$, the set of all mappings from $A$ to $B$.
\end{proposition}

\proof
Define a mapping $\gamma : B^A \to \mathcal{P}(A \times B)$ by letting 
\[\gamma(f) = (a,b)\in A \times B : \mbox{ $b = f(a)$ for some $a \in A$ } \}\;.\] 
Let $f,\, g \in B^A$ with $\gamma(f) = \gamma(g)$ and let $a \in A$. Then $(a,f(a)) \in \gamma(f)$ and so $(a,f(a)) = (a',g(a')$ for 
some $a' \in A$, since $\gamma(f) = \gamma(g)$. 
Thus $a = a'$ and $f(a) = g(a') = g(a)$ and so $f(a) = g(a)$ for all $a \in A$, i.e.,
$f = g$. This shows that $\gamma$ is injective and hence by Propositions~\ref{prop_fs_1}, 
\ref{prop_fs_2} (1) and \ref{prop_fs_3} $B^A$ is finite. \eop

Note that in this proof the only property of $B^A$ that has been used is that mappings are determined by their values, meaning that if $f,\,g \in B^A$ with
$f(a) = g(a)$ for all $a \in A$ then $f = g$.

Recall that if $A$ and $B$ are finite sets then $G \subset\mathcal{P}(A \times B)$ is  an $A \times B$-graph if for each $a \in A$ there exists a unique
$b \in B$ with $(a,b) \in G$. In particular, if $f : A \to B$ is a mapping then the set 
$\Gamma_f = \{ (a,b) \in A \times B : \mbox{ $b = f(a)$ for some $a \in A$} \}$
is an $A\times B$-graph.

\begin{proposition}\label{prop_fs_333}
Let $A$ and $B$ be finite sets. Then for each $A \times B$-graph $G$  there exists a unique mapping $f : A \to B$ such that 
$G = \Gamma_f$.
\end{proposition}

\proof The uniqueness is clear since mappings are determined by their values . For the existence consider the the subset $\mathcal{S}$ of $\mathcal{P}(A)$
consisting of those $C \subset A$ having the property that for each $C \times B$-graph $G$ there exists a mapping $f : C \to B$ with $\Gamma_f = G$. 
Then $\varnothing \in \mathcal{S}$, since $\varnothing$ is an $\varnothing \times B$-graph, there is a mapping $p : \varnothing \to B$ and $\varnothing = \Gamma_p$.

Thus let $C \in \proper{\mathcal{S}}$, $c \in A \setminus \{c\}$, put $C' = C \cup \{c\}$ and let $G$ be a $C' \times B$-graph.
Then $G' = G \cap \mathcal{P} (C \times B)$ is a $C \times B$-graph and so there exists a mapping $h : C \to B$ with $G' = \Gamma_h$, since $C \in \mathcal{S}$. Extend $h$ to a mapping
$h' : C'= C \cup \{c\} \to B$ by letting $h'(c) = d$, where $d$ is the unique element of $B$ with $(c,b) \in G$. Then $\Gamma_{h'} = G$, which shows that
$C' \in \mathcal{S}$. Thus $\mathcal{S}$ is an inductive $A$-system and hence $A \in \mathcal{S}$. Therefore for each $A \times B$-graph $G$ there exists a mapping $f : A \to B$ such that $G = \Gamma_f$. \eop
\begin{proposition}\label{prop_fs_77} (Existence of a factor mapping).
Let $A$, $B$ and $C$ be finite sets, let  $f : C\to A$ be a surjective mapping and let $g : C\to B$ be a mapping. 
Then there exists a mapping $h : A \to B$ with $g = h\circ f$ if and only if $g(c) = g(c')$ whenever $c,\,c' \in C$
with  $f(c) = f(c')$. Moreover, if $h$ exists then it is unique.
\end{proposition}

\proof Suppose first that there exists $f : A \to B$ with  $g = h\circ f$. If $c,\,c' \in C$ with $f(c) = f(c')$ then $g(c) = h(f(c)) = h(f(c')) = g(c')$
and so $g(c) = g(c')$ whenever $f(c) = f(c')$. Moreover, $h(f(c)) = g(c)$ for each $c \in C$  and $f$ is surjective and hence $h$ is uniquely determined by
$f$ and $g$. 

Now suppose that $g(c) = g(c')$ whenever $c,\,c' \in C$ are such that $f(c) =  f(c')$. Let
\[G = \{(a,b) \in A \times B: \mbox{there exists $c \in C$ with $a = f(c)$ and $b = g(c)$}\}\;.\]
Let $a\in A$; then $a = f(c)$ for some $c \in C$, since $f$ is surjective and then $(a,b) \in G$ with $b = g(c)$. For each $a \in A$ 
there thus exists at least one $b \in B$ with $(a,b) \in G$. But if also $(a,b') \in G$ then there exists $c' \in C$ with $a = f(c')$ and $b' = g(c')$. In
particular $f(c) = f(c')$ and so $b = g(c) = g(c') = b'$, i.e., $b = b'$. Hence for each $a \in A$ there exists a unique $b \in B$ with $(a,b) \in G$, which
shows that $G$ is an $A \times B$-graph. Therefore by Proposition~\ref{prop_fs_333} there exists a unique mapping $h : A \to B$ such that
$G = \Gamma_h$. This means that $G =  \{(a,b) \in A \times B : h(a) = b \}$. 
Let $c \in C$; then $(f(c),h(f(c))) \in G$. But also $(f(c),g(c)) \in G$ and there is a unique $b \in B$ such that $(f(c),b) \in G$. Hence
$g(c) = h(f(c))$ and therefore $g = h \circ f$. \eop

It is not really necessary to use Proposition~\ref{prop_fs_333} here: Suppose that $g(c) = g(c')$ whenever $c,\,c' \in C$ with $f(c) = f(c')$. 
For each $a \in A$ let $G_a = \{c \in C : f(c) = a \}$. Thus $G_a \ne \varnothing$, since $f$ is surjective and if $a \ne a'$ then $G_a \cap G_{a'} = \varnothing$.
Let $\mathcal{E} = \{ E \in \mathcal{P}(C) : \mbox{$E = G_a$ for some $a \in A$}\}$ and define $r : A \to \mathcal{E}$ by 
$r(a) = G_a$ for each $a \in A$. Hence $r$ is a bijection. Now if $c ,\,c'\in r(a)$ then $f(c) = f(c')$ and so $g(c) = g(c')$.
There is thus a mapping $q : \mathcal{E} \to B$ such that $q(r(a)) =g(c)$, where $c$ is any element in $r(a)$ and note that $g(c)$ does not depend on which
element of $r(a)$ is used. Define $h : A \to B$ by $h =q\circ r$ and so $h\circ f = q \circ r \circ f$. 
Let $c \in C$; then $a = f(c) \in A$ and thus $r(a) = G_a \in \mathcal{E}$. Hence $q(G_a) = q(r(a)) = g(c)$, since $c \in r(a)$,
i.e., $(h \circ f)(c) = g(c)$, which shows that $h\circ f = g$. 

\begin{proposition}\label{prop_fs_222}
Let $A$ and $B$ be finite sets and let $f : A \to B$ be a bijection. Then there exists a unique mapping $f^{-1} : B \to A$ such that
$f^{-1} \circ f = \id_A$ and $f \circ f^{-1} = \id_B$. Moreover, $f^{-1}$ is a  bijection.

\end{proposition}
 
\proof 
Let $\mathcal{S}$ denote the set of subsets $C$ of $A$ for which there exists a unique mapping $f^{-1}_{|C} : f(C) \to C$ 
such that $f_{|C} \circ f^{-1}_{|C} = \id_{f(C)} $
and $f^{-1}_{|C} \circ f_{|C)} = \id_C$
 and so in particular $\varnothing \in \mathcal{S}$. Let $C \in  \proper{\mathcal{S}}$ and $a \in A \setminus C$; put
$C' = C \cup \{a\}$. We have the unique mapping $f^{-1}_{|C} :f(C) \to C$ such that $f_{|C} \circ f^{-1}_{|C} = \id_{f(C)} $
and $f^{-1}_{|C} \circ f_{|C)} = \id_C$
and can define $f^{-1}_{|C'} : f(C') \to C'$ by letting $f^{-1}_{|C'}(d) = f^{-1}_{|C}(d)$ if
$d \in f(C)$ and putting $f^{-1}_{|C'}(f(c) = c$. Then $f^{-1}_{|C'} : f(C') \to C'$ is the unique mapping
such that $f_{|C'} \circ f^{-1}_{|C'} = \id_{f(C')} $ and $f^{-1}_{|C'} \circ f_{|C')} = \id_{C'}$
and hence $C \cup \{c\} \in \mathcal{S}$. Thus $\mathcal{S}$ is an inductive $A$-system, and therefore $\mathcal{S} = \mathcal{P(A)}$, since $A$ is finite.
In particular, $A \in \mathcal{S}$. 
This shows that there exists a unique mapping $f^{-1} : B \to A$ such that $f^{-1} \circ f = \id_A$ and $f \circ f^{-1} = \id_B$. 
It is clear that $f^{-1}$ is a  bijection.
\eop

Proposition~\ref{prop_fs_222} also follows from Proposition~\ref{prop_fs_333}. If $f : A \to B$ is a bijection then for each $b \in B$ there exists a 
unique $a\in A$ with $(a,b) \in \Gamma_f$. Therefore the set
$G = \{ (b,a) \in B \times A : (a,b) \in \Gamma_f \}$ is a $B \times A$-graph and thus by Proposition~\ref{prop_fs_333} there exists a unique mapping
$g : B \to A$ with $\Gamma_g = G$ and $g$ is a bijection. Let $a \in A$; then $(a,f(a) \in \Gamma_f$ and so $(f(a),a) \in \Gamma_g$. Thus $g(f(a)) = a$, i.e., 
$g \circ f = \id_A$. Let $b \in B$; then $(b,g(b)) \in \Gamma_g$ and so $(g(b),b) \in \Gamma_f$. Thus $f(g(b)) = b$, i.e., $f \circ g = \id_B$.
Now $g(f(a)) = a$ for all $a \in A$ and $f$ is a bijection and hence $g$ is uniquely determined by $f$.This shows that $g = f^{-1}$. 
\medskip

The next result holds for arbitrary sets $E$ and $F$, the second statement then being the  
Cantor-Bernstein-Schr\"oder theorem. (The first statement only holds in general assuming the axiom of choice.)
As can be seen, the proofs for finite sets are trivial in comparison to those for the general case.

\begin{theorem}\label{theorem_fs_4}
Let $A$ and $B$ be finite sets. Then either there exists an injective mapping $f : A \to B$ or an injective 
mapping $g : B \to A$. Moreover, if there exists both an injective mapping $f : A \to B$ and an injective mapping
$g : B \to A$ then $A \approx B$.
\end{theorem}

\proof
Let $\mathcal{S}$ be the set consisting of those $C \in \mathcal{P}(A)$ for which there either there exists an injective mapping 
$p : C \to B$ or an injective mapping $q : B \to C$. Then $\varnothing \in \mathcal{S}$, since the only mapping 
$p : \varnothing \to B$ is injective. Let $C \in \proper{\mathcal{S}}$ and let $a \in A \setminus C$. There are two cases:

($\alpha$)\enskip 
There exists an injective mapping $q : B \to C$.  Then $q$ is still injective when considered as a mapping from 
$B$ to $C \cup \{a\}$.

($\beta$)\enskip 
There exists an injective mapping $p : C \to B$. If $p$ is not surjective then it can be extended to an injective 
mapping $p' : C \cup \{a\} \to B$ (with $p'(a)$ chosen to be any element in $B \setminus f(C)$). On the other 
hand, if $p$ is surjective (and hence a bijection) then the inverse mapping $p^{-1} : B \to C$ given in Proposition~\ref{prop_fs_222} is injective and 
so is still injective when considered as a mapping from $B$ to $C \cup \{a\}$. 

This shows that $B \cup \{a\} \in \mathcal{S}$ and thus that $\mathcal{S}$ is an inductive $A$-system. Hence 
$A \in \mathcal{S}$ and so there either exists an injective mapping $f : A \to B$ or an injective mapping $g : B \to A$.

Suppose there exists both an injective mapping $f : A \to B$ and an injective mapping $g : B \to A$. Then 
$f \circ g : B \to B$ is an injective mapping, which by Theorem~\ref{theorem_fs_1} is bijective. In particular 
$f$ is surjective and therefore bijective, i.e., $A \approx B$. 
\eop

For sets $E$ and $F$ we write $E \preceq F$ if there exists an injective mapping 
$f : E \to F$ and $A\prec B$ if $A\preceq B$ but $A \not\approx B$.  Theorem~\ref{theorem_fs_4} thus states that
if $A$ and $B$ are finite sets then  exactly one of the three statements $A \prec B$, $ B\prec A$ and $A \approx B$ holds.

The following result collects together some useful technical properties:

\begin{proposition}\label{prop_fs_112}

(1)\enskip Let $A$ and $B$ be finite sets with $B \preceq A$. Then there exists a subset $B'$ of $A$ with $B' \approx B$.

(2)\enskip Let $A$ and $B$ be finite sets. Then there exists either a subset $B'$ of $A$ with $B' \approx B$ or a subset $A'$ 
of $B$ with $A' \approx A$.
Moreover, if $A \not\approx B$ then there exists either a proper subset $B'$ of $A$ with $B' \approx B$ or a proper subset $A'$ 
of $B$ with $A' \approx A$.

(3)\enskip Let $A'$ and $B'$ be finite sets. Then there exist finite sets $A$ and $B$ with $A \approx A'$ and $B \approx B'$ and either $B \subset A$ or 
$A \subset B$.

(4)\enskip  Let $A$ and $E$ be sets with $A$ finite. Then there exists a set $A'$ disjoint from $E$ with $A \approx A'$.

(5)\enskip Let $A_1$ and $A_2$ be finite sets. Then there exist disjoint sets $B_1$ and $B_2$ with $B_1 \approx A_1$ and $B_2 \approx A_2$.
\end{proposition}

\proof (1)\enskip
There exists an injective mapping $g : B \to A$. Put $B' = g(B)$; then $B' \subset A$ with $B' \approx B$ 
(since $g$ as a mapping from $B$ to $B'$ is a bijection). 

(2)\enskip By Theorem~\ref{theorem_fs_4} either $A \preceq B$ or $B \preceq A$.
Thus by Proposition~\ref{prop_fs_112} (1)there exists either a subset $B'$ of $A$ with $B' \approx B$ or a subset $A'$ 
of $B$ with $A' \approx A$. The final statement now follows from Theorem~\ref{theorem_fs_2}.

(3)\enskip This follows directly from (2). 

(4)\enskip
 Let $\mathcal{S}$ be the set of subsets $B \in \mathcal{P}(A)$ for which there exists a set $B'$ disjoint from $E$ with $B \approx B'$, and so
$\varnothing \in \mathcal{S}$. Thus let $B \in \proper{\mathcal{S}}$ and $a \in A \setminus B$. Let $B'$ be disjoint from $E$ with $B \approx B'$.
By Lemma~\ref{lemma_intro_3} there exists an element $b$ not in $E \cup B'$. Then $B' \cup \{b\}$ is disjoint from $E$ and $B \cup \{a\} \approx B' \cup\{b\}$
and hence $B \cup \{a\} \in \mathcal{S}$. This shows that $\mathcal{S}$ is an inductive-$A$-system. Therefore $A \in \mathcal{S}$, i.e., there exists a set $A'$ 
disjoint from $E$ with $A \approx A'$.

(5)\enskip This follows directly from (4). 
\eop

The next result can be seen as a version of the axiom of choice for finite sets.

\begin{proposition}\label{prop_fs_113}
Let $A$ and $A'$ be finite sets and $f : A \to A'$ be a surjective mapping. Then there exists $C \subset A$
such that the restriction $f_{|C} : C \to A'$ is a bijection.
\end{proposition}

\proof
Let $\mathcal{S}$ be the set of subsets $B$ of $A$ such that if $p : B \to B'$ is a surjective mapping then there exists $D \subset B$
such that the restriction $p_{|D} : D \to B'$ is a bijection. Then $\varnothing \in \mathcal{S}$ and so let $B \in \proper{\mathcal{S}}$,
$b \in A \setminus B$ and $p: B \cup \{b\} \to B'$ be a surjective mapping. Then $p_{|B} : B \to B' \setminus \{p(b)\}$ is
surjective and so there exists $D \subset B$ such that $p_{|D} : D \to B' \setminus\{p(b)\}$ is a bijection, since $B \in \mathcal{S}$.
Put $D' = D \cup \{b\}$; then $D' \subset B \cup \{b\}$ and $p_{|D'} D' \to B'$ is a bijection.
Thus $B \cup \{b\} \in \mathcal{S}$ and hence $\mathcal{S}$ is an inductive $A$-system.
This shows that if $f : A \to A'$ is a surjective mapping then there exists $C \subset A$
such that the restriction $f_{|C} : C \to A'$ is a bijection. \eop

\begin{proposition}\label{prop_fs_111}
Let $A$ and $B$ be finite sets. 

(1)\enskip (Existence of a right-inverse).
Let $f : A \to B$ be surjective. Then there exists a mapping $g : B \to A$ such that $f \circ g =\id_B$. (The mapping
$g$ is clearly injective.) 

(2) \enskip (Existence of a left-inverse). Let $f :A \to B$ be injective. Then there exists a mapping $g : B \to A$ such that $g \circ f =\id_A$. (The mapping $g$ is clearly surjective.)
\end{proposition}

\proof (1)\enskip
By Proposition~\ref{prop_fs_113} there exists $C \subset A$ such that $f_{|C} : C \to B$ is a bijection and by Proposition~\ref{prop_fs_222}  there is the
inverse mapping $h : B \to C$ which we can consider as a mapping $g :B \to A$. Clearly $f \circ g = \id_B$. 

(2)\enskip Let We can assume that $A$ is non-empty and so choose an element $a \in A$. Let $f' : A \to f(A)$ be the restriction of $f$ to $f(A)$. Then $f'$ is a bijection and by Proposition~\ref{prop_fs_222}  there is the
inverse mapping $g' : f(A)\to A$. Extend $g'$ to a mapping $g : B \to A$ by letting $g(b) = a$ for all $b \notin f(A)$. 
Clearly $g \circ f = \id_A$. 
\eop

Proposition~\ref{prop_fs_111}~(2) holds for sets in general. However, Proposition~\ref{prop_fs_111}~(1) only holds for sets in general assuming the axiom of choice
 and it is in fact equivalent to the axiom of choice.

The next result is a kind of cancellation law for finite sets.

\begin{proposition}\label{prop_fs_6e} (1) If $A_1,\, A_2$ and $A$ are disjoint finite sets with $A_1 \cup A \approx A_2 \cup A$ then $A_1 \approx A_2$.

(2) If $A_1,\, A_2$ and $A$ are non-empty finite sets with $A_1 \times A \approx A_2 \times A$ then $A_1 \approx A_2$.
\end{proposition}

\proof 
(1) \enskip Suppose $A_1\not\approx A_2$. Then by Theorem~\ref{theorem_fs_4} (and without loss of generality) we can assume that there exists an injective
mapping $h : A_1 \to A_2$ which is not surjective. This mapping $h$ can be extended to a mapping $h' :A_1 \cup A \to A_2 \cup A$ by putting $h'(a) = a$ for all
$a\in A$. Then $h'$ is injective  but not surjective and hence $A_1 \cup A \not\approx A_2 \cup A$.

(2)\enskip  Again suppose $A_1\not\approx A_2$ and as in (1) can assume that there exists an injective
mapping $h : A_1 \to A_2$ which is not surjective. Let $h' :A_1 \times A \to A_2 \times A$ be the mapping given by $h'(a_1,h) = (h(a_1),h)$  for all
$a_1 \in A_1$, $a \in A$. Then $h'$ is injective but not surjective and hence $A_1 \times A \not\approx A_2 \times A$. \eop

\medskip

\begin{proposition}\label{prop_fs_88}
Let $E$ be an infinite set (i.e., $E$ is not finite. Then for each finite set $A$ there exists a finite subset $C$ of $E$ with $C \approx A$.
\end{proposition}

\proof Note that if $C$ is a finite subset of $E$ then $C \ne E$ and $C \cup \{c\}$ is also a finite subset of $E$ for each $c \in E \setminus C$. 
Let $A$ be finite set and let $\mathcal{S}$ denote the set of subsets $B$ of $A$ for which there exists a finite subset $C$ of $E$ with $C \approx B$.
Clearly $\varnothing \in \mathcal{S}$, thus consider $B \in \proper{\mathcal{S}}$, let $a \in A \setminus B$ and put $B' = B \cup \{a\}$. 
Then there exists a finite subset
$C$ of $E$ with $C \approx B$, since $B \in \mathcal{S}$. Also $C \ne B$ and so choose $d \in E \setminus C$. Thus $C' = C \cup\{d\}$ is a finite subset of $E$ with 
$C' \approx B'$, which shows that $B' \in \mathcal{S}$. Hence $\mathcal{S}$ is an inductive $A$-system and so $A \in\mathcal{S}$. \eop
\medskip

Recall that for each set $E$ the set of non-empty subsets of $E$ is denoted by $\mathcal{P}_0(E)$.

\begin{proposition}\label{prop_fs_888}
Let $A$ be a non-empty finite set and let $\approx_A$ be the restriction of the equivalence relation $\approx$ to $\mathcal{P}_0(A)$. Let
$\mathcal{E}(A)$ be the set of equivalence classes.
Then $\mathcal{E}(A) \approx A$.
\end{proposition}

\proof For each $B \in \mathcal{P}_0(A)$ let $\approx_B$ be the restriction of $\approx$ to $\mathcal{P}_0(B)$ and let
$\mathcal{E}(B)$ be the set of equivalence classes. Also for $D \in\mathcal{P}_0(B)$ let $[D]_B$ be the element of $\mathcal{E}(B)$ containing $D$.

Put $\mathcal{S} = \varnothing \cup  \{B \in \mathcal{P}_0(A) : \mathcal{E}(B) \approx B\}$. Consider $B \in \proper{\mathcal{S}}$ with $B \ne \varnothing$,
let $a \in A \setminus B$ and put $C = B \cup \{a\}$. Since $B \in \mathcal{S}$ there exists a bijective mapping $\alpha : B \to \mathcal{E}(B)$ and we
extend $\alpha$ to a mapping $\alpha' : C \to \mathcal{E}(C)$ by letting $\alpha'(a) = [C]_C$ and note that $[C]_C$
consists of the single element $\{C\}$.
Now $C \not\approx D$ for all $D \subset B$ and so $[C]_C \notin \mathcal{E}(B)$. Thus $\alpha'$ is injective. But $\alpha'$ is also surjective:
By definition $\alpha'(a) = [C]_C$ and so consider $k \in \mathcal{E}(C)$ with $k \ne [C]_C$ and let $D \in k$. Then $D$ is a proper subset of $C$
and thus there exists $D' \subset B$ with $D' \approx D$. Since $\alpha$ is surjective there exists $d \in B$ with $\alpha(d) = [D']_B$ It follows that 
$\alpha(d)= [D]_C$. Hence $\alpha'$ is bijective
which implies that $C \in \mathcal{S}$. This shows that
$\mathcal{S}$ is an inductive $A$-system and thus $A \in \mathcal{S}$, i.e., $\mathcal{E}(A) \approx A$. \eop

\begin{theorem}\label{theorem_fs_888}
Let $A$ and $B$ be finite sets. 

(1)\enskip Let $A \prec B$; then there is no surjective mapping $f : A \to B$.

(2)\enskip Let $B \prec A$; then there is no injective mapping $f : A \to B$.

\end{theorem}

\proof (1)\enskip Suppose $f :A \to B$ is surjective; then by Proposition~\ref{prop_fs_111} (1) there exists an injective mapping $g : B \to A$ and so
$B \preceq A$, which by Theorem~\ref{theorem_fs_4} is not the case.

(2)\enskip This is exactly the same, making use of Proposition~\ref{prop_fs_111} (2). \eop

Theorem~\ref {theorem_fs_888} (2) is equivalent to the usual formulation of the Pigeonhole Principle. This principle, introduced, by Dirichlet in 1834 as the
Shubfachprinzip, states that if $n$ objects are placed in $m$ containers and $m < n$ then there is at least one container that contains more than one
element.

Some of the results in the following sections involve partial orders and for their proofs Proposition~\ref{prop_fs_6} below will be needed. 

A \definition{partial order} on a set $E$ is a subset $\le$ of $E \times E$ such that $e \le e$ for all $e \in E$,
$e_1 \le e_2$ and $e_2 \le e_1$ both hold if and  only if $e_1 =  e_2$, and
$e_1 \le e_3$ holds whenever $e_1 \le e_2$ and $e_2 \le e_3$ for some $e_2 \in E$,
and where as usual $e_1 \le e_2$ is written instead of $(e_1,e_2)  \in \le$.  We also write $e_1 < e_2$ if $e_1 \le e_2$ but $e_1 \ne e_2$.

A \definition{partially ordered set} or \definition{poset} is a pair $(E,\le)$ consisting of a set $E$ and a partial order $\le$ on $E$.

If $(E,\le)$ is a poset and $D$ is a non-empty subset of $E$ then $d \in D$ is said to be a 
\definition{maximal} resp.\ \definition{minimal element of $D$}
if $d$ itself is the only element $e \in D$ with $d \le e$ resp.\ with $e \le d$.

\begin{proposition}\label{prop_fs_6}
If $(E,\le)$ is  a poset then every non-empty finite subset of $E$ possesses both a maximal and a minimal element.
\end{proposition}

\proof
Let $A$ be a non-empty finite subset of $E$ and let
$\mathcal{S}$ be the set consisting of the empty set $\varnothing$ together with those non-empty $B \in \mathcal{P}(A)$ which possess a maximal element.
By definition $\varnothing \in \mathcal{S}$. Let $B \in \proper{\mathcal{S}}$ and $a \in A \setminus B$.
We want to show that $B' = B \cup \{a\} \in \mathcal{S}$, and this holds trivially if $B = \varnothing$, since then $a$ is the
only element in $B'$. We can thus suppose that $B \ne \varnothing$, in which case  $B$ has a maximal element $b$.
If $b \le a$ then $a$ is a maximal element of $B'$ (since if $a \le c$ then $b \le c$, thus $b = c$ and so $a = c$).
On the other hand, if $b \le a$ does not hold then $b$ is still a maximal element of $B'$. In both cases $B'$ possesses a maximal
element and hence $B \cup \{a\} = B' \in \mathcal{S}$. This shows $\mathcal{S}$ is an inductive $A$-system and hence 
$A \in \mathcal{S}$, i.e., $A$ possesses a maximal element. Essentially the same  proof also shows that $A$ possesses a minimal element. 
\eop

A partial order $\le$ on $E$ is a \definition{total order} if for all $e_1,\,e_2 \in E$ either $e_1 \le e_2$ or $e_2 \le e_1$ and then $(E,\le)$ is called
a  \definition{totally ordered set}.
If $(E,\le)$ is a totally ordered set and $D$ is a non-empty subset of $E$ then a maximal element $d$ of $D$ is a
\definition{maximum} element, i.e., $e \le d$ for all $e \in D$. Moreover, if a
maximum element exists then it is unique.
In the same way, a minimal element $d$ of $D$ is then a \definition{minimum} element, i.e., $d \le e$ for all $e \in D$,
and if a minimum element exists then it is unique.

If $(E,\le)$ is a totally ordered set and $e_1 \le e_2$ then as usual we denote the set $\{ e \in E: e_1  \le e \le e_2\}$ by $[e_1,e_2]$,
the set $\{ e \in E: e_1  < e <e_2\}$ by $(e_1,e_2)$ and similarly for the sets $(e_1,e_2]$ and $[e_1,e_2)$.

If $(E,\le)$ is a totally ordered set then
Proposition~\ref{prop_fs_6} implies that every non-empty finite subset of $E$ possesses both a unique maximum and a unique minimum element.

\begin{lemma}\label{lemma_fs_7}
For each finite set $A$ there exists a totally ordered set $(E,\le)$ with $E \approx A$.
\end{lemma}

\proof This is clear: Let $A$ be a finite  set for which there exists a totally ordered set $(E,\le)$ with $E \approx A$ and let $a \notin A$. Let $e$
be an element not in $E$. We  can extend $\le$ to  a total order $\le'$ on $E' = E \cup \{e\}$ by defining  $e$ to be the maximum element in
$(E',\le')$.  Then $(E',\le')$ is a totally ordered set with $E' \approx A \cup \{a\}$. 
\eop

We end the section with a result which can be used as a replacement for certain kinds of proofs by induction. Let us start by describing
the type of situation which is involved here.

Suppose we are working with a set-up in which each finite set comes equipped with some additional structure, so  we are dealing with
pairs $(A,\mathcal{T})$, where $\mathcal{T}$ is the structure associated with the finite set $A$. For example, we might be interested
in finite partially ordered sets and in this case $\mathcal{T}$ would be a partial order defined on $A$.
It will be the case that for each pair $(A,\mathcal{T})$ with $A \ne \varnothing$ and for each $a \in A$ there is an induced structure
$\mathcal{T}_a$ on $A \setminus \{a\}$, resulting in a new pair $(A \setminus \{a\},\mathcal{T}_a)$.
Now we would like to show that each such pair $(A,\mathcal{T})$ has a certain property and
a standard approach to tackling this kind of problem is to proceed by induction as follows:
For each $n \in \Nat$ let $\prop(n)$ be the statement that the property holds for all pairs $(A,\mathcal{T})$ with $|A| \le n$
(where $|A|$ is the cardinality of $A$). It is usually clear that $\prop(0)$ holds,
thus take $n \in \Nat \setminus \{0\}$ and assume $\prop(n-1)$ holds.  Then in order to verify that $\prop(n)$ holds it is enough to show
that the property holds for each pair $(A,\mathcal{T})$ with $|A| = n$. Let $(A,\mathcal{T})$ be such a pair;
then for each $a \in A$ the pair $(A \setminus \{a\},\mathcal{T}_a)$ will have the property, since $|A \setminus \{a\}| < n$.
But the crucial step is to choose the element $a \in A$ in a way that allows us to deduce that $(A,\mathcal{T})$ has the property
from the fact that $(A \setminus \{a\},\mathcal{T}_a)$ does, and the correct choice of $a$ will clearly depend very much on the structure $\mathcal{T}$ and
the property involved.
 
Of course, the approach outlined above requires the natural numbers. However, this can be avoided with the help of the result which follows.
As an example, it will applied in Section~\ref{posets} to give a proof of Dilworth's decomposition theorem.

\begin{proposition}\label{prop_fs_7}
Let $A$ be a finite set and $\mathcal{S}$ be a subset of $\mathcal{P}(A)$ containing $\varnothing$. Suppose that
each non-empty subset $F$ of $A$ contains an element $s_F$ such that $F \in \mathcal{S}$ 
whenever $\mathcal{P}(F \setminus \{s_F\}) \subset \mathcal{S}$.
Then $\mathcal{S} = \mathcal{P}(A)$.
\end{proposition}

\proof 
Let $\mathcal{S}_*$ consist of those $B \in \mathcal{P}(A)$ for which $\{ E \in \mathcal{P}(A) : \mbox{$E \preceq B$} \} \subset \mathcal{S}$.
Thus $B \in \mathcal{S}_*$ if and only if $E \in \mathcal{S}$ for each $E \subset A$ with $E \preceq B$.
We will show that  $\mathcal{S}_*$ is an inductive $A$-system. It then follows that $\mathcal{S}_* = \mathcal{P}(A)$ and hence that
$\mathcal{S} = \mathcal{P}(A)$, since $\mathcal{S}_* \subset \mathcal{S}$.

To start with it is clear that $\varnothing \in \mathcal{S}_*$, since $E \preceq \varnothing$ is only possible with $E = \varnothing$
and $\varnothing \in \mathcal{S}$.
Thus let $B \in \proper{\mathcal{S}_*}$ and $a \in A \setminus B$; we want to show that
$B \cup \{a\} \in \mathcal{S}_*$ i.e., 
to show that if $F \in \mathcal{P}(A)$ with $F \preceq B \cup \{a\}$ then $F \in \mathcal{S}$.
If $F \not\approx B \cup \{a\}$ then $F \preceq B$, and in this case $F \in \mathcal{S}$, since $B \in \mathcal{S}_*$.
On the other hand, if $F \approx B \cup \{a\}$ then
$F \setminus \{s_F\} \approx B$ and $B \in \mathcal{S}_*$, and 
so in particular $\mathcal{P}(F \setminus \{s_F\}) \subset \mathcal{S}$.
Hence also $F \in \mathcal{S}$, i.e., $\mathcal{S}$ contains every subset $F$ of $A$ with
$F \preceq B \cup \{a\}$.
Therefore $B \cup \{a\} \in \mathcal{S}_*$, which shows that $\mathcal{S}_*$ is an inductive $A$-system. \eop


\startsection{Permutations}

\label{perm}
Let $E$ be a set; ;then$\,\Self{E}$ will denote the set of all mappings $f : E \to E$ of $E$ into itself, considered as a monoid with
functional composition $\circ$ as monoid operation and $\id_E$ as identity element. 

thus  $\Selfb{E}$ is a submonoid of $\Self{E}$ and it is a group.

If $A$ is finite then the elements of $\Selfb{A}$ are often referred to as \definition{permutations}.

An element $\tau$ of $\Selfb{E}$ will be called an \definition{$E$-transposition}, or just a \definition{transposition} when it is clear which set $E$ is
involved, if there exist $b,\,c \in E$ with $b \ne c$ such that
\[ 
   \tau(x) = \left\{ \begin{array}{cl}
                  c &\ \mbox{if}\ x = b\;,\\
                  b &\ \mbox{if}\ x = c\;,\\
                  x   &\  \mbox{otherwise}\;.\\
\end{array} \right. 
\] 
This transposition will be denoted by $\tau_{b,c}$, or by $\tau^E_{b,c}$ when the set $E$ cannot be determined from the context.
Each transposition is its own inverse.

Denote by $\twogroup$ the multiplicative group $\{+,-\}$ with ${+} \cdot {+} = {-} \cdot {-} = {+}$ and ${-} \cdot {+} = {+} \cdot {-} = {-}$.
For each element $s \in \twogroup$ the other element will be denoted by $-s$.

Let $E$ be a set; we call a mapping $\sigma_E : \Selfb{E} \to \twogroup$ an \definition{$E$-signature} if $\sigma_E(\id_E) = {+}\,$ and
$\sigma_E(\tau \circ f) = -\sigma_E(f)$ for each $f \in \Selfb{E}$ and for each $E$-transposition $\tau$. 
In particular, it then follows that $\sigma_E(\tau) = {-}\,$ for each $E$-transposition $\tau$.

\begin{theorem}\label{theorem_perm_1}
For each finite set $A$ there exists a unique $A$-signature $\sigma : \Selfb{A} \to \twogroup$.
Moreover, $\sigma(f \circ g) = \sigma(f) \cdot \sigma(g)$ for all $f,\,g \in \Selfb{A}$ and hence
$\sigma$ is a group homomorphism.
\end{theorem}

\proof
We first need some preparation and start by noting some of the standard identities concerning the composition of
transpositions which will be needed.

\begin{lemma}\label{lemma_perm_1}
Let $p,\,q,\,r,\,s$ be elements of some set $E$ with $p \ne q$ and $r \ne s$. Then

(1)\enskip $\tau_{r,s} \circ \tau_{p,q} = \tau_{p,q} \circ \tau_{r,s}$ if the elements $p,\,q,\,r,\,s$ are all different.

(2)\enskip $\tau_{r,s} \circ \tau_{p,q} = \tau_{q,s} \circ \tau_{r,s}$ if $p = r$ and $q \ne s$.

(3)\enskip $\tau_{r,s} \circ \tau_{p,q} = \id_E = \tau_{p,q} \circ \tau_{r,s}$ if $p = r$ and $q = s$.
\end{lemma}

\proof 
Just check what happens to the elements $p,\,q,\,r,\,s$. (All other elements in $E$ remain fixed.)
\eop

In what follows let $B$ be a set and $a \notin B$; put $C = B \cup \{a\}$. Denote by $\Selfb{C}^a$ the subgroup of $\Selfb{C}$ consisting
of those mappings $f \in \Selfb{C}$ with $f(a) = a$.
For each $f \in \Selfb{C}^a$ denote 
by $f_B$ the restriction of $f$ to $B$, considered as an element of $\Selfb{B}$. The mapping $\varphi_B : \Selfb{C}^a \to \Selfb{B}$ with
$\varphi_B(f) = f_B$ for each $f \in \Selfb{C}^a$ is clearly a group isomorphism.

Going in the other direction, if $g \in \Selfb{B}$ then denote by $g^a$ the extension of $g$ to $C$ with
$g^a(a) = a$; thus $g^a \in \Selfb{C}^a$ and $(g^a)_B = g$. The mapping $\psi^a : \Selfb{B} \to \Selfb{C}^a$ with $\psi^a(g) = g^a$
for each $g \in \Selfb{B}$ is the inverse of the isomorphism $\varphi_B$.

Note that if $b,\,b' \in B$ with $b \ne b'$ then $\tau^{C}_{b,b'} \in \Selfb{C}^a$ and $\varphi_B(\tau^{C}_{b,b'})  = \tau^B_{b,b'}$.

Consider $f \in \Selfb{C} \setminus \Selfb{C}^a$; then $b = f(a) \in B$ and $\tau^{C}_{a,b} \circ f \in \Selfb{C}^a$. 
Since $\tau^{C}_{a,b}$ is its own inverse we have $f = \tau^{C}_{a,b} \circ (\tau^{C}_{a,b} \circ f)$ which shows that each mapping in
$\Selfb{C} \setminus \Selfb{C}^a$ can be written as the composition of a mapping in $\Selfb{C}^a$ with a $C$-transposition.
The transposition can also be chosen to be on the other side of a mapping from $\Selfb{C}^a$:
There exists a unique element $b' \in B$ with $f(b') = a$, then $f \circ \tau^{C}_{a,b'} \in \Selfb{C}^a$ and
$f = (f \circ \tau^{C}_{a,b'}) \circ \tau^{C}_{a,b'}$.

\begin{lemma}\label{lemma_perm_2}
Let $\sigma_{C} : \Selfb{C} \to \twogroup$ be a $C$-signature and define 
$\sigma_B : \Selfb{B} \to \twogroup$ by $\sigma_B = \sigma_{C} \circ \psi^a$. Then $\sigma_B$ is a $B$-signature. Moreover, 
\begin{evlist}{15pt}{5pt}
\item[]  $\sigma_{C}(f) = \sigma_B(\varphi_B(f))$ for each $f \in \Selfb{C}^a$,

\item[]  
$\sigma_{C}(f) = -\sigma_B(\varphi_B(\tau^{C}_{a,b} \circ f))$ for each $f \in \Selfb{C} \setminus \Selfb{C}^a$,  
where $b = f(a) \in B$ (and as \phantom{xxxxxxxxxxxxxxxxxxxxxxxxxxx}above $\tau^{C}_{a,b} \circ f \in \Selfb{C}^a$).
\end{evlist}
In particular, $\sigma_{C}$ is uniquely determined by $\sigma_B$.
\end{lemma}

\proof
We have $\sigma_B(\id_B) = \sigma_{C}(\psi^a(\id_B)) = \sigma_{C}(\id_{C}) = +$ and
\[
\sigma_B(\tau \circ f) = \sigma_{C}(\psi^a(\tau \circ f)) = \sigma_{C}(\psi^a(\tau) \circ \psi^a(f)) = 
-\sigma_{C}(\psi^a(f)) = -\sigma_B(f)  
\]
for each $f \in \Selfb{B}$ and for each $B$-transposition $\tau$, since $\psi^a(\tau)$ is a $C$-transposition.
This shows that $\sigma_B$ is a $B$-signature. If $f \in \Selfb{C}^a$ then $f = \psi^a(\varphi_B(f))$ and hence
$\sigma_{C}(f) = \sigma_{C}(\psi^a(\varphi_B(f))) = \sigma_B(\varphi_B(f))$.
If $f \in \Selfb{C} \setminus \Selfb{C}^a$ then $f = \tau^{C}_{a,b} \circ (\tau^{C}_{a,b} \circ f)$ and $\tau^{C}_{a,b} \circ f \in \Selfb{C}^a$. 
Thus $\sigma_{C}(f) = -\sigma_{C}(\tau^{C}_{a,b} \circ f) = -\sigma_B(\varphi_B(\tau^{C}_{a,b} \circ f))$. \eop

\begin{lemma}\label{lemma_perm_3}
Let $\sigma_B : \Selfb{B} \to \twogroup$ be a $B$-signature and let $\sigma_{C} : \Selfb{C} \to \twogroup$ be the mapping
given by
\begin{evlist}{15pt}{5pt}
\item[]  $\sigma_{C}(f) = \sigma_B(\varphi_B(f))$ for each $f \in \Selfb{C}^a$,

\item[]  
$\sigma_{C}(f) = -\sigma_B(\varphi_B(\tau^{C}_{a,b} \circ f))$ for each $f \in \Selfb{C} \setminus \Selfb{C}^a$,  
where $b = f(a)$.
\end{evlist}
Then $\sigma_{C}$ is a $C$-signature.
\end{lemma}

\proof
To start with $\sigma_{C}(\id_{C}) = \sigma_B(\varphi_B(\id_{C})) = \sigma_B(\id_B) = +$. 
Thus let $f \in \Selfb{C}$ and let $c,\,d \in C$ with $c \ne d$; we must show that 
$\sigma_{C}(\tau^{C}_{c,d} \circ f) = -\sigma_{C}(f)$.  

Suppose first that both $c$ and $d$ lie in $B$. There are three cases:

(1)\enskip $f \in \Selfb{C}^a$.
In this  case $\tau^{C}_{c,d} \circ f \in \Selfb{C}^a$ and hence
\begin{eqnarray*}
\sigma_{C}(\tau^{C}_{c,d} \circ f) = \sigma_B(\varphi_B(\tau^{C}_{c,d} \circ f)) &=& \sigma_B(\varphi_B(\tau^{C}_{c,d}) \circ \varphi_B(f))\\
&=& \sigma_B(\tau^B_{c,d} \circ \varphi_B(f)) = -\sigma_B(\varphi_B(f)) = -\sigma_{C}(f)\;.
\end{eqnarray*}

(2)\enskip
$f \in \Selfb{C} \setminus \Selfb{C}^a$ with $b =f(a) \notin \{c,d\}$.
In this case we have $\tau^{C}_{c,d} \circ f \in \Selfb{C} \setminus \Selfb{C}^a$ with 
$(\tau^{C}_{c,d} \circ f)(a) = b$ and by Lemma~\ref{lemma_perm_1} $\tau^{C}_{a,b} \circ \tau^{C}_{c,d} = \tau^{C}_{c,d} \circ \tau^{C}_{a,b}$, 
since the elements $a,\,b,\,c,\,d$ are all different. Hence
\begin{eqnarray*}
\sigma_{C}(\tau^{C}_{c,d} \circ f) &=& 
-\sigma_B(\varphi_B(\tau^{C}_{a,b} \circ \tau^{C}_{c,d} \circ f)) = -\sigma_B(\varphi_B(\tau^{C}_{c,d} \circ \tau^{C}_{a,b} \circ f)) \\
&=& -\sigma_B(\varphi_B(\tau^{C}_{c,d}) \circ \varphi_B(\tau^{C}_{a,b} \circ f)) 
= -\sigma_B(\tau^B_{c,d} \circ \varphi_B(\tau^{C}_{a,b} \circ f)) \\
&=& \sigma_B(\varphi_B(\tau^{C}_{a,b} \circ f)) = -\sigma_{C}(f)\;. 
\end{eqnarray*}

(3)\enskip 
$f \in \Selfb{C} \setminus \Selfb{C}^a$ with $b =f(a) \in \{c,d\}$, and
without loss of generality assume $b = d$.
In this case $\tau^{C}_{c,d} \circ f \in \Selfb{C} \setminus \Selfb{C}^a$ with 
$(\tau^{C}_{c,d} \circ f)(a) = c$ and by Lemma~\ref{lemma_perm_1}~(2) $\tau^{C}_{c,d} \circ \tau^{C}_{a,b} = \tau^{C}_{a,c} \circ \tau^{C}_{c,d}$. 
Hence
\begin{eqnarray*}
\sigma_{C}(\tau^{C}_{c,d} \circ f) &=& 
 -\sigma_B(\varphi_B(\tau^{C}_{a,c} \circ \tau^{C}_{c,d} \circ f)) = -\sigma_B(\varphi_B(\tau^{C}_{c,d} \circ \tau^{C}_{a,b} \circ f)) \\
&=& -\sigma_B(\varphi_B(\tau^{C}_{c,d}) \circ \varphi_B(\tau^{C}_{a,b} \circ f)) 
= -\sigma_B(\tau^B_{c,d} \circ \varphi_B(\tau^{C}_{a,b} \circ f)) \\
&=& \sigma_B(\varphi_B(\tau^{C}_{a,b} \circ f)) = -\sigma_{C}(f)\;. 
\end{eqnarray*}

This deals with the cases when both $c$ and $d$ lie in $B$. 
Suppose now then that one of $c$ and $d$ is equal to $a$, and without loss of generality it can be assumed that $c = a$ (and so
$d \in B$). There are the same three cases as above:

(1)\enskip $f \in \Selfb{C}^a$.
Here $\tau^{C}_{c,d} \circ f \in \Selfb{C} \setminus \Selfb{C}^a$ with $(\tau^{C}_{c,d} \circ f)(a) = d$ and thus
\[
\sigma_{C}(\tau^{C}_{c,d} \circ f) = 
-\sigma_B(\varphi_B(\tau^{C}_{a,d} \circ \tau^{C}_{c,d} \circ f))  = -\sigma_B(\varphi_B(f)) 
= -\sigma_{C}(f)\;,
\] 
since $\tau^{C}_{a,d} \circ \tau^{C}_{c,d} = \tau^{C}_{c,d} \circ \tau^{C}_{c,d} = \id_C$.

(2)\enskip $f \in \Selfb{C} \setminus \Selfb{C}^a$ with $b = f(a) \ne d$.
Here $\tau^{C}_{c,d} \circ f \in \Selfb{C} \setminus \Selfb{C}^a$ with 
$(\tau^{C}_{c,d} \circ f)(a) = b$ and by Lemma~\ref{lemma_perm_1}~(2) $\tau^{C}_{a,b} \circ \tau^{C}_{c,d} = \tau^{C}_{b,d} \circ \tau^{C}_{a,b}$. Hence 
\begin{eqnarray*}
\sigma_{C}(\tau^{C}_{c,d} \circ f) &=& 
-\sigma_B(\varphi_B(\tau^{C}_{a,b} \circ \tau^{C}_{c,d} \circ f)) 
 = -\sigma_B(\varphi_B(\tau^{C}_{b,d} \circ \tau^{C}_{a,b} \circ f)) \\
&=& -\sigma_B(\varphi_B(\tau^{C}_{b,d}) \circ \varphi_B(\tau^{C}_{a,b} \circ f)) 
= -\sigma_B(\tau^B_{b,d} \circ \varphi_B(\tau^{C}_{a,b} \circ f)) \\
&=& \sigma_B(\varphi_B(\tau^{C}_{a,b} \circ f)) = -\sigma_{C}(f)\;. 
\end{eqnarray*}

(3)\enskip $f \in \Selfb{C} \setminus \Selfb{C}^a$ with $b = f(a) = d$, and hence $(a,b) = (c,d)$. 
Here $\tau^{C}_{c,d} \circ f \in \Selfb{C}^a$, since $(\tau^{C}_{c,d} \circ f)(a) = a$ and therefore 
\[
\sigma_{C}(\tau^{C}_{c,d} \circ f) = \sigma_B(\varphi_B(\tau^{C}_{c,d} \circ f)) = \sigma_B(\varphi_B(\tau^{C}_{a,b} \circ f)) = -\sigma_{C}(f)\:.
\]

This deals with the cases where one of $c$ and $d$ is equal to $a$, and so all of the possibilities have now been
exhausted. Thus  $\sigma_{C}$ is a $C$-signature. \eop

The first statement in Theorem~\ref{theorem_perm_1} follows directly from Lemmas \ref{lemma_perm_2} and \ref{lemma_perm_3}:
Let $A$ be a finite set and let $\mathcal{S}$ be the set consisting of those $B \in \mathcal{P}(A)$ for which there exists
a unique $B$-signature.
Then $\varnothing \in \mathcal{S}$, since $\Selfb{\varnothing} = \{\id_\varnothing\}$
(and there are no $\varnothing$-transpositions).
Now let $B \in \proper{\mathcal{S}}$, let $a \in A \setminus B$ and put $C = B \cup \{a\}$. Then by Lemma~\ref{lemma_perm_3} there exists a
$C$-signature $\sigma_{C}$ which is the unique $C$-signature, since by Lemma~\ref{lemma_perm_2} 
it is uniquely determined by the unique $B$-signature $\sigma_B$. Thus $B \cup \{a\} \in \mathcal{S}$.
Hence $\mathcal{S}$ is an inductive $A$-system and so $A \in \mathcal{S}$. This shows that for each finite set $A$ 
there is a unique $A$-signature $\sigma : \Selfb{A} \to \twogroup$.

In order to show that the second statement in Theorem~\ref{theorem_perm_1} holds (i.e., that the unique $A$-signature is group
homomorphism) we need the following fact:

\begin{proposition}\label{prop_perm_1}
For each finite set $A$ the group $\Selfb{A}$ is the least submonoid of $\Self{A}$ containing
the $A$-transpositions.
\end{proposition}

\proof
Let $A$ be a finite set and let $\mathcal{S}$ be the set consisting of those $B \in \mathcal{P}(A)$ for which
$\Selfb{B}$ is the least submonoid of $\Self{B}$ containing
the $B$-transpositions. Then $\varnothing \in \mathcal{S}$, since $\Selfb{\varnothing} = \Self{\varnothing} = \{\id_\varnothing\}$
(and there are no $\varnothing$-transpositions).

Now let $B \in \proper{\mathcal{S}}$ and $a \in A \setminus B$; put $C = B \cup \{a\}$ and consider any
submonoid $M$ of $\Self{C}$ containing the $C$-transpositions.
Then $\Selfb{C}^a \cap M$ is a submonoid of $T_{C}$ containing all $C$-transpositions of the form $\tau^{C}_{b,c}$ with $b,\,c \in B$
and hence $\varphi_B(\Selfb{C}^a \cap M)$ is a submonoid of $\Self{B}$ containing all the $B$-transpositions. Thus
$\Selfb{B} \subset  \varphi_B(\Selfb{C}^a \cap M)$  (since $B \in \mathcal{S}$) and it follows that $\Selfb{C}^a \subset M$.
But we have seen that each element of $\Selfb{C} \setminus \Selfb{C}^a$  
can be written in the form  $\tau \circ f$ with $\tau$ a 
$C$-transposition and $f \in \Selfb{C}^a$ and hence also $\Selfb{C} \setminus \Selfb{C}^a \subset M$.  
This shows $\Selfb{C} = (\Selfb{C} \setminus \Selfb{C}^a) \cup \Selfb{C}^a \subset M$, i.e., that $B \cup \{a\} \in \mathcal{S}$.

Hence $\mathcal{S}$ is an inductive $A$-system and so $A \in \mathcal{S}$. For each finite set $A$ the group $\Selfb{B}$ is thus the 
least submonoid of $\Self{B}$ containing the $B$-transpositions.
\eop

We also need the following standard fact:

\begin{lemma}\label{lemma_perm_4}
Let $(M,\bullet,e)$ be a monoid and let $T \subset M$.
If $Q$ is any subset of $M$ containing $e$ such that $t \bullet q \in Q$ for all $q \in Q$ and
all $t \in T$ then $\langle T \rangle \subset Q$, where $\langle T \rangle$ denotes the least submonoid of $M$ containing $T$.
\end{lemma}

\proof Let $N = \{a \in M : \mbox{$a \bullet q \in Q$ for all $q \in Q$}\}$; then clearly $e \in N$ and if
$a_1,\, a_2 \in N$ then $(a_1 \bullet a_2) \bullet q = a_1 \bullet (a_2 \bullet q) \in B$ for all $q \in Q$, i.e., $a_1 \bullet a_2 \in N$.
Thus $N$ is a submonoid of $M$ and by assumption $T \subset N$; hence
$\langle T \rangle \subset N$. But $N \subset Q$, since $e \in Q$, and therefore $\langle T \rangle \subset B$.
\eop

Let $A$ be a finite set and consider the unique $A$-signature $\sigma_A : \Selfb{A} \to \twogroup$.
Let $T$ be the set of $A$-transpositions and
\[ 
Q = \{ f \in \Selfb{A} :  \mbox{$\sigma_A(f \circ g) = \sigma_A(f) \cdot \sigma_A(g)$ for all $g \in \Selfb{A}$}\}\;;
\]
thus in particular $\id_A \in Q$. If $\tau \in T$ and $f \in Q$ then
\[
  \sigma_A((\tau \circ f) \circ g) = \sigma_A(\tau \circ (f \circ g)) = -\sigma_A(f \circ g) 
  = -\sigma_A(f) \cdot \sigma_A(g) = \sigma_A(\tau \circ f) \cdot \sigma_A(g)
\]
for all $g \in \Selfb{A}$ and therefore $\tau \circ f \in Q$. Hence by Lemma~\ref{lemma_perm_4} $\langle T \rangle \subset Q$.
But by Proposition~\ref{prop_perm_1} $\langle T \rangle = \Selfb{A}$ and so $Q = \Selfb{A}$. This shows that the unique $A$-signature $\sigma$ is a group
homomorphism, which completes the proof of Theorem~\ref{theorem_perm_1}. \eop


\startsection{Binomial coefficients}

\label{binom}

For each finite set $A$ and each $B \subset A$ denote the set $ \{ C \in \mathcal{P}(A) : C\approx B \}$  by  $A \,\Delta\, B$ and so 
$A\,\Delta\, B \in \mathcal{P}(\mathcal{P}(A))$. 
The set $A \,\Delta\, B$ plays the role of a \definition{binomial coefficient}: If $|A| = n$ (with $|A|$ the usual cardinality of the set $A$)  
and $|B|= k$ then $|A\,\Delta\,B | = {n\choose k}$.
In this section we establish results which correspond to some of the usual identities for binomial coefficients.

If $A$, $B$ and $C$ are finite sets then we write $C \approx A \amalg B$ if there exist disjoint sets $A'$ and $B'$ with $A \approx A'$, $B \approx B'$ and
$C\approx A' \cup B'$. (It is clear that whether this is the case does not depend on the choice of $A'$ and $B'$.)

The following theorem corresponds to the identity 
\[{{n+1}\choose{k+1}} = {n\choose{k+1}} + {n\choose k}\]

which is used to generate Pascal's triangle.

\begin{theorem}\label{theorem_binom_1}
Let $A$ be a finite set and $B$ be a proper subset of $A$. Let $a$ be an element not in $A$ and let $b \in A \setminus B$. 
Then 
\[(A \cup \{a\})\,\Delta\,(B \cup \{a\})\approx (A\,\Delta\, B) \amalg (A \,\Delta\, (B\cup \{b\})).\] 
\end{theorem}

\proof
Put $D = (A \cup \{a\})\,\Delta\,(B \cup \{a\})$,  $ D_a = \{ C \in (A \cup \{a\})\,\Delta\,(B \cup \{a\}) : a \in C\}$ 
and $D_b = \{ C \in (A \cup \{a\})\,\Delta\,(B \cup \{a\}) : a \notin C\} $, so $D$ is the disjoint union of $D_a$ and $D_b$.
Let $C \in (A \cup \{a\})\,\Delta\,(B \cup \{a\})$. Then $C \subset A \cup \{a\}$ with $C \approx B \cup \{a\}$. 
If $C \in D_a$ then $C' = C \setminus \{a\} \subset A$
and $C' \approx B$ and in this case $C' \in A\,\Delta\,B$. 
If $C\in D_b$ then $C \subset A$ and $ C  \approx B \cup \{b\}$ and in this case  $C  \in A\,\Delta\,(B \cup\{b\})$.
There is thus a mapping $\alpha : D_a \to A\,\Delta \,B$ given by $\alpha(C) = C \setminus \{a\}$ and a mapping  $\beta : D_b \to A\,\Delta\,(B \cup \{b\})$ given by 
$\beta(C) = C$.
Let $C\in  A\,\Delta\, B$. Then $C \subset A$ with $C \approx B$ and so $C \cup \{a\} \subset A \cup \{a\}$ with $C \cup \{a\} \approx B \cup \{a\}$. 
Hence  $C \cup \{a\}\in D_a$ and therefore there is a mapping 
$\alpha' : A \,\Delta\,(B \to D_a$ given by $\alpha'(C) = C \cup \{a\}$.
Now let $C \in A \,\Delta\, (B\cup \{b\})$. 
Then $C \subset A$ with  $C \approx (B\cup \{b\})$ and so $C \cup \{a\} \subset A \cup \{a\}$ with $C \approx (B\cup \{a\})$. Hence  
$C \in D_b$ and therefore there 
is a mapping $\beta' : A \,\Delta\,(B\cup \{b\} \to D_b$ given by $\beta'(C) = C$.

It follows immediately that $\alpha'$ is the inverse of $\alpha$ and that $\beta'$ is the inverse of $\beta$.
Thus the mappings $\alpha$ and $\beta$ are both bijections.
Therefore $D_a \approx A\,\Delta\,B$ and $D_b \approx A\,\Delta \,B \cup \{b\}$ which shows that

\[(A \cup \{a\})\,\Delta\,(B \cup \{a\}) = D = D_a \cup D_b \approx (A\,\Delta\, B) \amalg (A \,\Delta\, (B\cup \{b\}))\;.\ \eop\]

If $C$ is a finite set then, as in Section~\ref{perm}, let $\Selfb{C}$ denote the set(group) of bijections $h :C \to C$. If $B$ is a subset of a finite set $A$
then let $I_{B,A}$ denote the set of injective mappings $k : B \to A$. Note that $I_{A,A} = \Selfb{A}$, $I_{\varnothing,A} =\{\varnothing\}$ and $I_{\{a\},A} \approx A$
for each $a \in A$. 

\begin{theorem}\label{theorem_binom_2}
Let $B$ be a subset of a finite set $A$. Then
$I_{B,A} \approx (A\,\Delta\, B) \times \Selfb{B}$.
\end{theorem}

\proof
Let $u : I_{B,A} \to (A\,\Delta\,B$ be the mapping with $u(k) = k(B)$. Then $u$ is surjective and so by Proposition~\ref{prop_fs_111} (1) there exists a mapping 
$v : (A\,\Delta\,B) \to I_{B,A}$ with $u(v)(C) = C$ for all $C \in (A\,\Delta\,B)$ (and $v$ is injective). 
   
Let $k \in I_{B,A}$; then $k(B) \in A\,\Delta\,B$ and therefore  there exists a bijective mapping $t_k : B \to k(B)$. (Note that $t_k$ is not unique unless $B$ is
empty or contains only one element.)
If $k_1,\,k_2 \in I_{B,A}$ with $k_1(B) = k_2(B)$ and $h = (t_{k_2})^{-1} \circ t_{k_1}$ then $h \in \Selfb{B}$ and
$k_2 \circ h = k_1$. 
On the other hand, if there exists $h \in \Selfb{B}$ with $k_2 \circ h = k_1$ then  $k_1(B) = k_2(B)$.
Therefore $k_1(B) = k_2(B)$ if and only if there exists $h\in \Selfb{B}$ such that $k_2 \circ h = k_1$. Note that if there exists if $h\in \Selfb{B}$ such that 
$k_2 \circ h = k_1$ then $h$ is unique.

Now if $C \in A\,\Delta\,B$ then $v(C)$ is the unique element in $I_{B,A}$ with $v(C)(B) = C$. In particular, it follows for each $k \in I_{B,A}$ that
$v(k(B))$ is the unique element of $I_{B,A}$ with $v(k(B))(B) = k(B)$.
For each $k \in I_{B,A}$ there
thus exists a unique element $s_k \in \Selfb{B}$ such that $k = v(k(B)) \circ s_k$. Define $G : I_{B,A} \to (A\,\Delta\, B) \times \Selfb{B}$ by letting
$G(k) = (k(B), s_k)$. Let $(C,h) \in  (A\,\Delta\, B) \times \Selfb{B}$ and put $k= v(C) \circ h$. Then $k \in I_{B,A}$ with 
$k(B) = (v(C) \circ h)(B) = (v(C)(B) = C$ and
$s_k =h$, since $s_k$ is uniquely determined by the requirement that $k = v(k(B)) \circ s_k$. This shows that $G$ is surjective.
Now let $j,\,k \in I_{B,A}$ with $G(j) = G(k)$. Then $j(B) = k(B)$ and $s_j = s_k$ and hence $j = v(j(B)) \circ s_j = v(k(B)) \circ s_k = k$.This shows that $G$ is
injective and therefore $G$ is a bijection, i.e., $I_{B,A} \approx (A\,\Delta\, B) \times \Selfb{B}$. \eop

\begin{proposition}\label{prop_binom_2}
Let $B$ be a non-empty subset of a finite set $A$, let $b \in B$ and  put $B' = B \cup \{b\}$. Then
$I_{B'A} \approx I_{B,A} \times (A \setminus B)$.
\end{proposition}

\proof Let $r:  I_{B',A} \to I_{B,A}$ be the restriction mapping. Then $r$ is surjective, since $B \ne A$. For each $j \in I_{B,A}$ 

let $p(j) = \{k \in I_{B',A} : r(k) = j \}$. 
Now $j(B) \approx B$ and for each $c \in A \setminus j(B)$ there is a unique $k \in p(j)$ with $k(b) = c$. Hence $p(j) \approx A \setminus B$ and so let
$ s_j : p(j) \to A\setminus B$ be a bijection.  
Define $t : I_{B',A} \to (A \setminus B)$ by $t(k)$ = $s_j(k(b))$, where $j = r(k)$.   
Now define $R: I_{B'A} \to I_{B,A} \times (A \setminus B)$ by $R(k) = (r(k),t(k))$ for each $k \in I_{B',A}$. Let $k_1,\,k_2 \in I_{B',A}$ with $R(k_1) = R(k_2)$.
Then $r(k_1) = r(k_2)$ and $t(k_1) = t(k_2)$. Thus $p(r(k_1)) = p(r(k_2))$ and $s_j((k_1(b)) = s_j(k_2(b))$, where $j = r(k_1) = r(k_2)$. Hence $k_1(b) = k_2(b)$
and so $k_1 = k_2$, i.e., $R$ is injective. Now let $(j,C) \in I_{B,A} \times (A \setminus B)$ and let $k = (s_j)^{-1}(C)$. 
It follows that $R(k) = (r(k),t(k)) = (j,C)$ and so $R$ is surjective. Therefore $R$ is a bijection and hence $I_{B'A} \approx I_{B,A} \times (A \setminus B)$. \eop

The choice of the bijections $s_j$, $j \in I_{B,A}$ in the above proof can be made more explicit with help of Proposition~\ref{prop_fs_111}: 
Let $\Delta$ be the set of all subsets $C$  of $I_{B',A}$ with $C \approx A \setminus B$ and let $\Lambda$ be the set of all bijections $q : C \to A \setminus B$ with 
$C \in \Delta$. Define $u : \Lambda \to \Delta$ by letting $u(q)$ be the domain of $q$. Then $u$ is surjective and so by Proposition~\ref{prop_fs_111} (1)there exists a mapping
$v :\Delta \to \Lambda$ with $u \circ  v = \id_\Delta$. For each $j \in I_{B,A}$ put $s_j =v(p(j))$. 

\begin{proposition}\label{prop_binom_1}
Let $A$ be a finite set, let $a\notin A$ and  put $A' = A \cup \{a\}$. Then
$\Selfb{A'} \approx \Selfb{A} \times A'$.
\end{proposition}

\proof Let $p \in \Selfb{A'}$ and suppose that $p(a)\ne a$. Then there exists a unique element $c_p\in A$ with $p(c_p) = a$. Define $\lambda[p] :A \to A$ by letting 
$\lambda[p](d) = p(d)$ if $p(d) \in A$ and $\lambda[p](c_p) = p(a)$, and thus $\lambda[p] \in \Selfb{A}$.
If $p(a) = a$ then let $\lambda[p]$ be the restriction of $p$ to $A$. Define $\Lambda : \Selfb{A'} \to \Selfb{A} \times A'$ by putting  
$\Lambda(p) = (\lambda[p],p(a))$ for each $p \in \Selfb{A'}$.

Let $p,\, q \in \Selfb{A'}$ with $\Lambda(p)) = \Lambda(q)$. Then $\lambda[p] = \lambda[q]$ and $p(a) = q(a)$. Assume first that $b = p(a) \ne a$ and put
$r = \lambda[p]$. Then $r(c_p)$ and $r(c_q)$ are both equal to $b$ and hence $c_p = c_q$ since $r$ is a bijection. But if $d \ne c_p$ then $r(d) =p(d) = q(d)$
and it follows that $p = q$. If $p(a) = a$ then$ \lambda[p]$ is the restriction of $p$ to $A$ and $ \lambda[q]$ is the restriction of $q$ to $A$ and it again
follows that $p = q$. This shows that $\Lambda$ is injective.

Now let $q\in \Selfb{A}$ and $b\in A$; then there exists a unique $b_q \in A$ with $q(b_q) = b$.
Define $\omega[q,b]: A' \to A'$ by $\omega[q,b](c) =  q(c)$ if $c \in A \setminus \{b_q\}$, $\omega[q,b](b_q) = a$ and $\omega[q,b](a) = b$, so 
$\omega[q,b] \in \Selfb{A'}$. Then 
\[(\Lambda(\omega[q,b]) = (\lambda[\omega[q,b]],\omega[q,b](a)) =(\lambda[\omega[q,b]],b) = (\lambda[r],b) = (s,b)\,,\]
where $r= \omega[q,b]$ and $s = \lambda[r]$. Thus $r(c) = q(c)$ if $c \in A \setminus\{b_q\}$, $r(b_q) = a$ and $r(a) = b$ (and where $q(b_q) = b$). 
Also $s(d) = r(d)$ if $s(d) \in A \setminus \{c_r\}$ and $s(c_r) = r(a)$ (and  where $r(c_r) = a$). Since $r(c_a) = r(b_q) = a$ and $r$ is a bijection it follows 
that $b_q = c_r$. Therefore $s(d) = r(d) = q(d)$ for all $d \in A \setminus \{c_r\}$ and $s(c_r) = r(a) = b = q(b_q) = q(c_r)$ and hence $s = q$, i.e.,
$\Lambda(\omega[q,b]) = (q,b)$. Moreover, it is clear that $\Lambda(q') = (q,a)$ for each $q \in \Selfb{A}$, where $q'$ is the extension of $q$ to $\Selfb{A'}$
with $q'(a) = a$. This shows that $\Lambda$ is surjective and hence it is a bijection. In particular, $\Selfb{A'} \approx \Selfb{A} \times A'$. \eop

\begin{lemma}\label{lemma_binom_1}
Let $B$ be a non-empty subset of a finite set $A$ , let $ b\in B$ and  put $B' = B \cup \{b\}$. 
Then $I_{B',A} \times \Selfb{A\setminus B'}\approx I_{B,A} \times \Selfb{A\setminus B}$.
\end{lemma}

\proof 
By Proposition~\ref{prop_binom_2}  $I_{B',A} \approx I_{B,A} \times (A \setminus B$) and
by Proposition~\ref{prop_binom_1} it follows that  $\Selfb{A\setminus B} \approx \Selfb{A \setminus B'} \times (A \setminus B)$.
Therefore 
\[I_{B',A} \times \Selfb{A\setminus B'}\approx I_{B,A} \times (A\setminus B)\times \Selfb{A \setminus B'} \approx I_{B,A} \times \Selfb{A \setminus B}\;.\ \eop\]

\begin{theorem}\label{theorem_binom_5}
Let $B$ be a subset of a finite set $A$. Then $\Selfb{A} \approx I_{B,A} \times \Selfb{A \setminus B}$.
\end{theorem}

\proof Let $\mathcal{S}$ denote the set of subsets $B$ of $A$ for which $\Selfb{A} \approx I_{B,A} \times \Selfb{A \setminus B}$.
Then $\varnothing \in \mathcal{S}$, since $I_{\varnothing,A} = \{\varnothing\}$. Let $B\in \proper{\mathcal{S}}$ and $b \in A \setminus B$; put $B' = B \cup \{b\}$.
Then $\Selfb{A} \approx I_{B,A} \times \Selfb{A \setminus B}$ and so by Lemma~\ref{lemma_binom_1} $\Selfb{A} \approx I_{B',A} \times \Selfb{A \setminus B'}$.
Hence $B' \in \mathcal{S}$, which shows that $\mathcal{S}$ is an inductive $A$-system. Therefore $\mathcal{S} = \mathcal{P}(A)$ and in particular 
$A \in \mathcal{S}$, i.e., $\Selfb{A} \approx I_{B,A} \times \Selfb{A \setminus B}$. \eop

If $B$ is a subset of a finite set $A$ then by Theorem~\ref{theorem_binom_5} $\Selfb{A} \approx I_{B,A} \times \Selfb{A \setminus B}$
and so there exists a bijective mapping $h:\Selfb{A} \to I_{B,A} \times \Selfb{A \setminus B}$. However, there does not seem to be a natural candidate for
the mapping $h$.

The following theorem corresponds to the usual expression for binomial coefficients:
\[{n\choose m}  = \frac{n!}{m!\cdot (n-m)!}\] 

\begin{theorem}\label{theorem_binom_6}
Let $B$ be a subset of a finite set $A$. Then 
\[\Selfb{A} \times (A\, \Delta \,B) \approx \Selfb{B} \times \Selfb{A \setminus B}\;.\]
\end{theorem}

\proof By Theorem~\ref{theorem_binom_2} $I_{B,A} \approx (A\,\Delta\, B) \times \Selfb{B}$ 
by Theorem~\ref{theorem_binom_5} $\Selfb{A} \approx I_{B,A} \times \Selfb{A \setminus B}$. Therefore
$\Selfb{A} \approx I_{B,A} \times \Selfb{A \setminus B} \approx (A\,\Delta\,B) \times \Selfb{B} \times \Selfb{A\setminus B}$. \eop

The remark following Theorem~\ref{theorem_binom_5} also applies here:
There does not seem to be a natural candidate for a bijective mapping
$h: \Selfb{A} \times (A\, \Delta \,B) \to \Selfb{B} \times \Selfb{A \setminus B}$.

Theorem~\ref{theorem_binom_7} below corresponds to the following well-known identity for binomial coefficients

\[{n\choose m} {m\choose k} = {n\choose k} {{n- k}\choose {m - k}}\]

which holds for all $0 \le k \le m \le n$.

\begin{theorem}\label{theorem_binom_7}
Let $A$, $B$ and $C$ be finite sets with $C \subset B \subset A$. Then
\[(A \,\Delta\, B) \times( B\, \Delta\, C) \approx (A\,\Delta\,C) \times ((A \setminus C)\,\Delta\,(B \setminus C))\;.  \] 
\end{theorem}

\proof
By Theorem~\ref{theorem_binom_6} we have
\[\Selfb{A} \times \Selfb{B} \times (A\,\Delta\,B) \times (B\,\Delta\,C) \approx \Selfb{B}\times \Selfb{A \setminus B} \times \Selfb{C} \times\Selfb{B\setminus C}\]
\[\mbox{and} \quad\Selfb{A} \times \Selfb{A \setminus C}\times (A\,\Delta\,C) \times ((A \setminus C)\Delta\,(B\setminus C)) \approx 
 \Selfb{C}\times \Selfb{A \setminus C} \times \Selfb{B \setminus C}\times \Selfb{A \setminus B}\]
and therefore by Proposition~\ref{prop_fs_6e} \textit{(cancellation law for finite sets)} 
\[\Selfb{A} \times (A\,\Delta\,B) \times (B\,\Delta\,C) \approx \Selfb{A \setminus B} \times \Selfb{C} \times\Selfb{B\setminus C}
\approx \Selfb{C}\times \Selfb{B \setminus C}\times \Selfb{A \setminus B}\]
\[\mbox{and}\quad\Selfb{A} \times (A\,\Delta\,C) \times ((A \setminus C)\Delta\,(B\setminus C)) \approx 
 \Selfb{C}\times \Selfb{B \setminus C}\times \Selfb{A \setminus B}\]
and hence
\[\Selfb{A} \times (A\,\Delta\,B) \times (B\,\Delta\,C) \approx \Selfb{A}\times (A\,\Delta\,C) \times ((A \setminus C)\,\Delta\,(B \setminus C))\;.  \] 
Again making use of Proposition~\ref{prop_fs_6e} it follows that
\[(A \,\Delta\, B) \times( B\, \Delta\, C) \approx (A\,\Delta\,C) \times ((A \setminus C)\,\Delta\,(B \setminus C))\;. \ \eop \] 

Theorem~\ref{theorem_binom_7} implies there is a bijective mapping 
\[h : (A \,\Delta\, B) \times( B\, \Delta\, C) \to (A\,\Delta\,C) \times ((A \setminus C)\,\Delta\,(B \setminus C))\] 
but once  again there does not seem to be a natural candidate for this mapping. 
It is worth noting that in many text-books 
the following simple combinatorial argument is often used to justify the identity
 \[{n\choose m} {m\choose k} = {n\choose k} {{n- k}\choose {m - k}}\;.\]

The left-hand side is the number of ways of first choosing m objects from a set of n objects and
then choosing from these  m objects a subset of k objects.
But this is the same as first choosing k objects from the set of n objects and then choosing m-k objects from the remaining n-k objects, which is the right-hand
side.

The combinatorial argument would seem to suggest how a bijective mapping $h$ could be defined but it is not clear  how to implement this.


\startsection{Dilworth's decomposition theorem}

\label{posets}

In this section we prove Dilworth's decomposition theorem \cite{dilworth} by modifying a proof due to
Galvin \cite{galvin} to work with the present treatment of finite sets.

It is well-known that Dilworth's theorem can be used to provide straightforward  proofs of further important combinatorial results such as the theorems  
of K\"onig, Menger, K\"onig-Egev\'ary and Hall. (See, for example \cite{reich}.

Recall that the set of non-empty subsets of a set $E$ will be denoted by $\mathcal{P}_0(E)$ and that by a 
\definition{partition} of $E$ we mean  a subset $\mathcal{Q}$ of $\mathcal{P}_0(E)$ such that for each $e \in E$ there exists a unique 
$Q \in \mathcal{Q}$ such that $e \in Q$. In particular, different elements in a partition have to be disjoint. 
The only partition of the empty set $\varnothing$ is the empty set $\varnothing = \mathcal{P}_0(\varnothing)$.
If $A$ is finite then by Propositions \ref{prop_fs_3} and \ref{prop_intro_1} any partition of $A$ is also finite.
To each partition $\mathcal{Q}$ of a set $E$ there is the \definition{evaluation map} $i_\mathcal{Q} : E \to \mathcal{Q}$,
where $i_\mathcal{Q}(e)$ is the unique element in $\mathcal{Q}$ containing $e$. If $D \subset E$ then 
the restriction of $i_\mathcal{Q}$ to $D$ will be denoted by $i^D_\mathcal{Q}$.

Recall that a \definition{partial order} on a set $E$ is a mapping $\le :E \times E\to \Bool$ such that $e \le e$ for all $e \in E$,
$e_1 \le e_2$ and $e_2 \le e_1$ both hold if and  only if $e_1 =  e_2$, and
$e_1 \le e_3$ holds whenever $e_1 \le e_2$ and $e_2 \le e_3$ for some $e_2 \in E$,
and where as usual $e_1 \le e_2$ is written instead of $\le(e_1,e_2) = \class{T}$. 
A \definition{partially ordered set} (or \definition{poset}) is a pair $(E,\le)$ consisting of a set $E$ and a partial order $\le$ on $E$.
A finite poset $(A,\le)$ is a poset $(A,\le)$ with $A$ a finite set.

If $(E,\le)$ is a poset and $D$ a non-empty subset of $E$ then $d \in D$ is said to be a \definition{maximal element of $D$}
if $d$ itself is the only element $e \in D$ with $d \le e$. By Proposition~\ref{prop_fs_6} every non-empty finite subset of $E$ possesses a maximal element.

Let $(E,\le)$ be a poset.
A subset $C$ of $E$ is called a \definition{chain} if any two elements in $C$ are comparable, i.e., if
$c \le c'$ or $c' \le c$ for all $c,\,c' \in C$. 
If a chain possesses a maximal element then this is unique, and so will be referred
to as \definition{the} maximal element. 

A partition $\mathcal{C}$ of a subset $F$ of $E$ will be called a \definition{chain-partition of $F$} if each
element in $\mathcal{C}$ is a chain.
A subset $D$ of $E$ is called an \definition{antichain} if no two distinct elements in $D$ are comparable, i.e., if
neither $d \le d'$ nor $d' \le d$ holds whenever $d,\,d' \in D$ with $d \ne d'$. If $D \subset F$ then we say that $D$ is an
\definition{antichain in $F$}.

\begin{lemma}\label{lemma_posets_1}
Let $(E,\le)$ be a poset and $F \subset E$.
If $\mathcal{C}$ is a chain-partition of $F$ and $D$ is an antichain in $F$ then $D \preceq \mathcal{C}$.

In particular, if $\mathcal{C}$ is a chain-partition of $E$ and $D$ is an antichain then $D \preceq \mathcal{C}$.
\end{lemma}

\proof Each chain $C \in \mathcal{C}$ can contain at most one element of $D$ and hence the restricted evaluation
mapping $i^D_\mathcal{C} : D \to \mathcal{C}$ is injective.
\eop

The following important result is Dilworth's decomposition theorem \cite{dilworth}. As stated at the beginning of the section, 
the proof presented here is due to Galvin \cite{galvin}.

\begin{theorem}\label{theorem_posets_1}
Let $(A,\le)$ be a finite  poset. 
Then there exists a chain-partition $\mathcal{C}$ of $A$ and an antichain $D$
such that $D \approx \mathcal{C}$.
\end{theorem}

\proof We first need some preparation, and throughout the proof assume that $(A,\le)$ is a finite poset. 

Let us say that a subset $B$ of $A$ is \definition{regular} if there exists a chain-partition $\mathcal{C}$ of $B$ and an antichain 
$D$ in $B$ such that $D \approx \mathcal{C}$. We thus need to show that $A$ itself is regular.

Let $B$ be a regular subset of $A$.
If $\mathcal{C}$ and $\mathcal{C}'$ are chain-partitions of $B$
and $D$ and $D'$ are antichains in $B$ with $D \approx \mathcal{C}$ and $D' \approx \mathcal{C}'$ then 
$D \approx D' \approx \mathcal{C} \approx \mathcal{C}'$, since by Lemma~\ref{lemma_posets_1} 
$D \preceq \mathcal{C}' \approx D'$ and $D' \preceq \mathcal{C} \approx D$
and hence by Theorem~\ref{theorem_fs_4} $D \approx D'$.
We call $\mathcal{C}$ and $\mathcal{C}'$ \definition{minimal chain-partitions of $B$} and
$D$ and $D'$ \definition{maximal antichains in $B$}.
If $\mathcal{C}$ is a minimal chain-partition of $B$ then a chain-partition $\mathcal{C}'$ of $B$ 
is also minimal if and only if $\mathcal{C}' \approx \mathcal{C}$.
In the same way, if $D$ is a maximal antichain in $B$ then an antichain $D'$ in $B$ 
is also maximal if and only if $D' \approx D$.
If $\mathcal{C}$ is any minimal chain-partition of $B$ and $D$ any maximal antichain in $B$ then
the restricted evaluation mapping $i^D_\mathcal{C} : D \to \mathcal{C}$ is a bijection.

\begin{lemma}\label{lemma_posets_2}
Let $E$ be a subset of $A$ such that every subset of $E$ is regular.
Then there exists a maximal antichain $D_*$ in $E$ and for each $d \in D_*$ a minimal chain-partition $\mathcal{C}_d$ of $E$ such 
that the chain in $\mathcal{C}_d$ containing $d$ has $d$ as its maximal element.
\end{lemma}

\proof
This holds trivially if $E = \varnothing$ and so we can assume that $E$ is non-empty.
Denote the (non-empty) set of maximal antichains in $E$ by $\mathcal{D}$ and let $\Delta$ be the union of all the sets in $\mathcal{D}$, i.e.,
$\Delta = \{ e \in E : \mbox{ $e \in D$ for some $D \in \mathcal{D}$} \}$.

Now fix an arbitrary minimal chain-partition $\mathcal{C}$ of $E$.
For each $C \in \mathcal{C}$ the set $\Delta \cap C$ is non-empty (since it contains an element from each set in $\mathcal{D}$), thus let
$m_C$ be the maximal element in $\Delta \cap C$. Finally let $D_*$ be the set consisting of the elements 
$m_C$, $C \in \mathcal{C}$.

We show that $D_*$ is a maximal antichain in $E$:
Let $C_1,\,C_2 \in \mathcal{C}$ with $C_1 \ne C_2$. Then there exist 
$D_1,\,D_2 \in \mathcal{D}$ such that $m_{C_1}$ is the unique element in $C_1 \cap D_1$
and $m_{C_2}$ the unique element in $C_2 \cap D_2$. Let $b_{12}$ be the unique element in $C_1 \cap D_2$
and $b_{21}$ the unique element in $C_2 \cap D_1$. Then $b_{12} \le m_{C_1}$, since 
$b_{12} \in \Delta \cap C_1$ and $m_{C_1}$ is the maximal element of this set, and in the same way
$b_{21} \le m_{C_2}$. Now if $m_{C_1} \le m_{C_2}$ then it would follow that $b_{12} \le m_{C_2}$, which is not possible
since $b_{12}$ and $m_{C_2}$ are distinct elements of the antichain $D_2$.
The same argument shows also that $m_{C_2} \le m_{C_1}$ is not possible, and hence $D_*$ is an antichain.
Clearly $D_* \subset E$ and $D_* \approx \mathcal{C}$, since each chain in $\mathcal{C}$ contains exactly one element of $D_*$.
Therefore $D_*$ is a maximal antichain in $E$. 

Now for each $d \in D_*$ we obtain a new minimal chain-partition $\mathcal{C}_d$ of $E$. Let $C_d$ be the chain in $\mathcal{C}$ containing
$d$ and put $C'_d = \{ c \in C_d : c \le d \}$, thus $C'_d$ is a non-empty chain. 
By assumption the set $E_d = E \setminus C'_d$ is regular, 
and clearly $D_* \setminus \{d\}$ is an antichain in $E_d$. Suppose there exists an antichain
$D'$ in $E_d$ with $D_* \setminus \{d\} \prec D'$, i.e., with $D_* \preceq D'$. Then there is a subset $D''$ of $D'$ with
$D'' \approx D_*$ and $D''$ is antichain in $E_d$ and thus also an antichain in $E$, i.e., $D''$ is a maximal antichain in $E$.
But this is not possible, since any maximal antichain in $E$ intersects $C'_d$.
It follows that $D_* \setminus \{d\}$ is a maximal antichain in $E_d$. 
Let $\mathcal{C}'_d$ be a minimal chain-partition of $E_d$; thus $\mathcal{C}'_d \approx D_* \setminus \{d\}$.
Finally $\mathcal{C}_d = \mathcal{C}'_d \cup \{C'_d\}$ is a chain-partition of $E$
and $\mathcal{C}_d \approx \mathcal{C}'_d \cup \{C'_d\} \approx (D_* \setminus \{d\}) \cup \{d\} = D_*$, i.e., 
$\mathcal{C}_d$ is a maximal chain-partition of $E$. Moreover, $d$ is the maximal element of the chain 
$C'_d$ in $\mathcal{C}_d$. \eop

Let $\mathcal{S}$ be the set of regular subsets of $A$, and
in particular $\varnothing \in \mathcal{S}$, since in this case  the only chain-partition and antichain  are empty. 
We will show that $\mathcal{S} = \mathcal{P}(A)$ by applying Proposition~\ref{prop_fs_7}.
To do this we must show that each non-empty subset $F$ of $A$ contains an element $s_F$ such that $F \in \mathcal{S}$
whenever $\mathcal{P}(F \setminus \{s_F\}) \subset \mathcal{S}$, i.e., such that $F$ is regular
whenever each subset of $F \setminus \{s_F\}$ is regular.

Thus let $F$ be a non-empty subset of $A$ and take $s_F$ to be a maximal element of $F$ (whose existence is guaranteed by
Proposition~\ref{prop_fs_6}). Suppose that each subset of $E = F \setminus \{s_F\}$ is regular.
Then by Lemma~\ref{lemma_posets_2}
there exists a maximal antichain $D_*$ in $E$ and for each $d \in D_*$ a minimal chain-partition $\mathcal{C}_d$ of $E$ such 
that the chain in $\mathcal{C}_d$ containing $d$ has $d$ as its maximal element.

Now if $D_* \cup \{s_F\}$ is an antichain then it is immediate that $F$ is regular, since
$\mathcal{C}_d \cup \{s_F\}$ is a chain-partition of $F$ for any $d \in D_*$ and
$\mathcal{C}_d \cup \{s_F\} \approx D_* \cup \{s_F\}$.
Thus suppose that $D_* \cup \{s_F\}$ is not an antichain, and so there exists $d \in D_*$ with $d \le s_F$ (since $s_F$ is a maximal element
of $F$). Then $C_d \cup \{s_F\}$ is a chain, since $d$ is the maximal element in $C_d$ and $d \le s_F$, which means
$(\mathcal{C}_d \setminus C_d) \cup (C_d \cup \{s_F\})$ is a chain-partition of $F$. But
$(\mathcal{C}_d \setminus C_d) \cup (C_d \cup \{s_F\}) \approx \mathcal{C}_d \approx D_*$ and $D_*$ is also an antichain in $F$. 
It again follows that $F$ is regular.

Proposition~\ref{prop_fs_7} now implies that $\mathcal{S} = \mathcal{P}(A)$.
In particular $A \in \mathcal{S}$ and so $A$ is regular. The proof of
Theorem~\ref{theorem_posets_1} is complete. \eop

We end the section by applying Dilworth's theorem to give a proof of Hall's theorem on the existence of a system of distinct representatives \cite{hall}.

Let $A$ and $B$ be non-empty finite sets with $A \preceq B$ and let $f : A \to \mathcal{P}(B)$ be a mapping. 
Define $f^*:\mathcal{P}(A) \to \mathcal{P}(B)$ by letting $f^*C) = \bigcup_{c \in C} f(c)$ for each $C \subset A$.
A mapping $s: A \to B$ is a called a \definition{system of distinct representatives} if $s$ is injective and $s(a) \in f(a)$ for all $a \in A$. The following is a
famous theorem of Hall \cite{hall}:

\begin{theorem}\label{theorem_posets_2}
A  system of distinct representatives exists if and only if $C \preceq f^*(C)$ for all $C \subset A$.
\end{theorem}

\proof Assume first that there exists a system of distinct representatives $s : A \to B$ and let $C \subset A$. Then the restriction
$s_C : C \to B$ of $s$ to $C$ is injective and $s(C) \subset f^*(C)$. Hence $C \preceq f^*(C)$.

Now for the converse, so $C \preceq f(^*C)$ for each $C \subset A$. We can clearly assume that $A$ and $B$ are disjoint; put $E = A \cup B$. Define a partial order $\le$ on $E$ as follows: If $e_1, e_2\in E$
then $e_1 \le e_2$ if and only if either $e_1 = e_2$ or $e_1 \in A$, $e_2 \in B$ and $e_2 \in f(e_1)$. If $e_1, e_2, e_3 \in E$ then $e_1 < e_2$ and $e_2 < e_3$ cannot
both hold (since $A$ and $B$ are disjoint). Thus $\le$ is transitive, which shows $\le$ is a partial order.
 The chains in $E$ are exactly the singleton sets $\{e\}$, $e \in E$, together with the two-element sets of the form $\{a,b\}$ with $a \in A$ and
$b \in f(a)$. Note that $\{a\} \preceq f^*(\{a\})$ for each $a \in A$ and so $f^*(a) \ne \varnothing$, which implies that $a$  is contained in at least one 
two-element chain. Let $\mathcal{C}$ be a chain partition of $E$ and let $A_0$ be the subset of $A$ consisting of those $a \in A$ which are contained in a two-element chain in $\mathcal{C}$.
For each $a \in A_0$ denote the other element in the chain by $s'(a)$ This defines a mapping $s' :A_0 \to B$, which is injective (since the chains in $\mathcal{C}$
are disjoint) and $s'(a) \in f(a)$ for all $a \in A_0$. In particular $s'$  will be a system of distinct representatives if and only if $A_0 = A$.
Moreover, the chains in $\mathcal{C}$ are exactly the two-element chains $\{ a,s'(a)\}$, $a \in A_0$, together
with the singleton sets $\{a\}$, $a \in A \setminus A_0$, and the singleton sets $\{b\}$, $b \in B \setminus s'(A_0)$. It follows that
$\mathcal{C} \approx A \ \cup (B \setminus s'(A_0)) = E \setminus s'(A_0) \approx E \setminus A_0$ and hence $A_0  = A$ if and only if $\mathcal{C} \approx B$.

Let $D$ be an antichain in $E$ and put $J = D \cap A$, $K = J\cap B$, so $J$ and $K$ are disjoint and $D = J \cup K$. If $j \in J$  and $k \in f(j)$ then $j < k$ and 
so $k \notin D$. Thus $f^*(J) \cap D = \varnothing$ and hence $K \subset B \setminus f^*(J)$.
Now $J \preceq f^*(J)$.  Let $J'$ be a subset of $B$ with $J' \approx J$. Then $J' \preceq f^*(J)$ and so $B \setminus f^*(J) \preceq B \setminus J'$. 
Therefore $K \preceq B \setminus J'$. Hence $D = J \cup K \preceq J \cup (B \setminus J')\approx B $.  This shows that $D \preceq B$ for each antichain $D$ in $E$.
But $B$ itself is clearly an antichain and hence a maximal antichain in $E$. By Dilworth's theorem there thus exists a chain partition $\mathcal{C}$ with
$\mathcal{C} \approx B$ and we have seen above that the mapping $s':A_0 =A \to B$ is then a system of distinct representatives. \eop 


\startsection{Enumerators}
\label{enums}

In this section we give a further characterisation of a set being finite. This
can be seen as having something to do with enumerating the elements in the set. 
Let us begin with a very informal discussion.
Suppose we want to determine whether a given set  $E$ is finite or not. We could do this by marking the elements in $E$ one at a time and seeing if all
the elements can be marked in finitely many steps (whatever that means). At each stage of this process let us take a snapshot of the elements which 
have already been marked. This results in a subset $\mathcal{U}$ of $\mathcal{P}(E)$ whose elements are exactly these
snapshots; each $U \in \mathcal{U}$ is a subset of $E$ specifying the elements of $E$ which have already been marked at some stage in the process.
The empty set $\varnothing$ is in $\mathcal{U}$ because we take a snapshot before marking the first element.

The following definition will be employed to help make the above more precise. First some notation: A 
subset $\mathcal{U}$ of $\mathcal{P}(E)$ is called an \definition{$E$-selector} if $\varnothing \in \mathcal{U}$ and  for each 
$U \in \proper{\mathcal{U}}$ there exists  a unique element $e \in E \setminus U$ such that $U \cup \{e\} \in \mathcal{U}$. 
The set $\mathcal{U}$ of snapshots should thus be an $E$-selector. 
Moreover, if $E$ is finite then in the final snapshot all the elements of $E$ will
have been marked and so $\mathcal{U}$ should contain $E$.

This suggests that a necessary condition for a set $E$ to be finite is that there should exist an $E$-selector containing $E$.
However, something is missing here since
the condition as it stands is more-or-less vacuous since 
if $\mathcal{U}$ is any $E$-selector then $\mathcal{U} \cup \{E\}$ is also
an $E$-selector which contains $E$.
Thus the condition  is just that there should exist an $E$-selector and it is clear that this does not ensure that the set $E$ is finite since, for
example, $\{ [n] : n \in \Nat \}$ is an $\Nat$-selector for the infinite set $\Nat$.

To see what is missing we need another definition.
If $\mathcal{U}$ is an $E$-selector then a subset $\mathcal{V}$ of $\mathcal{U}$ is said to be \definition{invariant} if
$\varnothing \in \mathcal{V}$ and if $V \cup \{e\} \in \mathcal{V}$ for all $V \in \proper{\mathcal{V}}$, where 
$e$ is the unique element of $E \setminus V$ with $V \cup \{e\} \in \mathcal{U}$. 
In other words, a subset $\mathcal{V}$ of $\mathcal{U}$ is invariant if and only if it is itself an $E$-selector.
An $E$-selector  $\mathcal{U}$ is said to be \definition{minimal} if the only invariant subset of $\,\mathcal{U}$ is $\,\mathcal{U}$ itself.
Note that any $E$-selector contains a unique subset which is a minimal $E$-selector and which can be obtained by taking the intersection of all its 
invariant subsets.

Now consider the $E$-selector $\mathcal{U}$ described above whose elements are exactly the snapshots and let $\mathcal{V}$ be any invariant subset of $\mathcal{U}$. 
Then the process of taking the snapshots produces elements of $\mathcal{U}$ which start with the empty set $\varnothing$, which is in $\mathcal{V}$, and then,
given that the current snapshot is an element of $\mathcal{V}$, will produce a new snapshot which is also an element of $\mathcal{V}$. (This follows from the
definition of $\mathcal{V}$ being invariant.) The  process of taking snapshots can therefore only produce elements of $\mathcal{U}$ which lie in
$\mathcal{V}$. But if $E$ is finite then this process produces all the elements of $\mathcal{U}$, and so in this case we must conclude that 
$\mathcal{V} = \mathcal{U}$. This means that $\mathcal{U}$ must be a minimal $E$-selector. This suggests that
a  necessary condition for a set $E$ to be finite is that there should exist a minimal $E$-selector containing $E$ and 
such a minimal $E$-selector containing $E$ will now be called an \definition{$E$-enumerator}.

It turns out this necessary condition is also sufficient: In Theorem~\ref{theorem_enums_21} it will be shown that
an $E$-enumerator exists  if and only if $E$ is finite.

\begin{lemma}\label{lemma_enums_24}
If $A$ is finite then every $A$-selector contains $A$.
\end{lemma}
\proof 
An $A$-selector $\mathcal{U}$ is a non-empty subset of $\mathcal{P}(A)$ and thus by  Proposition~\ref{prop_intro_3} it contains a maximal 
element $U^*$. 
But each $U \in \proper{\mathcal{U}}$ is not maximal, since there exists an element $a \in A \setminus U$ with $U \cup \{a\} \in \mathcal{U}$. 
Hence $U^* = A$, and so $\mathcal{U}$ contains $A$.
\eop

Note that if $A$ is finite then by Lemma~\ref{lemma_enums_24} an $A$-selector is an $A$-enumerator if and only if it is minimal.

\begin{lemma}\label{lemma_enums_25}
If $A$ is finite then there exists an $A$-selector.
\end{lemma}

\proof
Let $\mathcal{S} = \{ B \in \mathcal{P}(A) : \mbox{there exists a $B$-selector} \}$. In particular $\varnothing \in \mathcal{S}$, since 
$\{\varnothing\}$ is a $\varnothing$-selector.
Let $B \in \proper{\mathcal{S}}$ and let $a \in A \setminus B$; put $B' = B \cup \{a\}$. By assumption there exists
a $B$-selector $\mathcal{U}$ which by Lemma~\ref{lemma_enums_24} contains $B$. It follows
that $\mathcal{U} \cup \{B'\}$ is a $B'$-selector and thus  $B' \in \mathcal{S}$.
This shows $\mathcal{S}$ is an inductive $A$-system and hence $\mathcal{S} = \mathcal{P}(A)$, since $A$ is finite. In particular
$A \in \mathcal{S}$ and so there exists an $A$-selector. 
\eop

\begin{theorem}\label{theorem_enums_21}
A set $E$ is finite if and only if there exists an $E$-enumerator.  
\end{theorem}

\proof
Suppose first that there exists an $E$-enumerator $\mathcal{U}$ and let $\mathcal{S}$ be an inductive $E$-system. 
Then $\mathcal{U} \cap \mathcal{S}$ is an invariant subset of $\mathcal{U}$ and therefore
$\mathcal{U} \cap \mathcal{S} = \mathcal{U}$, since $\mathcal{U}$ is minimal.
Thus $\mathcal{U} \subset \mathcal{S}$ and in particular $E \in \mathcal{S}$.
This shows that each inductive $E$-system contains $E$ and hence by Lemma~\ref{lemma_intro_2} $E$ is finite.
Suppose conversely that $E$ is finite. By Lemma~\ref{lemma_enums_24} there exists an $E$-selector and therefore there exists
a minimal $E$-selector which by Lemma~\ref{lemma_intro_2}
contains $E$. This shows that an $E$-enumerator $\mathcal{U}$ exists. \eop

A subset $\mathcal{U}$ of $\mathcal{P}(E)$ is said to be  \definition{totally ordered} if for all $E_1,\,E_2 \in \mathcal{U}$
 either $E_1 \subset E_2$ or $E_2 \subset E_1$. Note  that a set $E$ is finite if and only if there exists a totally ordered  $E$-enumerator.  
This follows from the fact that the $E$-selector obtained in Lemma~\ref{lemma_enums_25} is totally ordered.
In fact
in Theorem~\ref{theorem_enums_1} it will be shown that if $A$ is finite then any $A$-enumerator is automatically totally ordered.

Note that $\mathcal{U}^\varnothing = \{\varnothing\}$ is the single $\varnothing$-enumerator.

\begin{theorem}\label{theorem_enums_1} 
If $A$ is a finite set then an $A$-selector is minimal if and only if it is totally ordered and thus it is an $A$-enumerator if and only if it is totally ordered.
\end{theorem}

The proof requires some preparation.
Throughout the section $A$ always denotes a finite set. 

For each $E$-selector $\mathcal{U}$ let $\esuc{\mathcal{U}} : \proper{\mathcal{U}} \to E$ 
and $\ssuc{\mathcal{U}} : \proper{\mathcal{U}}  \to \mathcal{U} \setminus \{\varnothing\}$ be the mappings with
$\esuc{\mathcal{U}}(U) = e$ and $\ssuc{\mathcal{U}}(U) = U \cup \{e\}$, where $e$ is the unique 
element in $E \setminus U$ such that $U \cup \{e\} \in \mathcal{U}$. 

\begin{lemma}\label{lemma_enums_1} 
Let $\,\mathcal{U}$ be a totally ordered  $A$-selector and let $U, \,U' \in \mathcal{U}$. Then:

(1)\enskip $U'$ is a proper subset of $U$ if and only if $\ssuc{\mathcal{U}}(U') \subset U$.

(2)\enskip If $U \in \proper{\mathcal{U}}$ then $U'$ is a subset of $U$ if and only if it is a proper subset of $\ssuc{\mathcal{U}}(U)$.
\end{lemma}

\proof
(1)\enskip If 
$\ssuc{\mathcal{U}}(U') \subset U$ 
then $U'$ is a proper subset of $U$, since $\esuc{\mathcal{U}}(U') \notin U'$.
Conversely, suppose $U'$ is a proper subset of $U$. Then there is an injective mapping 
$i : \ssuc{\mathcal{U}}(U') \to U$. Thus if $U \subset \ssuc{\mathcal{U}}(U')$ then by Preposition~\ref{prop_fs_4} $U = \ssuc{\mathcal{U}}(U')$ and in particular 
$\ssuc{\mathcal{U}}(U') \subset U$. But either $U \subset \ssuc{\mathcal{U}}(U')$ or $\ssuc{\mathcal{U}}(U') \subset U$, and thus in both cases 
$\ssuc{\mathcal{U}}(U') \subset U$.

(2)\enskip If $U' \subset U$ then $U'$ is a proper subset of $\ssuc{\mathcal{U}}(U)$, since 
$\esuc{\mathcal{U}}(U) \notin U$. Conversely, suppose $U'$ is a proper subset of $\ssuc{\mathcal{U}}(U)$. 
Then there exists an injective mapping $i : U' \to U$. Thus if $U \subset U'$ then $U = U'$ and in particular 
$U' \subset U$. But either $U \subset U'$ or $U' \subset U$, and thus in both cases $U' \subset U$. \eop

\begin{lemma}\label{lemma_enums_2}   
Let $\,\mathcal{U}$ be a  $B$-enumerator, where  $B$ is a non-empty finite set. Then $u_0 = \esuc{\mathcal{U}}(\varnothing)\in U$ for all
$U \in \mathcal{U} \setminus \{\varnothing\}$ and 
$\mathcal{U}_0 = \{ U \setminus \{u_0\} : \mbox{$U \in \mathcal{U} \setminus \{\varnothing\}$} \}$ is a $(B \setminus \{u_0\})$-enumerator.
Moreover, $\,\mathcal{U}$ is totally ordered if and only if $\,\mathcal{U}_0$ is.
\end{lemma}

\proof The set $\mathcal{V} = \{\varnothing\} \cup \{ U \in \mathcal{U} : u_0 \in U \}$ contains $\varnothing$ and if
$U \in \proper{\mathcal{V}}$ then either $U = \varnothing$, 
in which case $u_0 \in \{u_0\} = \ssuc{\mathcal{U}}(U)$, or $U \ne \varnothing$, 
and then $u_0 \in U \subset \ssuc{\mathcal{U}}(U)$. 
Thus $\mathcal{V}$ is an invariant subset of $\mathcal{U}$  and so $\mathcal{V} = \mathcal{U}$, 
since $\mathcal{U}$ is minimal. Hence $u_0 \in U$ for all $U \in \mathcal{U} \setminus \{\varnothing\}$. Now $\esuc{\mathcal{U}}(U) \ne u_0$ whenever 
$U \in \mathcal{U} \setminus \{\varnothing\}$ (since $\esuc{\mathcal{U}}(U) \in B \setminus U$ and $u_0 \in U$), and so
$\mathcal{U}_0$ is a $(B \setminus \{u_0\})$-selector. Moreover, if $\mathcal{V}_0$ is an invariant
subset of $\mathcal{U}_0$
then $\mathcal{V}_0' = \{\varnothing\} \cup \{ U \cup \{u_0\} : U \in \mathcal{V}_0 \}$ is an invariant subset of 
$\mathcal{U}$. Therefore $\mathcal{V}_0' = \mathcal{U}$, since $\mathcal{U}$ is minimal, which implies that 
$\mathcal{V}_0 = \mathcal{U}_0$. This shows that $\mathcal{U}_0$ is a $(B \setminus \{u_0\})$-enumerator.
Finally, it is clear that $\mathcal{U}$ is totally ordered if and only if $\mathcal{U}_0$ is.
\eop

\

\textit{Proof of Theorem~\ref{theorem_enums_1}\ }
Assume first there exists an $A$-enumerator  which is not totally ordered. Then the subset
$\mathcal{S}$ of $\mathcal{P}(A)$ consisting of those subsets $B$ for which there exists a $B$-enumerator  which is not totally
ordered is non-empty. Hence by Proposition~\ref{prop_intro_2} $\mathcal{S}$ contains a minimal element
$B$ and $B$ is non-empty since the only $\varnothing$-enumerator  is trivially totally ordered. There thus exists a
$B$-enumerator  $\mathcal{U}$ which is not totally ordered and then Lemma~\ref{lemma_enums_2} implies that $B \setminus \{u_0\} \in \mathcal{S}$,
where $u_0 = \esuc{\mathcal{U}}(\varnothing)$. This contradicts the minimality of $B$ and therefore the assumption that
there exists an $A$-enumerator  which is not totally ordered is false. In other words, each $A$-enumerator is totally ordered.

For the converse let $\mathcal{U}$  be a totally ordered $A$-selector and 
suppose there exists an invariant proper subset $\mathcal{V}$ of  $\mathcal{U}$.
Put $\mathcal{S} = \mathcal{U} \setminus \mathcal{V}$. Then $\mathcal{S}$
is non-empty and hence by Proposition~\ref{prop_intro_2} it contains a minimal element $U_*$, and $U_* \ne \varnothing$, since $\varnothing \notin \mathcal{S}$.
Put $\mathcal{R} = \{ V \in \mathcal{V} : V \subset U_* \}$; then $\mathcal{R}$ is non-empty (since it contains $\varnothing$) and therefore by
Proposition~\ref{prop_intro_3} $\mathcal{R}$ contains a maximal element $U^*$.
Thus $U^* \subset U_*$, and in fact $U^*$ is a proper subset of $U_*$, since $U^* \in \mathcal{V}$ and $U_* \notin \mathcal{V}$.
Hence by Lemma~\ref{lemma_enums_1}~(1) $\ssuc{\mathcal{U}}(V^*) \subset U_*$. 
But $\ssuc{\mathcal{U}}(U^*) \in \mathcal{V}$, since $\mathcal{V}$ is invariant, and so $\ssuc{\mathcal{U}}(U^*) \in \mathcal{R}$. 
However, this contradicts the maximality of $U^*$ and we conclude that
$\mathcal{V} = \mathcal{U}$. Therefore $\mathcal{U}$ is minimal, i.e., $\mathcal{U}$ 
is an $A$-enumerator.

This completes the proof of Theorem~\ref{theorem_enums_1}.
 \eop

If $\mathcal{U}$ is an $A$-enumerator then for each $U \in \mathcal{U}$ the set $\mathcal{U} \cap \mathcal{P}(U)$
will be denoted by $\mathcal{U}_U$. 
Note that, as far as the definition of $\mathcal{U}_U^p$ is concerned,
$\mathcal{U}_U$ is considered here to be a subset of $\mathcal{P}(U)$ and so
$\mathcal{U}_U^p  =\{ U' \in \mathcal{U} : \mbox{$U'$ is a proper subset of $U$} \}$. 
If $U' \in \mathcal{U}_U^p$ then by Lemma~\ref{lemma_enums_1}~(1) $\ssuc{\mathcal{U}}(U') \in \mathcal{U}_U$ and so 
$\mathcal{U}_U$ is a $U$-selector. But  $\mathcal{U}_U$ is clearly totally ordered 
and therefore by Theorem~\ref{theorem_enums_1}
it is in fact a $U$-enumerator.
If $U \in \proper{\mathcal{U}}$ and $U^* = \ssuc{\mathcal{U}}(U)$ 
then by Lemma~\ref{lemma_enums_1}~(2) $\,\proper{\mathcal{U}_{U^*}} = \mathcal{U}_U$.

\begin{lemma}\label{lemma_enums_4}
If $\,\mathcal{U}$ is an $A$-enumerator then for all $U \in \mathcal{U}$
\[ U = \{ a \in A : \mbox{$a = \esuc{\mathcal{U}}(U')$ for some $U' \in \mathcal{U}_U^p$}\}\;.\]
\end{lemma}

\proof 
Let $\mathcal{V}$ be the set consisting of those elements
$U \in \mathcal{U}$ for which the statement above holds, i.e., for which
$U = \{ a \in A : \mbox{$a = \esuc{\mathcal{U}}(U')$ for some $U' \in \mathcal{U}_U^p$}\}$,
and hence in particular $\varnothing \in \mathcal{V}$.
Let $U \in \proper{\mathcal{V}}$ and put $U_* = \ssuc{\mathcal{U}}(U)$; then by Lemma~\ref{lemma_enums_1}~(2)
\begin{eqnarray*}
 \ssuc{\mathcal{U}}(U) &=& U \cup \{\esuc{\mathcal{U}}(U)\} 
=  \{ a \in A : \mbox{$a = \esuc{\mathcal{U}}(U')$ for some $U' \in \mathcal{U}_U$ } \} \\ 
  &=&  \{ a \in A : \mbox{$a = \esuc{\mathcal{U}}(U')$ for some $U' \in \proper{\mathcal{U}_{U_*}}$ } \}
\end{eqnarray*}
and hence $\ssuc{\mathcal{U}}(U) = U_* \in \mathcal{V}$. Thus $\mathcal{V}$ is an invariant subset of $\mathcal{U}$, and so 
$\mathcal{V} = \mathcal{U}$. \eop

\begin{proposition}\label{prop_enums_1} 
If $\,\mathcal{U}$ is an $A$-enumerator then the mappings
$\ssuc{\mathcal{U}} : \proper{\mathcal{U}} \to \mathcal{U} \setminus \{\varnothing\}$ 
and $\esuc{\mathcal{U}} : \proper{\mathcal{U}}  \to A$ are  both bijections.
In particular, if $\,\mathcal{U}$ and $\mathcal{V}$ are  $A$-enumerators then $\mathcal{U} \approx \mathcal{V}$. 
(This means, somewhat imprecisely, that any $A$-enumerator contains one more element than $A$.)
\end{proposition}

\proof
Let 
$\mathcal{V} = \{\varnothing\} \cup \{ U \in \mathcal{U} :  \mbox{there exists $U' \in \proper{\mathcal{U}}$ such that $U = \ssuc{\mathcal{U}}(U')$} \}$. 
Then $\varnothing \in \mathcal{V}$ and $\ssuc{\mathcal{U}}(U)$ is trivially an element of $\mathcal{V}$ for all 
$U \in \proper{\mathcal{U}}$, and so in particular for all $U \in \proper{\mathcal{V}}$. Therefore $\mathcal{V}$ is an invariant subset of 
$\mathcal{U}$ and so $\mathcal{V} = \mathcal{U}$. This shows that the mapping $\ssuc{\mathcal{U}}$ 
is surjective.
Now consider the mapping $\ssuc{\mathcal{U}}' : \mathcal{U} \to \mathcal{U}$ with $\ssuc{\mathcal{U}}'(A) = \varnothing$ and
$\ssuc{\mathcal{U}}'(U) = \ssuc{\mathcal{U}}(U)$ whenever $U \in \proper{\mathcal{U}}$. Then $\ssuc{\mathcal{U}}'$ is surjective, since 
$\ssuc{\mathcal{U}}$ is, and hence by 
Theorem~\ref{theorem_fs_1} $\ssuc{\mathcal{U}}'$ is bijective, since $\mathcal{U}$ is finite. It follows that $\ssuc{\mathcal{U}}$ is also bijective.

Now to the mapping $\esuc{\mathcal{U}}$.
Let $U,\,U' \in \proper{\mathcal{U}}$ with $U \ne U'$. By Theorem~\ref{theorem_enums_1} $\mathcal{U}$ is totally ordered, thus either 
$U \subset U'$ or $U' \subset U$ and so without loss
of generality assume that $U' \subset U$. Therefore $U'$ is a proper subset of $U$, hence by Lemma~\ref{lemma_enums_1}~(1)
$\ssuc{\mathcal{U}}(U') \subset U$ and in particular $\esuc{\mathcal{U}}(U') \in U$.
But $\esuc{\mathcal{U}}(U) \notin U$, which implies that $\esuc{\mathcal{U}}(U) \ne \esuc{\mathcal{U}}(U')$. This shows the mapping 
$\esuc{\mathcal{U}}$ is injective. 
Moreover, by Lemma~\ref{lemma_enums_4} (with $U = A$)
$A = \{ a \in A : \mbox{$a = \esuc{\mathcal{U}}(U)$ for some $U \in \proper{\mathcal{U}}$} \}$,
and thus the mapping $\esuc{\mathcal{U}}$ is surjective.
\eop

\begin{lemma}\label{lemma_enums_3}  
If $B \subset A$ then there exists an $A$-enumerator $\mathcal{U}$ with $B \in \mathcal{U}$.
\end{lemma}

\proof 
Let $\mathcal{V}$ be a $B$-enumerator and $\mathcal{V}'$ be an $(A \setminus B)$-enumerator. Then by Theorem~\ref{theorem_enums_1}
$\mathcal{U} = \mathcal{V} \cup \{ B \cup C : C \in \mathcal{V}' \setminus \{\varnothing\} \}$ is clearly a totally ordered
$A$-selector containing $B$   and thus by Theorem~\ref{theorem_enums_1} it is an $A$-enumerator containing $B$. 
\eop

In what follows $B$ is always a finite set.

If $\,\mathcal{U}$ is an $A$-enumerator and  $\,\mathcal{V}$ a $B$-enumerator then
a  mapping $\pi : \mathcal{U} \to \mathcal{V}$ is called a \definition{homomorphism} 
if $\pi(\varnothing) = \varnothing$, $\pi(\proper{\mathcal{U}}) \subset \proper{\mathcal{V}}$ and 
$\pi(\ssuc{\mathcal{U}}(U)) = \ssuc{\mathcal{V}}(\pi(U))$ for all $U \in \proper{\mathcal{U}}$. 

\begin{proposition}\label{prop_enums_2}  
If $\pi : \mathcal{U} \to \mathcal{V}$ is a homomorphism from an $A$-enumerator $\mathcal{U}$ to a $B$-enumerator $\mathcal{V}$ 
then $\pi(U) \approx U$ for all $U \in \mathcal{U}$ and $\pi$ maps $\,\mathcal{U}$ bijectively
onto $\,\mathcal{V}_{\pi(A)}$.
\end{proposition}

\proof 
Let $\mathcal{U}_0$ denote the set consisting of those $U \in \mathcal{U}$ for which  $\pi(U') \approx U'$ for all $U' \in \mathcal{U}_U$
and for which $\pi$ maps 
$\mathcal{U}_U$ bijectively onto $\mathcal{V}_{\pi(U)}$. Clearly $\varnothing \in \mathcal{U}_0$ since $\pi(\varnothing) = \varnothing$ and 
$\mathcal{U}_\varnothing = \mathcal{V}_{\pi(\varnothing)} = \{\varnothing\}$.

Consider $U \in \proper{\mathcal{U}_0}$ and put $U^* = \ssuc{\mathcal{U}}(U)$; by Lemma~\ref{lemma_enums_1}~(2) 
$\mathcal{U}_{U^*} = \mathcal{U}_U \cup \{U^*\}$. 
Now $\pi(U^*) = \pi(\ssuc{\mathcal{U}}(U)) = \ssuc{\mathcal{V}}(\pi(U)) = \pi(U) \cup \{c\}$, where $c \notin \pi(U)$,
$U^* = U \cup \{d\}$, where $d \notin U$ and $\pi(U) \approx U$, since $U \in \mathcal{U}_U$. It follows that $\pi(U^*) \approx U^*$, and hence
$\pi(U') \approx U'$ for all $U' \in \mathcal{U}_{U^*}$.
Moreover, $\mathcal{V}_{\pi({U^*})} = \mathcal{V}_{\ssuc{\mathcal{V}}(\pi(U))} = \mathcal{V}_{\pi(U)} \cup \{\ssuc{\mathcal{V}}(\pi(U))\}$,
$\pi(U^*) = \ssuc{\mathcal{V}}(\pi(U))$ and $\pi$ maps $\mathcal{U}_U$ bijectively onto $\mathcal{V}_{\pi(U)}$. It follows that
$\pi$ maps $\mathcal{U}_{U_*}$ bijectively onto $\mathcal{V}_{\pi({U^*})}$. This shows $U^* = \ssuc{\mathcal{U}}(U) \in \mathcal{U}_0$, and so
$\mathcal{U}_0$ is an invariant subset of $\mathcal{U}$. 
Thus  $\mathcal{U}_0 = \mathcal{U}$ and then by Lemma~\ref{lemma_enums_24}
$A \in \mathcal{U}_0$, i.e., $\pi(U) \approx U$ for all $U \in \mathcal{U}$ and $\pi$ maps $\mathcal{U}$ bijectively
onto $\mathcal{V}_{\pi(A)}$. 
\eop

If $\pi : \mathcal{U} \to \mathcal{V}$ is a homomorphism as above then by Proposition~\ref{prop_enums_1}
$\pi(A) \approx A$. But $\pi(A)$ is a subset of $B$, hence $\pi(A) \preceq B$ and thus $A \preceq B$.
This  necessary condition for the existence of a homomorphism is also sufficient:

\begin{proposition}\label{prop_enums_3} 
If $A \preceq B$, $\,\mathcal{U}$ is an $A$-enumerator and  $\mathcal{V}$ a $B$-enumerator
then there exists a unique homomorphism $\pi : \mathcal{U} \to \mathcal{V}$.
\end{proposition}

\proof
Let $\mathcal{U}_0$ denote the set consisting of those $U \in \mathcal{U}$ for which there exists a homomorphism $\pi_U : \mathcal{U}_U \to \mathcal{V}$. 
Clearly $\varnothing \in \mathcal{U}_0$ since $\mathcal{U}_\varnothing = \{\varnothing\}$ and $\proper{\mathcal{U}_\varnothing} = \varnothing$. 

Consider $U \in \proper{\mathcal{U}_0}$ and let $\pi_U : \mathcal{U}_U \to \mathcal{V}$ be a homomorphism.
Now $A \preceq B$  and $U$ is a proper subset of $A$ and hence $U \not\approx B$; it follows that 
$\pi_U(U) \ne B$, since by Proposition~\ref{prop_enums_1} $\pi_U(U) \approx U$. 
Let $U^* = \ssuc{\mathcal{U}}(U)$; by Lemma~\ref{lemma_enums_1}~(2) $\mathcal{U}_{U^*} = \mathcal{U}_U \cup \{U^*\}$ and so
we can define $\pi_{U^*} : \mathcal{U}_{U^*} \to \mathcal{V}$  by putting $\pi_{U^*}(U') = \pi_U(U')$ if $U' \in \mathcal{U}_U$ and letting
$\pi_{U^*}(U^*) = \ssuc{\mathcal{V}}(\pi_U(U))$ (recalling that $\pi_U(U) \ne B$). 
If $U' \in \proper{\mathcal{U}_U}$ then $U' \in \mathcal{U}_U$ and $\ssuc{\mathcal{U}}(U') \in \mathcal{U}_U$ and thus
\[\pi_{U^*}(\ssuc{\mathcal{U}}(U')) = \pi_U(\ssuc{\mathcal{U}}(U')) = \ssuc{\mathcal{V}}(\pi_U(U')) = \ssuc{\mathcal{V}}(\pi_{U^*}(U'))\;.\]
Also $\pi_{U^*}(\ssuc{\mathcal{U}}(U)) = \pi_{U^*}(U^*) = \ssuc{\mathcal{V}}(\pi_U(U)) = \ssuc{\mathcal{V}}(\pi_{U^*}(U))$ and
$\mathcal{U}_U = \proper{\mathcal{U}_{U^*}}$, and thus $\pi_{U^*}(\ssuc{\mathcal{U}}(U')) = \ssuc{\mathcal{V}}(\pi_{U^*}(U'))$ 
for all $U' \in \proper{\mathcal{U}_{U^*}}$. Hence $\pi_{U^*}$ is a homomorphism and so $\ssuc{\mathcal{U}}(U) = U^* \in \mathcal{U}_0$.
This shows $\mathcal{U}_0$ is an invariant subset of $\mathcal{U}$. It follows that $\mathcal{U}_0 = \mathcal{U}$ and then by Lemma~\ref{lemma_enums_24}
$B \in \mathcal{U}_0$, which means that $\pi = \pi_B : \mathcal{U} \to \mathcal{V}$ is a homomorphism.

It remains to consider the uniqueness. Let $\pi' : \mathcal{U} \to \mathcal{V}$ be any homomorphism and put
$\mathcal{U}_0 = \{ U \in \mathcal{U} : \pi'(U) = \pi(U) \}$. Clearly $\varnothing \in \mathcal{U}_0$ and
if $U \in \proper{\mathcal{U}_0}$ then
$\pi'(\ssuc{\mathcal{U}}(U)) = \ssuc{\mathcal{V}}(\pi'(U)) = \ssuc{\mathcal{V}}(\pi(U)) = \pi(\ssuc{\mathcal{U}}(U))$, i.e.,
$\ssuc{\mathcal{U}}(U) \in \mathcal{U}_0$. Thus $\mathcal{U}_0$ is an invariant subset of $\mathcal{U}$, and so $\mathcal{U}_0 = \mathcal{U}$. 
This shows that $\pi' = \pi$, i.e., there is a unique homomorphism $\pi :  \mathcal{U} \to \mathcal{V}$.
\eop

\begin{theorem}\label{theorem_enums_2} 
If $A \approx B$, $\,\mathcal{U}$ is an $A$-enumerator and  $\,\mathcal{V}$ is a $B$-enumerator
then there exists a unique homomorphism $\pi : \mathcal{U} \to \mathcal{V}$ and $\pi$ maps $\,\mathcal{U}$ bijectively onto $\,\mathcal{V}$. 
Moreover, $\pi(A) = B$.
\end{theorem}

\proof 
This follows from Propositions \ref{prop_enums_2} and \ref{prop_enums_3}. (Note that $\pi(A)$ is a subset of $B$ with $\pi(A) \approx A$ and $A \approx B$ 
and thus with $\pi(A) \approx B$. Hence by Theorem~\ref{theorem_fs_2} $\pi(A) = B$.)
\eop

An important special case of Theorem~\ref{theorem_enums_2} is when there is a second  $A$-enumerator $\mathcal{U}'$. There then exists
a unique homomorphism $\pi : \mathcal{U} \to \mathcal{U}'$, $\pi$ maps $\mathcal{U}$ bijectively onto $\mathcal{U}'$ and $\pi(A) = A$. 

Let $\mathcal{U}$ be an $A$-enumerator. 
By Proposition~\ref{prop_enums_1} the mapping $\esuc{\mathcal{U}} : \proper{\mathcal{U}} \to A$ is a  bijection and so there is a unique
binary relation $\le$ on $A$ such that
$\esuc{\mathcal{U}}(U) \le \esuc{\mathcal{U}}(U')$ holds for $U,\, U' \in \proper{\mathcal{U}}$ if and only if $U \subset U'$.
More explicitly, this means that $a \le a'$ if and only if $\esuc{\mathcal{U}}^{-1}(a) \subset \esuc{\mathcal{U}}^{-1}(a')$,
where $\esuc{\mathcal{U}}^{-1} : A \to \proper{\mathcal{U}}$ is the inverse of the mapping $\esuc{\mathcal{U}}$.
It is clear that $\le$ is a total order and it will be called the \definition{total order associated with $\mathcal{U}$}.

\begin{proposition}\label{prop_enums_4}
Let $\le$ be a total order on $A$ and put $L_a = \{ a' \in A : a' < a \}$ for each $a \in A$ (where as usual $a' < a$  means that both $a' \le a$ and 
$a' \ne a$ hold). Then
$\,\mathcal{U} = \{ U \in \mathcal{P}(A) : \mbox{$U = L_a$ for some $a \in A$} \} \cup \{A\}$
is an $A$-enumerator with $\esuc{\mathcal{U}}(L_a) = a$ for each $a \in A$. 
\end{proposition}

\proof 
We assume that $A$ is non-empty, since the result holds trivially when $A = \varnothing$.
By Proposition~\ref{prop_fs_6} the non-empty set $A$ contains a unique $\le$-minimum element $a_0$ and then $L_{a_0} = \varnothing$, which shows that
$\varnothing \in \mathcal{U}$. Let $U \in \proper{\mathcal{U}}$ and let $a \in A$ be such that $U = L_a$; put $U' = U \cup \{a\}$. If
$U' = A$ then $a$ is trivially the unique element in $A \setminus U$ with $U \cup \{a\} \in \mathcal{U}$, so consider the case with
$U' \ne A$.
Then by Proposition~\ref{prop_fs_6} the non-empty set $A \setminus U'$ contains a unique $\le$-minimum element $a'$. 
Now $U' = \{ b \in A : b \le a \}$, 
thus $A \setminus U' = \{ b \in A : a < b \}$ and so $U' \subset L_{a'}$. But $a < c < a'$ for each $c \in L_{a'} \setminus U'$ and hence 
$L_{a'} \setminus U' = \varnothing$, since $a'$ is the $\le$-minimum element in $\{ b \in A : a < b \}$. 
It follows that $U' = L_{a'}$, i.e., $U \cup \{a\} \in \mathcal{U}$.
Suppose $U \cup \{b\} \in \mathcal{U}$ for some other $b \in A \setminus U$. Then $U \cup \{b\} = L_{b'}$ for some $b' \in A$ (since $U \cup \{b\} \ne A$)
and then $a \le b < a'$, since $L_a \subset L_{a'}$ and $b \in L_{a'} \setminus L_a$. But this implies
$a \in L_{a'}$ and hence $b = a$, since $a \notin L_a$. Thus $a$ is the unique element in $A \setminus U$ such that $U \cup \{a\} \in \mathcal{U}$,
This shows that $\mathcal{U}$ is an $A$-selector (since it is clearly totally ordered)
 and that $\esuc{\mathcal{U}}(L_a) = a$ for each $a \in A$. Moreover,  $\mathcal{U}$ Hence by Theorem~\ref{theorem_enums_1}
$\mathcal{U}$ is an $A$-enumerator. \eop

If $\le$ is a total order on $A$ then the $A$-enumerator $\mathcal{U}$ given in Proposition~\ref{prop_enums_4} will be called the
\definition{$A$-enumerator associated with $\le$}.

\begin{theorem}\label{theorem_enums_3}
(1)\enskip
If $\le$ is the total order associated with an $A$-enumerator $\mathcal{U}$ then $\mathcal{U}$ is the $A$-enumerator associated with $\le$.

(2)\enskip
If $\mathcal{U}$ is the $A$-enumerator associated with a total order $\le$ on $A$ then $\le$ is the total order associated with $\mathcal{U}$.

\end{theorem}

\proof 
(1)\enskip 
Let $\le$ be the total order associated with the $A$-enumerator $\mathcal{U}$ and let $\mathcal{V}$ be the $A$-enumerator associated with $\le$.
Thus if $U,\,U' \in \proper {\mathcal{U}}$ then $\esuc{\mathcal{U}}(U) \le \esuc{\mathcal{U}}(U')$ if and only if $U \subset U'$.
Now
$\,\mathcal{V} = \{ V \in \mathcal{P}(A) : \mbox{$U = L_a$ for some $a \in A$} \} \cup \{A\}$, 
 where $L_a = \{ a' \in A : a' < a \}$ for each $a \in A$. Let $U \in \proper{\mathcal{U}}$ and put
$a = \esuc{\mathcal{U}}(U)$. Then 
\[L_a = \{ a' \in A : a' < a \} = \{ b \in A : \mbox{$b = \esuc{\mathcal{U}}(U')$ for some $U' \in \mathcal{U}_U^p$}\}\]
and therefore by Lemma~\ref{lemma_enums_4} $L_a  = U$. It follows that $\mathcal{V} = \mathcal{U}$.

(2)\enskip
Let $\mathcal{U}$ be the $A$-enumerator associated with the total order $\le$  
and let  $\le'$ be the total order associated with $\mathcal{U}$.
Thus if $U,\,U' \in \proper {\mathcal{U}}$ then there exist $a,\,a' \in  A$ with $U = L_a$ and $U' = L_{a'}$
and by Proposition~\ref{prop_enums_4} $\esuc{\mathcal{U}}(U) = a$ and $\esuc{\mathcal{U}}(U') = a'$.
Thus $a \le' a'$ if and only if $U \subset U'$ which means that $a \le' a'$ if and only if
$\{ b \in A: b < a \} \subset \{ b \in A: b < a' \}$.
Suppose $a \le a'$. If $b < a$ then $b < a'$ and so $\{ b \in A: b < a \} \subset \{ b \in A: b < a' \}$.
Conversely, suppose $\{ b \in A: b < a \} \subset \{ b \in A: b < a' \}$.
Now either $a \le a'$ or $a' \le a$ and if $a' \le a$ then $\{ b \in A: b < a' \} \subset \{ b \in A: b < a \}$ from which it follows that
$\{ b \in A: b < a \} =\{ b \in A: b < a' \}$ and this is only possible if $a = a'$.
Therefore ${\le'} = {\le}$. \eop


\startsection{Iterators and assignments}

\label{iterators}

In this section we introduce what will be called an assignment of finite sets in a triple $\,\mathbf{I} = (X, f, x_0)$, where $X$ is some class of
objects, $f : X \to X$ is a mapping of the class $X$ into itself and $x_0$ is an object of $X$.  Such a triple will be called
an \definition{iterator}. The results of this section will
will be applied in Section~\ref{ordinals} to define the finite ordinals. In this case $(X,f,x_0) = (\fin,\vonn, \varnothing)$, where 
$\vonn : \fin \to \fin$ is the mapping given by $\vonn(A) = A \cup\{A\}$  for each finite set $A$.

We assume that the usual set-theoretic operations are also valid for classes. If $X$ and $Y$ are arbitrary classes there is a class $X \cup Y$ (their union), 
a class  $X \times Y$ (their product), a class $Y ^X$ (the class of all mappings from $X$ to $Y$)
and the power class $\mathcal{P}(X)$ of $X$ ( the class of all subclasses of $X$). 
The union, product and power class are defined as for sets and mappings are defined either as graphs or similarly to how mappings were defined in the finite case.

The archetypal example of an iterator whose first component is  a set is $(\Nat,\textsf{s},0)$, where the successor mapping $\textsf{s} : \Nat \to \Nat$ is 
given by $\textsf{s}(n) = n + 1$ for each $n \in \Nat$. However, we will also be dealing with examples in which the first component is a proper class or a finite set.
 In what follows let us fix an iterator $\mathbf{I} = (X,f, x_0)$.

A mapping $\val : \fin \to X$ will be called  an \definition{assignment of finite sets in $\mathbf{I}$} or, 
when it is clear what $\mathbf{I}$ is, simply  an 
\definition{assignment of finite sets}
if $\val(\varnothing) = x_0$ and 
\[\val(A \cup \{a\}) = f(\val(A))\] 
for each finite set $A$ and each element $a \notin A$.

For each finite set $A$ denote the cardinality of $A$ by $|A|$ (with $|A|$ defined as usual in terms of  $\Nat$). Then it is clear that the mapping
$|\cdot| : \fin \to \Nat$  defines an assignment of finite sets in $(\Nat,\textsf{s},0)$ and that this is the 
unique  such assignment. 

Let $X$ be a class and let $\lambda : \fin \to X$ be a mapping. For each $A \in \fin$ let $\lambda_A : \mathcal{P}(A) \to X$ be the restriction of
$\lambda$ to $\mathcal{P}(A)$. Then the family of mappings
$\{ \lambda_A :A \in \fin\}$, is \definition{compatible} in the sense that whenever $A,\, B \in \fin$ with $B \subset A$ then $\lambda_B$ is the restriction of
$\lambda_A$ to $\mathcal{P}(B)$, i.e., $\lambda_B(C) = \lambda_A(C)$ for all $C \in \mathcal{P}(B)$.

\begin{lemma}\label{lemma_iterators_111}
Let $\{ \lambda_A : A \in \fin\}$, be a compatible family  and define a mapping $\lambda : \fin \to X$ 
by setting $\lambda(A) = \lambda_A(A)$ for each $A \in \fin$. Then $\lambda_A$ is the restriction of $\lambda$ to $\mathcal{P}(A)$ for each $A \in \fin$. 
\end{lemma}

\proof  This is clear, since if  $ B \in \mathcal{P}(A)$ then $\lambda(B) = \lambda_B(B) = \lambda_A(B)$. \eop

In general, we assume that the following statements are valid for mappings between classes:

(1)\enskip Mappings are determined by their values. This implies in particular that the mapping $\lambda : \fin \to X$ defined in terms of a compatible family, 
is unique.

(2)\enskip For each class $X$ there is a mapping $\id_X : X \to X$ satisfying $\id_X(x) = x$ for all $x \in X$.

(3)\enskip If $X$, $Y$ and $Z$ are classes and $f : X \to Y$ and $g :Y \to Z$ are mappings then there is a mapping $g\circ f : X \to Z$ (their composition) satisfying
$(g \circ f)(x) = g(f(x))$ for all $x \in X$.

We will also need the following which corresponds to Proposition~\ref{prop_fs_77}:

\begin{proposition}\label{prop_iterators_77}
Let $X$, $Y$ and $Z$ be classes, let  $f : Z\to X$ be a surjective mapping and let $g : Z\to Y$ be a mapping. 
Then there exists a mapping $h : X \to Y$ with $g = h\circ f$ if and only if $g(z) = g(z')$ whenever $z,\,z' \in Z$
with  $f(z) = f(z')$. Moreover, if $h$ exists then it is unique.
\end{proposition}

\proof Assuming a one-to-one correspondence between mappings and graphs then  the proof of
Proposition~\ref{prop_fs_77} can be used. Without this assumption we can proceed exactly as in the alternative proof of Proposition~\ref{prop_fs_77}.
 Suppose first that there exists $f : X \to Y$ with  $g = h\circ f$. If $z,\,z' \in Z$ with $f(z) = f(z')$ then $g(z) = h(f(z)) = h(f(z')) = g(z')$
and so $g(z) = g(z')$ whenever $z,\,z' \in Z$ with $(z) = f(z')$. Moreover, $h(f(z)) = g(z)$ for each $z \in Z$  and $f$ is surjective and hence $h$ is uniquely determined by $f$ and $g$. 

Suppose conversely that $g(z) = g(z')$ whenever $z,\,z' \in Z$ with $f(z) = f(z')$. 
For each $x \in X$ let $G_x = \{ z \in Z : f(z) = x \}$. Thus $G_x \ne \varnothing$, since $f$ is surjective and if $x \ne x'$ then $G_x \cap G_{x'} = \varnothing$.
Let $\mathcal{E} = \{ E \in \mathcal{P}(Z) : \mbox{$E = G_x$ for some $x \in X$}\}$ and define $r : X \to \mathcal{E}$ by 
$r(x) = G_x$ for each $x \in X$. Hence $r$ is a bijection. Now if $z ,\,z'\in r(x)$ then $f(z) = f(z')$ and so $g(z) = g(z')$.
There is thus a mapping $q : \mathcal{E} \to Y$ such that $q(r(x)) =g(z)$, where $z$ is any element in $r(x)$ and note that $g(z)$ does not depend on which
element of $r(x)$ is used. Define $h : X \to Y$ by $h =q\circ r$ and so $h\circ f = q \circ r \circ f$. 
Let $z \in Z$; then $x = f(z) \in X$ and thus $r(x) = G_x \in \mathcal{E}$. Hence $q(G_x) = q(r(x)) = g(z)$, since $z \in r(x)$,
i.e., $(h \circ f)(z) = g(z)$, which shows that $h\circ f = g$. \eop

\begin{proposition}\label{prop_iterators_777}
Let $X$ and $Y$ be classes and let $f : Y \to X$ be a bijection. Then there exists the inverse mapping $f^{-1} : X \to Y$. This is the unique mapping
$g : X \to Y$ satisfying $g \circ f = \id_Y$ and $f \circ g = \id_X$.
\end{proposition}

\proof If we again assume that there is a one-to-one correspondence between mappings and graphs then then the proof of Proposition~\ref{prop_fs_222} based on
Proposition~\ref{prop_fs_333} can be used.
If not then we can proceed as follows: 
Let $f : Y \to X$ be a bijection. Set $Z = Y$ and so $f : Z \to X$ is a bijection. Also put $p : Z \to Y = \id_Y$.
If $z,\,z' \in Z$ with $f(z) = f(z')$ then $z= z'$, since $f$ is a bijection and so $p(z) =p(z')$. Thus by Proposition~\ref{prop_iterators_77} there exists a 
unique mapping $g : X\to Y$ such that $p = g \circ f$, i.e., with $g \circ f = \id_Y$.
Repeating the above construction with the bijection $g : X \to Y$ there exists a bijection $f' : Y \to X$ such that $f' \circ g =\id_X$.
Then $f' =f \circ g \circ \circ f' = f$, i.e., $f' = f$. This shows that $g = f^{-1}$.  \eop

\begin{theorem}\label{theorem_iterators_1}
(1)\enskip There exists a unique assignment $\val$ of finite sets in $\mathbf{I}$. 

(2) If $A$ and $B$ are finite sets with $A \approx B$ then $\val(A) = \val(B)$. 
\end{theorem}

\proof
Let $A$ be a finite set; 
then a mapping $\val_A : \mathcal{P}(A) \to X$ will be called an \definition{$A$-assignment} if 
$\val_A(\varnothing) = x_0$ and $\val_A(B \cup \{a\}) = f(\val_A(B))$ for each proper subset $B$ of $A$ and each 
$a \in A \setminus B$.

\begin{lemma}\label{lemma_iterators_1}
For each finite set $A$ there exists a unique $A$-assignment.
\end{lemma}

\proof 
Let $A$ be a finite set and let $\mathcal{S}$ be the set consisting of those $B \in \mathcal{P}(A)$ for which there exists a
unique $B$-assignment.
Then $\varnothing \in \mathcal{S}$, since the mapping $\val_\varnothing : \mathcal{P}(\varnothing) \to X$ with 
$\val_\varnothing(\varnothing) = x_0$ is clearly the unique $\varnothing$-assignment.

Let $B \in \proper{\mathcal{S}}$ with unique $B$-assignment $\val_B$, and let $a \in A \setminus B$; put 
$B' = B \cup \{a\}$. Now $\mathcal{P}(B')$ is 
the disjoint union of the sets $\mathcal{P}(B)$ and $\{ C \cup \{a\} : C \subset B\}$ and so we can define 
a mapping $\val_{B'} : \mathcal{P}(B') \to X$ by letting $\val_{B'}(C) = \val_B(C)$ and $\val_{B'}(C\cup \{a\}) = f(\val_B(C))$ 
for each $C \subset B$. Then $\val_{B'}(\varnothing) = \val_B(\varnothing) = x_0$, and so consider
$C' \subset B'$ and $b \in B' \setminus C'$. There are three cases:

The first is with $C' \subset B$ and $b \in B \setminus C'$ and here
\[\val_{B'}(C' \cup \{b\}) = \val_B(C' \cup \{b\}) = f(\val_B(C')) = f(\val_{B'}(C'))\;.\]

The second is with $C' \subset B$ and $b = a$. In this case
\[ \val_{B'}(C' \cup \{b\}) = \val_{B'}(C' \cup \{a\}) = f(\val_B(C')) = f(\val_{B'}(C'))\;. \]

The final case is with $C' = C \cup \{a\}$ for some $C \subset B$ and $b \in B \setminus C$, and here
\begin{eqnarray*}
 \val_{B'}(C' \cup \{b\}) &=& \val_{B'}(C \cup \{a\} \cup \{b\}) = f(\val_B(C \cup \{b\}))\\ 
&=& f(f(\val_B(C))) =  f(\val_{B'}(C \cup \{a\})) =  f(\val_{B'}(C'))\;. 
\end{eqnarray*}
In all three cases $\val_{B'}(C' \cup \{b\}) = f(\val_{B'}(C')$, which shows
$\val_{B'}$ is a $B'$-assignment.

Now let $\val'_{B'}$ be an arbitrary $B'$-assignment. In
particular $\val'_{B'}(C \cup \{b\}) = f(\val'_{B'}(C))$ for all $C \subset B$ and all $b \in B \setminus C$, and 
from  the uniqueness of the $B$-assignment $\val_B$ it follows that $\val'_{B'}(C) = \val_B(C)$ and thus also that
\[\val'_{B'}(C \cup \{a\}) = f(\val'_{B'}(C)) = f(\val_{B'}(C)) = \val_{B'}(C \cup \{a\})\]
 for all $C \subset B$, i.e., 
$\val'_{B'} = \val_{B'}$. Hence $B \cup \{a\} \in \mathcal{S}$.

Therefore $\mathcal{S}$ is an inductive $A$-system and thus $A \in \mathcal{S}$. This shows there exists a
unique $A$-assignment. \eop

\begin{lemma}\label{lemma_iterators_2}
If $A,\,B \in \fin$ with $B \subset A$; then the unique $B$-assignment $\val_B$ is the restriction of 
the unique $A$-assignment $\val_A$ to $\mathcal{P}(B)$.
\end{lemma}

\proof 
This follows immediately from the uniqueness of $\val_B$.
\eop

Lemma~\ref{lemma_iterators_2} shows that the family $\{\val_A : A \in \fin \}$ is compatible and therefore there there exists a unique mapping
$\val : \fin \to X$ such that $\val_A$ is the restriction of $\val$ to $\mathcal{P}(A)$ for each $A \in \fin$.
In particular, $\val(A) = \val_A(A)$ for each $A \in \fin$.
Thus $\val(\varnothing) = \val_\varnothing(\varnothing) = x_0$  and if $A \in \fin$ and $a \notin A$ then by Lemma~\ref{lemma_iterators_2} 
\[ 
\val(A\cup \{a\}) = \val_{A\cup \{a\}}(A\cup \{a\}) = f(\val_{A\cup \{a\}}(A)) = f(\val_A(A)) = f(\val(A))\;.  
\]
Hence $\val$ is an assignment of finite sets in $\mathbf{I}$.
For the uniqueness consider an arbitrary assignment $\val'$ of finite sets in $\mathbf{I}$. Then for each $A  \in \fin$ the restriction of $\val'$ to
$\mathcal{P}(A)$ is an $A$-iterator and thus equal to $\val_A$. It follows that $\val' = \val$.
This shows that there is a unique assignment $\val$ of finite sets in $\mathbf{I}$.

(2)\enskip We must show that if $A$ and $B$ are finite sets with $A \approx B$ then $\val(A) = \val(B)$. Let $A$ be a finite set and 
$\mathcal{S}$ be the set consisting of those $C \in \mathcal{P}(A)$ for which $\val(C) = \val(B)$ whenever $B$ is a finite set with $B \approx C$.
Then $\varnothing \in \mathcal{S}$, since $B \approx \varnothing$ if and only if $B = \varnothing$.
Consider $C \in \proper{\mathcal{S}}$ and $a \in A \setminus C$, and let 
$B$ be a finite set with $B \approx C \cup \{a\}$; thus $B \ne \varnothing$, so let $b \in B$. Then $B' = B \setminus \{b\} \approx C$
and hence $\val(B') = \val(C)$. Thus
$\val(B) = \val(B' \cup \{b\}) = f(\val(B')) = f(\val(C)) = \val(C \cup \{a\})$.
This shows that $C \cup \{a\} \in \mathcal{S}$. Therefore $\mathcal{S}$ is an inductive $A$-system and so
$A \in \mathcal{S}$, i.e., $\val(A) = \val(B)$ whenever $A \approx  B$.
This completes the proof of Theorem~\ref{theorem_iterators_1}. \eop

Consider the equivalence relation $\approx$ on $\fin$ and denote by $\fin_{/\approx}$ the corresponding class of equivalence classes. 
By Theorem~\ref{theorem_iterators_1} (2) there is then an induced mapping 
$\val_{/\approx} : \fin_{/\approx} \to X$.

For what follows it is necessary to determine the range of the assignment $\val$, this being the subclass 
$\,X_0 = \{ x \in X : \mbox{$x = \val(A)$ for some finite set $A$} \}$ of $X$. 
A subclass $Y$ of $X$ is said to be \definition{$f$-invariant} if $f(y) \in Y$ for all $y \in Y$. 
The next result shows that $X_0$ is the least $f$-invariant subclass of $X$ containing $x_0$.

\begin{lemma}\label{lemma_iterators_3}
$X_0$ is an $f$-invariant subclass of $X$
containing $x_0$. Moreover, if $X'$ is any $f$-invariant subclass of $X$ containing $x_0$ then $X_0 \subset X'$.
\end{lemma} 

\proof Clearly  $x_0 \in X_0$ since $x_0 = \val(\varnothing)$. Thus let $x \in X_0$, and so there exists a finite set $A$ with 
$x = \val(A)$. By Lemma~\ref{lemma_intro_3} there exists an element not in $A$;
it then follows that $\val(A \cup \{a\}) = f(\val(A)) = f(x)$, which implies that $f(x) \in X_0$. Hence 
$X_0$ is $f$-invariant. 
Now let $X'$ be any $f$-invariant subclass of $X$ containing $x_0$.
Let $A$ be a finite set and let $\mathcal{S} = \{ B \in \mathcal{P}(A) : \val(B) \in X' \}$.
Then $\varnothing \in \mathcal{S}$ since $\val(\varnothing) = x_0 \in X'$.
Consider $B \in \proper{\mathcal{S}}$ (and so $\val(B) \in X'$) and let $a \in A \setminus B$. 
Then $\val(B \cup \{a\}) = f(\val(B)) \in X'$, since $X'$ is $f$-invariant, and hence $B \cup \{a\} \in \mathcal{S}$.
Thus $\mathcal{S}$ is an inductive $A$-system and hence $A \in \mathcal{S}$, i.e., $\val(A) \in X'$.
This shows that $\val(A) \in X'$ for each finite set $A$  and it follows that $X_0 \subset X'$.
\eop

The iterator $\,\mathbf{I}$ is said to be \definition{minimal} if the only 
$f$-invariant subclass of $X$ containing $x_0$ is $X$ itself, thus $\,\mathbf{I}$ is minimal if and only if
$X_0 = X$. In particular, it is easy to see that the Principle of Mathematical Induction is exactly the requirement that the iterator 
$(\Nat,\textsf{s},0)$ be minimal.

Note that the iterator $\,\mathbf{I}_0 = (X_0,f_0,x_0)$ is always minimal, where $f_0 : X_0 \to X_0$
is the restriction of $ f $ to $X_0$.

For a minimal iterator Lemma~\ref{lemma_iterators_3} takes the form:

\begin{proposition}\label{lemma_iterators_4}
An iterator $\,\mathbf{I}$ is minimal if and only if the mapping $\val : \fin \to X$ is surjective, and thus if and only if
the induced mapping $\val_{/\approx} : \fin_{/\approx} \to X$ is surjective. 
\end{proposition} 

\proof 
This is a special case of Lemma~\ref{lemma_iterators_3}.
\eop

From now on we will make use of Lemma~\ref{lemma_intro_3} (guaranteeing the existence of an element $a$ not in a set $A$) without referring explicitly to this result.

\begin{proposition}\label{prop_iterators_1}
Suppose $\,\mathbf{I}$ is minimal; then $\{x_0\} \cup f(X) = X$. Thus for each $x \ne x_0$ there exists an $x' \in X$ with
$x = f(x')$. Moreover, the mapping $f$ is surjective if and only if $x_0 \in f(X)$.
\end{proposition}

\proof
For a general iterator the subclass $(\{x_0\} \cup f(X))$ is always $f$-invariant and contains $x_0$. Thus, since 
$\,\mathbf{I}$ is minimal it follows that $\{x_0\} \cup f(X) = X$. \eop

The iterator $\,\mathbf{I}$ will be called \definition{regular} if $B_1 \approx B_2$ whenever $B_1$ and $B_2$ are finite sets with 
$\val(B_1) = \val(B_2)$. Thus
the iterator $\,\mathbf{I}$ is regular if and only if the induced mapping $ \val_{/\approx} : \fin_{/\approx}  \to X$ is injective.

Note that if $\,\mathbf{I}$ is regular then so is  the minimal iterator $\,\mathbf{I}_0 = (X_0,f_0,x_0)$.

$\,\mathbf{I}$ will be called a \definition{Peano iterator} if it is minimal
and \standard, where \definition{\standard} means that the mapping $f$ is injective and $x_0 \notin f(X)$.
The Peano axioms thus require $(\Nat,\textsf{s},0)$ to be a Peano iterator. If $\,\mathbf{I}$ is a Peano iterator then $f$ is injective and so
Proposition~\ref{prop_iterators_1} implies that for each $x \ne x_0$ there exists a unique $x' \in X$ with $x = f(x')$. If $\,\mathbf{I} = (X,f,x_0)$
is a Peano iterator then by Theorem~\ref{theorem_fs_1} $X$ cannot be a finite set.

Theorem~\ref{theorem_iterators_2} below states that a minimal iterator is regular if and only if it is a Peano iterator. 
This will be applied to prove the recursion theorem for Peano iterators.

\begin{theorem}\label{theorem_iterators_2}
A minimal iterator $\,\mathbf{I}$ is regular if and only if it is \standard, i.e., if and only if it is a Peano iterator. 
Thus an iterator $\,\mathbf{I}$ is a Peano iterator if and only if the induced mapping $\val_{/\approx} : \fin_{/\approx} \to X$ is a bijection.
\end{theorem}

\proof
Assume first that $\,\mathbf{I}$ is \standard.
Let $A$ be a finite set and  let $\mathcal{S}$ denote the set of subsets $B$ of $A$ such that $B \approx B'$
whenever $B' \subset A$ with $\val(B') = \val(B)$.
Let $B \subset A$ with $B \ne \varnothing$, let $b \in B$ and put $B' = B \setminus \{b\}$. It then follows that
$\val(B) = \val(B' \cup \{b\}) = f(\val(B'))$, and so $\val(B) \ne x_0$, since $x_0 \notin f(X)$. Thus 
$\val(B) \ne \val(\varnothing)$, which shows that $\varnothing \in \mathcal{S}$, since $\varnothing \approx B$ if and 
only if $B = \varnothing$.

Let $B \in \proper{\mathcal{S}}$ and let $a \in A \setminus B$. Consider 
$B' \subset A$ with $\val(B') = \val(B \cup \{a\})$; then $\val(B') = f(\val(B)) \in f(X)$, hence $\val(B') \ne x_0$ and so 
$B' \ne \varnothing$. Let $b \in B'$ and put $C = B' \setminus \{b\}$; then
$f(\val(C)) = \val(C \cup \{b\}) = \val(B') = f(\val(B))$ 
and thus $\val(C) = \val(B)$, since $f$ is injective, and it follows that $C \approx B$, since $B \in \mathcal{S}$.
But $B' = C \cup \{b\}$ with $b \notin C$, $a \notin B$ and $C \approx B$, and therefore 
$B' = C \cup \{b\} \approx B \cup \{a\}$. Hence $B \cup \{a\} \in \mathcal{S}$,
which shows that $\mathcal{S}$ is an inductive $A$-system and thus that $\mathcal{S} = \mathcal{P}(A)$.
This implies that if $B_1,\,B_2$ are subsets of $A$ with $\val(B_1) = \val(B_2)$ then $B_1 \approx B_2$.

Now let $B_1$ and $B_2$ be arbitrary finite sets with $\val(B_1) = \val(B_2)$. Applying the above with $A = B_1 \cup B_2$ then shows that
$B_1 \approx B_2$. Thus $\,\mathbf{I}$ is regular. 

For the converse we assume $\,\mathbf{I}$ is not \standard{} and show this implies it is not regular.
Suppose first that $x_0 = f(x)$ for some $x \in X$.
By Lemma~\ref{lemma_iterators_4} there exists a finite set $A$ with $x = \val(A)$ and there exists some element $a$ not
in $A$. Then $A \cup \{a\} \not\approx \varnothing$ but $\val(A \cup \{a\}) = f(\val(A)) = f(x) = x_0 = \val(\varnothing)$.
Thus $(X,f,x_0)$ is not regular. 
Suppose now that $f$ is not injective and so there exist $x,\, x' \in X$ with $x \ne x'$ and $f(x) = f(x')$. By Lemma~\ref{lemma_iterators_4} there
exist finite sets $A$ and $B$ with $x = \val(A)$ and $x' = \val(B)$ and by Theorem~\ref{theorem_fs_4} and Proposition~\ref{prop_fs_112} (1)
we can assume that $B \subset A$. Thus $B$ is a proper subset of $A$,
since $\val(A) = x \ne x' = \val(B)$. Let $a \notin A$; then $B \cup \{a\}$ is a proper subset of
$A \cup \{a\}$ and so by Theorem~\ref{theorem_fs_2} $B \cup \{a\} \not\approx A \cup \{a\}$.
But $\val(B \cup \{a\}) = f(\val(B)) = f(x') = f(x) = f(\val(A)) = \val(A \cup \{a\})$,
which again shows $\mathbf{I}$ is not regular.
\eop

Here is the recursion theorem (which first appeared in Dedekind \cite{dedekind}).

\begin{theorem}\label{theorem_iterators_3}
If $\,\mathbf{I}$ is a Peano iterator then  for each iterator $\,\mathbf{J} = (Y,g,y_0)$ there exists a unique mapping
$ \pi: X \to Y$ with $\pi(x_0) = y_0$ such that $\pi \circ f = g \circ \pi$.
\end{theorem}

\proof 
As before let $\val$ be the assignment of finite sets in $\,\mathbf{I}$ and denote the assignment of finite sets in $\,\mathbf{J}$ by $\val'$.
If $A,\,B \in \fin$ with $\val(A) = \val(B)$ then by Theorem~\ref{theorem_iterators_2} $A \approx B$
and therefore by Theorem~\ref{theorem_iterators_1} (2) $\val'(A) = \val'(B)$. Moreover, by Proposition~\ref{lemma_iterators_4}
$\val$ is surjective and thus by Proposition~\ref{prop_iterators_77} there exists a unique factor mapping $\pi : X \to Y$ such that $\pi(\val(A)) = \val'(A)$ for each $A \in \fin$. In particular, 
 $\pi(x_0) = \pi(\val(\varnothing)) = \val'(\varnothing) = y_0$. Let $x \in X$; as above there exists a finite set $A$ 
with $x = \val(A)$, and there exists an element $a$ not contained in $A$. Hence 
\begin{eqnarray*}
  \pi(f(x)) = \pi(f(\val(A))) &=& \pi(\val(A \cup \{a\}))\\ 
 &=& \val'(A \cup \{a\}) = g(\val'(A)) = g(\pi(\val(A))) = g(\pi(x))
\end{eqnarray*}
and this shows that $\pi \circ f = g \circ \pi$. 

The proof of the uniqueness only uses the fact that $\,\mathbf{I}$ is 
minimal: Let $\pi' : X \to Y$ be a further mapping with $\pi'(x_0) = y_0$ and such that $\pi' \circ f = g \circ \pi'$ and let
$X' = \{ x \in X : \pi(x) = \pi'(x) \}$. Then $x_0 \in X'$, since $\pi(x_0) = y_0 = \pi'(x_0)$, and if $x \in X'$ then
$\pi'(f(x)) = g(\pi'(x)) = (\pi(x)) = \pi(f(x))$, i.e., $f(x) \in X'$. Thus $X'$ is an $f$-invariant subclass of 
$X$ containing $x_0$ and so $X' = X$, i.e., $\pi' = \pi$.
\eop

\begin{theorem}\label{theorem_iterators_4}
Let $\,\mathbf{I}$ be minimal; then the class $X$ is  a finite set if and only if $\,\mathbf{I}$ is not regular.
\end{theorem}

\proof
Suppose first that $X$ is a finite set.
Since $\,\mathbf{I}$ is minimal Proposition~\ref{prop_iterators_1} states that $f$ is surjective if and only if $x_0 \in f(X)$, and since $X$ 
is a finite set 
Theorem~\ref{theorem_fs_1} implies  $f$ is surjective if and only if it is injective. Therefore either $x_0 \in f(X)$ or $f$ is not injective,
which means that $\mathbf{I}$ is not \standard. It thus follows from Theorem~\ref{theorem_iterators_2} that $\,\mathbf{I}$ is not
regular. This can also be shown directly without using Theorem~\ref{theorem_iterators_2}:
Assume first that $x_0 = f(x)$ for some $x \in X$.
By Lemma~\ref{lemma_iterators_4} there exists a finite set $A$ with $x = \val(A)$; let $a$ be some element not
in $A$. Then $A \cup \{a\} \not\approx \varnothing$ but $\val(A \cup \{a\}) = f(\val(A)) = f(x) = x_0 = \val(\varnothing)$.
Thus $\mathbf{I}$ is not regular. 

Assume now that $f$ is not injective and so there exist $x,\, x' \in X$ with $x \ne x'$ and $f(x) = f(x')$. By Lemma~\ref{lemma_iterators_4} there
exist finite sets $A$ and $B$ with $x = \val(A)$ and $x' = \val(B)$ and by Theorem~\ref{theorem_fs_4} and Proposition~\ref{prop_fs_112} (1)
we can assume that $B \subset A$. Thus $B$ is a proper subset of $A$,
since $\val(A) = x \ne x' = \val(B)$. Let $a \notin A$; then $B \cup \{a\}$ is a proper subset of
$A \cup \{a\}$ and so by Theorem~\ref{theorem_fs_2} $B \cup \{a\} \not\approx A \cup \{a\}$.
But $\val(B \cup \{a\}) = f(\val(B)) = f(x') = f(x) = f(\val(A)) = \val(A \cup \{a\})$,
which again shows $(X,f,x_0)$ is not regular.

Suppose conversely that $\mathbf{I}$ is not regular, so there exist finite sets 
$A$ and $A'$ with $\val(A) = \val(A')$ and $A \not\approx A'$. Then by Proposition~\ref{prop_fs_112} (2) and Theorem~\ref{prop_iterators_1}  
there exist such subsets $A$ and $A'$ with $A'$ a proper subset of $A$.
We show that for each finite set $B$ there exists $C \subset A$ with $\val(C) = \val(B)$.
By Lemma~\ref{lemma_iterators_4} it then follows that the mapping $\val_A: \mathcal{P}(A)\to X$ with $\val_A(B) = \val(B)$ for each $B \subset A$ is surjective,
and hence by the remark following Proposition~\ref{prop_fs_2}~(2) that $X$ is a finite set, since by Proposition~\ref{prop_fs_3} $\mathcal{P}(A)$ is finite.

Thus let $B$ be a finite set; by Theorem~\ref{theorem_fs_4} and Proposition~\ref{prop_fs_112} (1) there exists a finite set $D$ with 
$D \approx B$ and either $D \subset A$ or $A \subset D$, and by Theorem~\ref{theorem_iterators_1} (2) $\val(D) = \val(B)$.
If $D \subset A$ then $C = D$ is the required subset of $A$. It remains to show that if $D$ is a finite set with $A \subset D$ then there
exists $C \subset A$ with $\val(C) = \val(D)$.

Thus let $D$ be a finite set with $A \subset D$. Put $D' = D \setminus A$ and let $\mathcal{S}$ be the set consisting of those $E \in \mathcal{P}(D')$ for which
there exists $C \subset A$ with $\val(C) = \val(A \cup E)$. In particular $\varnothing \in \mathcal{S}$. Consider $E \in \proper{\mathcal{S}}$ and so
$\val(C) = \val(A \cup E)$ for some $C \subset A$, let $b \in D' \setminus E$. If $C$ is a proper subset of $A$ and $a \in A \setminus C$ then
$C \cup \{a\} \subset A$ and $\val(C \cup \{a\}) = f(\val(C)) = f(\val(A \cup E)) = \val(A \cup (E \cup \{b\}))$. 
On the other hand, if $C = A$ and $a \in A \setminus A'$ then $A' \cup \{a\} \subset A$ and
 \[\val(A \cup (E \cup \{b\})) = f(\val(A \cup E)) = f(\val(C)) = f(\val(A)) = f(\val(A')) = \val(A' \cup \{a\})\;.\]
Thus $E \cup \{b\} \in \mathcal{S}$, which shows $\mathcal{S}$ is an inductive $D'$-system. Therefore $D' \in \mathcal{S}$ and hence
there exists $C \subset A$ with $\val(C) = \val(A \cup D') = \val(D)$.
\eop

Theorems~\ref{theorem_iterators_2} and \ref{theorem_iterators_4} imply that for a minimal iterator $\mathbf{I} = (X,f,x_0)$ there are two mutually exclusive possibilities: Either $\,\mathbf{I}$ is  a Peano iterator or $X$ is a finite set.

\begin{proposition}\label{prop_iterators_2}
Let $\,\mathbf{I}$ be minimal. If $x_0 \in f(X)$ then $X$ is a finite set and the mapping $f$ is bijective.
\end{proposition}

\proof
Exactly as in the proof above the fact that $x_0 \in f(X)$ implies $\,\mathbf{I}$ is not regular, and thus by Theorem~\ref{theorem_iterators_4}
$X$ is a finite set. Moreover, by Proposition~\ref{prop_iterators_1} $f$ is surjective, since $x_0 \in f(X)$, and therefore by
Theorem~\ref{theorem_fs_1} $f$ is bijective, since $X$ is a finite set.
\eop

We now give an example of a Peano iterator which is defined without making use of the natural numbers or any other infinite set. 
Let $h : \fin \to \fin$ be the mapping with $h(A) = \mathcal{P}(A)$ for all $A \in \fin$ and so we have an iterator $\,\mathbf{H'} = (\fin,h,\varnothing)$.
Also let $\,\mathbf{H} = (H,h,\varnothing)$ be the corresponding minimal iterator, thus $H$ is the minimal $h$-invariant subclass of $\fin$ containing $\varnothing$
and we denote the restriction  of $h$ to $H \to H$ again by $h$.

\begin{proposition}\label{prop_iterators_222}
$\,\mathbf{H}$ is a Peano iterator and $V \subset h(V)$ for all $V \in H$. Moreover, the sets in $H$ are transitive, where a set $A$ is \definition{transitive} if $x \subset A$ whenever $x \in A$. 
\end{proposition}

\proof If $A \ne A'$ then $\mathcal{P}(A) \ne \mathcal{P}(A')$ and hence $h$ is injective. Also $\mathcal{P}(A) \ne \varnothing$ for any set $A$ and so
$h(V) \ne \varnothing$ for all $V \in H$. This shows that $\,\mathbf{H}$  is a Peano iterator.

Note that $\mathcal{P}(E) \subset \mathcal{P}(\mathcal{P}(E))$for any set $E$
since if $F \in \mathcal{P}(E)$ then $F \subset E$ and thus $F$ is a subset of the subset $E$ of $E$;
hence $F \in\mathcal{P}(\mathcal{P}(E))$. Also, if $V \in H \setminus \{\varnothing\}$ then by Proposition~\ref{prop_iterators_1} 
$V = h_0(V') = \mathcal{P}(V')$ for some $V' \in H$. It follows that $V \subset h(V)$ for all $V \in H \setminus \{\varnothing\}$ and therefore
$V \subset h(V)$ for all $V \in H$, since $\varnothing \subset h(\varnothing)$ holds trivially. 

The transitivity of the sets in $V$ follows from the minimality of
$\,\mathbf{H}$ and the fact that if $A$ is transitive then so is its power set $\mathcal{P}(A)$: Let $B \in \mathcal{P}(A)$, i.e., $B \subset A$
and let $b \in B$. Then $b \in A$ and hence $b\subset A$, since $A$ is transitive, and so $b \in \mathcal{P}(A)$. Thus $b \in \mathcal{P}(A)$ for all $b \in B$,
i.e., $B \subset \mathcal{P}(A)$, which shows $\mathcal{P}(A)$ is transitive. \eop

Denote by $V_\omega$ the union of the sets in $H$. Thus a finite set $A$ is an element of $V_\omega$ if there exists $V \in H$ with $A \in V$ and
then $A \subset V$, since $V$ is transitive. On the other hand, if $A \subset V$ for some $V \in H$ then
$A \in \mathcal{P}(V) = h_0(V)$ and so $A \in V_\omega$. Thus the elements in $V_\omega$ are also exactly  the subsets of the sets in $H$. The elements of $V_\omega$
 are called \definition{hereditarily finite sets}. It is easy to see that $V_\omega$ cannot be a finite set and thus if the negation of the  axiom of infinity is assumed then $V_\omega$ must be a proper class.

Note that, although the sets in $H$ are finite, they rapidly become extremely large. Let us index the sets in $H$ using the natural
numbers with $V_0 = \varnothing$ and $V_{n+1} = h_0(V_n) = \mathcal{P}(V_n)$ for all $n \in \Nat$. Then $|V_0| = 0$ and $|V_{n+1}| = 2^{|V_n|}$ for all $n \in \Nat$.
Thus $|V_0| = 0$, $|V_1| = 1$, $|V_2| = 2$, $|V_3| = 4$, $|V_4| = 16$, $|V_5| = 65336$ and $|V_6| = 2^{65336}$.

Before going any further we need to be more explicit about the structure preserving mappings between iterators. In the following let
 $\,\mathbf{I} =(X,f,x_0)$ and 
$\,\mathbf{J}= (Y,g,y_0)$ be iterators; a mapping $\mu : X \to Y$ is called a \definition{morphism} from 
$\,\mathbf{I}$ to $\,\mathbf{J}$ if $\mu(x_0) = y_0$ and $g\circ \mu = \mu \circ f$.
 This will also be expressed by 
saying that $\mu : \mathbf{I} \to \mathbf{J}$ is a morphism. The recursion theorem thus states that if $\,\mathbf{I}$ is a Peano iterator then for each iterator
$\,\mathbf{J}$ there exists a unique morphism $ : \,\mathbf{I} \to \,\mathbf{J}$.

\begin{lemma}\label{lemma_iterators_5}
(1)\enskip The identity mapping $\id_X$ is a morphism from $\,\mathbf{I}$ to $\,\mathbf{I}$.  

(2)\enskip Let $\,\mathbf{K} = (Z,h,z_0)$ be a further iterator. If $\mu : \mathbf{I} \to \mathbf{J}$ and $\nu : \mathbf{J} \to \mathbf{K}$
 are morphisms then $\nu \circ \mu$ is a morphism from $\,\mathbf{I}$ to $\,\mathbf{K}$. 
\end{lemma}

\proof 
(1)\enskip
This is clear, since $\id_X(x_0) = x_0$ and $f \circ \id_X = f = \id_X \circ f$. 

(2)\enskip This follows since $(\nu \circ \mu)(x_0) = \nu(\mu(x_0)) = \nu(y_0) = z_0$ and
\[
  \beta\circ (\nu \circ \mu) =  (\beta\circ \nu) \circ \mu =  (\nu \circ g) \circ \mu 
=  \nu \circ (g \circ \mu) =  \nu \circ (\mu \circ f) = (\nu \circ \mu) \circ f\;. \eop\]

If $\mu : \mathbf{I} \to \mathbf{J}$ is a morphism then clearly $\mu \circ \id_X = \pi = \id_Y \circ \mu$, and 
if $\mu,\,\nu$ and $\tau$ are morphisms for which the compositions are defined then
$(\tau \circ \nu) \circ \mu = \tau \circ (\nu \circ \mu)$. 

\begin{lemma}\label{lemma_iterators_15}
(1)\enskip If $\,\mathbf{I}$ is minimal then there is at most one morphism $\mu : \mathbf{I} \to \mathbf{J}$.

(2)\enskip If $\,\mathbf{J}$ is minimal and $\mu : \mathbf{I} \to \mathbf{J}$ is a morphism  then $\mu$ is surjective.

(3)\enskip If $\,\mathbf{I}$ is minimal and $\mu :\mathbf{I} \to \mathbf{J}$ is a morphism  then $\mu(X) = Y_0$, where $Y_0$ is the least
$g$-invariant subclass of $Y$ containing $y_0$. In particular, if $\mu$ is surjective then $\,\mathbf{J}$ is minimal.
\end{lemma}

\proof (1) \enskip Let $\mu,\,\mu' :\,\mathbf{I}\to\,\mathbf{J}$ be morphisms and let 
$X_0 = \{ x \in X : \mu(x) = \mu'(x) \}$.
Then $x_0 \in X_0$, since $\mu(x_0) = \nu'(x_0) = y_0$ and if $x \in X_0$ then 
\[\mu(f(x)) = g(\mu(x)) = g(\mu'(x)) = \mu'(f(x))\]
and therefore $f(x) \in X_0$. Hence $X_0$ is an $f$-invariant subclass of $X$ containing $x_0$ and so 
$X_0 = X$,
i.e., $\mu = \mu'$. 

(2) \enskip Let $Y_0 = \{ \mu(x) : x \in X \}$. Then $y_0 = \mu(x_0) \in Y_0$ and if  
$h = \mu(x) \in Y_0$ then $g(h) = g(\mu(x)) =\mu(g(x)) \in Y_0$. 
Thus $Y_0$ is a $g$-invariant subclass of $Y$ containing $y_0$ and hence $Y_0 = Y$. This shows $\mu$ is surjective. 

(3)\enskip Let $X_0 = \{ x \in X : \mu(x) \in Y_0 \}$. Then $x_0 \in X_0$, since 
$\mu(x_0) = y_0 \in Y_0$, and if $x \in X_0$ then $\mu(f(x)) = g(\mu(x)) \in X_0$, 
since $Y_0$ is $g$-invariant. Therefore $X_0$ is a $f$-invariant
subclass of $X$ containing $x_0$ and hence $X_0 =X$, i.e., $\mu(X) \subset Y_0$. Now since
$\mu(X) \subset Y_0$ we can consider $\mu$ as a morphism $\mu_0$ from $\,\mathbf{I}$ to $\,\mathbf{J}_0$, 
where  $\,\mathbf{J}_0$ is the corresponding minimal iterator, and by (2) $\mu_0$ is surjective. But this implies that $\mu(X) = Y_0$. \eop 

The iterators $\,\mathbf{I}$ and $\,\mathbf{J}$ are said to be 
\definition{isomorphic} if there exists a  morphism $\mu : \mathbf{I} \to \mathbf{J} $ and a morphism $\nu : \mathbf{J} \to \mathbf{I}$ such that 
$\nu \circ \mu = \id_{X}$ and  $\mu \circ \nu = \id_{Y}$ and the morphism $\mu$ is then said to be an isomorphism. In particular, the mappings $\mu$ and $\nu$ are then both bijections.   

\begin{lemma}\label{lemma_iterators_16}
If $\mu : \mathbf{I} \to \mathbf{J}$ is a morphism and the mapping $\mu : X \to Y$ is a bijection then the inverse mapping
$\mu^{-1} : Y \to X$ is a morphism from  $\,\mathbf{J}$ to $\,\mathbf{I}$ and so $\,\mathbf{I}$ and $\,\mathbf{J}$ are isomorphic.
\end{lemma}

\proof
We have $g = g\circ \mu \circ \mu^{-1} = \mu \circ f \circ \mu^{-1}$ and so
$\mu^{-1} \circ g =\mu^{-1} \circ \mu \circ f \circ \mu^{-1} = f \circ \mu^{-1}$.
Thus $\mu^{-1}$ is a morphism from  $\,\mathbf{J}$ to $\,\mathbf{I}$, since also $\mu^{-1}(y_0) = x_0$. \eop

\begin{theorem}\label{theorem_iterators_7}
An iterator which is isomorphic to a Peano iterator is itself a Peano iterator.
\end{theorem}
\proof Let $\,\mathbf{I} = (X,f,x_0)$ be a Peano iterator and $\mu : \mathbf{I} \to \mathbf{J}$ be an isomorphism with $\,\mathbf{J} = (Y,g,y_0)$.
Thus $\mu : X \to Y$ is a bijection with $\mu(x_0) = y_0$ and $g \circ \mu = \mu \circ f$. Hence $g = \mu \circ f \circ \mu^{-1}$ and so
$g$ is injective, since $f$ is injective. Also $f(x) \ne x_0$ for all $x \in X$ and thus $(\mu \circ f )(x) \ne y_0$ for all $x \in X$. Therefore
$(g \circ \mu)(x) \ne g_0$ for all $x\in X$ which implies that $g(y) \ne y_0$ for all $y \in Y$. This shows that $\,\mathbf{J}$ is $\Nat$-like. Hence by 
Theorem~\ref{theorem_iterators_2} $\,\mathbf{J}$ is a Peano iterator, since by Lemma~\ref{lemma_iterators_15} (3) $\,\mathbf{J}$ is minimal. \eop

The iterator $\,\mathbf{I}$ is said to be \definition{initial} if for each iterator $\,\mathbf{J}$ there 
is a unique morphism from $\,\mathbf{I}$ to $\,\mathbf{J}$. Theorem~\ref{theorem_iterators_3} (the recursion theorem) thus states that a Peano iterator is initial.

\begin{lemma}\label{lemma_iterators_6}
(Let  $\,\mathbf{I}$ be initial  and $\pi : \mathbf{I} \to \mathbf{J}$ be the unique morphism.

1(1)\enskip If $\,\mathbf{J}$ is initial then $\pi$ is an isomorphism and so $\,\mathbf{I}$ and $\,\mathbf{J}$ are isomorphic. This says that any two initial
iterators are isomorphic and so  in particular, any two
Peano iterators are isomorphic.

(2)\enskip If $\pi$ is an isomorphism then $\,\mathbf{J}$ is initial.
\end{lemma}

\proof (1)\enskip There exists a unique morphism $\tau : \mathbf{J} \to \mathbf{I}$ (since $\,\mathbf{J}$ is initial). Thus  
$\tau \circ \pi : \mathbf{I} \to \mathbf{I}$ is a morphism. But $\id_{X} : \mathbf{I} \to \mathbf{I}$ is also a morphism and there is a unique morphism
from $\,\mathbf{I}$ to $\,\mathbf{I}$ (since $\,\mathbf{I}$ is initial) and hence $\tau \circ \pi = \id_{X}$. In the same way $\pi \circ \tau = \id_{Y}$.
Therefore $\pi$ is an isomorphism and so $\,\mathbf{I}$ and $\,\mathbf{J}$ are isomorphic. 

(2)\enskip Let $\,\mathbf{K}$ be an iterator and $\mu : \,\mathbf{I} \to \mathbf{K}$ be  the unique morphism. Then $\mu \circ \pi^{-1}$ is a morphism from 
$\,\mathbf{J}$ to $\,\mathbf{K}$. If $\nu : \mathbf{J}\to \mathbf{K}$ is any morphism then $\nu \circ \pi$ is a morphism from $\,\mathbf{I}$ to $\,\mathbf{K}$
and thus $\nu \circ \pi = \mu$. Hence $\nu = \mu \circ \pi^{-1}$ and so there is a unique morphism from $\,\mathbf{J}$ to $\,\mathbf{K}$, which shows that
$\,\mathbf{J}$ is initial. \eop

By Proposition~\ref{prop_iterators_222} there exists a
Peano (and thus an initial iterator) $\,\mathbf{H}$ which is defined without making use of an infinite set and by Lemma~\ref{lemma_iterators_6}
$\,\mathbf{H}$ is, up to isomorphism, the unique initial iterator.  

In Section~\ref{ordinals} we will exhibit another initial iterator
$\mathbf{O} = (O,\vonn,\varnothing)$ which is also defined only in terms of finite sets. The elements of
$O$ are the finite ordinals. 

The following result of Lawvere \cite{lawvere} shows that the converse of the recursion theorem holds. 

\begin{theorem}\label{theorem_iterators_6}
An initial iterator $\,\mathbf{I} = (X,f,x_0)$ is a Peano iterator.
\end{theorem}

\proof
We first show that $\,\mathbf{I}$ is minimal, and then that it is \standard.

\begin{lemma}\label{lemma_iterators_7}
The initial iterator $\,\mathbf{I}$ is minimal.
\end{lemma}

\proof Let $\,X_0 = \{ x \in X : \mbox{$x = \val(A)$ for some finite set $A$} \}$ of $X$,
let $f_0 : X_0 \to X_0$ be the restriction of $f$ to $X_0$, 
Then the iterator $\,\mathbf{I}_0 = (X_0,f_0,x_0)$ is minimal
and the inclusion mapping $\mathrm{inc} : X_0 \to X$ defines a morphism from $\,\mathbf{I}_0$ to 
$\,\mathbf{I}$. Let $\mu : \mathbf{I} \to \mathbf{I}_0$ be the unique morphism; then
$\mathrm{inc} \circ \mu = \id_X$, since by Lemma~\ref{lemma_iterators_5} $\mathrm{inc} \circ \mu$ and $\id_X$ 
are both morphisms from $\,\mathbf{I}$ to $\,\mathbf{I}$ (and there is only one such morphism, since $\,\mathbf{I}$ is 
initial). In particular, $\mathrm{inc}$ is surjective, which implies that $X_0 = X$, i.e., $\,\mathbf{I}$ is minimal.
\eop

\begin{lemma}\label{lemma_iterators_8}
The initial iterator $\,\mathbf{I}$ is \standard.
\end{lemma}

\proof
Let $\diamond$ be an element not contained in $X$, put $X_\diamond = X \cup \{\diamond\}$ and define a mapping 
$f_\diamond : X_\diamond \to X_\diamond$ by putting $f_\diamond(x) = f(x)$ for $x \in X$ and 
$f_\diamond(\diamond) = x_0$; thus $\,\mathbf{I}_\diamond = (X_\diamond,f_\diamond,\diamond)$ is an iterator. Since $\,\mathbf{I}$ is 
initial there exists a unique morphism $\mu : \mathbf{I} \to \mathbf{I}_\diamond$. Consider the subclass 
$X' = \{ x \in X : f_\diamond(\mu(x)) = x \}$; then $x_0 \in X'$, since
$f_\diamond(\mu(x_0)) = f_\diamond(\diamond) = x_0$ and if $x \in X'$ then $f_\diamond(\mu(x)) = x$ and so
\[  f_\diamond(\mu(f(x))) = f_\diamond(f_\diamond(\mu(x))) = f_\diamond(x) = f(x)\;, \]
i.e., $f(x) \in X'$. Thus $X'$ is a $f$-invariant subclass of $X$ containing $x_0$ and hence $X' = X$, since by 
Lemma~\ref{lemma_iterators_7} $\mathbf{I}$ is minimal. Thus $\mu(f(x)) = f_\diamond(\mu(x)) = x$ for all $x \in X$, which 
implies that $f$ is injective. Moreover, $x_0 \notin f(X)$, since
\[\nu(f(x)) = f_{\diamond}(\mu(x)) \ne \diamond =\mu(x_0))\]
 for all $x\in X$. Hence $\mathbf{I}$ is \standard.
\eop

This completes the proof of Theorem~\ref{theorem_iterators_6}.
\eop

We now give another proof that the
definition of a finite set being used here is equivalent to the usual one. The usual definition of $A$ being 
finite is that there exists $n \in \Nat$ and a 
bijective mapping $h : [n] \to A$, where $[n] = \{0,1,\ldots,n-1\}$ for 
$n \in \Nat \setminus \{0\}$ and $[0] = \varnothing$. Moreover, if $h : [n] \to A$ is a bijective mapping then $n$ is the 
cardinality of $A$, i.e., $n = |A|$, and so $A \approx [\,|A|\,]$ for each finite set $A$. The problem here is 
to assign a meaning to the expression $\{0,1,\ldots,n-1\}$, and one way to do this is to make use of the fact 
that $\{0,1,\ldots,n\} =  \{ m \in \Nat : m < n \}$ for all $n \in \Nat$. A similar approach works with any Peano iterator 
$\,\mathbf{I} = (X,f,x_0)$, since there is a natural way to define a total order $\le$ on $X$ which corresponds to the usual total order on $\Nat$.

Let $\,\mathbf{I} = (X,f,x_0)$ be an iterator. A total order $\le$ on $X$ will be called \definition{compatible with $\,\mathbf{I}$} if
$x \le f(x)$ for all $x \in X$.

\begin{lemma}\label{lemma_iterators_199}
Let $\,\mathbf{I} = (X,f,x_0)$ be a Peano iterator with $\val$ the assignment of finite sets in $\,\mathbf{I}$ and $\le$ be a total order on $X$ compatible
with $\,\mathbf{I}$. Let $x,\,y \in X$ and $A,\,B \in \fin$ with $x = \val(A)$ and $y = \val(B)$. Then $y \le x$ if and only if $B \preceq A$.
\end{lemma}

\proof We first show that if $B \preceq A $ then $y \le x$ and by by Proposition~\ref{prop_fs_112} (4) and Theorem~\ref{theorem_iterators_1} (2) we can assume, without loss of generality, that $B \subset A$. Put
$C = A\setminus B$ and let $\mathcal{S} = \{ D \in \mathcal{P}(C) : y \le \val(B \cup D) \}$; in particular, $\varnothing \in \mathcal{S}$. Thus let
$D \in \proper{\mathcal(S)}$ and $d \in C \setminus D$. Then 
$y \le \val(B \cup D) \le f(\val(B \cup D)) = \val(B \cup D \cup \{d\})$
and hence $D \cup \{d\} \in \mathcal{S}$. This shows that $\mathcal{S}$ is an inductive-$C$-system and so $C \in \mathcal{S}$, i.e., 
$y \le \val(B \cup C) = \val(A) = x$. We next assume that $y \le x$ and show that then $B \preceq A$. Now by Theorem~\ref{theorem_fs_4} either
$B\preceq A$ or $A \preceq B$. If $B \preceq A$ then by the above $y \le x$. On the other hand, if $A \preceq B$ then by the above $x \le y$ and so
$x = y$, since by assumption $y \le x$. Thus $A \approx B$, because $\,\mathbf{I}$ is regular, and in particular $B \preceq A$.
\eop

\begin{theorem}\label{theorem_iterators_199}
Let $\,\mathbf{I} = (X,f,x_0)$ be a Peano iterator.
Then there exists a unique total order $\le$ on $X$ compatible with $\,\mathbf{I}$ and the following hold:

(1)\enskip $x_0 \le x$ and $x <f(x)$ for all $x \in X$, where as usual we write $y < x$ if $y \le x$ but $y \ne x$.

(2)\enskip If $y \le x$ then $f(y) \le f(x)$.

(3)\enskip Let $x,\,y \in X$. Then $y \le f(x)$ if and only if $y \le x$ or $y = f(x)$. Thus $y < f(x)$ if and only if $y \le x$.

(4)\enskip For each $x \in X$ let $L_x = \{ y \in X : y < x \}$. Then $L_{x_0} = \varnothing$ and
$L_{f(x)}$ is the disjoint union of $L_x$ and the singleton set $\{x\}$  for each $x \in X$. Moreover, $L_{f(x)}$ is also the disjoint union of $\{x_0\}$ and $f(L_x)$
for each $x \in X$.

(5) \enskip The subclass  $L_x$ is a finite set for each $x \in X$.

(6) \enskip Each non-empty subclass $Y$ of $X$ contains a minimum element, i.e., an element $x$ with $x \le y$ for all $y \in Y$.
\end{theorem}

\proof\enskip Let $\val$ be the assignment of finite sets in $\,\mathbf{I}$. Define  a binary relation $\le$ on $X$ as follows: 
If $x,\,y \in X$ then $y \le x$ if and only if there exist
$A,\,B \in \fin$ with $x = \val(A)$, $y = \val(B)$ and such that $ B\preceq A$.

Let $x \in X$; then there  exists $A \in \fin$ with $x = \val(A)$, since $\,\mathbf{I}$ is minimal. If also $x = \val(A')$ then
$A \approx A'$, since $\,\mathbf{I}$ is regular. Now let $x,\,y \in X$ and $A,\,A',\,B,\,B'$ be finite sets with 
$x = \val(A) = \val(A')$ and $y = \val(B) = \val(B')$, and so $A \approx A'$ and $B \approx B'$. Hence if $B \preceq A$ then also $B' \preceq A'$. Thus if
$y \le x$ then $B \preceq A$ for all $A,\,B \in \fin$ with $x = \val(A)$ and $y = \val(B)$.

Clearly $x \le x$ for all $x \in X$ since $A \preceq A$ for all $A \in \fin$. Next let $x,\,y,\,z \in X$ with $x \le y$ and $y \le z$; then
$x \le z$, since if $A,\,B,\,C \in \fin$ with $A \preceq B$ and $B \preceq C$ then $A \preceq C$. Moreover, if $x,\,y \in X$ with $x \le y$ and $y \le x$
and $A,\,B \in \fin$ are such that $x = \val(A)$ and $y = \val(B)$ then $B \preceq A$ and $A \preceq B$. Then by Theorem~\ref{theorem_fs_4} 
$A \approx B$ and hence by Theorem~\ref{theorem_iterators_1}~(2) $x = y$. This shows that $\le$ is a partial order on $X$. Finally, if $A,\,B \in \fin$ then
by Theorem~\ref{theorem_fs_4} either $A \preceq B$ or $B \preceq A$ and thus if $x,\,y \in X$ then either $x \le y$ or $y \le x$. Hence $\le$ is a total order.

(1)\enskip Let $x \in X$ and $A \in \fin$ with $x = \val(A)$; let $a$ be an element not in $A$. Then $\val(A \cup \{a\}) = f(x)$ and $A \preceq A \cup \{a\}$
and therefore $x \le f(x)$. But if $x = f(x)$ then $f(x) = f(f(x))$, which is not possible since $f$ is injective. Thus $x < f(x)$.
Moreover, $x_0 \le x$ since $x_0 = \val(\varnothing)$ and $\varnothing \preceq A$ for all $A \in \fin$. In particular, $x \le f(x)$ for all
$x \in X$ and so $\le$ is compatible with $\,\mathbf{I}$. Moreover, by Lemma~\ref{lemma_iterators_199} $\le$ is the unique total order on $X$ compatible
with $\,\mathbf{I}$.

(2)\enskip Let $x,\,y \in X$ with $y \le x$ and let $A,\,B \in \fin$ with $x = \val(A)$ and $y = \val(B)$, and so $B \preceq A$. Choose an element $C$ not
in $A \cup B$. Thus $B \cup \{c\} \preceq A \cup \{c\}$ and $f(y) = \val(B \cup \{c\})$, $f(x) = \val(A \cup \{c\})$. Therefore $f(y) \le f(x)$.

(3)\enskip Let $x,\, y \in X$. If $y \le x$ then $y \le f(x)$ since $x \le f(x)$ and if $y = f(x)$ then clearly $y \le f(x)$. Thus if $y \le x$ or $y = f(x)$
then $y \le f(x)$.  Suppose conversely that $y \le f(x)$ and let $A,\,C \in \fin$ with $y = \val(B)$ and $f(x) = \val(C)$, and so $B \preceq C$. Let $p : B \to C$ 
be an injective mapping. If $B \approx C$
then $y = f(x)$ and thus assume this is not the case. Then $p$ is not surjective. Choose $a \in  C \setminus p(B)$ and put $A = C \setminus \{a\}$. Then $p$,
considered as a mapping from  $B$ to $A$ is injective. Hence $B \preceq A$ which means that $y \le \val(A)$. But $f(\val(A)) = \val(C) = f(x)$ and $f$ is injective,
which implies that  $\val(A) = x$. Therefore if $y \le f(x)$ then either $y \le x$ or $y = f(x)$.

(4)\enskip Clearly $L_{x_0} = \varnothing$, since $x_0 \le x$ for all $x \in X$.
Now $L_{f(x)} = \{ y \in X : y < f(x) \}$ and so by (3) $L_{f(x)}$ is the disjoint union of $L_x$ and $\{x\}$. Moreover,
\[f(L_x) = \{ y \in X \setminus \{x_0\} : \mbox{ $ y = f(z)$ with $ z < x$}\} = \{ y \in X  \setminus \{x_0\} : y \le x\}\]
and hence $fL_x) \cup \{x_0\} = \{y \in X : y \le x\} =L_x \cup \{x\} = L_{f(x)}$. The union $f(L_x) \cup \{x_0\}$ is disjoint since
$x_0 \notin f(X)$.

(5)\enskip Let $X_0 = \{ x \in X : \mbox{$L_x$ is a finite set} \}$. Then by (4) $X_0$  contains $x_0$ and is $f$-invariant. Thus 
$X_0 = X$, since $\,\mathbf{I}$ is minimal. 

(6) \enskip Let $z \in Y$; by (5)  $L_z \cup \{z\}$ is then  is a non-empty finite totally ordered set and hence by Proposition~\ref{prop_fs_6} it 
contains a minimal element $x$. Thus $x \le y$ for all $y \in Y$. \eop

The following result is needed when dealing with lists in Section~\ref{lists}.

\begin{proposition}\label{prop_iterators_399}
Let $\,\mathbf{I} = (X,f,x_0)$ be a Peano iterator with $\le$ the unique total order
on $X$ compatible with $\,\mathbf{I}$. For each $x\in X$ put $X_x  = \{ y \in X :x \le y \}$. Then: 

(1)\enskip $X_x$ is the least $f$-invariant subclass of $X$ 
containing $x$.

(2)\enskip $\,\mathbf{I}_x = (X_x,f,x)$ is a Peano iterator.

(3)\enskip Let $\pi_x :\mathbf{I} \to \mathbf{I}_x$ be the unique morphism, which is an isomorphism, since both $\,\mathbf{I}$ and $\,\mathbf{I}_x$ are Peano iterators. 
If $y_1,\,y_2 \in X$ then $y_1 \le y_2$ if and only if $\pi_x(y_1) \le \pi_x(y_2)$.

(4)\enskip If $y_1,\,y_2  \in X$ with $y_1 < y_2$ then 
\[\pi_x(\{ z \in X : y_1 \le z < y_2 \}) = \{ z' \in X : \pi_x(y_1) \le z' < \pi_x(y_2) \}\;.\] 

(5)\enskip If $y \in X$ then $L_{\pi_x(y)}$ is the disjoint union of  $ L_x$ and  $\pi_x(L_y)$.  

\end{proposition}

\proof (1)\enskip $X_x$ is an $f$-invariant subclass of $X$ containing $x$. Suppose there is an $f$-invariant subclass $X'$ of $X$ containing $x$ which is a proper
subclass of $X_x$. Then $X_x \setminus X'$ contains a minimum element $y$ and $y > x$ since $x \in X'$. Thus there exists $y' \in X_x$ with $y = f(y')$ and then $y' < y$. Hence $y' \in X'$. But then $y = f(y') \in X'$, since $X'$ is $f$-invariant. This contradiction shows that $X_x$ is the least $f$-invariant subclass of $X$
containing $x$.

(2)\enskip By (1) $\mathbf{I}_x$ is minimal. If $y \in f(X_x)$ then $y = f(y')$ for some $y' \in X_x$ and so $y > x$. Thus $x \notin f(X_x)$.
Hence $\,\mathbf{I}_x$ is a Peano iterator, since $f$ is also injective.

(3)\enskip Let $A \in \fin$ with $\val(A) = x$, where $\val$ is the assignment of finite sets in $\,\mathbf{I}$.
We show that if  $y \in X$ and $B \in \fin$ with $B \cap A = \varnothing$ and $\val(B) = y$ then $\pi_x(y) = \val(A \cup B)$.
Let $\mathcal{S}$ be the set of subsets $C$ of $B$ for which
$\pi_x(\val(C)) = \val(A \cup C)$ and so $\varnothing \in \mathcal{S}$. Let $D \in \proper{\mathcal{S}}$ and $d \in C \setminus D$; put $D' = D \cup \{d\}$.
Then $\pi_x(\val(D')) = f(\pi_x(\val(D))) = f(\val(A \cup D)) = \val(A \cup D')$ and hence $D' \in \mathcal{S}$. This hows that $\mathcal{S}$ is an
inductive-$B$-system and so $B \in \mathcal{S}$, i.e., $\pi_x(y) = \pi_x(\val(B)) = \val(A \cup B)$. 

Now let $y_1,\,y_2 \in X$ with $y_1 \le y_2$ and let $B_1,\,B_2 \in \fin$ with $y_1 = \val(B_1)$ and $y_1 = \val(B_1)$. Then by Lemma~\ref{lemma_iterators_199}
$B_1 \preceq B_2$. As usual we can choose $B_1$ and $B_2$ with $B_1 \subset B_2$  and also so that $A \cap B_2 = \varnothing$, where $A$ is again
such that $x = \val(A)$. Then $\pi_x(y_1) = \val(A \cup B_1)$ and $\pi_x(y_2) = \val(A \cup B_2)$ and hence by Lemma~\ref{lemma_iterators_199}
$\pi_x(y_1) \le \pi_x(y_2)$, since $A \cup B_1 \subset A \cup B_2)$. Suppose conversely that $\pi_x(y_1) \le \pi_x(y_2)$. Now either $y_1 \le y_2$ or
$y_2 \le y_1$ and if $y_2 \le y_1$ then $\pi_x(y_2) \le \pi_x(y_1)$ and hence $\pi_x(y_1) = \pi_x(y_2)$, in which case $y_1 = y_2$, since $\pi_x$ is injective. 
Thus in both cases $y_1 \le y_2$.

(4)\enskip We must show that $L = R$, where $L= \pi_x(J)$, $J = \{ z\in X : y_1 \le z < y_2 \}$ and $R = \{ z' \in X : \pi_x(y_1) \le z' < \pi(y_2) \}$. 
If $z \in J$ then by (3) $\pi_1(y_1) \le \pi_x(z) < \pi_x(y_2)$ and thus
$\pi_x(J) = L \subset R$. Let $z' \in R$. Now $\pi_x(X) = X_x$, since $\pi_x : X \to X_x$ is a bijection, and so $z' \in \pi_x(X)$, since $\pi_x(y_1) \in X_x$ and $z' \ge \pi_x(y_1)$.
Hence there exists $z\in X$ such that $z' = \pi_x(z)$ and then by (3) $z \in J$ and so
$z' \in L$. This shows that $L = R$. 

(5)\enskip By (4) $\pi_x(L_y) = \{ z' \in X : \pi_x(x_0) \le z' < \pi_x(y) \} = \{ z' \in X : x \le z' < \pi_x(y) \}$ 
and thus $L_{\pi_x(y)} =  L_x \cup\pi_x(L_y)$. It is clear that $L_x$ and  $\pi_x(L_y)$ are disjoint.    
\eop

\begin{theorem}\label{theorem_iterators_99}
Let $\,\mathbf{I} = (X,f,x_0)$ be a Peano iterator with $\le$ the unique total order
on $X$ compatible with $\,\mathbf{I}$. Define an iterator $\,\mathbf{I}_\le = (X_\le,f_\le,\varnothing)$
by letting $X_\le = \{L_x : x \in X\}$ and with $f_\le : X_\le \to X_\le$ given by $f_\le(L_x) = L_{f(x)}$ for all $x \in X$. Also define $\pi_\le : X \to X_\le$
by $\pi_\le(x) = L_x$ for all $x \in X$. 
Then:

(1)\enskip $\pi_\le : \mathbf{I} \to \mathbf{I}_\le$ is an isomorphism.

(2)\enskip $\,\mathbf{I}_\le$ is a Peano iterator.

(3)\enskip $f_\le(L_x)$ is the disjoint union of $\{x_0\}$ and $f(L_x)$ for each $x \in X$.

(4)\enskip The sets in $X_\le$ are totally ordered by inclusion and inclusion is the unique total order on $X_\le$ compatible with $\,\mathbf{I}_\le$.

(5)\enskip Let $\val_\le : \fin \to X_\le$ be the assignment of finite sets in $\,\mathbf{I}_\le$.
Then $\val_\le(A) \approx A$ for all $A \in \fin$.
 
\end{theorem}

\proof (1)\enskip  $\pi_\le(x_0) = \varnothing$ and $f_\le(\pi_\le(x)) = L_{f(x)} = \pi_\le(f(x))$ for all $x\in X$ which means that
$f_\le \circ \pi_\le = \pi_\le \circ f$. Thus $\pi_\le : \mathbf{I} \to \mathbf{I}_\le$ is a morphism. But  $\pi_\le : X \to X_\le$ is clearly a bijection
and hence by Lemma~\ref{lemma_iterators_16} $\pi_\le$ is an isomorphism.

(2)\enskip This follows immediately from (1) and Theorem~\ref{theorem_iterators_7}.

(3)\enskip This follows from Theorem~\ref{theorem_iterators_199} (4).

(4)\enskip This is clear. If $x,\,y \in X$ then either $x \le y$ or $y \le x$. If $x \le y$ then $L_x \subset L_y$ and if $y \le x$ then $L_y \subset L_x$.
Moreover, $L_x \subset L_{f(x)} = f_\le({L_x})$ and so inclusion is compatible with $\,\mathbf{I}_\le$.

(5)\enskip Let $A$ be a finite set and let $\mathcal{S} = \{ B \subset A : \val_\le(B) \approx B \}$ and so $\varnothing \in \mathcal{S}$, since
$\val_\le(\varnothing) = \varnothing$. Thus let $B \in \proper{\mathcal{S}}$, $a \in A \setminus B$ and put $B' = B \cup \{a\}$. Then 
 $\val_\le(B') = f_\le(\val_\le(B))$. If $\val_\le(B) = L_x$ 
then $\val_\le(B')  = f_\le(L_x) = L_{f(x)}$, which is the disjoint union of $L_x$ and the singleton set $\{x\}$.
Thus $\val_\le(B')$ is the disjoint union of $\val_\le(B)$ and a singleton set. 
But $\val_\le(B) \approx B$, since $B \in \mathcal{S}$ and hence $\val_\le(B') \approx B'$, i.e., $B' \in \mathcal{S}$. This shows that
$\mathcal{S}$ is an inductive $A$-system and so $\val_\le(A) \approx A$. \eop

We call $\,\mathbf{I}_\le$ the \definition{initial segment iterator associated with $\,\mathbf{I}$}.

If Theorem~\ref{theorem_iterators_99} is applied to the Peano iterator $(\Nat,\mathsf{s},0)$ then it is easy to see that
$L_n = \{0,1,\ldots, n-1\}$ for each $n \in \Nat$.
\medskip
 
Theorem~\ref{theorem_iterators_199} will now be employed  to obtain a result which is sometimes called \definition{The Definition by Induction Theorem}.

\begin{theorem}\label{theorem_iterators_11}
Let $\,\mathbf{I} = (X,f,x_0)$ be a Peano iterator, let $Z$ be a class, $\beta_0 \in Z$ and $\alpha : X \times Z \to Z$ be a mapping. 
Then there is a unique mapping $\pi : X \to Z$ with $\pi(x_0) = \beta_0$ such that $\pi(f(x)) = \alpha(x,\pi(x))$ for all $x \in X$. 
\end{theorem}

\proof The notation is as in Theorem~\ref{theorem_iterators_199}. Put $\Lambda_x = \{y \in X : y \le x \}$ for each $x \in X$. Then
Theorem~\ref{theorem_iterators_199} (4) implies that
$\Lambda_{x_0} =\{x_0\}$ and $f(\Lambda_x)$ is the disjoint union of $\Lambda_x$ and $\{f(x)\}$ for all $x \in X$.
Let $X_0$ be the subclass of $X$ consisting of those $x \in X$ for 
which there exist a unique mapping $\pi_x : \Lambda_x \to Z$ with 
$\pi_x(x_0) = \beta_0$ and
$\pi_x(f(x') = \alpha(x',\pi_x(x')$ for all 
$x' \in \Lambda_x$.
Since $\Lambda_{x_0} =\{x_0\}$ 
we have to define $\pi_{x_0} : \{x_0\} \to Z$ by $\pi_{x_0}(x_0) = \beta_0$, which 
shows that $x_0 \in X_0$. Thus let $x \in  X_0$ with unique mapping 
$\pi_x : \Lambda_x \to Z$. 
Since $f(\Lambda_x)$ is the disjoint union of $\Lambda_x$ and $\{f(x)\}$ we have to 
define $\pi_{f(x)} : f(\Lambda_x) \to Z$ by letting 
$\pi_{f(x)}(x') = \pi_x(x')$ if $x' \in \Lambda_x$ and 
$\pi_{f(x)}(f(x)) = \alpha(x,\pi_x(x))$.
It follows that $f(x) \in X_0$ and hence $X_0$ is an $f$-invariant subclass of $X$ containing $x_0$. 
Therefore $X_0 = X$.
Now define $\pi : X \to Z$ with by letting $\pi(x) = \pi_x(x)$ for all 
$x \in X$.
Then $\pi(x_0)= \beta_0$ and if $x \in X$ then
\begin{eqnarray*}
\pi(f(x)) =\pi_{f(x)}(f(x)) &=& \alpha(x,\pi_{\pi(x)}(x))\\
 &=& \alpha(x,\pi_x(x)) = \alpha(x,\pi(x))
\end{eqnarray*}
Finally, if $\pi'$ is another mapping satisfying the conditions of the theorem then it easy to see that 
$\{ x \in X : \pi(x) = \pi(x')\}$ is a 
$f$-invariant subclass of $X$ containing $x_0$
and hence $X_0 = X$, i.e., $\pi' = \pi$. Therefore the mapping $\pi$ is unique. \eop

The following is a more elaborate version of the previous theorem:

\begin{theorem}\label{theorem_iterators_111}
Let $\,\mathbf{I} = (X,f,x_0)$ be a Peano iterator, let $Y$ and $Z$ be classes and let 
$\beta : Y \to Z$
and $\alpha : X \times Y \times Z \to Z$ be mappings. Then there is a unique mapping
$\pi : X \times Y \to Z$ with $\pi(x_0,y) = \beta(y)$ for all $y \in Y$ such that
$\pi(f(x),y) = \alpha(x,y,\pi(x,y))$ 
for all $x \in X$, $y \in Y$.
\end{theorem}

\proof The notation is as in Theorem~\ref{theorem_iterators_199}. Put $\Lambda_x = \{y \in X : y \le x \}$ for each $x \in X$. Then
Theorem~\ref{theorem_iterators_199} (4) implies that
$\Lambda_{x_0} =\{x_0\}$ and $f(\Lambda_x)$ is the disjoint union of $\Lambda_x$ and $\{f(x)\}$ for all $x \in X$.
Let $X_0$ be the subclass of $X$ consisting of those $x \in X$ for 
which there exist a unique mapping $\pi_x : \Lambda_x \times Y \to Z$ with 
$\pi_x(x_0,y) = \beta(y)$ for all $y \in Y$ and
$\pi_x(f(x'),y) = \alpha(x',y,\pi_x(x',y))$ for all 
$x' \in \Lambda_x$, $y \in Y$. 
Since $\Lambda_{x_0} =\{x_0\}$ 
we have to define $\pi_{x_0} : \{x_0\} \times Y \to Z$ by $\pi_{x_0}(x_0,y) = \beta(y)$, which 
shows that $x_0 \in X_0$. Thus let $x \in  X_0$ with unique mapping 
$\pi_x : \Lambda_x \times Y \to Z$. 
Since $f(\Lambda_x)$ is the disjoint union of $\Lambda_x$ and $\{f(x)\}$ we have to 
define $\pi_{f(x)} : f(\Lambda_x) \times Y \to Z$ by letting 
$\pi_{f(x)}(x',y) = \pi_x(x',y)$ if $x' \in \Lambda_x$ and 
$\pi_{f(x)}(f(x),y) = \alpha(x,y,\pi_x(x,y))$ for all  $y \in Y$.
It follows that $f(x) \in X_0$ and hence $X_0$ is an $f$-invariant subclass of $X$ containing $x_0$. 
Therefore $X_0 = X$.
Now define $\pi : X \times Y \to Z$ with by letting $\pi(x,y) = \pi_x(x,y)$ for all 
$x \in X$, $y \in Y$.
Then $\pi(x_0,y)= \pi_{x_0}(x_0,y) = \beta(y)$ for all $y \in Y$ and if $x \in X$ then
\begin{eqnarray*}
\pi(f(x),y) =\pi_{f(x)}(f(x),y) &=& \alpha(x,y,\pi_{\pi(x)}(x,y))\\
 &=& \alpha(x,y,\pi_x(x,y)) = \alpha(x,y,\pi(x,y))
\end{eqnarray*}
 for all $y \in Y$. Finally, if $\pi'$ is another mapping satisfying the conditions of the theorem then it easy to see that 
$\{ x \in X : \pi(x,y) = \pi(x,y) \mbox{ for all $y \in Y$} \}$ is a 
$f$-invariant subclass of $X$ containing $x_0$
and hence $X_0 = X$, i.e., $\pi' = \pi$. Therefore the mapping $\pi$ is unique. \eop

We next look at the relationship between enumerators and iterators and in what follows let 
$\,\mathbf{I} = (X,f,x_0)$
 be a fixed iterator, where for simplicity we assume that the first component $X$ is a set. Let $\val$ be the assignment of finite sets in $\,\mathbf{I}$.
For the finite set $A$ the element $\val(A)$ can be thought of as the analogue of the cardinality of $A$
for the iterator $\,\mathbf{I}$.
Now the cardinality $|A|$ of a finite set $A$ can also be determined by counting or enumerating its elements
and the analogue of this procedure can also be carried out in the iterator $\,\mathbf{I}$.
To explain what this means let us return to the informal discussion presented at the beginning of Section~\ref{enums}. There
we determined whether $A$ is finite or not by marking the elements in $A$ one at a time and seeing if all
the elements can be marked in finitely many steps. At each stage of this process we took a snapshot of the elements which 
have already been marked, which resulted in the $A$-enumerator $\mathcal{U}$ whose elements are exactly the
snapshots. The set $A$ being finite meant that $A \in \mathcal{U}$.

Now suppose that each act of marking an element of $A$ is registered with the iterator $\,\mathbf{I}$. Each such act produces an element of $X$ 
which can be considered as the current state of the registering process.
Before the first element of $A$ has been marked the current state is $x_0$.
If at some stage the current state is $x$ then marking the next element of $A$ changes the current state to $f(x)$.

This registering process can be regarded as a mapping $\alpha_{\mathcal{U}} : \mathcal{U} \to X$, where $\alpha_{\mathcal{U}}(U)$ gives the current state
when the elements in the subset $U$ have been marked. The above interpretation then requires that $\alpha_{\mathcal{U}}(\varnothing) = x_0$ and
$\alpha_{\mathcal{U}}(\ssuc{\mathcal{U}}(U)) = f(\alpha_{\mathcal{U}}(U))$ for all $U \in \proper{\mathcal{U}}$. 
In Proposition~\ref{prop_iterators_501} it is shown that there is a unique mapping $\alpha_{\mathcal{U}}$ satisfying these requirements.

Since $A$ is finite  the registering process ends when all the elements in $A$ have been marked and the final current state is then the element 
$\alpha_{\mathcal{U}}(A)$ of $X$.
Now if the analogy with the iterator $(\Nat,\textsf{s},0)$ and the cardinality $|A|$ is valid then we would expect that
$\alpha_{\mathcal{U}}(A) = \val(A)$ holds for each finite set. In fact this does hold, as is shown in Theorem~\ref{theorem_iterators_502}.

Theorem~\ref{theorem_iterators_502} implies that
$\alpha_{\mathcal{U}}(A)$ does not depend on the $A$-enumerator $\mathcal{U}$. This is the fundamental reason why counting makes sense:
It does not matter in which order the elements in a finite set are counted; the same number always comes out in the end.

\begin{proposition}\label{prop_iterators_501}  
Let $\,\mathcal{U}$ be an $A$-enumerator.
Then there exists a unique mapping $\alpha_{\mathcal{U}} : \mathcal{U} \to X$ with $\alpha_{\mathcal{U}}(\varnothing) = x_0$ 
such that $\alpha_{\mathcal{U}}(\ssuc{\mathcal{U}}(U)) = f(\alpha_{\mathcal{U}}(U))$ for all $U \in \proper{\mathcal{U}}$. 
\end{proposition}

\proof
This is essentially the same as the proof of Theorem~\ref{theorem_enums_2}.
Let $\mathcal{U}_0$ denote the set consisting of those $U \in \mathcal{U}$ for which there exists a mapping $\alpha_U : \mathcal{U}_U \to X$
with $\alpha_U(\varnothing) = x_0$ and such that $\alpha_U(\ssuc{\mathcal{U}}(U')) = f(\alpha_U(U'))$ for all $U' \in \proper{\mathcal{U}_U}$. 
Clearly $\varnothing \in \mathcal{U}_0$ (since $\mathcal{U}_\varnothing = \{\varnothing\}$).

Consider $U \in \proper{\mathcal{U}_0}$ and let $\alpha_U : \mathcal{U}_U \to X$ be a mapping
with $\alpha_U(\varnothing) = x_0$ and such that $\alpha_U(\ssuc{\mathcal{U}}(U')) = f(\alpha_U(U'))$ for all $U' \in \proper{\mathcal{U}_U}$. 
Write $U^*$ for $\ssuc{\mathcal{U}}(U)$; by Lemma~\ref{lemma_enums_1}~(2) $\mathcal{U}_{U^*} = \mathcal{U}_U \cup \{U^*\}$ and so
we can define $\alpha_{U^*} : \mathcal{U}_{U^*} \to X$  by putting $\alpha_{U^*}(U') = \alpha_U(U')$ if $U' \in \mathcal{U}_U$ and letting
$\alpha_{U^*}(U^*) = f(\alpha_U(U))$.
If $U' \in \proper{\mathcal{U}_U}$ then $U' \in \mathcal{U}_U$ and $\ssuc{\mathcal{U}}(U') \in \mathcal{U}_U$ and thus
\[\alpha_{U^*}(\ssuc{\mathcal{U}}(U')) = \alpha_{U}(\ssuc{\mathcal{U}}(U')) = 
f(\ssuc{\mathcal{U}}(U')) = f(\alpha_{\mathcal{U^*}}(U'))\;.\]
Also $\alpha_{U^*}(\ssuc{\mathcal{U}}(U)) = \alpha_{U^*}(U^*) = f(\alpha_U(U)) = f(\alpha_{U^*}(U)))$ and
$\mathcal{U}_U = \proper{\mathcal{U}_{U^*}}$, and so $\alpha_{U^*}(\ssuc{\mathcal{U}}(U')) = f(\alpha_{U^*}(U'))$ 
for all $U' \in \proper{\mathcal{U}_{U^*}}$. Hence $\ssuc{\mathcal{U}}(U) = U^* \in \mathcal{U}_0$.
This shows $\mathcal{U}_0$ is an invariant subset of $\mathcal{U}$. Thus $\mathcal{U}_0 = \mathcal{U}$ and then by Lemma~\ref{lemma_enums_24}
$A \in \mathcal{U}_0$, which means there
exists a mapping $\alpha_{\mathcal{U}} : \mathcal{U} \to X$ with $\alpha_{\mathcal{U}}(\varnothing) = x_0$ and
such that $\alpha_{\mathcal{U}}(\ssuc{\mathcal{U}}(U)) = f(\alpha_{\mathcal{U}}(U))$ for all $U \in \proper{\mathcal{U}}$. 

It remains to consider the uniqueness. Let $\alpha'_{\mathcal{U}} : \mathcal{U} \to X$ be any mapping with
$\alpha'_{\mathcal{U}}(\varnothing) = x_0$ and
such that $\alpha'_{\mathcal{U}}(\ssuc{\mathcal{U}}(U)) = f(\alpha'_{\mathcal{U}}(U))$ for all $U \in \proper{\mathcal{U}}$ and consider the set
$\mathcal{U}_0 = \{ U \in \mathcal{U} : \alpha'_{\mathcal{U}}(U) = \alpha_{\mathcal{U}}(U) \}$. Clearly $\varnothing \in \mathcal{U}_0$ and
if $U \in \proper{\mathcal{U}_0}$ then
$\alpha'_{\mathcal{U}}(\ssuc{\mathcal{U}}(U)) = f(\alpha'_{\mathcal{U}}(U)) = f(\alpha_{\mathcal{U}}(U)) = \alpha_{\mathcal{U}}(\ssuc{\mathcal{U}}(U))$, i.e.,
$\ssuc{\mathcal{U}}(U) \in \mathcal{U}_0$. Thus $\,\mathcal{U}_0$ is an invariant subset of $\,\mathcal{U}$, and so $\,\mathcal{U}_0 = \mathcal{U}$. 
This shows that $\alpha'_{\mathcal{U}} = \alpha_{\mathcal{U}}$.
\eop

The mapping $\alpha_{\mathcal{U}} : \mathcal{U} \to X$ in Proposition~\ref{prop_iterators_501} will be referred to as the
\definition{$\mathcal{U}$-valuation in $\,\mathbf{I}$}, or just as the \definition{$\mathcal{U}$-valuation} if it is clear which iterator
is involved. 

The uniqueness of the $\mathcal{U}$-valuation implies that for each $U \in \mathcal{U}$ the $\mathcal{U}_U$-valuation 
$\alpha_{\mathcal{U}_U}$ is just the restriction of $\alpha_{\mathcal{U}}$ 
to $\mathcal{U}_U$, i.e., $\alpha_{\mathcal{U}_U} : \mathcal{U}_U \to X$ is the mapping with $\alpha_{\mathcal{U}_U}(U') = \alpha_{\mathcal{U}}(U')$ for all 
$U' \in \mathcal{U}_U$.

\begin{proposition}\label{prop_iterators_502}  
Let $\mathcal{U}$ be an $A$-enumerator and $\mathcal{V}$ a $B$-enumerator.
Suppose $A \preceq B$ and so by Proposition~\ref{prop_enums_3} there exists a unique homomorphism $\pi : \mathcal{U} \to \mathcal{V}$.
Then $\alpha_{\mathcal{U}} = \alpha_{\mathcal{V}} \circ \pi$ (with $\alpha_{\mathcal{U}}$ the $\mathcal{U}$-valuation and $\alpha_{\mathcal{V}}$ the 
$\mathcal{V}$-valuation).
\end{proposition}

\proof 
The mapping $\alpha_{\mathcal{V}} \circ \pi : \mathcal{U} \to X$ is also a $\mathcal{U}$-valuation since
\[(\alpha_{\mathcal{V}} \circ \pi)(\ssuc{\mathcal{U}}(U)) = \alpha_{\mathcal{V}}(\pi(\ssuc{\mathcal{U}}(U))) 
= \alpha_{\mathcal{V}}(\ssuc{\mathcal{V}}(\pi(U))) = f(\alpha_{\mathcal{V}}(\pi(U))) = f((\alpha_{\mathcal{V}} \circ \pi)(U))\]
for all $U \in \proper{\mathcal{U}}$,
and $(\alpha_{\mathcal{V}} \circ \pi)(\varnothing) = \alpha_{\mathcal{V}}(\pi(\varnothing)) = \alpha_{\mathcal{V}}(\varnothing) = \varnothing$.
By the uniqueness of the $\mathcal{U}$-valuation it therefore follows that $\alpha_{\mathcal{U}} = \alpha_{\mathcal{V}} \circ \pi$.  
\eop

\begin{proposition}\label{prop_iterators_503}
If $A \approx B$ then $\alpha_{\mathcal{U}}(A) = \alpha_{\mathcal{V}}(B)$
for each $A$-enumerator $\,\mathcal{U}$ and each $B$-enumerator $\,\mathcal{V}$.
In particular, $\alpha_{\mathcal{U}}(A) = \alpha_{\mathcal{U}'}(A)$ for all $A$-enumerators $\,\mathcal{U}$ and $\,\mathcal{U}'$.
\end{proposition}

\proof
By Theorem~\ref{theorem_enums_2} there exists a unique homomorphism $\pi : \mathcal{U} \to \mathcal{V}$ which maps $\mathcal{U}$ bijectively onto 
$\mathcal{V}$, and  $\pi(A) = B$. Also, by Proposition~\ref{prop_iterators_502}  
$\alpha_{\mathcal{U}} = \alpha_{\mathcal{V}} \circ \pi$ and therefore
$\alpha_{\mathcal{U}}(A) = \alpha_{\mathcal{V}}(\pi(A)) = \alpha_{\mathcal{V}}(B)$. 
\eop

\begin{theorem}\label{theorem_iterators_502}  
For each $A$-enumerator $\,\mathcal{U}$ we have $\alpha_{\mathcal{U}}(A) = \val(A)$.
\end{theorem}

\proof 
If $\,\mathcal{U}$ and $\,\mathcal{U}'$ are $A$-enumerators then by Proposition~\ref{prop_iterators_503}
$\alpha_{\mathcal{U}}(A) = \alpha_{\mathcal{U}'}(A)$ and hence the element $\alpha_{\mathcal{U}}(A)$ of $X$ does not depend on which $A$-enumerator
$\,\mathcal{U}$ is used. Denote this element by $\val_*(A)$. By the uniqueness of the valuation it is enough to show that the
assignment $A \mapsto \val_*(A)$ is a valuation.
Now if $\mathcal{U}^\varnothing$ is the unique $\varnothing$-enumerator then $\mathcal{U}^\varnothing(\varnothing) = x_0$ and so
$\val_*(\varnothing) = x_0$. Thus consider a finite set $A$ and let $a$ be an element with $a \notin A$;
put $A' = A \cup \{a\}$. By Lemma~\ref{lemma_enums_3} there exists an $A'$-enumerator $\mathcal{U}'$ with $A \in \mathcal{U}'$ and then
$A' = \ssuc{\mathcal{U}'}(A)$, since $a$ is the only element in $A' \setminus A$. Moreover $\mathcal{U}'_A$ is an $A$-enumerator and $\alpha_{\mathcal{U}'_A}$ 
is the restriction of $\alpha_{\mathcal{U}'}$ to $\mathcal{U}'_A$ and so
\[\val_*(A') = \alpha_{\mathcal{U}'}(A') = \alpha_{\mathcal{U}'}(\ssuc{\mathcal{U}'}(A)) =  f(\alpha_{\mathcal{U}'}(A)) =  f(\alpha_{\mathcal{U}'_A}(A))
= f(\val_*(A)).\]
Hence $\val_*$ is a valuation and therefore $\val_* = \val$, i.e., $\alpha_{\mathcal{U}}(A) = \val(A)$ for each finite set $A$ and each $A$-enumerator $\,\mathcal{U}$.
\eop

Consider the case with $A \approx B$.
Let $\,\mathcal{U}$ be an $A$-enumerator and $\,\mathcal{V}$ a $B$-enumerator. 
By Proposition~\ref{prop_iterators_503} $\alpha_{\mathcal{U}}(A) = \alpha_{\mathcal{V}}(B)$
and hence by Theorem~\ref{theorem_iterators_502} 
\[\val(A) = \alpha_{\mathcal{U}}(A) = \alpha_{\mathcal{V}}(B) = \val(B)\ .\] 

We end the section with a couple of remarks about the case when $\mathbf{I} = (X,f,x_0)$ is a Peano iterator with $X$ an infinite set.
A set $E$ is defined to be \definition{Dedekind-infinite} if there exists an injective mapping $h : E \to E$ which is not surjective,
and so by Theorem~\ref{theorem_fs_1} a Dedekind-infinite set is infinite,i.e., it is not finite. The converse also holds (i.e., every infinite set is 
Dedekind-infinite) provided a suitable form of the axiom of choice is assumed. In models without the axiom of choice there can exist infinite sets which are
Dedekind-finite.
If $(X,f,x_0)$ is a Peano iterator then the set $X$ is Dedekind-infinite. Conversely, if $E$ is a Dedekind-infinite set and
$h : E \to E$ is injective but not surjective and $e_0 \in E \setminus h(E)$ then 
$(E_0,y_0,e_0)$ is a Peano iterator, where $E_0$ is the least $h$-invariant subset of $E$ containing $e_0$ 
and $y_0 :E_0 \to E_0$ is the restriction of $h$ to $E_0$.
Thus a Peano iterator whose first component is a set exists if and only if there exists a Dedekind-infinite set.

The following somewhat strange result can be found in \cite{wr} and in \cite{tarski}.

\begin{proposition}\label{prop_iterators_11}
Suppose that there exists a Dedekind-infinite set (for example, the set of natural numbers $\Nat$ is a such a set). Then a set $A$ is finite if and only if 
$\mathcal{P}(\mathcal{P}(A))$ is Dedekind-finite. 
 \end{proposition}

\proof Let $ \mathbf{I} = (X,f,x_0)$ be a Peano iterator (which exists since a Dedekind-infinite set exists)
and let $\val$ be the assignment of finite sets in $\mathbf{I}$. Now let $E$ be an infinite set and consider the mapping $h : \mathcal{P}(E) \to X$ given by 
$h(A) = \val(A)$  if $A$ is finite, and $h(Y) = x_0$ if $Y$ is infinite. If $A$ is finite then by Proposition~\ref{prop_fs_88} there exists a finite subset $C$ of 
$E$ with $C \approx A$ and thus by  Lemma~\ref{lemma_iterators_4} $h$ is surjective. Hence there exists an injective 
mapping $g : \mathcal{P}(X) \to \mathcal{P}(\mathcal{P}(E))$. (If $f : Y \to Z$ is any surjective  mapping then the mapping $g : \mathcal{P}(Z) \to \mathcal{P}(Y)$
given by $g(E) =f^{-1}(E)$ for each $E\in \mathcal{P}(Z)$  is injective.)
There is  then an injective mapping $\alpha : X \to \mathcal{P}(\mathcal{P}(E))$ given by
$\alpha(x) = g(\{x\})$ for all $x \in X$. Therefore by Lemma~\ref{lemma_iterates_11} below $\mathcal{P}(\mathcal{P}(E))$ is Dedekind-infinite. On the other hand,
if $E$ is finite then by Proposition~\ref{prop_fs_3} $\mathcal{P}(E)$ and hence $\mathcal{P}(\mathcal{P}(E))$ is finite, and so by Theorem~\ref{theorem_fs_1}
$\mathcal{P}(\mathcal{P}(E))$ is Dedekind-finite. Thus if there exists a
Dedekind-infinite set then a set $A$ is finite if and only if $\mathcal{P}(\mathcal{P}(A))$ is Dedekind-finite. \eop

Note that the assumption about the Dedekind-infinite set is only needed to show that if $Y$ is infinite then $\mathcal{P}(\mathcal{P}(Y))$ is Dedekind-infinite. 
But if $Y$ is infinite then in fact it can be shown that the class $O$ of finite ordinals is actually a set (and thus a Dedekind-infinite set)
and so the hypothesis is not actually needed.

\begin{lemma}\label{lemma_iterates_11}
A set containing a Dedekind-infinite set is itself Dedekind-infinite.
\end{lemma}

\proof Let $E$ contain a Dedekind-infinite set $F$ and so there exists an injective mapping $h : F \to F$ which is not surjective. Define $g : E \to E$ by letting 
$g(x) = h(x)$ if $x \in F$ and $g(x) = x $ if $x\in E \setminus F$. Then $g$ is injective but not surjective and hence $E$ is Dedekind-infinite. \eop


\startsection{Finite ordinals}

\label{ordinals}

In this section  we study finite ordinals using the standard approach introduced by von Neumann \cite{vonn}. As in Section~\ref{iterators} 
we denote the class of all finite sets by $\,\fin$. Let $\vonn : \fin \to \fin$ be the mapping given by $\vonn(A) = A \cup \{A\}$ 
for each finite set $A$. Note that $\vonn(A)$ is either equal to $A$ or to the disjoint union of $A$ and the singleton set $\{A\}$ and in the latter case
we say that $A$ is \definition{$\sigma$-regular}.

\begin{lemma}\label{lemma_ordinal_31}
(1)\enskip If $A \in \fin$ is transitive then so is $\vonn(A)$. (Recall that a set $B$ is transitive if $b \subset B$  for all $b \in B$.)

(2)\enskip Let $A \in \fin$; then $A$ is $\sigma$-regular if and only if $a \ne A$ for all $a \in A$.

(3)\enskip Let $A \in\fin$ be transitive. Then $A$ is $\sigma$-regular if and only if $a$ is a proper subset of $A$ for all $a \in A$.

(4) \enskip If $A \in \fin$ is transitive and $\sigma$-regular then $\vonn(A)$ is $\sigma$-regular.
 \end{lemma}

\proof
(1) \enskip Let $A$ be transitive and $b \in \vonn(A)$. Then either $b \in A$, in which case $b \subset A \subset \vonn(A)$,
or $b = A$, in which case $b \subset \vonn(A)$. Thus $\vonn(A)$ is transitive.  

(2)\enskip Let $A \in \fin$; then $A$ is not $\sigma$-regular if and only if
$\vonn(A) = A \cup \{A\} = A$ and this holds if and only if $\{A\} \subset A$, which in turn holds if and only if 
$A \in A$. Also $A \in A$ if and only if there exists $a \in A$ with $a = A$.
Thus $A$ is $\sigma$-regular if and only if $a \ne A$ for all $a \in  A$.

(3)\enskip This follows from (2), since if $A$ is transitive then $a \ne A$ if and only $a$ is a proper subset of $A$.

(4)\enskip Let $A \in \fin$ be transitive and $\sigma$-regular. Then by (1) $\vonn(A)$ is transitive and by (3) $a$ is a proper subset of $A$ for each $a \in A$. Thus
$a$ is also a proper subset of $\vonn(A)$ for all $a \in A$ and $A$ is a proper subset of $\vonn(A)$, since $A$ is $\sigma$-regular. Hence $a'$ is a proper subset of
$\vonn(A)$ for all $a' \in \vonn(A)$ and so by (3) $\vonn(A)$ is $\sigma$-regular. \eop

 If we iterate the operation $\vonn$ starting with the empty set and label the resulting sets using the natural numbers then we obtain the following:

 $0 =\varnothing$,\\
 $1 = \vonn(0) = 0 \cup\{0\} = \varnothing \cup \{0\} = \{0\}$,\\ 
 $2 = \vonn(1) = 1 \cup \{1\} = \{0\} \cup \{1\} = \{0,1\}$,\\
 $3 = \vonn(2) = 2 \cup\{2\} = \{0,1\} \cup \{2\} = \{0,1,2\}$,\\
 $4 = \vonn(3) = 3 \cup\{3\} = \{0,1,2\} \cup \{3\} = \{0,1,2,3\}$,\\
 $5 = \vonn(4) = 4 \cup \{4\} = \{0,1,2,3\} \cup \{4\} = \{0,1,2,3,4\}$,

 $n+1 =  \vonn(n) = n \cup \{n\} = \{0,1,2,\ldots,n-1\} \cup\{n\} = \{0,1,2,\ldots,n\}$ .
 \bigskip

By Lemma~\ref{lemma_ordinal_31} (1) and (4) the sets $\vonn(n)$, $n \in \Nat$, are all transitive and $\sigma$-regular.

 Denote by $\mathbf{O}'$ the iterator $(\fin,\vonn,\varnothing)$.
 Then by Theorem~\ref{theorem_iterators_1} there exists a unique assignment $\ordass$
 of finite sets in $\mathbf{O}'$. Thus $\ordass : \fin \to \fin$  is the unique mapping with 
 $\ordass(\varnothing) = \varnothing$ and such that
 \[\ordass(A \cup \{a\}) = \vonn(\ordass(A))\] 
 for each finite set $A$ and each element $a \notin A$.
 Moreover, if $A$ and $B$ are finite sets with $A \approx B$ then $\ordass(A) = \ordass(B)$.

 \begin{theorem}\label{theorem_ordinal_1}
 For each finite set $A$ we have $\ordass(A) \approx A$, and thus  $\ordass(A) = \ordass(B)$ if and only if $A \approx B$. Moreover,
$\ordass(A)$ is transitive  and $\sigma$-regular for each finite set $A$.
 \end{theorem}

\proof Let $A$ be a finite set and let 
\[\mathcal{S} = \{ B \subset A : \mbox{$ \ordass(B)$ is transitive and $\sigma$-regular with $\ordass(B) \approx B$} \}\;.\] 
In particular $\varnothing \in \mathcal{S}$, since 
 $\ordass(\varnothing) = \varnothing$. Let $B \in \proper{\mathcal{S}}$, $a \in A \setminus B$ and put $B' = B \cup \{a\}$.
 Then $\ordass(B)$ is transitive and $\sigma$-regular with $\ordass(B) \approx B$ and $\ordass(B') = \vonn(\ordass(B))$. Hence by Lemma~\ref{lemma_ordinal_31} $\ordass(B')$ is transitive and $\sigma$-regular 
and $\ordass(B') =\vonn(\ordass(B)) \approx B'$. Thus $B' \in \mathcal{S}$  and so $\mathcal{S}$ is an inductive $A$-system.
 Therefore $A \in \mathcal{S}$, i.e., $\ordass(A)$ is transitive and $\sigma$-regular with  $\ordass(A) \approx A$. \eop

It follows from Theorem~\ref{theorem_ordinal_1} that $\ordass(\ordass(A)) = \ordass(A)$ for each finite set $A$.

Let $\,O = \{ B \in \fin : \mbox{$B = \ordass(A)$ for some finite set $A$} \}$. The elements of $O$ will be called 
\definition{finite ordinals}. Thus for  each finite set $A$  there exists  a unique  $o \in O$  with $o = \ordass(A)$.
By Lemma~\ref{lemma_iterators_3} $O$ is the least $\vonn$-invariant subclass of $\fin$ containing $\varnothing$.
We denote the restriction of $\vonn$ to a mapping $O \to O$  again by $\vonn$.
Thus $\mathbf{O} = (O,\vonn,\varnothing)$ is a minimal iterator.

We do not assume that $O$ is a set but if it were then it would not be 
finite. (If $O$ were finite then $o' = \vonn(\ordass(O))$ would be a finite ordinal, but $o' \not\approx o$ for each $o \in O$.)   

\begin{theorem}\label{theorem_ordinal_2}
$\,\mathbf{O}$ is a Peano iterator. Therefore by Theorem~\ref{theorem_iterators_3} (the recursion theorem) it follows that
for each iterator $ \mathbf{J} = (H,\delta,h_0)$ there exists a unique mapping $\pi : O \to H$ with $\pi(\varnothing) = h_0$ such that $\pi \circ \vonn = \delta \circ \pi$.
\end{theorem}

\proof  By Theorem~\ref{theorem_ordinal_1} $\,\mathbf{O}'$ is regular and hence also $\,\mathbf{O}$ is regular. Therefore by Theorem~\ref{theorem_iterators_2}
$\,\mathbf{O}$ is a Peano iterator. 
\eop

The iterator $\mathbf{O}$ is obtained only making use of finite sets. If the negation of the axiom of infinity  is assumed then
$(\Nat,\mathsf{s},0)$ does not exist. However, the Peano iterator 
$\mathbf{O}$ does exist and in this case $O$ is not a set. 
If $(\Nat,\mathsf{s},0)$ does exist then by Lemma~\ref{lemma_iterators_6} it is isomorphic to $\,\mathbf{O}$. Thus $\,\mathbf{O}$ can be considered
as a particular version of $(\Nat,\mathsf{s},0)$ and we can use the usual notation for the elements of $\Nat$ to denote the elements of $O$.  
In Section~\ref{am} we show how the arithmetic operations of addition, multiplication and exponentiation can be introduced for any minimal iterator. In
particular, this can applied to the iterator $\,\mathbf{O}$.

\begin{proposition}\label{prop_ordinal_1}
For each finite set $A$
\[\ordass(A) = \{ o \in O : o = \ordass(A') \mbox{ for some proper subset $A'$ of $A$} \}.\]
\end{proposition}

\proof Let $A$ be a finite  set and denote by $\mathcal{S}$ the set of subsets $B$ of $A$  for which
\[\ordass(B) = \{ o \in O: o = \ordass(B') \mbox{ for some proper subset $B'$ of $B$} \}.\]

In particular $\varnothing \in \mathcal{S}$. Thus consider $B \in \proper{\mathcal{S}}$
and $a \in A \setminus B$. Then 
\[\ordass(B\cup\{a\}) = \vonn(\ordass(B)) = \ordass(B) \cup \{\ordass(B)\} = \{ o \in O : o = \ordass(B') \mbox{ for some $B' \subset B$} \}.\]
But if $C$ is a proper subset of $B \cup \{a\}$ then $C \approx C'$ for some $C' \subset C$ and then by 
Theorem~\ref{theorem_iterators_1} $\ordass(C) = \ordass(C')$. It follows that
\[\ordass(B \cup\{a\}) = \{ o \in O : o = \ordass(B') \mbox{ for some proper subset $B'$ of $B \cup \{a\}$} \}\]
and thus $B \cup \{a\} \in \mathcal{S}$. Hence $\mathcal{S}$ is an inductive $A$-system and so $A \in \mathcal{S}$.  Therefore
\[\ordass(A) = \{ o \in O : o = \ordass(A') \mbox{ for some proper subset $A'$ of $A$} \}.\ \eop\]
\begin{proposition}\label{prop_ordinal_2}
Let $A$ be a finite set and $B\subset A$; then $\ordass(B) \subset \ordass(A)$.
\end{proposition}

\proof This follows immediately from Proposition~\ref{prop_ordinal_1}. \eop

Let $o \in O$ and let $a$ be an element not in $o$. Then $\ordass(o \cup \{a\}) = \vonn(o)$.

\begin{proposition}\label{prop_ordinal_3}

For each $o \in O$ we have
\[o = \{ o' \in O : \mbox{$o'$ is a proper subset of $o$} \},\]
\[\vonn(o) = \{ o' \in O :\mbox{ $o'$ is a subset of $o$} \}.\]
\end{proposition}

\proof This follows from Proposition~\ref{prop_ordinal_1}. \eop

A set $E$ is said to be \definition{totally ordered with respect to set membership} if, whenever $e_1$ and $e_2$ are distinct elements of $E$ then exactly one of 
$e_1\in e_2$ and $e_2\in e_1$ holds.  By Proposition~\ref{prop_ordinal_3} each finite ordinal $\alpha$ is totally ordered with respect to set inclusion and, moreover,
each element of $\alpha$ is a subset of $\alpha$ and so $\alpha$ is a transitive set.

In the general (non-finite case) the usual definition of an ordinal is as a set having these two properties \cite{vonn}.

\begin{proposition}\label{prop_ordinal_11}
Let $o,\, o'\in O$ with $o \ne o'$.  Then either $o$ is a proper subset of $o'$ or $o'$ is a proper subset of
$o$.
\end{proposition}

\proof By Theorem~\ref{theorem_ordinal_1} $o \not\approx o'$ and so by Proposition~\ref{prop_fs_112} (2) there either exists a proper subset $B$ of $o$
with $B \approx o'$ or there exists a proper subset $B'$ of $o'$ with $B' \approx o$. Suppose the former holds. Then $\ordass(B) = o'$
and therefore by Proposition~\ref{prop_ordinal_1}
\[o' = \ordass(B) = \{ b : b = \ordass(B') \mbox{ for some proper subset $B'$ of $B$} \},\]
\[o = \ordass(o) = \{ b : b = \ordass(B') \mbox{ for some proper subset $B'$ of $o$} \}.\]

It follows that $o'$ is a proper subset of $o$. If the latter holds then, in the same way, $o$ is a proper subset of $o$'. \eop

If $o, \,o' \in O$ then we write $o' \le o$ if $o' \subset o$. By Proposition~\ref{prop_ordinal_11} $\le$
defines a total order on $O$. The total order $\le$ is compatible with $\,\mathbf{O}_0$, since $o \subset \vonn(o)$  for all $o \in O$. Thus by
Theorem~\ref{theorem_iterators_199} $\le$ is the unique total order on $O$ compatible with $\,\mathbf{O}_0$. In particular, it follows from
Theorem~\ref{theorem_iterators_199} (5) that each non-empty subclass $O'$ of $O$ contains a minimum element, i.e., an element $o$ with $o \le o'$ for all 
$o' \in O'$.

The following induction principle for finite ordinals corresponds to Theorem~\ref{theorem_fs_3}.

\begin{proposition}\label{prop_ordinal_4}
Let $\,\prop$ be a statement about finite ordinals. Suppose
$\,\prop(0)$ holds 
and that $\,\prop\vonn(o))$ holds whenever $\,\prop(o)$ holds for  $o \in O$.
Then $\,\prop$ is a property of finite ordinals, i.e., $\prop(o)$ holds for every  $o \in O$.
\end{proposition}

\proof Let $A$ be a finite set and put $\mathcal{S} = \{ B \in \mathcal{P}(A) : \mbox{$\,\prop(\ordass(B))$ holds } \}$.
Then  $\,\prop(\varnothing) = \,\prop(0)$ holds, so let $B \in \proper{\mathcal{S}}$ and $a \in A \setminus B$.
Then $\ordass(b \cup \{a\}) =  \vonn(\ordass(B))$ and thus $B \cup \{a\} \in \mathcal{S}$. Therefore\ $\mathcal{S}$ is an inductive $A$-system and so
$A \in \mathcal{S}$. Let $o \in O$; then $o$ is a finite set and hence applying the above with $A = o$ shows that
$\,\prop(o)$ holds. \eop
\medskip

We next look at a further Peano iterator $\,\mathbf{U}$. This has nothing to do with ordinals, except that
the iterator $\,\mathbf{O}$ is involved in its definition. What it has in common with the iterator $\,\mathbf{O}$ and the iterator
$\,\mathbf{H}$ introduced in Section~\ref{iterators}  is that it is defined 
'absolutely' and is in fact constructed solely from operations performed on the empty set $\varnothing$. First consider the iterator 
$\,\mathbf{U}' = (\fin,\alpha,\varnothing)$, where  the  mapping $\alpha : \fin \to \fin$ is given by 
$\alpha'(A) = \{A\}$ for each finite set $A$. The recursion theorem  for the Peano iterator $\,\mathbf{O}$ applied to the iterator 
$\,\mathbf{U}'$ implies 
there exists a unique morphism  $\eta : \,\mathbf{O} \to \mathbf{U}'$. Thus $\eta : O \to \fin$ is the unique mapping  with
$\eta(0) = \varnothing$ such that $\eta(\vonn(o)) = \alpha'(\eta(o))$ for each $o \in O$. 
If we iterate the operation $\alpha'$ starting with the empty set then we obtain the following:

$\eta(0) = \varnothing$,\\
$\eta(1) = \alpha(\eta(0)) = \{\varnothing\}$,\\
$\eta(2) = \alpha(\eta(1)) =\{\{\varnothing\}\}$,\\
$\eta(3) = \alpha(\eta(2)) = \{\{\{ \varnothing\}\}\} $,\\
$\eta(4) = \alpha(\eta(3)) = \{\{\{\{\varnothing \}\}\}\}\}$,\\ 
$\vdots$
$\eta(n) = \{\{\{\{\{\{ \cdots \{ \varnothing\} \cdots\}\}\}\}\}\}$.

\medskip

 Let $\varphi$ be the assignment of finite sets in $\,\mathbf{U}'$, 
thus $\varphi : \fin \to \fin$ is the unique mapping with $\varphi(\varnothing) = \varnothing$ such that $\varphi(A \cup \{a\}) = \alpha(\varphi(A))$ for each
finite set $A$ and each $a \notin A$, and by Theorem~\ref{theorem_iterators_1} $\varphi(A) = \varphi(B)$ whenever $A \approx B$.
Let $U = \{ v \in \fin : \mbox{$v = \varphi(A)$ for some finite set $A$} \}$, so by Lemma~\ref{lemma_iterators_3} $U$ is  the least
$\alpha'$-invariant subclass of $\fin$ containing $\varnothing$; 
the restriction of $\alpha$ to a mapping
$U \to U$  will again be denoted by $\alpha$. Thus $\,\mathbf{U} = (U, \alpha, \varnothing)$ is a minimal iterator.

\begin{lemma}\label{lemma_ordinal_14} For each finite set $A$ we have $\varphi(A) = \eta(\ordass(A))$.
\end{lemma}

\proof
Let $A$ be a finite set and let $\mathcal{S} = \{ B \in \mathcal{P}(B) : \varphi(B) = \eta(\ordass(B)) \}$ and in particular $\varnothing \in \mathcal{S}$,
since $\varphi(\varnothing) = \eta(\ordass(\varnothing)) = \varnothing$. Thus let $B \in \mathcal{S}$ and $b \in A \setminus B$. Then
$\varphi(B \cup\{b\}) =\alpha(\varphi(B)) = \alpha(\eta(\ordass(B)))= \eta(\vonn(\ordass(B))) =\eta(\ordass(B \cup\{b\})$  and hence $B \cup\{b\} \in \mathcal{S}$
.This shows that $\mathcal{S}$ is an inductive $A$-system and therefore $A \in \mathcal{S}$, i.e., $\varphi(A) = \eta(\ordass(A))$. \eop

\begin{lemma}\label{lemma_ordinal_5} The mapping $\eta : O \to \fin$ maps $O$ bijectively onto $U$.
\end{lemma}

\proof By Lemma~\ref{lemma_iterators_15} $\eta(O) =U$ and so it remains to show that $\eta$ is injective. Suppose this is not the case and let
\[O_0 = \{ o \in O : \mbox{there exists $o' \in O$ with $o < o'$ and $\eta(o) = \eta(o')$} \}\;.\]
Thus $O_0$ is non-empty and so it contains a minimum element $o_0$ with $o_0 < o$ for all $o \in O_0$, and since
$o_0 \in O_0$ there exists $o_1 \in O_0$ with $o_0 < o_1$ and $\eta(o_0) = \eta(o_1)$. Now if $p \in O \setminus \{0\}$ then by
Proposition~\ref{prop_iterators_1} there exists a unique $q \in O$ with $p = \vonn(q)$ and then  $\eta(p) = \eta(\vonn(q)) = \alpha(\eta(q)) = \{\eta(q)\}$.
Thus if $p \in O \setminus \{0\}$ then there exists a unique $q \in O$ with $\eta(p) = \{\eta(q)\}$. In particular, $\eta(p) \ne \varnothing$ if 
$p \ne 0$ and so $o_1 \ne 0$, since $\eta(0) = \varnothing$ and $o_2 \ne 0$.
There thus exist unique $q_1,\,q_2 \in O$ with $\eta(o_1) = \{\eta(q_1)\}$ and $\eta(o_2) = \{\eta(q_2)\}$. Then $\{\eta(q_1)\} = \{\eta(q_2)\}$ and so 
$\eta(q_1) = \eta(q_2)$. But $q_1 < o_1$ and $q_1 < q_2$, which contradicts the minimality of $o_1$. Therefore $\eta$ is injective. \eop

\begin{proposition}\label{prop_ordinal_15}
The iterator $\,\mathbf{U}$ is a Peano iterator.
\end{proposition}

\proof By Lemma~\ref{lemma_ordinal_5} and Lemma~\ref{lemma_iterators_16} the morphism $\eta$ is an isomorphism and therefore by Theorem~\ref{theorem_iterators_7}
$\,\mathbf{U}$ is a Peano iterator. \eop

Except for being a Peano iterator the iterator $\,\mathbf{U}$ has none of the properties enjoyed by $\,\mathbf{O}$. It corresponds to the perhaps most primitive
method of counting by representing the number $n$ with something like $n$ marks, in this case the empty set enclosed in $n$ braces.
 
We can improve the situation somewhat by considering the finite section iterator $\,\mathbf{U}_\le = (U_\le, \alpha_\le,\varnothing)$ associated with $\,\mathbf{U}$
which was introduced in Theorem~\ref{theorem_iterators_99}. Here $\le$ is the unique total order on $U$ compatible with $\,\mathbf{U}$.
 
Thus $U_\le = \{ L_B : B \in U \}$ and $\alpha_\le(L_B) = L_{\alpha(B)}$ for all $B \in U$.
Define $\pi_\le : U \to U_\le$
by $\pi_\le(B) = L_B$ for all $B \in U$. 

Then Theorem~\ref{theorem_iterators_99} states that

(1)\enskip $\pi_\le : \mathbf{U} \to \mathbf{U}_\le$ is an isomorphism.

(2)\enskip $\,\mathbf{U}_\le$ is a Peano iterator.

(3)\enskip $f_\le(L_B)$ is the disjoint union of $L_B$ and the singleton set $\{B\}$ for each $B \in U$.
Moreover, $f_\le(L_B)$ is also the disjoint union of $\{x_0\}$ and $f(L_B)$ for each $B \in U$.

(4)\enskip The sets in $U_\le$ are totally ordered by inclusion and inclusion is the unique total order on $U_\le$ compatible with $\,\mathbf{U}_\le$.

(5)\enskip Let $\val_\le : \fin \to U_\le$ be the assignment of finite sets in $\,\mathbf{U}_\le$.
Then $\val_\le(A) \approx A$ for all $A \in \fin$.

\medskip

We end the section by considering a situation which is somewhat more general than that occurring with the iterator $\,\mathbf{O}$.
Let $\,\mathbf{J} = (T,h,\varnothing)$ be a minimal iterator with $T$ a subclass of $\fin$. We call $\,\mathbf{J}$ an \definition{ordinal iterator}
if for each $B \in T$ there exists an element $b \notin B$ such that $h(B) = B \cup \{b\}$. Note that if it is not assumed that $\,\mathbf{J}$ is minimal
then the associated minimal iterator $\,\mathbf{J}_0$ will be an ordinal iterator. The archetypal example of an ordinal iterator is of course 
$\,\mathbf{O}$. 
Moreover, if $\,\mathbf{I}$ is a Peano iterator and $\,\mathbf{I}_\le$ is the initial segment iterator associated with $\,\mathbf{I}$
then Theorem~\ref{theorem_iterators_199} (4) shows that  $\,\mathbf{I}_\le$ 
will be an ordinal iterator.
In what follows 
let $\,\mathbf{J} = (T,h,\varnothing)$ be an ordinal iterator and let $\tau : \fin \to T$ be the evaluation of finite sets in $\,\mathbf{J}$. 

\begin{theorem}\label{theorem_ordinal_3}
(1)\enskip $\tau(A) \approx A$ for all $A \in \fin$ and therefore $\tau(A) = \tau(A')$ if and only if $A \approx A'$. In particular
$\tau(\tau(A)) = \tau(A)$ for all $A \in \fin$ and $\tau(B) = B$ for all $B \in T$ (since $\tau$ is surjective).

(2)\enskip $\,\mathbf{J}$ is  a Peano iterator.

(3)\enskip If $A,\,A' \in \fin$ with $A \subset A'$ then $\tau(A) \subset \tau(A')$.

(4)\enskip If $A,\,A' \in \fin$ with $A \preceq A'$ then $\tau(A) \subset \tau(A')$.

(5)\enskip For all $B,\,B' \in T$ either $B \subset B'$ or $B' \subset B$. Thus $T$ is totally ordered by inclusion. Moreover, inclusion is the unique total
order on $T$ compatible with $\,\mathbf{J}$.
\end{theorem}

\proof(1) \enskip This is the same as the  proof of Theorem~\ref{theorem_ordinal_1}.

(2) \enskip This is the same as the  proof of Theorem~\ref{theorem_ordinal_2}.

(3) \enskip It is enough to show that if $A,\,B \in \fin$ are disjoint then $\tau(A) \subset \tau(A \cup B)$.
Let $\mathcal{S} = \{ C \in \mathcal{P}(B) :\tau(A) \subset \tau(A \cup C) \}$. Clearly $\varnothing \in \mathcal{S}$, so let
$C \in \proper{\mathcal{S}}$ and $c \in B \setminus C$, put $C' = C \cup \{c\}$.
 Then $\tau(A) \subset \tau(A \cup C)$ and hence also $\tau(A) \subset \tau(A \cup C')$, since $\tau(A \cup C') = h(\tau(A \cup C))$ and
$\tau(A \cup C) \subset h(\tau(A \cup C))$. Thus $C \cup \{c\} \in \mathcal{S}$, which implies $\mathcal{S}$ is an inductive $B$-system.
Therefore $B \in \mathcal{S}$, i.e., $\tau(A) \subset \tau(A \cup B)$.
 
(4) \enskip This follows from (3), Theorem~\ref{theorem_fs_4} and Proposition~\ref{prop_fs_112} (4).

(5) \enskip By Theorem~\ref{theorem_fs_4} either $B \preceq B'$ or $B' \preceq B$ and $\tau(B) = B$ and $\tau(B') = B'$. Thus by (4) either $B \subset B'$ or
$B' \subset B$. Moreover, $B \subset h(B)$ for all $B \in T$ and thus inclusion is the unique total order on $T$ compatible with $\,\mathbf{J}$.
\eop

We have already noted that  if $\,\mathbf{I}$ is a Peano iterator then the finite segment iterator $\,\mathbf{I}_\le$ associated with $\,\mathbf{I}$ is an
ordinal iterator.
 We now show that the construction in Theorem~\ref{theorem_iterators_199} can be reversed.
Again let $\,\mathbf{J} = (T,h,\varnothing)$ be an ordinal iterator.
Let $T^\dag$ be the class consisting of all elements $c$ for which there exists $B \in T$ such that $h(B) = B \cup \{c\}$.
Define $\gamma: T \to T^\dag$ by letting $\gamma(B) = c$, where $h(B) = B \cup \{c\}$. Thus $h(B) = B \cup \{\gamma(B)\}$  for all $B \in T$.

\begin{theorem}\label{theorem_ordinal_4}
(1)\enskip
The mapping $\gamma : T \to T^\dag$ is a bijection.

(2)\enskip Define  a mapping $h^\dag: T^\dag \to T^\dag$ by $h^\dag = \gamma \circ h \circ \gamma^{-1}$ and an iterator by $\,\mathbf{J}^\dag = (T^\dag,h^\dag,t^\dag_0)$, where $t^\dag_0 = \gamma(\varnothing)$ 
and so $t^\dag_0$ is the single element in $h(\varnothing)$.  Then $\gamma : \,\mathbf{J} \to \mathbf{J}^\dag$ is an isomorphism.

(3)\enskip $\,\mathbf{J}^\dag$ is a Peano iterator.

(4)\enskip $h(\gamma^{-1}(t)) = \gamma^{-1}(t) \cup \{t\}$ for all $t \in T^\dag$.

(5)\enskip  Define $\le$ on $T^\dag$ by letting $t' \le t$ if and only if
$\gamma^{-1}(t') \subset \gamma^{-1}(t)$. 
Then $\le$ is the unique a total order on $T^\dag$      compatible with $\,\mathbf{J}^\dag$.  

(6)\enskip $L_t =  \gamma^{-1}(t)$ for all $t \in T^\dag$.

(7)\enskip $h(L_t) = L_t \cup \{t\}$ for all $t \in T^\dag$.

\end{theorem}

\proof (1) \enskip
Let $B,\,B' \in T$ with $B \ne B'$; then by Theorem~\ref{theorem_ordinal_3}~(5) and without loss of generality we can assume $B$ is a proper subset of $B'$.
Then $h(B) \subset B'$; thus $\gamma(B) \in B'$ and $\gamma(B') \notin B'$ and so $\gamma(B) \ne \gamma(B')$. This shows that $\gamma$ is injective, and since $\gamma$ is
clearly surjective it follows that $\gamma$ is a bijection.

(2) \enskip
We have $t^\dag_0 =\gamma(\varnothing)$ and $h^\dag \circ \gamma = \gamma \circ h$ and therefore
$\gamma: \,\mathbf{J} \to \,\mathbf{J}^\dag$ is a morphism. But $\gamma$ is a bijection and hence by Lemma~\ref{lemma_iterators_16} $\gamma$ is an isomorphism.

(3)\enskip
It now follows from (2) and Theorem~\ref{theorem_iterators_7} that $\,\mathbf{J}^\dag$ is a Peano iterator. 

(4)\enskip  If $t \in T^\dag$ and $B = \gamma^{-1}(t)$ then
$h(\gamma^{-1}(t)) = h(B) = B \cup \{\gamma(B)\} = \gamma^{-1}(t) \cup \{t\}$.

(5)\enskip $\le$ is a total order on $T^\dag$ since inclusion defines a total order on $T$.
Let $t \in T^\dag$; then 
$\gamma^{-1}(t) \subset  h(\gamma^{-1}(t)) = \gamma^{-1}(h^\dag(t))$,
 since inclusion is the unique total order on $T$ compatible with $\,\mathbf{J}$ and by definition
$h^\dag(t) = \gamma \circ h \circ \gamma^{-1}(t)$.
Thus $t \le h^\dag(t)$ and hence $\le$ is the unique total order on $T^\dag$ compatible with $\,\mathbf{J}^\dag$. 

(6)\enskip Let $T^\dag_0 = \{ t \in T^\dag :  L_t =  \gamma^{-1}(t)\}$. In particular
$t^\dag_0 \in T^\dag_0$, since $L_{t^\dag_0}$ and $\gamma^{-1}(t^\dag_0)$ are both empty. 
Let $t \in T^\dag_0$; then by (4) $h(\gamma^{-1}(t)) = \gamma^{-1}(t)) \cup \{t\}$  and so
$\gamma^{-1}(h^\dag(t)) = h(\gamma^{-1}(t)) = \gamma^{-1}(t) \cup \{t\}$. Moreover, 
Theorem~\ref{theorem_iterators_199} (4) implies that
$L_{h^\dag(t)} = L_t \cup \{t\}$ and hence $L_{h^\dag(t)} = \gamma^{-1}(h^\dag(t))$,
i.e., $h^\dag(t) \in T^\dag_0$. It follows that $T^\dag_0 = T^\dag$, since $\,\mathbf{J}$ is minimal and therefore $L_t = \gamma^{-1}(t)$ for all $t \in T^\dag$.

(7)\enskip By (4) and (6) we have  $h(L_t) = h(\gamma^{-1}(t)) = \gamma^{-1}(t) \cup \{t\} = L_t \cup \{t\}$. \eop

We call the iterator $\,\mathbf{J}^\dag$ the \definition{dual iterator} to the ordinal iterator $\,\mathbf{J}$.

Consider the canonical ordinal iterator $\mathbf{O} = (O,\vonn,\varnothing)$ and let $(O^\dag$  
be the class consisting of all elements $c$ for which there exists $o \in O$ such that $\vonn(o) = o \cup \{c\}$.
Then $O^\dag = O$, since $\vonn(o) = o \cup \{o\}$ for each $o \in O$. Moreover, if $\gamma' : O \to O^\dag= O$ is the mapping
corresponding to the mapping $\gamma : T \to T^\dag$ in Theorem~\ref{theorem_ordinal_4} then clearly  $\gamma = \id_O$. Therefore
$\,\mathbf{O}^\dag = \,\mathbf{O}$ and so the iterator $\,\mathbf{O}$ is its own dual.

\begin{proposition}\label{prop_ordinal_22}
Let $\,\mathbf{I} = (X,f,x_0)$ be a Peano iterator and let $\,\mathbf{I}_\le$ be the finite segment iterator associated with $\,\mathbf{I}$.
Also let $\,\mathbf{I}_\le^\dag$ be the dual  iterator  to the ordinal iterator $\,\mathbf{I}_\le$. 
Then $\,\mathbf{I}_\le^\dag = \,\mathbf{I}$.
\end{proposition}

\proof 
Let $\le$ be the unique total order on $X$ compatible with $\,\mathbf{I}$. We thus have  $\,\mathbf{I}_\le = (X_\le,f_\le, \varnothing)$, where
$X_\le = \{ L_x : x \in X \}$ and $f_\le(L_x) = L_{f(x)} = L_x \cup \{x\}$ for all $x \in X$. 
Put $\,\mathbf{I}^\dag_\le = (Y,g,y_0)$. Then  $Y$ is the class of all elements $c$ for which there exists $L_x \in X_\le$ such that $f_\le(L_x) = L_x \cup \{c\}$. But
$f_\le(L_x) = L_x \cup \{x\}$ and hence $Y = X$. Moreover, if $\beta : X_\le \to Y = X$ is the mapping corresponding to the mapping 
$\gamma : T \to T^\dag$ in Theorem~\ref{theorem_ordinal_4} then $\beta(L_x) = x$  for all $x \in X$. Also $g : X \to X$ is the mapping 
$\beta \circ f_\le \circ \beta^{-1} = f$, and $y_0$ is the single element in $f_\le(\varnothing)$ which is $x_0$. This shows that $\,\mathbf{I}_\le^{*}  = \,\mathbf{I}$. 
\eop

\begin{proposition}\label{prop_ordinal_23}
Let $\,\mathbf{J} = (T,h,\varnothing)$ be an ordinal iterator and let $\,\mathbf{J}^\dag = (T^\dag,h^\dag,t^\dag_0)$ be the dual iterator to $\,\mathbf{J}$.
Also let $\,\mathbf{J}^\dag_\le$ be the finite segment iterator associated with $\,\mathbf{J}^\dag$.
Then $\,\mathbf{J} = \mathbf{J}^\dag_\le$.
\end{proposition}

\proof Let $\le$ be the unique total order on $T^\dag$ compatible with $\,\mathbf{J}^\dag$. and let 
$(Y,g,\varnothing)$ be the components of the iterator $\,\mathbf{J}^\dag_\le$. Then by definition $Y = \{ L_t : t \in T^\dag \}$ and hence by
Theorem~\ref{theorem_ordinal_4} (6) $Y=  \{ \gamma^{-1} (t) : t \in T^\dag \} = T$.
Also by Theorem~\ref{theorem_iterators_199} (4) $g(L_t) = L_t \cup \{t\}$ and by Theorem~\ref{theorem_ordinal_4} (7) 
$h(L_t) = L_t \cup \{t\} = g(t)$. Therefore $g = h$ which shows that $\,\mathbf{J} = \,\mathbf{J}^\dag_\le$. \eop


\startsection{Finite minimal iterators}

\label{finiteiter}

Theorems~\ref{theorem_iterators_2} and \ref{theorem_iterators_4} imply that for a minimal iterator $\mathbf{I} = (X,f,x_0)$ there are two mutually 
exclusive possibilities: Either $\,\mathbf{I}$ is  a Peano iterator or $X$ is a finite set. In this section we deal with case in which $X$
is a finite set.

Thus in what follows let $\,\mathbf{I} = (X,f,x_0)$ be a minimal iterator with $X$ a finite set. 
For each $x \in X$ let $X_x$ be the least $f$-invariant subset of $X$ containing $x$ and let $f_x : X_x \to X_x$ be the restriction of $f$ to $X_x$. Thus
$\,\mathbf{I}_x = (X_x,f_x,x)$ is a minimal iterator. Also let $\val_x$ be the unique assignment of finite sets in $\,\mathbf{I}_x$; put $\val = \val_{x_0}$. 
By Proposition~\ref{prop_iterators_1} and Theorem~\ref{theorem_fs_1} $f_x$ is a bijection if and only if $x \in f_x(X_x)$.

The iterator $\,\mathbf{I}$ can be considered as a
finite dynamical system with the dynamics given by the mapping $f : X \to X$ and with initial state $x_0$. 
Now it is an elementary fact that the mapping $f$ is then eventually periodic. What this means 
can best be described by  by fixing  a 
Peano iterator
$\,\mathbf{I} = (N,s,0)$ and to reduce the clutter we prefer to use $0$ instead of $n_0$ for the third component of $\,\mathbf{I}$. Let $\le$ be the unique total order on $N$. Let $x \in X$; then applying the recursion theorem for the iterator $(N,s,,0)$
to the iterator $\,\mathbf{I}_x$ there exists a unique mapping $\pi_x : N \to X_x$ with $\pi_x(0) = x$ such that $\pi_x(s(n)) = f_x(\pi_x(n))$ for all 
$n \in N$. We use the standard notation and write $f^n(x)$ instead of $\pi_x(n)$. Thus $f^0(x) = x$ and $f^{s(n)}(x) = f(f^n(x))$ for all $n \in N$.
The element $f^n(x)$ can be thought of as  the $n$-th iterate of $f$ when the initial state is $x$.

An element $x\in X$ is periodic if $x = f^n(x)$ for some $n \ne 0$ and if $n$ is the least such index then  the set $ \{ f^k(x) : 0 \le k < n \}$ is the 
corresponding periodic cycle. The elementary fact about $\,\mathbf{I}$ states that there is a unique periodic  cycle and that there exists $m \in N$ such that 
$f^n(x_0)$ is periodic for all $n \ge m$. Note that the uniqueness of the periodic cycle only holds because $\,\mathbf{I}$ is minimal.

For each $x \in X$ let $X'_x$  be the set consisting of those $y \in X$ such that $y = f^n(x)$ for some $n \in N$. 

\begin{lemma}\label{lemma_finiteiter_23}
Let $x \in X$. Then:

(1)\enskip   $X'_x= X_x$. 
 
(2)\enskip $x$ is periodic if and only if $x \in f(X_x)$ and thus $x$ is periodic if and only if $f_x$ is a bijection. 

(3)\enskip  $X_{f(x)} = f(X_x)$ and thus $x$ is periodic if and only if $x \in X_{f(x)}$.
\end{lemma}

\proof (1)\enskip Clearly $X'_x$ is $f$-invariant  and contains $x$ and so $X_x \subset X'_x$. 

Let $N_0 = \{ n \in N: f^n(x) \in X_x \}$.Then $0 \in N_0$ and $N_0$ is $s$-invariant, since if
$ y  = f^n(x) \in X_x$ then $f^{s(n)} (x) = f(y) \in X_x$. Thus $N_0 = N$, which shows that $X'_x \subset X_x$.

(2)\enskip Suppose first that $x$ is periodic and so there exists $n \ne 0$ with $x = f^n(x)$. There exists a unique $m \in N$ with $n = s(m)$ and then
$x = f(f^m(x))$. But by (1) $f^m(x) \in X_x$ and hence $x \in f(X_x)$. Suppose conversely that $x \in f(X_x)$. Then by (1) there exists $n \in N$ such that 
$x = f(y)$ with $y = f^n(x)$. But then $x = f^m(x)$ with $m =s(n)$ , and $m \ne 0$. Thus $x$ is periodic.

(3)\enskip 
$f(X_x)$ is $f$-invariant and contains $f(x)$ and so $X_{f(x)} \subset f(X_x)$. Let $Y_x = \{ y \in X_x : f(y) \in X_{f(x)}\}$. Then $Y_x$ is $f$-invariant and contains 
$x$ and hence $Y_x  = X_x$. Thus also $f(X_x) \subset X_{f(x)}$. \eop

By Lemma~\ref{lemma_finiteiter_23} we do not need to employ the Peano iterator $\mathbf{I}$. We define $x$ to be periodic if $x \in f(X_x)$, which is 
equivalent to requiring that $f_x$ be a bijection. Then $x$ will also be periodic if and only if $x \in X_{f(x)}$. In the following the Peano
iterator $\mathbf{I}$ no longer appears.

\begin{theorem}\label{theorem_finiteiter_1}
(1) \enskip Let $X_P = \{ x \in X : \mbox{ $x$ is periodic} \}$.
Then $X_P$ is non-empty and  $X_x = X_y$) for all $x,\,y \in X_P$.  Thus $f$ maps $X_P$ bijectively onto itself.

(2) \enskip Let $X_N =\{ x \in X : \mbox{ $x$ is not periodic} \}$ and suppose $X_N \ne \varnothing$.
Then $f$ is injective on $X_N$ and there exists a unique element $u \in X_N$ such that $f(u)$ is periodic. Moreover, there exists a unique element 
$v \in X_P$ such that $f(v) = f(u)$ and $u$ and $v$ are the unique elements of $X$ with $u \ne v$ such that $f(u) = f(v)$.
\end{theorem}

For the proof of (1) we need the following:

\begin{lemma}\label{lemma_finiteiter_2}
Let $x \in X$. Then:

(1)\enskip  $X_{f(x)} \subset X_x \subset \{x\} \cup X_{f(x)}$. 
Thus $X_{f(x)}$ is either $X_x$ or $X_x \setminus\{x\}$. If $X_{f(x)} = X_x$ then $x \in X_{f(x)}$ and so
$x$ is periodic. If $X_{f(x)} =  X_x \setminus \{x\}$ then $x$ is not periodic. Thus $x$ is periodic if and only if $X_{f(x)} = X_x$.

(2)\enskip If $x$ is periodic then so is $f(x)$.

(3)\enskip If $x$ is periodic then $y$ is periodic and $X_y = X_x$ for all $y \in X_x$.

(4)\enskip If $A$ is a finite set, $B \subset A$ and $x = \val(B)$ then $\val(A) \in X_x$.

(5)\enskip For each $x \in X$ we have $X_y \subset X_x$ for all $y \in X_x$.
\end{lemma}

\proof\enskip (1) $X_x$ is $f$-invariant and contains $f(x)$ and hence $X_{f(x)} \subset X_x$. Also $\{x\} \cup  X_{f(x)}$ is $f$-invariant and contains $x$ 
and therefore $X_x \subset \{x\} \cup X_{f(x)}$. 

(2)\enskip This follows from the final statement in (1).

(3)\enskip Let $\mathcal{S} = \{ y \in X : \mbox{ $y$ is periodic and $X_y = X_x$} \}$ then $x \in \mathcal{S}$ and if $y \in \mathcal{S}$ then by (1) $f(y)$ 
is periodic and so$X_{f(y)} =X_y = X_x$. Thus $f(y) \in \mathcal{S}$ and so $\mathcal{S}$ is $f$-invariant. Hence $\mathcal{S} \subset X_x$, i.e.,
$y$ is periodic and $X_y = X_x$ for all $y \in X_x$.

(4)\enskip Let $C = A \setminus B$ and put $\mathcal{S} = \{ D \subset C : \val(D) \in X_x \}$. Then $\varnothing \in \mathcal{S}$ and if $D \in\proper{\mathcal{S}}$
and $d \in C \setminus D$ then $\val(D \cup \{d\}) = f(\val(D)) \in \mathcal{S}$, since $\val(D) \in \mathcal{S}$ and $X_x$ is $f$-invariant.
Therefore $\mathcal{S}$ is an inductive $C$ system and thus $C \in \mathcal{S}$, i.e., $A = B \cup C \in X_x$.

(5) \enskip Let $X' = \{ y \in X : X_y \subset X_x\}$. Then $x \in X'$ and if $y \in X'$ (and so $X_y \subset X_x$) then
$X_{f(y)} \subset X_y \subset X_x$ and hence $f(y) \in X'$. Therefore $X_x \subset X'$, i.e., $X_y \subset X_x$ for all $y \in X_x$.
\eop

\it{Proof of Theorem~\ref{theorem_finiteiter_1} (1)} \enskip
 Let $\mathcal{S} = \{ X_y : y \in X \}$. Then by Proposition~\ref{prop_intro_2} there exists
$x\in X$ such that $X_x$ is a  minimal element of 
$\mathcal{S}$. But by Lemma~\ref{lemma_finiteiter_2}  $X_{f(x)} \subset X_x$ and so $X_{f(x)} = X_x$. Hence  by 
Lemma~\ref{lemma_finiteiter_2}  $x$ is periodic, i.e., $x \in X_P$. 
Thus by Lemma~\ref{lemma_finiteiter_2} $y$ is periodic and $X_y = X_x$ for all $y \in X_x$.

Let $x,\,y \in X_P$; by Lemma~\ref{lemma_iterators_4} there exist finite sets $A$ and $B$ with $x = \val(A)$ and $y = \val(B)$ and by 
Proposition~\ref{prop_fs_112} (3)
and Theorem~\ref{theorem_iterators_1} we can assume that $B\subset A$ or $A \subset B$ and without loss of generality assume that $B \subset A$ . 
Put $C = A \setminus B$ and let $\mathcal{S} = \{ D \subset C : 
\mbox{$\val(B \cup D)$ is periodic and $ X_{\val(B \cup D)} = X_y$}\}$ and so $\varnothing \in \mathcal{S}$. 
Thus let $D \in \proper{\mathcal{S}}$ and $d \in C\setminus D$. Then $\val(B \cup D \cup \{d\}) = f(z)$ where $z = \val(B \cup D))$
and $z$ is periodic and $X_z = X_y$, since $D \in \mathcal{S}$. Hence by Lemma~\ref{lemma_finiteiter_2}  $f(z)$ is periodic and $X_{f(z)} = X_z  = X_y$, 
which shows that $D \cup \{d\}  \in\mathcal{S}$.
Therefore $\mathcal{S}$ is an inductive $C$-system and so $C \in \mathcal{S}$, i.e., $X_x=X_y$.
\eop

For the proof of (2) we need the following:

\begin{lemma}\label{lemma_finiteiter_22} Let $s,\,t \in X$ with $s \ne t$ and $f(s) = f(t)$. Then $u = f(s) = f(t)$ is periodic. Moreover, one of $s$ and $t$ is
periodic. 
\end{lemma}

\proof By Lemma~\ref{lemma_iterators_4} there exist finite sets
 $B$ and $C$ such that $s = \val(B)$ and $t = \val(C)$ and by Theorem~\ref{theorem_iterators_1} and Proposition~\ref{prop_fs_112}  we can assume
without loss of generality that $B \subset C$, and so $B$ is a proper subset of $C$.  Let $d \notin C$ and put $B'= B \cup \{d\}$, $C' = C \cup \{d\}$.
Then $\val(B') = \val(C')= u$. Now let $a \in C \setminus B$ and put
$C'' = C' \setminus \{a\}$. But $u = \val(B')$ and $B' \subset C''$ and so Lemma~\ref{lemma_finiteiter_2}   implies that $\val(C'') \in X_u$ 
and then $u = f(\val(C'')) \in f(X_u)$. This shows that $u$ is periodic. Now $B$ is a proper subset of $C$ and so there exists $D \supset B'$ with
$D \approx C$. Thus $u = \val(B')$ and $t = \val(C) = \val(D)$ and $B' \subset D$ and so by Lemma~\ref{lemma_finiteiter_2} 
$t \in X_u$. Therefore by Lemma~\ref{lemma_finiteiter_2} $t$ is periodic.
\eop

\it{Proof of Theorem~\ref{theorem_finiteiter_1} (2)} \enskip
Lemma~\ref{lemma_finiteiter_22} implies that $f$ is injective on $X_N$. By Proposition~\ref{prop_intro_2} there exists
$u \in X_N$ such that $X_u$ is a minimal element of the set $\mathcal{S}_N = \{ X_y : y \in X_N \}$.  Then $u$ is not periodic and so by 
Lemma~\ref{lemma_finiteiter_2}  $X_{f(u)}$ is a
 proper subset of $X_u$ and hence $X_{f(u)} \notin \mathcal{S}_N$, i.e., $f(u)$ is periodic. Suppose there exist $u_1,\,u_2 \in X_N$ with $u_1 \ne u_2$
and such that $f(u_1)$ and $f(u_2)$ are both periodic. Then there exist finite sets $A_1,\,A_2$ with $u_i = \val(A_i)$ for $i = 1,\,2$ and as usual we can
assume that $A_1$ is a proper subset of $A_2$. Let $a \in A_2 \setminus A_1$ and so $A_1' = A_1 \cup \{a\}\subset A_2$. 
Then $\val(A_1') = f(\val(A_1) = f(u_1)$
is periodic and
by Lemma~\ref{lemma_finiteiter_2} 
$u_2 \in X_{f(u_1)}$. Hence again  using Lemma~\ref{lemma_finiteiter_2} $u_2$ would be periodic. This contradiction
shows that there is a unique $u \in X_N$ such that $z = f(u)$ is periodic.
 Thus $z \in f_z(X_z)$ and so there exists a finite set $A'$ such
that $z = f(\val_v(A'))$. Let $a$ be an element not in $A'$ and put $A = A' \cup \{a\}$. Then $A$ is a non-empty finite set and $z = \val_z(A)$. Let 
$v = \val_v(A')$. Then $v \in X_v$ and so $v$ is periodic and $f(v) = z = f(u)$. Moreover, $v$ is the unique element of $X_P$ with $f(v) = f(u)$, since
$f$ maps $X_P$ bijectively onto itself.  Finally, let $u',\,v' \in X$ with $u' \ne v'$ and $f(u') = f(v')$. Since $f$ is injective on $X_N$ and on $X_P$
one of these elements is in $X_N$ and the other in $X_P$. Label them so $u' \in X_N$ and $v' \in X_P$. Then by Lemma~\ref{lemma_finiteiter_22} $f(u') \in X_P$
and by the uniqueness of $u$ it follows that $u' = u$ and by the uniqueness of $v$ it follows that $v' = v$.
\eop 

We next consider the special case in which $X$ contains a fixed-point, i.e., an element $z$ with $f(z) = z$. If $z = x_0$ then $X = \{x_0\}$ and we assume that
this is not the case, and thus $z \ne x_0$. Since $z$ is periodic it follows from Theorem~\ref{theorem_finiteiter_1}~(1) that
$X_P = X_z =\{z\}$. Hence $X_N = X \setminus \{z\}$ and by Theorem~\ref{theorem_finiteiter_1}~(2) $f$ is injective on $X_N$ and there exists a unique
$w \in X_N$ such that $f(w) = z$.

Let $A$ be a finite set with $z = \val(A)$ and by Proposition~\ref{prop_intro_2} we can assume that $z \ne \val(B)$ for each proper subset $B$ of $A$.

\begin{lemma}\label{lemma_finiteiter_3} For each $x \in X$ there exists $B \subset A$ with $x = \val(B)$.
\end{lemma}

\proof Let $X_0 = \{x \in X : \mbox{ $x = \val(B)$ for some $B \subset A$}\}$, and thus $x_0 \in X_0$, since $x_0 = \val(\varnothing)$.
Let $x \in X_0$ with $x = \val(B)$. If $B$ is a proper subset of $A$ and $a \in A \setminus B$ then $B \cup \{a\} \subset A$ and 
$f(x) = \val(B \cup \{a\})$ and so $f(x) \in X_0$. But if $B = A$ then $x = z$ and so $f(z) = z \in X_0$. Thus $X_0$ is $f$-invariant and contains $x_0$
and hence $X_0 = X$. \eop
 
\begin{lemma}\label{lemma_finiteiter_4} If $B,\,B' \in \mathcal{P}(A)$ with $\val(B) = \val(B')$ then $B \approx B'$.
\end{lemma}

\proof Suppose there exist $B,\,B' \in \mathcal{P}(A)$ with $\val(B) = \val(B')$ and $B \not\approx B'$. Then by Proposition~\ref{prop_fs_112} (2) (and if necessary
exchanging the r\^oles of $B$ and $B'$)
there exist such $B,\,B'$ with $B' \subset B$,i.e., with $B'$ a proper subset of $B$. Let $\mathcal{S}$ be the subset of
$\mathcal{P}(A)$ consisting of those subsets $B$ which contain a proper subset $B'$ with $\val(B) = \val(B')$.
Thus $\mathcal{S}$ is non-empty and hence by Proposition~\ref{prop_intro_3} $\mathcal{S}$ contains a maximal element $C$; let $C'$ be a proper subset
of $C$ with $\val(C) = \val(C')$. Now $C \ne A$, since otherwise $C'$ would be a proper subset of $A$ with $\val(C') = z$.  Choose $a \in A \setminus C$;
then $C' \cup \{a\}$ is a proper subset of $C \cup \{a\}$. But
\[\val(C' \cup \{a\}) = f(\val(C')) = f(\val(C)) = \val(C \cup \{a\})\:,\]
 which contradicts the maximality of $C$. Therefore $B \approx B'$ whenever $B,\,B'$ are subsets of $A$ with $\val(B) = \val(B')$. \eop

\begin{lemma}\label{lemma_finiteiter_5} If $x_1\, x_2 \in X$ with $X_{x_1} = X_{x_2}$ then $x_1 = x_2$.
\end{lemma}

\proof Let $x_1,\,x_2 \in X$ with $x_1 \ne x_2$. Then by Lemma~\ref{lemma_finiteiter_3} and Proposition~\ref{prop_fs_112} (2) (and if 
necessary exchanging the r\^oles of $x_1$ and $x_2$) there exist $B_1,\, B_2 \in \mathcal{P}(A)$. 
with $x_k = \val(B_k)$ for $k = 1,\,2$ and such that $B_1$ is a proper subset of $B_2$. Let $b \in  B_2 \setminus B_1$; then $B_1' = B_1 \cup \{a\} \subset B_2$
and $\val(B_1') = f(x_1)$. Thus by Lemma~\ref{lemma_finiteiter_2}  $x_2 \in X_{f(x_1)}$ and so $X_{x_2} \subset X_{f(x_1)}$.  
If $x_1 = z$ then $X_{f(x_1)} = \{z\}$ and so $X_{x_2} = \{z\}$. But $x_2 \ne z$, and Lemma~\ref{lemma_finiteiter_2}  implies  that $X_{f(x_2)}$ would be
 a proper subset of $X_{x_2} = \{z\}$. Therefore this case cannot occur. Hence $x_1 \ne z$  and
by Lemma~\ref{lemma_finiteiter_2}  $X_{f(x_1)}$ is a proper subset of $X_{x_1}$. Therefore$X_{x_1} \ne X_{x_2}$. \eop

Define a relation $\le$ on $X$ by stipulating that $x_1 \le x_2$ if $X_{x_2} \subset X_{x_1}$. 
As usual, $x_1 < x_2$ means that $x_1 \le x_2$ but $x_1 \ne x_2$.

\begin{proposition}\label{prop_finiteiter_11} (1)\enskip
The relation $\le$ is  a total order on $X$ with $x_0 \le x \le z$ for all $x \in X$, Moreover, if $x \in X \setminus \{z\}$ then $x < f(x)$.

(2)\enskip If $x \le y \le f(x)$ then $y = x$ or $y = f(x)$.

(3)\enskip Let $x, y \in X$. if $x \le y$ if then $f(x) \le f(y)$. Moreover, if $x,\,y \in X \setminus \{z\}$ and $f(x) \le f(y)$ then $x \le y$.

(4)\enskip  Let $x, y \in X \setminus \{z\}$. Then $x < y$ if and only if $f(x) < f(y)$.
\end{proposition}

\proof (1)\enskip It is clear that $\le$ is transitive and if both $x_l \le x_2$ and $x_2 \le x_1$ hold then by Lemma~\ref{lemma_finiteiter_5} $x_1 = x_2$.
Let $x_1,\,x_2 \in X$ and as usual we can assume that there exist $B_1,\,B_2 \in \mathcal{P}(A)$ with $x_k = \val(B_k)$ for $k = 1,\,2$ and either $B_2 \subset B_1$
or $B_1 \subset B_2$. If $B_2 \subset B_1$
then by Lemma~\ref{lemma_finiteiter_2}  $x_1 \in X_{x_2}$ and hence $X_{x_1} \subset X_{x_2}$. 
If $B_1 \subset B_2$ then in the same way $X_{x_2} \subset X_{x_1}$. Therefore either $x_1 \le x_2$ or $x_2 \le x_1$ and this shows that $\le$ is a total order.
It is clear that $x_0 \le x \le z$ for all $x \in X$ If $x \in X \setminus \{z\}$ then by Lemma~\ref{lemma_finiteiter_2} 
$X_{f(x)}$ is a proper subset of $X_x$ and so $x < f(x)$. 

(2)\enskip If $x \le y \le f(x)$ then $X_y \subset X_x$ and $X_{f(x)} \subset X_y$ and hence
\[\{x\} \cup X_y \subset \{x\} \cup X_x = X_x = \{x\} \cup X_{f(x)} \subset  \{x\} \cup X_y\;.\] 
It follows that
$\{x\} \cup X_y = X_x = \{x\} \cup X_{f(x)} =  \{x\} \cup X_y$ and in particular
$\{x\} \cup X_y = X_x$.
Thus either $x\in X_y$, in which case $X_x = X_y$, or $x \notin X_y$, in which case $X_y = X_x \setminus \{x\}$ and then by Lemma~\ref{lemma_finiteiter_2}  
$X_y = X_{f(x)}$. Therefore by  Lemma~\ref{lemma_finiteiter_5} either $y = x$ or $y = f(x)$. 

(3)\enskip Note that for all $x \in X$ both $x \le z$ and $f(x) \le f(z)$ hold trivially. Thus we can assume that $y \in X \setminus \{z\}$.
Suppose $x \le y$ and so also $x \in X \setminus \{z\}$.
Then $X_y \subset X_x$ and $X_{f(y)} = X_y \setminus \{y\}$, $X_{f(x)} = X_x \setminus \{x\}$
and thus $X_{f(y)} \subset X_{f(x)}$ provided $x \notin X_y$. But if $x \in X_y$ then $X_x \subset X_y$ and so $X_x = X_y$. Hence either
$f(x) \le f(y)$ or $x = y$ and in both cases $f(x) \le f(y)$. Next let $x,\,y \in X \setminus \{z\}$ with $f(x) \le f(y)$. Then either
$x \le y$ or $y \le x$ and if $y \le x$ then by the above $f(y) \le f(x)$ and so $f(x) = f(y)$. But if $f(x) = f(y)$ then $x = y$, since $f$ is injective on
$X \setminus \{z\}$, and again $x \le y$.  

(4)\enskip If $x < y$ then by (3) $f(x) \le f(y)$ and $f(x) \ne f(y)$ since $x \ne y$ and $f$ is injective on $X \setminus \{z\}$. Hence $f(x) < f(y)$.
In the same way, if $f(x) < f(y)$ then by (3) $x \le y$. But $x \ne y$ since $f(x) \ne f(y)$ and thus $x < y$. 
\eop

The construction given above can be reversed. Let $(X,\le)$ be a finite totally ordered set  with least element $x_0$ and greatest element $z$. Define a 
mapping $f : X \to X$ by letting $f(x)$ be the least element in $\{ y \in X: y > x \}$ if $x \ne z$ and putting $f(z) = z$. 
Then $\,\mathbf{I} = (X,f,x_0)$ is a finite iterator with fixed point $z$.

\begin{proposition}\label{prop_finiteiter_1} The iterator $\,\mathbf{I}$ is minimal.
\end{proposition}
 
\proof
Let  $X_0$ be an $f$-invariant subset of $X$ containing $x_0$ and suppose $X_0 \ne X$. Let $u$ be the least element in $X \setminus X_0$ and so 
$u \ne x_0$. Let $v$ be the greatest element in $\{ y \in X : y < u \}$. Then $v < u$, hence $v \in X_0$ and so $f(v) \in X_0$. But this is not possible, since
  $f(v) = u$. Therefore $X_0 = X$ which shows that $\,\mathbf{I}$ is minimal. \eop

Let $G$ be a finite group, with the product of $a$ and $b$ in $G$ denoted just by $ab$, with identity element $1$ and with $a^{-1}$ the inverse of $a$.
For each $a \in G$ let $n_a  : G \to G$ be given by $n_a(b) = ab$ for each $b \in G$. Then there is the iterator $\,\mathbf{I}_a = (G,n_a,1)$.
Also let $G_a$ be the least $n_a$-invariant subset of $G$ containing $1$ and $m_a : G_a \to G_a$ be the restriction of $n_a$ to $G_a$, so
$(G_a,m_a,1)$ is a minimal iterator.

\begin{proposition}\label{prop_finiteiter_2} For each $a\in G$ the mapping $m_a$ is a bijection and so $a$ is periodic.
\end{proposition}

\proof Suppose $a$ is not periodic.
Then by Theorem~\ref{theorem_finiteiter_1} (2)
there exists a unique non-periodic element $u \in G_a$ such that $m_a(u)$ is periodic and a unique periodic element $v \in G_a$ with $m_a(v) = m_a(u)$, 
i.e., $av = au$. But then $v = u$, which contradicts the fact that $u$ is not periodic and $v$ is periodic. Hence $a$ is periodic and $m_a$ is a bijection.
\eop

\begin{proposition}\label{prop_finiteiter_3} Let $H$ be a subset of $G$ containing $1$ and such that $ab \in H$ for all $a,\,b \in H$. Then $H$ is a subgroup of
$G$.
\end{proposition}

\proof Let $a \in H$; then  $H$ is an $m_a$-invariant subset of $G_a$ containing $1$ and so $G_a \subset H$. Now by Proposition~\ref{prop_finiteiter_2} 
$m_a : G_a \to G_a$ is a bijection and $1 \in G_a$ and so there exists $b \in G_a$ with $m_a(b) = ab  = 1$. Thus $b \in H$ and $b = a^{-1}$. 
This shows that $H$ is a subgroup of $G$. \eop


\startsection{Addition and multiplication}

\label{am}

In this section we show how an addition and a multiplication can be defined for any minimal iterator. 
These operations are associative and commutative and can be specified by the rules (a0), (a1), (m0) and (m1) 
below, which are usually employed when defining the operations on $\Nat$ via the Peano axioms.

Note that,even if we do not assume the existence of an infinite set, we can apply the results of this section to the Peano iterator $\,\mathbf{O}$.

In the following let $\,\mathbf{I} = (X,f,x_0)$ be a minimal iterator  with $\val$ the assignment of finite sets in $\,\mathbf{I}$.
We first state the main results (Theorems \ref{theorem_am_1} and 
\ref{theorem_am_2}) and then develop the machinery required to prove them. In the following section we give 
alternative proofs for these theorems.
  
\begin{theorem}\label{theorem_am_1}
There exists a unique binary operation $\oplus$ on $X$ such that
\[ 
\val(A) \oplus \val(B) = \val(A \cup B) 
\]
whenever $A$ and $B$ are disjoint finite sets. This operation $\oplus$ is both associative and commutative, 
$x \oplus x_0 = x$ for all $x \in X$ and for all $x_1,\,x_2 \in X$ there is an $x \in X$ such that either 
$x_1 = x_2 \oplus x$ or $x_2 = x_1 \oplus x$. Moreover, $\oplus$ is the unique binary operation $\oplus$ on $X$ 
such that
\begin{evlist}{15pt}{6pt}
\item[(a0)]  $x \oplus x_0 = x$ for all $x \in X$. 

\item[(a1)]  $x \oplus f(x') = f(x \oplus x')$ for all $x,\,x' \in X$.
\end{evlist}
\end{theorem}

\begin{theorem}\label{theorem_am_2}
There exists a unique binary operation $\otimes$ on $X$ such that
\[ 
\val(A) \otimes \val(B) = \val(A\times B) 
\]
for all finite sets $A$ and $B$. This operation $\otimes$ is both associative and commutative, 
$x \otimes x_0 = x_0$ for all $x\in X$ and $x \otimes f(x_0) = x$ for all $x \in X$ with $ \ne x_0$ 
(and so $f(x_0)$ is a multiplicative identity element) and 
the distributive law holds for $\oplus$ and $\otimes$: 
\[x \otimes (x_1 \oplus x_2) = (x \otimes x_1) \oplus (x \otimes x_2)\]
for all $x,\,x_1,\,x_2 \in X$.
Moreover, $\otimes$ is the unique binary operation on $X$ such that 
\begin{evlist}{15pt}{6pt}
\item[(m0)]  $x \otimes x_0 = x_0$ for all $x \in X$. 

\item[(m1)]  $x \otimes f(x') = x \oplus (x \otimes x')$ for all $x,\,x' \in X$.
\end{evlist}
\end{theorem}

We now prepare for the proofs of Theorems~\ref{theorem_am_1} and Theorem~\ref{theorem_am_2} and first look at what is common to these two theorems. 
Let $D$ be a subset of $\fin \times \fin$ and let $\gamma : D \to \fin$ be a mapping. In Theorem~\ref{theorem_am_1} we will have
\[D = \{(A,B) \in \fin \times \fin : \mbox{$A$ and $B$ are disjoint}\}\] and $\gamma(A,B) = A \cup B$ and in Theorem~\ref{theorem_am_2}
$D = \fin \times \fin$ and $\gamma(A,B) = A \times B$. We assume that the mapping $\val' : D \to X \times X$ with $\val'(A,B) = (\val(A),\val(B))$ is 
surjective. By Lemma~\ref{lemma_iterators_4} this is clearly the case for Theorem~\ref{theorem_am_2} and for Theorem~\ref{theorem_am_1} it follows from the
next result.
\begin{lemma}\label{lemma_am_1}
For all $(x,\,x') \in X \times X$ there exist disjoint finite sets $A$ and $B$ with $(x,x') = (\val(A),\val(B))$.
\end{lemma}

\proof 
By Lemma~\ref{lemma_iterators_4} there exists $(C,D) \in \fin \times \fin$ with $(\val(C),\val(D)) = (x,x')$ and by
Proposition
~\ref{prop_fs_112} (4) there exist disjoint finite sets $A$ and $B$ with $A \approx C$
and $B \approx D$. Therefore by  Theorem~\ref{theorem_iterators_1} (2) the sets $A$ and $B$ are  disjoint with  $(x,x') = (\val(A),\val(B))$.
\eop

Theorem~\ref{theorem_am_1} and Theorem~\ref{theorem_am_2} state for the appropriate mapping $\gamma : D \to \fin$ there exists a binary operation 
$\odot$ on $X$ such that $\val(A) \odot \val(B) = \val(\gamma(A,B))$ for all $(A,B) \in D$. 

\begin{proposition}\label{prop_am_111}
Let $\gamma : D \to \fin$ be an arbitrary mapping for which the mapping $\val' : D : \to X\times X$ is surjective.
Then there exists a binary operation $\odot$ on $X$ such that 
$\val(A) \odot \val(B) = \val(\gamma(A,B))$ for all $(A,B) \in D$ if and only if 
\begin{evlist}{26pt}{6pt}
\item[($\heartsuit$)]
$\val(\gamma(A,B)) = \val(\gamma(A'.B'))$ 
whenever $(A,B)$ and $(A',B')$ are elements  of $D$ with $\omega'(A,B) = \omega'(A',B')$.
\end{evlist}

Moreover, if ($\heartsuit$) holds then $\odot$ is the unique binary operation $\odot'$ on  $X$ such that
$\val(A) \odot' \val(B) = \val((\gamma(A,B))$ for all $(A,B) \in D$.
\end{proposition}

\proof This is a special case of Proposition~\ref{prop_iterators_77}. 
If $\odot$ is written as a prefix operation then
the requirement on $\odot$ is that $\odot(\val'(A,B)) = \val(\gamma((A,B))$ for all $(A,B) \in D$ which in turn is the requirement that
$\odot \circ \val' = \val \circ \gamma$.
\eop

Let $D_\oplus = (A,B) \in \fin \times \fin : \mbox{ $A$ and $B$ are disjoint}\}$
and let $\gamma_\oplus : D_\oplus \to \fin$ be given by $\gamma_\oplus(A,B) = A \cup B$ for all $(A,B) \in D_\oplus$.
Also let $D_\otimes = \fin \times \fin$ and let $\gamma_\otimes : D_\otimes \to \fin$ be given by $\gamma_\otimes(A,B) = A \times B$ for all $(A,B) \in D_\otimes$.
We will establish the existence of the operations $\oplus$ and $\otimes$ in Theorems \ref{theorem_am_1}and \ref{theorem_am_2} by showing that
the mappings $\gamma_\oplus$ and $\gamma_\otimes$ satisfy condition ($\heartsuit$) and then applying Proposition~\ref{prop_am_111}. 

For a Peano iterator $\,\mathbf{I}$ this is not a problem. Consider $\gamma_\oplus$: If $(A,B),\, (A',B') \in D_\oplus$ with 
$\val'(A,B) = \val'(A',B')$ then by Theorem~\ref{theorem_iterators_2}  $A \approx A'$ and $B \approx B'$, from which it 
easily follows that $A \cup B \approx A'\cup B'$ and therefore by Theorem~\ref{theorem_iterators_1} (2) we have 
$\val(A \cup B) = \val(A'\cup B')$, i.e., $\val(\gamma_\oplus(A,B)) = \val(\gamma_\oplus(A',B'))$ and so $\gamma_\oplus$ satisfies condition ($\heartsuit$).
 Essentially he same proof also shows that $\gamma_\otimes$ satisfies ($\heartsuit$).

Once it is known that the  operation $\oplus$ exists then the remaining properties 
of $\oplus$ listed in Theorem~\ref{theorem_am_1} follow from the corresponding properties of the union operation 
$\cup$ (for example, that it is associative and commutative).

The following shows that $\gamma_\oplus$ satisfies condition ($\heartsuit$).

\begin{lemma}\label{lemma_am_2}
If $(A,B)$ and $(A',B')$ are elements of $D_\oplus$ with $\val'(A,B) = \val'(A',B')$ then $\val(A \cup B) = \val(A' \cup B'))$.
\end{lemma}

\proof 
Consider finite sets $A$ and $A'$ with $\val(A) = \val(A')$ and a finite set $B$ disjoint 
from $A$ and $A'$. Let $\mathcal{S} = \{ C \in \mathcal{P}(B) : \val(A \cup C) = \val(A' \cup C) \}$.
Then $\varnothing \in \mathcal{S}$, since $\val(A \cup \varnothing) = \val(A) = \val(A') = \val(A' \cup \varnothing)$.
Let $C \in \proper{\mathcal{S}}$ and $b \in B \setminus C$. Then
\begin{eqnarray*}
\val(A \cup (C \cup \{b\})) &=& \val((A \cup C) \cup \{b\}) =  f(\val(A \cup C))\\ 
&=& f(\val(A' \cup C)) = \val((A' \cup C) \cup \{b\}) = \val(A' \cup (C \cup \{b\})) 
\end{eqnarray*}
and hence $C \cup \{b\} \in \mathcal{S}$. Thus $\mathcal{S}$ is an inductive $B$-system and so $B \in \mathcal{S}$.
Therefore $\val(A \cup B) = \val(A' \cup B)$.

For any set $C$ and any element $d$ put $C_d = C \times \{d\}$ (and so $C_d \approx C$). 
Now let $(A,B),\, (A',B') \in \gamma_\oplus$ with $\val'(A,B) = \val'(A',B')$, and choose distinct elements $\triangleleft$ and 
$\triangleright$; then $\val(A \cup B ) = \val(A_\triangleleft \cup B_\triangleleft)$ (since 
$A \cup B \approx A_\triangleleft \cup B_\triangleleft$), 
$\val(A'_\triangleright \cup B'_\triangleright) = \val(A' \cup B')$ 
(since $A'_\triangleright \cup B'_\triangleright \approx A' \cup B'$),
$\val(A_\triangleleft) = \val(A'_\triangleright)$ (since 
$A_\triangleleft \approx A$ and  $A' \approx A'_\triangleright)$)
and $\val(B_\triangleleft) = \val(B'_\triangleright)$ (since
$B_\triangleleft \approx B$ and  $B' \approx B'_\triangleright$), which gives us the following data: 
\begin{evlist}{16pt}{6pt}
\item[--] $\val(A \cup B) = \val(A_\triangleleft \cup B_\triangleleft)$, 
\item[--] $\val(A_\triangleleft) = \val(A'_\triangleright)$ and $B_\triangleleft$ is disjoint from both $A_\triangleleft$ 
and $A'_\triangleright$,
\item[--] $\val(B_\triangleleft) = \val(B'_\triangleright)$ and $A'_\triangleright$ is disjoint from both $B_\triangleleft$ 
and $B'_\triangleright$,
\item[--] $\val(A'_\triangleright \cup B'_\triangleright) = \val(A' \cup B')$. 
\end{evlist}
Thus by two applications of the first part of the proof
\begin{eqnarray*} 
\val(A \cup B) &=& \val(A_\triangleleft \cup B_\triangleleft) = \val(A'_\triangleright \cup B_\triangleleft)\\
&=& \val(B_\triangleleft \cup A'_\triangleright) = \val(B'_\triangleright \cup A'_\triangleright)
= \val(A'_\triangleright \cup B'_\triangleright) = \val(A' \cup B')\;.\eop
\end{eqnarray*}

\textit{Proof of Theorem~\ref{theorem_am_1}:}\enskip
By Lemma~\ref{lemma_am_2}  
condition ($\heartsuit$) holds for the mapping $\gamma_\oplus$ and thus by Proposition~\ref{prop_am_111} there exists a binary operation $\oplus$ on $X$ such that
$\val(A) \oplus \val(B) = \val(A \cup B)$ whenever $A$ and $B$ are disjoint finite sets. 
Moreover, $\oplus$ is the unique operation with this property.
We show that $\oplus$ is associative and commutative: Let $x_1,\,x_2,\,x_3 \in X$; then by 
Lemma~\ref{lemma_iterators_4} there exist finite sets $A_1,\,A_2$ and $A_3$ with $x_1 = \val(A_1)$, $x_2 = \val(A_2)$ and 
$x_3 = \val(A_3)$ and by Proposition~\ref{prop_fs_112} (5) and Theorem~\ref{theorem_iterators_1} (2) we can assume that these sets are disjoint. Clearly
$(A_1 \cup A_2) \cup A_3 \approx A_1 \cup (A_2 \cup A_3)$ and therefore
\begin{eqnarray*}
(x_1 \oplus x_2) \oplus x_3 &=& (\val(A_1) \oplus \val(A_2)) \oplus \val(A_3) 
 = \val(A_1 \cup A_2) \oplus \val(A_3) \\
 &=& \val((A_1 \cup A_2) \cup A_3) = \val(A_1 \cup (A_2 \cup A_3)) \\
 &=& \val(A_1) \oplus \val(A_2 \cup A_3) = \val(A_1) \oplus (\val(A_2) \oplus \val(A_3))\\ 
&=& x_1 \oplus (x_2 \oplus x_3) \;.
\end{eqnarray*}
In the same way $\oplus$ is commutative. Let $x_1,\,x_2 \in X$; then, as above there exist  disjoint finite sets $A_1$ and $A_2$ with $x_1 = \val(A_1)$   and 
$x_2 = \val(A_2)$. Also clearly $A_1 \cup A_2 \approx A_2 \cup A_1$ and hence
\begin{eqnarray*}
x_1 \oplus x_2 &=& \val(A_1) \oplus \val(A_2)\\ 
&=& \val(A_1 \cup A_2) = \val(A_2 \cup A_1) = \val(A_2) \oplus \val(A_1) = x_2 \oplus x_1 \;.
\end{eqnarray*}
Moreover, if $x \in X$ and $A$ is a finite set with $x = \val(A)$ then
\[ 
x \oplus x_0 = \val(A) \oplus \val(\varnothing)  = \val(A \cup \varnothing) = \val(A) =  x\,,
\]
and so $x \oplus x_0 = x$ for all $x \in X$.

Let $x_1,\,x_2 \in X$, and so by Lemma~\ref{lemma_iterators_4} there exist finite sets $A$ and $B$ such that
$x_1 = \val(A)$ and $x_2 = \val(B)$. By Theorem~\ref{theorem_fs_4} there either exists an injective mapping 
$g : A \to B$ or an injective mapping $h : B \to A$. Assume the former holds and put $B' = g(A)$ and 
$C = B \setminus B'$. Then $B'$ and $C$ are disjoint and $B = B' \cup C$; moreover, $A \approx B'$ (since $g$ 
considered as a mapping from $A$ to $B'$ is a bijection) and so by Theorem~\ref{theorem_iterators_1} (2) $\val(A) = \val(B')$.
Thus, putting $x = \val(C)$, it follows that
$x_2 = \val(B) = \val(B' \cup C) = \val(B') \oplus \val(C) = \val(A) \oplus \val(C) = x_1 \oplus x$. On the other hand, if 
there exists an injective mapping $h : B \to A$ then the same argument shows that $x_1 = x_2 \oplus x$ for some 
$x \in X$.

Now to (a0) and (a1), and we have seen above that (a0) holds. Let $x,\,x' \in X$, so by Lemma~\ref{lemma_am_1} 
there exist disjoint finite sets $A$ and $B$ with $x = \val(A)$ and $x' = \val(B)$. Let $b \notin A \cup B$; then 
\begin{eqnarray*}
x \oplus f(x') &=& \val(A) \oplus f(\val(B)) = \val(A) \oplus \val(B \cup \{b\}) = \val(A \cup (B \cup \{b\})) \\
 &=& \val((A \cup B) \cup \{b\}) =  f(\val(A \cup B)) = f(\val(A) \oplus \val(B)) = f(x \oplus x') 
\end{eqnarray*}
and hence (a1) holds. If $\oplus'$ is another binary operation on $X$ satisfying (a0) and (a1) then it is easy 
to see that $X_0 = \{ x' \in X : x \oplus' x' = x \oplus x'\ \mbox{for all $x \in X$} \}$ is an $f$-invariant 
subclass of $X$ containing $x_0$. Hence $X_0 = X$, since $\,\mathbf{I}$  is minimal, which implies that 
${\oplus'} = {\oplus}$.

This completes the proof of Theorem~\ref{theorem_am_1}.
\eop

Theorem~\ref{theorem_am_2} will be dealt with in a similar manner. 
We obtain the operation $\otimes$ by showing that $\gamma_\otimes$ satisfies condition ($\heartsuit$).

As with the addition $\oplus$, once it is known that the operation $\otimes$ exists  
then the remaining properties of $\otimes$ listed in Theorem~\ref{theorem_am_2} follow from the corresponding 
properties of the cartesian product operation $\times$ (for example, that it is (modulo the relation $\approx$) associative and commutative) and from 
the relationship between $\cup$ and $\times$.

The following shows that $\gamma_\otimes$ satisfies condition ($\heartsuit$).

\begin{lemma}\label{lemma_am_3}
If $A$, $B$, $A'$, $B'$ are finite sets with $\val(A) = \val(A')$ and $\val(B) = \val(B')$
 then $\val(A \times B) = \val(A' \times  B')$.
\end{lemma}

\proof 
Consider finite sets $A$ and $A'$ with $\val(A) = \val(A')$ and let $B$ be any finite set. Put 
$\mathcal{S} = \{ C \in \mathcal{P}(B) : \val(A \times C) = \val(A' \times C) \}$.
Then $\varnothing \in \mathcal{S}$, since $A \times \varnothing = \varnothing = A' \times \varnothing$ and so 
$\val(A \times \varnothing) = \val(A' \times \varnothing)$. 
Let $C \in \proper{\mathcal{S}}$ and let $b \notin B \setminus C$. Then the sets $A \times C$ and 
$A \times \{b\}$ are disjoint and $A \times (C \cup \{b\}) = (A \times C) \cup (A \times \{b\})$. It follows 
that
\[ 
\val(A \times (C \cup \{b\})) = \val((A \times C) \cup (A \times \{b\}) = \val(A \times C) \oplus \val(A \times \{b\}) 
\]
and in the same way $\val(A' \times (C \cup \{b\})) = \val(A' \times C) \oplus \val(A' \times \{b\})$. Clearly 
$A \times \{b\} \approx A$ and so by Theorem~\ref{theorem_iterators_1} (2) $\val(A \times \{b\}) = \val(A)$, and in the same 
way $\val(A' \times \{b\}) = \val(A')$. Therefore
\begin{eqnarray*}
\val(A \times (C \cup \{b\})) &=& \val(A \times C) \oplus \val(A \times \{b\}) = \val(A' \times C) \oplus \val(A)\\ 
&=& \val(A' \times C) \oplus \val(A' \times \{b\}) = \val(A' \times (C \cup \{b\})) 
\end{eqnarray*}
and so $B \cup \{b\} \in \mathcal{S}$. Hence $\mathcal{S}$ is an inductive $B$-system and so $B \in \mathcal{S}$. 
Therefore $\val(A \times B) = \val(A' \times B)$.
Now let $A$, $B$, $A'$, $B'$ be finite sets with $\val(A) = \val(A')$ and $\val(B) = \val(B')$.
Then clearly we have
$A' \times B \approx B \times A'$ and $A' \times B' \approx B' \times A'$ and hence by 
Theorem~\ref{theorem_iterators_1} (2) $\val(A' \times B) = \val(B \times A')$ and $\val(A' \times B') = \val(B' \times A')$.
Hence by the first part 
\[ 
\val(A \times B) = \val(A' \times B) = \val(B \times A') = \val(B' \times A') = \val(A' \times B')\;.\eop
\]

\textit{Proof of Theorem~\ref{theorem_am_2}:}\enskip
By Lemma~\ref{lemma_am_3}  
condition ($\heartsuit$) holds for the mapping $\gamma_\otimes$ and thus by Proposition~\ref{prop_am_111} there exists a binary operation $\otimes$ on $X$ such that
$\val(A) \otimes \val(B) = \val(A \times B)$ whenever $A$ and $B$ are finite sets. 
Moreover, $\otimes$ is the unique operation with this property.

We show that $\otimes$ is associative and commutative: Let $x_1,\,x_2,\,x_3 \in X$; then by 
Lemma~\ref{lemma_iterators_4} there exists finite sets $A_1,\,A_2,\,A_3$ with $x_1 = \val(A_1)$, $x_2 = \val(A_2)$ and 
$x_3 = \val(A_3)$. Now it is easy to check that $(A_1 \times A_2) \times A_3 \approx A_1 \times (A_2 \times A_3)$
and so by Theorem~\ref{theorem_iterators_1} (2) $\val((A_1 \times A_2) \times A_3) = \val(A_1 \times (A_2 \times A_3))$.
Therefore 
\begin{eqnarray*}
(x_1 \otimes x_2) \otimes x_3 &=& (\val(A_1) \otimes \val(A_2)) \otimes \val(A_3) = \val(A_1 \times A_2) \otimes \val(A_3) \\
 &=& \val((A_1 \times A_2) \times A_3) = \val(A_1 \times (A_2 \times A_3)) \\
 &=& \val(A_1) \otimes \val(A_2 \times A_3) = \val(A_1) \otimes (\val(A_2) \otimes \val(A_3))\\ 
&=& x_1 \otimes (x_2 \otimes x_3) 
\end{eqnarray*}
which shows $\otimes$ is associative. Let $x_1,\,x_2 \in X$; by Lemma~\ref{lemma_iterators_4} there exist finite sets 
$A_1$ and $A_2$ with $x_1 = \val(A_1)$ and $x_2 = \val(A_2)$. Then by Theorem~\ref{theorem_iterators_1} (2) we have
$\val(A_1 \times A_2) = \val(A_2 \times A_1)$, since clearly $A_1 \times A_2 \approx A_2 \times A_1$. Thus
\begin{eqnarray*}
x_1 \otimes x_2 &=& \val(A_1) \otimes \val(A_2)\\ 
&=& \val(A_1 \times A_2) = \val(A_2 \times A_1) = \val(A_2) \otimes \val(A_1) = x_2 \otimes x_1 
\end{eqnarray*}
which shows that $\otimes$ is also commutative. 

Let $x \in X$, so by Lemma~\ref{lemma_iterators_4} there exists a finite set $A$ with $x = \val(A)$. Then
\[
x \otimes x_0 = \val(A) \otimes \val(\varnothing) = \val(A \times \varnothing) = \val(\varnothing) = x_0\;.
\]
Moreover, if $x \ne x_0$ then $A \ne \varnothing$, so if $a$ is any element then by Theorem~\ref{theorem_iterators_1} (2) $\val(A \times \{a\}) = \val(A)$, since 
$A \times \{a\} \approx A$, and hence
\begin{eqnarray*}
 x \otimes f(x_0) &=& \val(A) \otimes f(\val(\varnothing)) = \val(A) \otimes \val(\varnothing \cup \{a\})\\
&=& \val(A) \otimes \val(\{a\}) = \val(A \times \{a\}) = \val(A) = x \;.
\end{eqnarray*}
Thus $x \otimes x_0 = x_0$ for each $x \in X$ and $x \otimes f(x_0) = x$ for each $x \ne x_0$ (and note that the first statement is 
(m0)).

Now for the distributive law. Let $x,\,x_1,\,x_2 \in X$. There exists a finite set $A$ with $x = \val(A)$ and 
disjoint finite sets $B$ and $C$ with $x_1 = \val(B)$ and $x_2  = \val(C)$. Then $A \times (B \cup C)$ is the disjoint union of 
$A \times B$ and $A \times C$ and thus
\begin{eqnarray*}
(x \otimes x_1) \oplus (x \otimes x_2) &=&  (\val(A) \otimes \val(B)) \oplus (\val(A) \otimes \val(C))\\ 
&=&  \val(A \times B) \oplus \val(A \times C) =  \val((A \times B) \cup (A \times C))\\ 
&=&  \val(A \times (B \cup C)) = \val(A) \otimes \val(B \cup C)\\ 
&=& \val(A) \otimes (\val(B) \oplus \val(C)) =  x \otimes (x_1 \oplus x_2)\;.
\end{eqnarray*}

We have already seen that (m0) holds and, since $f(x_0)$ is an identity element, (m1) is a special case of the distributive 
law: Let $x,\,x' \in X$; then by (a0) and (a1) and since $\oplus$ is commutative it follows that
$f(x') = f(x' \oplus x_0) = x' \oplus f(x_0) = f(x_0) \oplus x'$, and hence 
$x \otimes f(x') = x \otimes (f(x_0) \oplus x')
= (x \otimes f(x_0)) \oplus (x \otimes x') = x \oplus (x \otimes x')$, which is (m1). Finally, if $\otimes'$ is 
another binary operation satisfying (m0) and (m1) then it is easy to see that 
$X_0 = \{ x' \in X : x \otimes' x' = x \otimes x'\ \mbox{for all $x \in X$} \}$
is a $f$-invariant subclass of $X$ containing $x_0$. Hence $X_0 = X$, since $\,\mathbf{I}$  is minimal, which implies 
that ${\otimes'} = {\otimes}$.
This completes the proof of Theorem~\ref{theorem_am_2}. \eop
 
We next give some results 
about the operation $\oplus$ for special cases of $\,\mathbf{I}$.

\begin{proposition}\label{prop_am_1}
If $f$ is injective then the cancellation law holds for $\oplus$ (meaning  that $x_1 = x_2$ whenever 
$x_1 \oplus x = x_2 \oplus x$ for some $x \in X$). In particular, $x \ne x \oplus x'$ for all $x,\, x' \in X$ with $x' \ne x_0$ 
(since $x = x \oplus x_0$).
\end{proposition}

\proof 
Let $x_1,\,x_2 \in X$ with $x_1 \ne x_2$ and let $X_0 = \{ x \in X : x_1 \oplus x \ne x_2 \oplus x \}$; then 
$x_0 \in X_0$, since by (a0) $x_1 \oplus x_0 = x_1 \ne x_2 = x_2 \oplus x_0$. Let $x \in X_0$, then by (a1), and 
since $f$ is injective, $x_1 \oplus f(x) = f(x_1 \oplus x) \ne f(x_2 \oplus x) = x_1 \oplus f(x)$, i.e.,
$f(x) \in X_0$. Thus $X_0$ is a $f$-invariant subset of $X$ containing $x_0$ and so $X_0 = X$, since
$\,\mathbf{I}$ is minimal. Hence if $x_1 \ne x_2$ then $x_1 \oplus x \ne x_2 \oplus x$ for all $x \in X$, which 
shows that the cancellation law holds for $X$. 
\eop

\begin{proposition}\label{prop_am_2}
If $x_0 \in f(X)$ (and so by Theorem~\ref{theorem_iterators_4} and Proposition~\ref{prop_iterators_2} $X$ is finite and 
$f$ is bijective) then $(X,\oplus,x_0)$ is an abelian group: For each $x \in X$ there exists $x' \in X$ such that 
$x \oplus x' = x_0$. Moreover, $X$ is the group generated by the element $f(x_0)$.
\end{proposition}

\proof 
 By Lemma~\ref{lemma_iterators_4} there exists a non-empty finite set $A$ with $\val(A) = x_0$ and then for each finite set $C$ there exists 
$B \subset A$ with $\val(B) = \val(C)$. By Lemma~\ref{lemma_iterators_4} there exists a finite set $C$ with $x = \val(C)$ and hence there also
exists a finite set $B \subset A$ with $x = \val(B)$. Put $B' = A \setminus B$ and let $x' = \val(B')$. 
Then $B$ and $B'$ are disjoint and hence 
$x \oplus x' = \val(B) \oplus \val(B') = \val(B \cup B') = \val(A) = x_0$.
Let $X_0$ be the least subgroup of $X$ containing $f(x_0)$. Then $x_0 \in X_0$ and if $x \in X_0$ then 
by (a0) and (a1) $f(x) =f(x \oplus x_0) = x \oplus f(x_0)$ and hence  $f(x) \in X_0$. Thus $X_0$ is a
$f$-invariant subset of $X$ containing $f(x_0)$ and so $X_0 = X$. Therefore $X$ is the group generated by  $f(x_0)$.
\eop

Until further notice let $\,\mathbf{I} = (X,f,x_0)$ be a Peano iterator and let $\le$ be the unique total order on $X$ compatible with $\,\mathbf{I}$ as introduced
in  Theorem~\ref{theorem_iterators_199}. Thus $\le$ is the unique total order on $X$ such that $x \le f(x)$ for all $x \in X$.
We now give the usual characterisation of $\le$ in terms of the addition $\oplus$.
 \begin{proposition}\label{prop_am_33}
Let $x,\,y \in X$; then $y \le x$ if and only if there exists $z \in X$ with $x = y \oplus z$.
\end{proposition}

\proof Let $x,\,y \in X$; then by definition $y \le x$ if and only if  there exist $A,\,B \in \fin$ with $x = \val(A)$, $y = \val(B)$ such that
$B \preceq A$. Let $y \le x$; then by Proposition~\ref{prop_fs_112} (4) and Theorem~\ref{theorem_iterators_1} (2) there exist $A,\, B \in \fin$ with
$x = \val(A)$, $y = \val(B)$ and $B \subset A$. Let $C = A \setminus B$ and put $z = \val(C)$. Then $A$ is the disjoint union of $B$ and $C$ and hence
$x = y \oplus z$. Suppose conversely that there exists $z \in X$ with $x = y \oplus z$. Let $B,\,C \in \fin$ be disjoint with $y = \val(B)$ and
$z = \val(C)$. Then $x = y \oplus z = \val(B \cup C)$ and $B \preceq B \cup C$ and thus $y \le x$. \eop

\begin{lemma}\label{lemma_am_5}
If $x_1,\,x_2 \in X$ then $x_1 \oplus x_2 = x_0$ if and only if $x_1 = x_2 = x_0$.
\end{lemma}

\proof Suppose that $x_1 \oplus x_2 = x_0$. Now there exist disjoint finite sets $A_1,\,A_2$ with $\val(A_1) = x_1$ and $\val(A_2) = x_2$ and then 
$x_1 \oplus x_2 = \val(A_1 \cup A_2) = \val(\varnothing)$. Therefore $A_1 \cup A_2 \approx \varnothing$, since $\,\mathbf{I}$ is a Peano iterator.
Thus $A_1 = A_2 = \varnothing$, i.e., $x_1 = x_2 = x_0$. The converse holds trivially
\eop

\begin{proposition}\label{prop_am_3}
For all $x_1,\,x_2 \in X$ with $x_1 \ne x_2$ exactly one of the following two
statements holds:

There exists  a unique $x \in X$ such that $x_1 = x_2 \oplus x$.

There exists a unique $x' \in X$ such that $x_2 = x_1 \oplus x'$.
\end{proposition}

\proof
Note that if there exists $x \in X$ with $x_1 = x_2 \oplus x$ then by Proposition~\ref{prop_am_1} $x$ is unique
and in the same way if there exists $x' \in X$ with $x_2 = x_1 \oplus x'$ then $x'$ is unique. Also Theorem~\ref{theorem_am_1} states that at least one 
of the statements holds.
Suppose there exists  both $x \in X$ with $x_1 = x_2 \oplus x$ and $x' \in X$ with  $x_2 = x_1 \oplus x'$.
Then $x_1 = x_2 \oplus x = (x_1 \oplus x') \oplus x = x_1 \oplus (x' \oplus x)$. Thus by Proposition~\ref{prop_am_1} $x' \oplus x = x_0$ and by
Lemma~\ref{lemma_am_5} it then follows that $x= x' = x'$. But this implies  $x_1 = x_2$, contradicting the assumption that $x_1 \ne x_2$. Therefore exactly
one of the statements holds. \eop


For a Peano iterator $\,\mathbf{I}$ the operations $\oplus$ and $\otimes$ are all that are needed to develop an elementary theory of numbers in $\,\mathbf{I}$
which corresponds exactly to the standard elementary theory of numbers in $\Nat$. The concepts and proofs can all be taken over without change from the
case with $\,\mathbf{I} = \Nat$.  Note in particular that the classical proof that there are infinitely many primes in $\Nat$ now shows that the class
of prime elements in $\,\mathbf{I}$ cannot be a finite set.

We end the section by looking at the operation of exponentiation. Here we have to be more careful: For example, 
$2 \cdot 2 \cdot 2 = 2$ in $\Int_3$ and so $2^3$ is not well-defined if the exponent $3$ is considered as an 
element of $\Int_3$ (since we would also have to have $2^0 = 1$). However, $2^3$ does make sense if $2$ is 
considered as an element of $\Int_3$ and the exponent $3$ as an element of $\Nat$.

In general we will see that if $(Y,g,y_0)$ is a Peano iterator then we can define an element of $X$ which is 
`$x$ to the power of $y$' for each $x \in X$ and each $y \in Y$ and this operation has the properties which might 
be expected. 

In what follows let $\,\mathbf{J} = (Y,g,y_0)$ be a Peano iterator with $\val'$ the assignment of finite sets in $\mathbf{J}$. (As 
before $ \,\mathbf{I} = (X,f,,x_0)$ is assumed to be minimal with $\val$ the assignment of finite sets in $\,\mathbf{I}$.)
Also let $\oplus$ and $\otimes$ be the operations given in Theorems~\ref{theorem_am_1} and \ref{theorem_am_2} for the iterator
$\,\mathbf{I}$.

\begin{theorem}\label{theorem_am_3}
There exists a unique operation ${\uparrow} : X \times Y \to X$ such that  
\[ 
\val(A) \uparrow \val'(B) = \val(A^B) 
\]
for all finite sets $A$ and $B$. This operation ${\uparrow}$ satisfies
\[     
x \uparrow (y_1 \oplus y_2) = (x \uparrow y_1) \otimes (x \uparrow y_2) 
\]
for all $x \in X$ and all $y_1,\,y_2 \in Y$ and
\[     
(x_1 \otimes x_2) \uparrow y = (x_1 \uparrow y) \otimes (x_2 \uparrow y) 
\]
for all $x_1,\,x_2 \in X$ and $y \in Y$. Moreover, ${\uparrow}$ is the unique operation such that
\begin{evlist}{15pt}{6pt}
\item[(e0)]  $x \uparrow y_0 = f(x_0)$ for all $x \in X$. 

\item[(e1)]  $x \uparrow g(y) = x \otimes (x \uparrow y)$ for all $x \in X$, $y \in Y$.
\end{evlist}
\end{theorem}

\begin{lemma}\label{lemma_am_4}
If $B,\,C$ are finite sets with $\val(B) = \val(C)$ then for all finite sets $A$ we have $\val(B^A) = \val(C^A)$. 
\end{lemma}

\proof
Let $B$ and $C$ be finite sets with $\val(B) = \val(C)$, let $A$ be a finite set and put 
$\mathcal{S} = \{ D \in \mathcal{P}(A) : \val(B^D) = \val(C^D) \}$. Then $\varnothing \in \mathcal{S}$, since
$\val(B^\varnothing) =\val(C^\varnothing)$. (For any set $X$ the set $X^\varnothing$ consists of the single element $\{\varnothing\}$.)
Let $D \in \proper{\mathcal{S}}$ and $a \in A \setminus D$. Now $\val(B^D) = \val(C^D)$ 
(since $D \in \mathcal{S}$) and $\val(B) = \val(C)$; therefore  by Lemma~\ref{lemma_am_3} and Theorem~\ref{theorem_iterators_1} (2) 
\[\val(B^{D \cup \{a\}}) = \val(B^D \times B)  = \val(C^D \times C) = \val(C^{D \cup \{a\}})\] 
(since $E^{D\cup \{a\}} \approx E^D \times E$ for each set $E$), and so $D \cup \{a\} \in \mathcal{S}$.
Hence $\mathcal{S}$ is an inductive $A$-system and so $A \in \mathcal{S}$. Thus $\val(B^A) = \val(C^A)$.
\eop

\textit{Remark:}\enskip
If $B,\, C$ are finite sets with $\val(B) = \val(C)$ then $\val(A^B) = \val(A^C)$ does not hold in general for a finite 
set $A$.

\medskip
\textit{Proof of Theorem~\ref{theorem_am_3}:}\enskip
Let $A_1,\,A_2,\,B_1,\,B_2$ be finite sets with $\val(A_1) = \val(A_2)$ and $\val'(B_1) = \val'(B_2)$; then by 
Lemma~\ref{lemma_am_4} $\val(A_1^{B_1}) = \val(A_2^{B_1})$ and by Theorem~\ref{theorem_iterators_2} $B_1 \approx B_2$. 
Since $B_1 \approx B_2$ it follows that $A_2^{B_1} \approx A_2^{B_2}$ and then by Theorem~\ref{theorem_iterators_1} (2)
$\val(A_2^{B_1}) = \val(A_2^{B_2})$. This shows that $\val(A_1^{B_1}) = \val(A_2^{B_2})$. Therefore by 
Lemma~\ref{lemma_iterators_4} we can define $x \uparrow y$ to be $\val(A^B)$, where $A$ and $B$ are any finite sets 
with $x = \val(A)$ and $y = \val'(B)$. Then $\val(A) \uparrow \val'(B) = \val(A^B)$ 
for all finite sets $A$ and $B$ and this requirement clearly determines ${\uparrow}$ uniquely.
 
Let $x \in X$ and $y_1,\,y_2 \in Y$; then by Lemma~\ref{lemma_am_1} there exists a disjoint pair $(B_1,B_2)$ 
with  $(y_1,y_2) = \val'(B_1,B_2)$ and by Lemma~\ref{lemma_iterators_4} there exists a finite set $A$ with $x = \val(A)$.
Moreover, it is easily checked that $A^{B_1 \cup B_2} \approx A^{B_1} \times A^{B_2}$ and thus by 
Theorem~\ref{theorem_iterators_1} (2)
\begin{eqnarray*}
x \uparrow (y_1 \oplus y_2)  
  &=& \val(A) \uparrow (\val'(B_1) \oplus \val(B_2))\\ 
&=& \val(A) \uparrow \val'(B_1 \cup B_2) = \val(A^{B_1 \cup B_2}) = \val(A^{B_1} \times A^{B_2})\\
&=& \val(A^{B_1}) \otimes \val(A^{B_2}) = (x \uparrow y_1) \otimes (x \uparrow y_2) \;.
\end{eqnarray*}
Now let $x_1,\,x_2 \in X$ and $y \in Y$. By Lemma~\ref{lemma_iterators_4} there exist finite sets $A_1,\,A_2$ and $B$ 
such that $x_1 = \val(A_1)$, $x_2 = \val(A_2)$ and $y = \val'(B)$ and $(A_1\times A_2)^B  \approx A_1^B \times A_2^B$. 
Thus by Theorem~\ref{theorem_iterators_1} (2)
\begin{eqnarray*}
(x_1 \otimes x_2) \uparrow y  &=& (\val(A_1) \otimes \val(A_2)) \uparrow \val'(B)\\ 
&=& \val(A_1\times A_2) \uparrow \val'(B) = \val((A_1\times A_2)^B) = \val(A_1^B \times A_2^B)\\
&=& \val(A_1^B) \otimes \val(A_2^B) = (x_1 \uparrow y) \otimes (x_2 \uparrow y) \;.
\end{eqnarray*}
It remains to consider the properties (e0) and (e1). Now for each finite set $A$ we have
$\val(A) \uparrow \val'(\varnothing) = \val(A^\varnothing)= \val(\{\varnothing\}) = f(x_0)$ and hence 
$x \uparrow y_0 = f(x_0)$ for each $x \in X$, i.e., (e0) holds. Let $A$ and $B$ be finite sets and let 
$b \notin B$. Then, since $A^{B\cup \{b\}} \approx A \times A^B$, it follows from Theorem~\ref{theorem_iterators_1} (2) that
\begin{eqnarray*}
\val(A) \uparrow g(\val'(B)) &=& \val(A) \uparrow \val'(B\cup \{b\}) = \val(A^{B\cup \{b\}}) = \val(A \times A^B)\\ 
&=& \val(A) \otimes \val(A^B) = \val(A) \otimes (\val(A) \uparrow \val'(B)) 
\end{eqnarray*}
and this shows 
$x \uparrow g(h) = x \otimes (x \uparrow y)$ for all $x \in X$, $y \in Y$, i.e., (e1) holds. Finally, if 
${\uparrow}'$ is another operation satisfying (e0) and (e1) then 
\[
Y_0 = \{ y \in Y : x \uparrow' y = x \uparrow y\ \mbox{for all $x \in X$} \}
\]
is a $g$-invariant subset of $Y$ containing $y_0$. Therefore $Y_0 = Y$, since $(Y,,y_0)$  is minimal, which 
implies that ${\uparrow'} = {\uparrow}$.
\eop


\startsection{Another take on addition and multiplication}
\label{mam}

In the following again let $\,\mathbf{I} =(X,f,x_0)$ be a minimal iterator with $\val$ the assignment of finite sets in $\,\mathbf{I}$. 
In this section we give alternative proofs for Theorems \ref{theorem_am_1} and \ref{theorem_am_2}. 

In Section~\ref{am} only the single assignment $\val$ was used. Here we make use of a family of assignments 
$\{ \val_x : x \in X \}$, which arise as follows: For each $x \in X$ there is the iterator $\mathbf{I}_x = (X,f,x)$ (which 
will usually not be minimal) and the assignment of finite sets in $\,\mathbf{I}_x$ will be denoted by $\val_x$. Thus 
$\val_x(\varnothing) = x$ and $\val_x(A \cup \{a\}) = f(\val_x(A))$ whenever $A$ is a 
finite set and $a \notin A$. In particular we have $\val = \val_{x_0}$. 
Now it is more convenient to repackage the information given by the assignments $\val_A$, $x \in X$, by introducing 
for each finite set $A$ the mapping $f_A : X \to X$ with $f_A(x) = \val_x(A)$ for all $x \in X$, and so 
$\val(A) = \val_{x_0}(A) = f_A(x_0)$.

Consider disjoint finite sets $A$ and $B$; then $\val(A \cup B)$ can be thought of as the element of $X$  obtained by iterating $f$ through the
elements of $A\cup B$ starting with $x_0$. This element can also be determined by first iterating $f$ through the elements of $B$ starting with $x_0$, giving
the result $\val(B)$ and then iterating $f$ through the elements of $A$, but starting with the element $\val(B)$ and not with $x_0$.
The result is $\val_z(A)$, where $z = \val(B)$, and $\val_z(A) = f_A(z) = f_A(\val(B)) = f_A(f_B(x_0)) = (f_A \circ f_B)(x_0)$, 
and so we would expect that $\val(A \cup B) =  (f_A \circ f_B)(x_0)$. But if $\oplus$ is the operation given by 
Theorem~\ref{theorem_am_1} then $\val(A \cup B) = \val(A) \oplus \val(B) = f_A(x_0) \oplus f_B(x_0)$, which suggests 
that the following should hold:
\begin{evlist}{26pt}{6pt}
\item[($\diamondsuit$)]\enskip
$f_A(x_0) \oplus f_B(x_0) = (f_A \circ f_B)(x_0)$ whenever $A$ and $B$ are disjoint finite sets.
\end{evlist}
It will be seen later that ($\diamondsuit$) does hold. What is perhaps more important, though, is that
($\diamondsuit$) can actually be used to define $\oplus$, as we now explain.

Denote by $\,\Self{X}$ 
the set of all mappings from $X$ to itself and so $A \mapsto f_A$ defines a mapping from $\fin$
to $\,\Self{X}$. 
 Then $(\Self{X},\circ,\id_X)$, where $\circ$ 
is functional composition and $\id_X : X \to X$ is the identity mapping, is a monoid. (A \definition{monoid} is 
any triple $(M,\bullet,e)$ consisting of a class $M$, an associative operation $\bullet$ on $M$ and an identity element 
$e \in M$ satisfying $a \bullet e = e \bullet a = a$ for all $a \in M$.) Lemma~\ref{lemma_mam_6} shows that
\[ 
M_f = \{ u \in \Self{X} : \mbox{$u = f_A$ for some finite set $A$} \}
\]
is a submonoid of $(\Self{X},\circ, \id_X)$, meaning that $\id_X \in M_f$ and $u_1 \circ u_2 \in M_f$ for all 
$u_1,\,u_2 \in M_f$, and that this submonoid is commutative, i.e., $u_1 \circ u_2 = u_2 \circ u_1$ for all 
$u_1,\,u_2 \in M_f$. (The monoid $(\Self{X},\circ,\id_X)$ itself is not commutative except when $X = \{x_0\}$.)

Let $\Phi_{x_0} : M_f \to X$ be the mapping with $\Phi_{x_0}(u) = u(x_0)$ for each $u \in M_f$, and so in 
particular $\Phi_{x_0}(f_A) = f_A(x_0) = \val(A)$ for each finite set $A$. Lemma~\ref{lemma_mam_7} will show that 
$\Phi_{x_0}$ is a bijection, and therefore there exists a unique operation $\oplus$ on $X$ such that
\begin{evlist}{26pt}{6pt}
\item[($\heartsuit$)]\enskip
$\Phi_{x_0}(u) \oplus \Phi_{x_0}(v) = \Phi_{x_0}(u \circ v)$ for all $u,\,v \in M_f$.
\end{evlist}
This is how $\oplus$ will be defined below. Note that if $A$ and $B$ are (not necessarily disjoint) finite sets 
then by ($\heartsuit$)
\[ 
f_A(x_0) \oplus f_B(x_0) = \Phi_{x_0}(f_A) \oplus \Phi_{x_0}(f_B) = \Phi_{x_0}(f_A \circ f_B) = (f_A \circ f_B)(x_0)
\]
and so in particular ($\diamondsuit$) holds.

We now give the details of the approach outlined above.

\begin{lemma}\label{lemma_mam_1}
The mapping $A \mapsto f_A$  is the unique mapping from $\fin$ to $\Self{X}$ with $f_\varnothing = \id_X$ such that $f_{A \cup \{a\}} = f \circ f_A$ whenever $A$ is a finite set and $a \notin A$. 
\end{lemma}

\proof 
We have $f_\varnothing(x) = \val_x(\varnothing) = x = \id_X(x)$ for all $x \in X$, and thus $f_\varnothing = \id_X$.
Moreover, if $A$ is a finite set and $a \notin A$ then
\[ 
f_{A \cup \{a\}}(x) = \val_x(A \cup \{a\}) = f(\val_x(A)) = f(f_A(x)) = (f \circ f_A)(x)
\]
for all $x \in X$ and hence $f_{A \cup \{a\}} = f \circ f_A$. Finally, consider a further mapping 
$A \mapsto f'_A$ with $f'_\varnothing = \id_X$ and such that $f'_{A \cup \{a\}} = f \circ f'_A$ whenever $A$ is a 
finite set and $a \notin A$. Let $A$ be a  finite set and put $\mathcal{S} = \{ B \in \mathcal{P}(A) : f'_B = f_B \}$.
Then $\varnothing \in \mathcal{S}$, since $f'_\varnothing = \id_X = f_\varnothing$.
Let $B \in \proper{\mathcal{S}}$ (and so $f'_B = f_B$) and let $a \in A \setminus B$. Then
$f'_{B \cup \{a\}} = f \circ f'_B = f \circ f_B = f_{B \cup \{a\}}$ and therefore $A \cup \{a\} \in \mathcal{S}$.
Thus $\mathcal{S}$ is an inductive $A$-system and so $A \in \mathcal{S}$.
Hence $f'_A = f_A$.
\eop

The mapping $A \mapsto f_A$ will be called the \definition{functional assignment of finite sets in $\,\mathbf{I}$}. In particular 
$f_{\{a\}} = f$ for each element $a$, since 
$f_{\{a\}} = f_{\varnothing \cup \{a\}} = f \circ f_\varnothing = f \circ \id_X = f$.

\begin{lemma}\label{lemma_mam_2}
$f \circ f_A = f_A \circ f$  for each finite set $A$.
\end{lemma}

\proof 
Let $A$ be a finite set and put $\mathcal{S} = \{ B \in \mathcal{P}(A) : f \circ f_B = f_B \circ f \}$. Then $\varnothing \in \mathcal{S}$ since
$f \circ f_\varnothing = f \circ \id_X = f = \id_X \circ f = f_\varnothing \circ f$.
Let $B \in \proper{\mathcal{S}}$ and $a \in A \setminus B$. Then 
$f \circ f_{B \cup \{a\}} = f \circ f \circ f_B = f \circ f_B \circ f = f_{B \cup \{a\}} \circ f$
and so $B \cup \{a\} \in \mathcal{S}$. Thus $\mathcal{S}$ is an inductive $A$-system and so $A \in \mathcal{S}$.
Hence $f \circ f_A = f_A \circ f$.
\eop

The next result establishes an important relationship between $\val$ and the functional assignment.

\begin{proposition}\label{prop_mam_1}
If $A$ and $B$ are finite sets then $f_A = f_B$ holds if and only if $\val(A) = \val(B)$.
\end{proposition}

\proof 
By definition $\val(C) = f_C(x_0)$ for each finite set $C$, and so $\val(A) = \val(B)$ whenever $f_A = f_B$. Suppose 
conversely that $\val(A) = \val(B)$ and consider the set $X_0 = \{ x \in X : f_A(x) = f_B(x) \}$. Then $X_0$ is 
$f$-invariant, since if $x \in X_0$ then by Lemma~\ref{lemma_mam_2} 
$f_A(f(x)) = f(f_A(x)) = f(f_B(x)) = f_B(f(x))$, i.e., $f(x) \in X_0$. Also $x_0 \in X_0$, since 
$f_A(x_0) = \val(A) = \val(B) = f_B(x_0)$. Hence $X_0 = X$, since $\,\mathbf{I}$ is minimal. This shows that 
$f_A(x) = f_B(x)$ for all $x \in X$, i.e., $f_A = f_B$. 
\eop

There is another way of obtaining the functional assignment: Consider the iterator $\,\mathbf{I}_* =(\Self{X},f_*,\id_X)$, where 
$f_* : \Self{X} \to \Self{X}$ is defined by $f_*(h) = f \circ v$ for all $v \in \Self{X}$.

\begin{lemma}\label{lemma_mam_3}
If $\val_*$ is the assignment of finite sets in $\,\mathbf{I}_*$ then $\val_*(A) = f_A$ for each finite set $A$.
\end{lemma}

\proof 
By definition $\val_*(\varnothing) = \id_X$ and $\val_*(A \cup \{a\}) = f_*(\val_*(A)) = f \circ \val_*(A)$ for each finite set $A$ 
and each $a \notin A$. Thus by the uniqueness in Lemma~\ref{lemma_mam_1} $\val_*(A) = f_A$ for each 
finite set $A$.
\eop

\begin{proposition}\label{prop_mam_2}
If $A$ and $B$ are finite sets with $A \approx B$ then $f_A = f_B$.
\end{proposition}

\proof 
This follows from Theorem~\ref{theorem_iterators_1} (applied to $\val_*$) and Lemma~\ref{lemma_mam_3}. 
\eop

\begin{lemma}\label{lemma_mam_4}
If $A$ and $B$ are disjoint finite sets then $f_{A \cup B} = f_A \circ f_B$.
\end{lemma}

\proof 
Let $A$ and $B$ be finite sets and let $\mathcal{S} = \{ C \in \mathcal{P}(A) : f_{C \cup B} = f_C \circ f_B \}$. 
Then $\varnothing \in \mathcal{S}$ since $f_{\varnothing \cup B} = f_B = \id_X \circ f_B = f_\varnothing \circ f_B$.
Let $C \in \proper{\mathcal{S}}$ and let $a \in A \setminus C$; then $f_{C \cup B} = f_C \circ f_B$ and hence by Lemma~\ref{lemma_mam_1}
\[
f_{(C \cup \{a\}) \cup B} = f_{(C \cup B) \cup \{a\}} = f \circ f_{C \cup B} =  f \circ f_C \circ f_B  
= f_{C \cup \{a\}} \circ f_B \;.
\]
This shows that $C \cup \{a\} \in \mathcal{S}$. Thus $\mathcal{S}$ is an inductive $A$-system and so $A \in \mathcal{S}$. 
Hence $f_{A \cup B} = f_A \circ f_B$.
\eop

\begin{lemma}\label{lemma_mam_5}
(1)\enskip
If $f$ is bijective then $f_A$ is bijective for each finite set $A$. 

(2)\enskip
If $f$ is injective then $f_A$ is also injective for each finite set $A$. 
\end{lemma}

\proof
(1)\enskip
Let $A$ be a finite set and put $\mathcal{S} = \{ B \in \mathcal{P}(A) : \mbox{$f_B$ is bijective} \}$.
Then $\varnothing \in \mathcal{S}$, since $f_{\varnothing} = \id_X$ is bijective.
Consider $B \in \proper{\mathcal{S}}$ and let $a \in A \setminus B$. Then $f_{B \cup \{a\}} = f \circ f_B$ and so $f_{B \cup \{a\}}$, 
as the composition of two bijective mappings, is itself bijective, i.e., $B \cup \{a\} \in \mathcal{S}$.
Thus $\mathcal{S}$ is an inductive $A$-system and so $A \in \mathcal{S}$. Hence $f_A$ is bijective.

(2)\enskip 
Just replace `bijective' by `injective' in (1).
\eop

As above let $M_f = \{ u \in \Self{X} : \mbox{$u = f_A$ for some finite set $A$} \}$. Then in particular 
$\id_X \in M_f$, since $\id_X = f_\varnothing$, and $f \in M_f$, since $f = f_{\{a\}}$ for each element $a$. 
Moreover, if $f$ is injective (resp.\ bijective) then by Lemma~\ref{lemma_mam_5} each element in $M_f$ is 
injective (resp.\ bijective).

\begin{lemma}\label{lemma_mam_6}
For all $u_1,\,u_2 \in M_f$ we have $u_1 \circ u_2 \in M_f$ and $u_1 \circ u_2 = u_2 \circ u_1$. (Since also
$\id_X \in M_f$ this means that $M_f$ is a commutative submonoid of the monoid $(\Self{X},\circ,\id_X)$.)
\end{lemma}

\proof 
Let $u_1,\,u_2 \in M_f$ and so there exist finite sets $A$ and $B$ with $u_1 = f_A$ and $u_2 = f_B$. There then 
exists a disjoint pair $(A',B')$ with $(A',B') \approx (A,B)$ and hence by Proposition~\ref{prop_mam_2} and 
Lemma~\ref{lemma_mam_4}
\[ 
u_1 \circ u_2 = f_A \circ f_B = f_{A'} \circ f_{B'} = f_{A' \cup B'} = f_{B' \cup A'}
= f_{B'} \circ f_{A'} = f_B \circ f_A = u_2 \circ u_1\;, 
\]
i.e., $u_1 \circ u_2 = u_2 \circ u_1$. Moreover, since $u_1 \circ u_2 = f_{A' \cup B'}$ and 
$f_{A' \cup B'} \in M_f$, this also shows that $u_1 \circ u_2 \in M_f$.
\eop

As above let $\Phi_{x_0} : M_f \to X$ be the mapping with $\Phi_{x_0}(u) = u(x_0)$ for all $u \in M_f$. Then 
$\Phi_{x_0}(\id_X) = x_0$ and $\Phi_{x_0}(f_A) = f_A(x_0) = \val(A)$ for each finite set $A$. An important property 
of $\Phi_{x_0}$ is that 
\begin{evlist}{26pt}{6pt}
\item[$\mathrm{(\sharp)}$]
$\ u(\Phi_{x_0}(v)) = \Phi_{x_0}(u \circ v)\,$ for all $u,\,v \in M_f$, 
\end{evlist}
which holds since $u(\Phi_{x_0}(v)) = u(v(x_0)) = (u \circ v)(x_0) = \Phi_{x_0}(u \circ v)$. The special case
of this with $u = f$ gives us $f(\Phi_{x_0}(v)) = \Phi_{x_0}(f \circ v)$ for all $v \in M_f$.

\begin{lemma}\label{lemma_mam_7}
The mapping $\Phi_{x_0}$ is a bijection.
\end{lemma}

\proof 
If $g \in X$ then by Lemma~\ref{lemma_iterators_4} there exists a finite set $A$ with $x = \val(A)$ and it follows that
$\Phi_{x_0}(f_A) = f_A(x_0) = \val(A) = x$. Thus $\Phi_{x_0}$ is surjective. Now let $u_1,\,u_2 \in M_f$ with 
$\Phi_{x_0}(u_1) = \Phi_{x_0}(u_2)$. By the definition of $M_f$ there exist finite sets $A$ and $B$ with 
$u_1 = f_A$ and $u_2 = f_B$, and hence 
\[ 
\val(A) = \Phi_{x_0}(f_A) = \Phi_{x_0}(u_1) = \Phi_{x_0}(u_2)= \Phi_{x_0}(f_B) = \val(B)\;.
\]
Therefore by Proposition~\ref{prop_mam_1} $f_A = f_B$, i.e., $u_1 = u_2$, which shows that $\Phi_{x_0}$ is also 
injective. 
\eop

\textit{Proof of Theorem~\ref{theorem_am_1}:}\enskip
Since $\Phi_{x_0} : M_f \to X$ is a bijection there clearly exists a unique binary relation $\oplus$ on $X$ such 
that
\[ 
\Phi_{x_0}(u_1) \oplus \Phi_{x_0}(u_2) = \Phi_{x_0}(u_1 \circ u_2)
\]
for all $u_1,\,u_2 \in M_f$. The operation $\oplus$ is associative since $\circ$ has this property: 
If $x_1,\,x_2,\,x_3 \in X$ and $u_1,\,u_2,\,u_3 \in M_f$ are such that $x_j = \Phi_{x_0}(u_j)$ for each $j$ then 
\begin{eqnarray*}
(x_1 \oplus x_2) \oplus x_3 &=& (\Phi_{x_0}(u_1) \oplus \Phi_{x_0}(u_2)) \oplus \Phi_{x_0}(u_3) \\
&=& \Phi_{x_0}(u_1 \circ u_2) \oplus \Phi_{x_0}(u_3) = \Phi_{x_0}( (u_1 \circ u_2) \circ u_3)\\ 
&=& \Phi_{x_0}( u_1 \circ (u_2 \circ u_3)) =  \Phi_{x_0}(u_1) \oplus \Phi_{x_0}(u_2 \circ u_3)\\
&=& \Phi_{x_0}(u_1) \oplus (\Phi_{x_0}(u_2) \oplus \Phi_{x_0}(u_3))
=  x_1 \oplus (x_2 \oplus x_3) \;.
\end{eqnarray*}
In the same way $\oplus$ is commutative, since by Lemma~\ref{lemma_mam_6} the restriction of $\circ$ to $M_f$ 
has this property: If $x_1,\,x_2 \in X$ and $u_1,\,u_2 \in M_f$ are such that $x_1 = \Phi_{x_0}(u_1)$ and 
$x_2 = \Phi_{x_0}(u_2)$ then $u_1 \circ u_2 = u_2 \circ u_1$ and so 
\begin{eqnarray*}
x_1 \oplus x_2  &=& \Phi_{x_0}(u_1) \oplus \Phi_{x_0}(u_2)\\ 
&=& \Phi_{x_0}(u_1 \circ u_2) = \Phi_{x_0}(u_2 \circ u_1) = \Phi_{x_0}(u_2) \oplus \Phi_{x_0}(u_1) = x_2 \oplus x_1\;.
\end{eqnarray*}
Moreover, if $x \in X$ and $u \in M_f$ is such that $x = \Phi_{x_0}(u)$ then
\[ 
x \oplus x_0 = \Phi_{x_0}(u) \oplus \Phi_{x_0}(\id_X) = \Phi_{x_0}(u \circ \id_X) = \Phi_{x_0}(u) =  x\,,
\]
and so $x \oplus x_0 = x$ for all $x \in X$.

Let $x_1,\,x_2 \in X$; we next show that for some $x \in X$ either $x_1 = x_2 \oplus x$ or $x_2 = x_1 \oplus x$. 
Let $u_1,\,u_2 \in M_f$ be such that $x_1 = \Phi_{x_0}(u_1)$ and $x_2 = \Phi_{x_0}(u_2)$ and let $A$ and $B$ be 
finite sets with $u_1 = f_A$ and $u_2 = f_B$. By Theorem~\ref{theorem_fs_4} there either exists an injective 
mapping $p : A \to B$ or an injective mapping $q : B \to A$. Assume the former holds and put $B' = p(A)$ and 
$C = B \setminus B'$. Then $B'$ and $C$ are disjoint and $B = B' \cup C$; moreover, $A \approx B'$ (since $p$ 
considered as a mapping from $A$ to $B'$ is a bijection) and so by Proposition~\ref{prop_mam_2} $f_A = f_{B'}$. 
Thus, putting $x = \Phi_{x_0}(f_C)$, it follows that
\begin{eqnarray*}
x_2 &=& \Phi_{x_0}(u_2) = \Phi_{x_0}(f_B) = \Phi_{x_0}(f_{B' \cup C}) = \Phi_{x_0}(f_{B'}\circ f_C)\\
&=& \Phi_{x_0}(f_{B'}) \oplus \Phi_{x_0}(f_C) = \Phi_{x_0}(f_A) \oplus \Phi_{x_0}(f_C)
= \Phi_{x_0}(u_1) \oplus x = x_1 \oplus x\:.
\end{eqnarray*}
On the other hand, if there exists an injective mapping $q : B \to A$ then the same argument shows there exists 
$x \in X$ with $x_1 = x_2 \oplus x$.

Now to (a0) and (a1), and we have seen above that (a0) holds. Let $x,\,x' \in X$ and let $u,\,u' \in M_f$ with 
$x = \Phi_{x_0}(u)$ and $x' = \Phi_{x_0}(u')$. Then by $\mathrm{(\sharp)}$
\begin{eqnarray*}
x \oplus f(x') &=& \Phi_{x_0}(u) \oplus f(\Phi_{x_0}(u')) = \Phi_{x_0}(u) \oplus \Phi_{x_0}(f \circ u'))\\
&=& \Phi_{x_0}(u \circ f \circ u') = \Phi_{x_0}(f \circ u \circ u') = f(\Phi_{x_0}(u \circ u'))\\
&=& f(\Phi_{x_0}(u) \oplus \Phi_{x_0}(u')) = f(x \oplus x')
\end{eqnarray*}
and so (a1) holds. If $\oplus'$ is another binary operation on $X$ satisfying (a0) and (a1) then it is easy to 
see that $X_0 = \{ x' \in X : x \oplus' x' = x \oplus x'\ \mbox{for all $x \in X$} \}$ is a $f$-invariant subclass 
of $X$ containing $x_0$. Hence $X_0 = X$, since $\mathbf{I}$  is minimal, which implies that ${\oplus'} = {\oplus}$.

Finally, if $A$ and $B$ are disjoint finite sets then
\[ 
\val(A) \oplus \val(B) = \Phi_{x_0}(f_A) \oplus \Phi_{x_0}(f_B) 
= \Phi_{x_0}(f_A \circ f_B) = \Phi_{x_0}(f_{A\cup B}) = \val(A \cup B)\;.
\] 
Moreover, $\oplus$ is uniquely determined by this requirement: Consider any binary operation $\oplus'$ on $X$ 
for which $\val(A) \oplus' \val(B) = \val(A \cup B)$ whenever $A$ and $B$ are disjoint finite sets. If $C$ and $D$ are 
any finite sets then there exists a disjoint pair $(A,B)$ with $(A,B) \approx (C,D)$ and hence by 
Theorem~\ref{theorem_iterators_1} 
\[ 
\val(C) \oplus' \val(D) = \val(A) \oplus' \val(B) = \val(A \cup B) = \val(A) \oplus \val(B) = \val(C) \oplus \val(D) 
\] 
and so by Lemma~\ref{lemma_iterators_4} ${\oplus'} = {\oplus}$. This completes the proof of 
Theorem~\ref{theorem_am_1}. 
\eop

Let $\oplus$ be the operation given in Theorem~\ref{theorem_am_1}. The theorem shows in particular that  
$(X,\oplus,x_0)$ is a commutative monoid. The next result generalises Propositions \ref{prop_am_1} and 
\ref{prop_am_2} and shows how properties of the mapping $f$ correspond to properties of the monoid $(X,\oplus,x_0)$.

\begin{proposition}\label{prop_mam_3}
(1)\enskip
$(X,\oplus,x_0)$ is a group if and only if $f$ is a bijection.

(2)\enskip
The cancellation law holds in $(X,\oplus,x_0)$ if and only if $f$ is injective.
\end{proposition}

\proof
We have the commutative monoid $(X,\oplus,x_0)$, and also the commutative monoid $(M_f,\circ,\id_X)$. Now the 
operation $\oplus$ was defined so that
\[ 
\Phi_{x_0}(u_1) \oplus \Phi_{x_0}(u_2) = \Phi_{x_0}(u_1 \circ u_2)
\]
for all $u_1,\,u_2 \in M_f$ and, since  $\Phi_{x_0}(\id_X) = x_0$ and $\Phi_{x_0}$ is a bijection, this means 
$\oplus$ was defined to make $\Phi_{x_0} : (M_f,\circ,\id_X) \to (X,\oplus,x_0)$ a monoid isomorphism. It follows 
that the cancellation law holds in $(X,\oplus,x_0)$ if and only if it holds in $(M_f,\circ,\id_X)$ and that 
$(X,\oplus,x_0)$ will be a group if an only if $(M_f,\circ,\id_X)$ is. It is thus enough to prove the statements 
in the proposition with $(X,\oplus,x_0)$ replaced by $(M_f,\circ,x_0)$.

(1)\enskip
We first show that $u^{-1} \in M_f$ whenever $u \in M_f$ is a bijection. This follows from the fact that 
$u^{-1}(x_0) \in X$ and $\Phi_{x_0}$ is surjective and so there exists $v \in M_f$ with 
$\Phi_{x_0}(v) = u^{-1}(x_0)$; thus by $\mathrm{(\sharp)}$
\[  
\Phi_{x_0}(u \circ v) = u(\Phi_{x_0}(v)) = u(u^{-1}(x_0)) = x_0 = \Phi_{x_0}(\id_X) 
\]
and therefore $u \circ v = \id_X$, since $\Phi_{x_0}$ is injective. Hence $u^{-1} = v \in M_f$. Now clearly $M_f$ 
is a group if and only if  each mapping $u \in M_f$ is a bijection and $u^{-1} \in M_f$, and we have just seen that
$u^{-1} \in M_f$ holds automatically whenever $u \in M_f$ is a bijection. Moreover, by 
Lemma~\ref{lemma_mam_5}~(1) each element of $M_f$ is a bijection if and only if $f$ is a bijection.

(2)\enskip
Suppose  the cancellation law holds in $(M_f,\circ,\id_X)$, and let $x_1,\,x_2 \in X$ with $f(x_1) = f(x_2)$. Then 
there exist $u_1,\,u_2 \in M_f$ with $\Phi_{x_0}(u_1) = x_1$ and $\Phi_{x_0}(u_2) = x_2$ (since $\Phi_{x_0}$ is 
surjective), and hence by $\mathrm{(\sharp)}$
\[   
\Phi_{x_0}(f\circ u_1) = f(\Phi_{x_0}(u_1)) = f(x_1) = f(x_2) = f(\Phi_{x_0}(u_2)) = \Phi_{x_0}(f\circ u_2)\;.
\] 
It follows that $f \circ u_1 = f \circ u_2$ (since $\Phi_{x_0}$ is injective) and so $u_1 = u_2$. In particular 
$x_1 = x_2$, which shows that $f$ is injective. The converse is immediate, since if $f$ is injective then by 
Lemma~\ref{lemma_mam_5}~(2) so is each $u \in M_f$ and hence $u_1 = u_2$ whenever $u \circ u_1 = u \circ u_2$.
\eop

We now begin the preparations for the proof of Theorem~\ref{theorem_am_2}.

For each $u \in M_f$ consider the iterator $\,\mathbf{I}^u_* =(\Self{X},u_*,\id_X)$, where 
$u_* : \Self{X} \to \Self{X}$ is defined by $u_*(v) = u \circ v$ for all $v \in \Self{X}$ and let
$\val^u_*$ be the assignment of finite sets in $\,\mathbf{I}^u_*$.
Thus  $\val^u_* : \fin \to \Self{X}$ is the unique mapping with $\val^u_*(\varnothing) = \id_X$ such that 
$\val^u_*(A \cup \{a\}) = u_*(\val^u_*(A)) = u \circ \val^u_*(A)$ for each finite set $A$ and each $a \notin A$. 

Now it is more convenient to write $u_A$ instead of $\val^u_*(A)$ (this being consistent with the previous notation for the special case with $u = f$. Thus $A \mapsto u_A$ is the unique mapping with $ u_\varnothing = \id_{X}$ such that
$u_{(A \cup \{a\})} = u \circ u_A$ for each finite set $A$ and each $a \notin A$.

\begin{lemma}\label{lemma_mam_8}
(1)\enskip $(f_B)_A = f_{B\times A}$ for all finite sets $A$ and $B$.

(2)\enskip $(f_B)_A = (f_A)_B$ for all finite sets $A$ and $B$.
\end{lemma}

\proof
(1)\enskip
Let $A$ and $B$ be finite sets and put $\mathcal{S} = \{ C \in \mathcal{P}(A) : (f_B)_C = f_{B\times C} \}$.
Then $\varnothing \in \mathcal{S}$ since $(f_B)_{\varnothing} = \id_X = f_\varnothing = f_{B \times \varnothing}$.
Let $C \in \proper{\mathcal{S}}$ (and so $(f_B)_C = f_{B\times C}$) and let $a \in A \setminus C$. 
Then $B \times (C \cup \{a\})$ is the disjoint union of the sets $B \times \{a\}$ and $B \times C$ and 
$B \times \{a\} \approx B$; thus by Proposition~\ref{prop_mam_2} and Lemma~\ref{lemma_mam_4}
\[  
(f_B)_{C \cup \{a\}} = f_B \circ (f_B)_C = f_B \circ f_{B\times C} = f_{B\times \{a\}} \circ f_{B\times C} 
= f_{(B\times \{a\}) \cup (B\times C)} = f_{B \times (C\cup \{a\})}
\]
and so $C \cup \{a\} \in \mathcal{S}$. Hence $\mathcal{S}$ is an inductive $A$-system and so $A \in \mathcal{S}$. 
Hence $(f_B)_A = f_{B\times A}$.

(2)\enskip
By Proposition~\ref{prop_mam_2} $f_{B\times A} = f_{A\times B}$, since clearly $B \times A \approx A \times B$, 
and therefore by (1) $(f_B)_A = f_{B\times A} = f_{A \times B} = (f_A)_B$. 
\eop

The next result is not needed in what follows, but it shows that ($\spadesuit$) in Section~\ref{am}
holds, and could thus be used instead of Lemma~\ref{lemma_am_3} in the previous proof of Theorem~\ref{theorem_am_2}.

\begin{lemma}\label{lemma_mam_9}
Let $(A,B)$ and $(A',B')$ be pairs of finite sets.

(1)\enskip
If $f_A = f_{A'}$ and $f_B = f_{B'}$ then $f_{A \times B} = f_{A' \times B'}$.

(2)\enskip
If $\val(A,B) = \val(A',B')$ then $\val(A \times B) = \val(A' \times B')$.
\end{lemma}

\proof 
(1)\enskip
By several applications of Lemma~\ref{lemma_mam_8} (1) and (2) we have
\[
f_{A \times B} = (f_B)_A = (f_{B'})_A = (f_A)_{B'} = (f_{A'})_{B'} = (f_{B'})_{A'} = f_{A' \times B'}\;.
\]

(2)\enskip
This follows immediately from (1) and Proposition~\ref{prop_mam_1}. \eop

\begin{lemma}\label{lemma_mam_10}
Let $v \in M_f$. Then:

(1)\enskip 
$v_A \in M_f$ for every finite set $A$.

(2)\enskip
If $A$ and $B$ are finite sets with $f_A = f_B$ then $v_A = v_B$.
\end{lemma}

\proof
(1)\enskip
Let $A$ be a finite set and put $\mathcal{S} = \{ B \in \mathcal{P}(A) : v_B \in M_f \}$.
Then $\varnothing \in \mathcal{S}$, since $v_{\varnothing} = \id_X \in M_f$.
Consider $B \in \proper{\mathcal{S}}$ (and so $v_B \in M_f$) and let $a \in A \setminus B$. 
Then $v_{B \cup \{a\}} = v \circ v_B \in M_f$, since $M_f$ is a submonoid of $(\Self{X},\circ,\id_X)$, and so 
$B \cup \{a\} \in \mathcal{S}$. Hence $\mathcal{S}$ is an inductive $A$-system and so $A \in \mathcal{S}$.
Thus $v_A \in M_f$.

(2)\enskip
There exists a finite set $C$ with $v = f_C$ and thus by Lemma~\ref{lemma_mam_8}~(2)
\[ 
v_A = (f_C)_A = (f_A)_C = (f_B)_C = (f_C)_B = v_B\;.\ \eop
\]

Let $v \in M_f$; then by Lemma~\ref{lemma_mam_10} there exists a unique mapping $\psi_v : M_f \to M_f$ such 
that $\psi_v(f_A) = v_A$ for each finite set $A$, and in particular $\psi_v(f) = v$ (since if $a$ is any element 
then $f = f_{\{a\}}$ and $v = v_{\{a\}}$). Moreover, if $v = f_B$ then by Lemma~\ref{lemma_mam_8}~(1) 
$\psi_v(f_A) = f_{A \times B}$. 

A mapping $\psi : M_f \to M_f$ is an \definition{endomorphism} (of the monoid $(M_f,\circ,\id_X)$) if 
$\psi(\id_X) = \id_X$ and $\psi(u_1 \circ u_2) = \psi(u_1) \circ \psi(u_2)$ for all $u_1,\,u_2 \in M_f$.

\begin{lemma}\label{lemma_mam_11}
(1)\enskip
$\psi_v$ is an endomorphism for each $v \in M_f$. 

(2)\enskip
$\psi_v(u) = \psi_u(v)$ for all $u,\,v \in M_f$.
\end{lemma}

\proof
(1)\enskip
If $u_1,\,u_2 \in M_f$ then there exist disjoint finite sets $A$ and $B$ with $u_1 = f_A$ and $u_2 = f_B$ and 
hence by Lemma~\ref{lemma_mam_4}
\[ 
\psi_v(u_1 \circ u_2) = \psi_v(f_A \circ f_B) = \psi_v(f_{A \cup B}) = v_{A \cup B} = 
  v_A \circ v_B = \psi_v(f_A) \circ \psi_v(f_B) 
\] 
(noting that proof of Lemma~\ref{lemma_mam_4} also shows that $v_{A \cup B} = v_A \circ v_B$.  Moreover, we have
$\psi_v(\id_X) = \psi_v(f_\varnothing) = v_\varnothing = \id_X$ and hence $\psi_v$ is an endomorphism.

(2)\enskip
Let $A,\,B$ be finite sets with $u = f_A$ and $v = f_B$. Then by Lemma~\ref{lemma_mam_8}~(2)
\[ 
\psi_v(u) = \psi_v(f_A) = v_A = (f_B)_A = (f_A)_B = u_B = \psi_u(f_B) = \psi_u(v)\;.\ \eop 
\]

\textit{Proof of Theorem~\ref{theorem_am_2}:}\enskip
Define a binary operation $\diamond : M_f \times M_f \to M_f$ by letting
\[ 
u \diamond v = \psi_u(v) 
\] 
for all $u,\,v \in M_f$. In particular, if $A$ and $B$ are finite sets then by Lemma~\ref{lemma_mam_8}~(1) 
$f_A \diamond f_B = \psi_{f_A}(f_B) = (f_A)_B = f_{A \times B}$ and therefore
\[ 
f_A \diamond f_B = f_{A \times B} 
\]
for all finite sets $A$ and $B$. Let $u,\,v,\,w \in M_f$ and $A$, $B$ and $C$ be finite sets with $u = f_A$, 
$v = f_B$ and $w = f_C$. Then clearly $A \times (B \times C) \approx (A \times B) \times C$, so by 
Proposition~\ref{prop_mam_2} $f_{A \times (B \times C)} = f_{(A \times B) \times C}$ and thus
\begin{eqnarray*}
u \diamond (v \diamond w) &=& f_A \diamond (f_B \diamond f_C) = f_A \diamond f_{B\times C}\\ 
&=& f_{A\times (B\times C)} = f_{(A\times B)\times C} = f_{A\times B} \diamond f_C 
= (f_A \diamond f_B) \diamond f_C = (u \diamond v) \diamond w\;.
\end{eqnarray*}
Hence $\diamond$ is associative. Moreover, Lemma~\ref{lemma_mam_11}~(2) shows that $\diamond$ is commutative, 
since $u \diamond v = \psi_u(v) = \psi_v(u) = v \diamond u$ for all $u,\,v \in M_f$. Also (with $a$ any element) 
$u \diamond f = u \diamond f_{\{a\}} = u_{\{a\}} = u$, i.e., $u \diamond f = u$ for all $u \in M_f$, and by 
Lemma~\ref{lemma_mam_11}~(1) $u \diamond \id_X = \id_X$ and 
$u \diamond (v_1 \circ v_2) = (u \diamond v_1) \circ (u \diamond v_2)$ for all $u,\,v_1,\,v_2 \in M_f$.

Now since $\Phi_{x_0} : M_f \to X$ is a bijection there clearly exists a unique binary relation $\otimes$ on 
$X$ such that
\[ 
\Phi_{x_0}(u_1) \otimes \Phi_{x_0}(u_2) = \Phi_{x_0}(u_1 \diamond u_2)
\]
for all $u_1,\,u_2 \in M_f$, and exactly as in the proof of Theorem~\ref{theorem_am_1} the operation $\otimes$ 
is associative and commutative since $\diamond$ has these properties. The same holds true of the distributive 
law: Let $x,\,x_1,\,x_2 \in X$, and $u,\,v_1,\,v_2 \in M_f$ be such that $x = \Phi_{x_0}(u)$, 
$x_1 = \Phi_{x_0}(v_1)$ and $x_2 = \Phi_{x_0}(v_2)$. Then
\begin{eqnarray*}
x \otimes (x_1 \oplus x_2) &=& \Phi_{x_0}(u) \otimes (\Phi_{x_0}(v_1) \oplus \Phi_{x_0}(v_2)) \\
&=& \Phi_{x_0}(u) \otimes \Phi_{x_0}(v_1 \circ v_2) = \Phi_{x_0}(u \diamond (v_1 \circ v_2)) \\
&=& \Phi_{x_0}((u \diamond v_1) \circ (u \diamond v_2))
= \Phi_{x_0}(u \diamond v_1) \oplus \Phi_{x_0}(u \diamond v_2)\\ 
&=& (\Phi_{x_0}(u) \otimes \Phi_{x_0}(v_1)) \oplus (\Phi_{x_0}(u) \otimes \Phi_{x_0}(v_2)) 
=  (x \otimes x_1) \oplus (x \otimes x_2)
\end{eqnarray*}

Next, if $x \in X$ and $u \in M_f$ is such that $x = \Phi_{x_0}(u)$ then
\begin{eqnarray*} 
&&x \otimes x_0 = \Phi_{x_0}(u) \otimes \Phi_{x_0}(\id_X) = \Phi_{x_0}(u \diamond \id_X) = \Phi_{x_0}(\id_X) 
=  x_0\,,\\
&&x \otimes f(x_0) = \Phi_{x_0}(u) \otimes \Phi_{x_0}(f) = \Phi_{x_0}(u \diamond f) = \Phi_{x_0}(u) =  x
\end{eqnarray*}
and so $x \otimes x_0 = x_0$ and $x \otimes f(x_0) = x$ for all $x \in X$.

We have already seen  (m0) holds and, since $f(x_0)$ is an identity element, (m1) is a special case of the distributive law: 
Let $x,\,x' \in X$; then by (a0) and (a1) and since $\oplus$ is commutative it follows that
$f(x') = f(x' \oplus x_0) = x' \oplus f(x_0) = f(x_0) \oplus x'$, and hence
$x \otimes f(x') = x \otimes (f(x_0) \oplus x')
= (x \otimes f(x_0)) \oplus (x \otimes x') = x \oplus (x \otimes x')$, which is (m1). Finally, if $\otimes'$ is 
another binary operation satisfying (m0) and (m1) then it is easy to see that 
$X_0 = \{ x' \in X : x \otimes' x' = x \otimes x'\ \mbox{for all $x \in X$} \}$
is a $f$-invariant subset of $X$ containing $x_0$. Hence $X_0 = X$, since $\,\mathbf{I}$ is minimal, which implies 
that ${\otimes'} = {\otimes}$.
\eop


\startsection{The generalised associative law}

\label{assoc}
Let $\bullet$ be a binary operation on a set $X$, written using infix notation, so $x_1 \bullet x_2$ denotes the product of $x_1$ and $x_2$. The large majority of 
such operations occurring in mathematics are \definition{associative}, meaning that $(x_1 \bullet x_2) \bullet x_3 = x_1 \bullet ( x_2 \bullet x_3)$ 
for all $x_1,\,x_2,\,x_3 \in X$. If $\bullet$ is associative and $x_1,\,x_2,\,\ldots,\,x_n \in X$ then the product $x_1 \bullet x_2 \bullet \cdots \bullet x_n$ 
is well-defined, meaning its value does not depend on the order in which the operations are carried out.

This result will be  established in the present section. We first define a particular order of carrying out the operations. This is the order in which at each
stage the product of the current first and second components are taken. For example, the product of the 6 components $x_1,\,x_2,\,x_3,\, x_4,\,x_5, \, x_6$
 evaluated using this order results in the value $\,\bullet(x_1,\ldots,x_6) = (((((x_1 \bullet x_2) \bullet x_3) \bullet x_4)  \bullet x_5)\bullet x_6)$.
In general, the corresponding product of $n$ terms will be denoted by $\bullet(x_1,\ldots,x_n)$.

Theorem~\ref{theorem_assoc_1} states that if $\bullet$ is associative then
\[\bullet(x_1,\ldots,x_m,x_{m+1},\ldots,x_n) = \alpha \bullet \beta\;,\]
where $\alpha = \bullet(x_1,\ldots,x_m)$ and $ \beta = \bullet(x_{m+1},\ldots,x_n)$. This is a weak form of the generalised associative law, although it is one 
which is often all that is needed.

Theorem~\ref{theorem_assoc_2} gives the general form of the generalised associative law and states that if $\bullet$ is associative then
$\bullet(x_1,\ldots,x_n) = \bullet_{\mathbf{R}}(x_1,\ldots,x_n)$ for each $\,\mathbf{R}$ from the set of prescriptions describing how the operations are carried out. 
The main task is to give a rigorous definition of this set. We do this using partitions of intervals of the form $\{ k \in \Int : m \le  k \le n \}$ in which each 
element in the partition is also an interval of this form.

Recall that by a \definition{partition} of a set $S$ we mean  a subset $\mathcal{Q}$ of $\mathcal{P}_0(S)$ such that for each $s \in S$ there exists a unique 
$Q \in \mathcal{Q}$ such that $s \in Q$. Thus, different elements in a partition of $S$ are disjoint and their union is $S$. 

Consider the product $((x_1 \bullet x_2) \bullet ((x_3 \bullet x_4) \bullet x_5))$. The order of operations involved here can described with the help of the 
following sequence of partitions of the set $\{1,2,3,4,5\}$:

$\{\{1\},\{2\},\{3\},\{4\},\{5\}\}$\\
$\{\{1\},\{2\},\{3,4\},\{5\}\}$\\
$\{\{1,2\},\{3,4\},\{5\}\}$\\
$\{\{1,2\},\{3,4,5\}\}$\\
$\{\{1,2,3,4,5\}\}$

For each of these partitions (except the last one) the next partition  is obtained by amalgamating two adjacent partitions. Corresponding to these partitions 
there  is a sequence of partial evaluations:

$\{\{x_1\},\{x_2\},\{x_3\},\{x_4\},\{x_5\}\}$\\
$\{\{x_1\},\{x_2\},\{(x_3\bullet x_4)\},\{x_5\}\}$\\
$\{\{(x_1\bullet x_2)\},\{(x_3\bullet x_4)\},\{x_5\}\}$\\
$\{\{(x_1 \bullet x_2)\},\{((x_3 \bullet x_4) \bullet x_5)\}\}$\\
$\{\{((x_1 \bullet x_2 )\bullet ((x_3 \bullet x_4) \bullet x_5))\}\}$

and the final expression is essentially the product we started with.

\medskip

Instead of using intervals of the form $\{ k \in \Int : m \le  k \le n \}$ we prefer to use a formulation in terms of a finite totally ordered set.

Thus let $(E,\le_E)$ be a non-empty finite  total ordered set with maximum element $z$ and define a mapping $f : E \to E$ by letting $f(e)$ be the least element in 
$\{ e' \in E: e' >_E e \}$ if $e \ne z$ and putting $f(z) = z$. Then $\,\mathbf{I} = (E,f,e_0)$, with $e_0$ the minimum element in $E$, is a finite iterator with 
fixed point $z$ and by Proposition~\ref{prop_finiteiter_1} $\,\mathbf{I}$ is minimal. The totally ordered set $(E,\le)$ is considered fixed in what follows, but 
recall that if $A$ is any finite non-empty set then by Lemma~\ref{lemma_fs_7} there exists a totally ordered set $(E,\le)$ with $E \approx A$.
 
For all $r,\,s \in E$ with $r \le_E s$ put $[r,s] = \{ e  \in E : r \le_E e \le_E s \}$. Sets of this form will be called \definition{intervals}. 
Also put $[r,s) = \{ e \in E : r \le_E e <_E s \}$.

\begin{lemma}\label{lemma_assoc_1} Let $[p,q]$ be an interval and let $B \subset [p,q]$. Suppose $p \in B$ and $f(s) \in B$ for all $s \in [p,q)$. Then $B = [p,q]$.
\end{lemma}

\proof Consider the restriction $\le_{[p,q]}$ of the total order $\le_E$ to the interval $[p,q]$ and define $f_q :[p,q] \to [p,q]$ by letting $f_q(r) = f(r)$ if 
$r \in [p,q)$ and putting $f_q(q) = q$. Then $\,\mathbf{I}_{[p,q]} = ([p,q],f_q,p)$ is the iterator given by Proposition~\ref{prop_finiteiter_1} for the totally 
ordered set $([p,q],\le_{[p,q]})$ and by Proposition~\ref{prop_finiteiter_1} $\,\mathbf{I}_{[p,q]}$ is minimal. But  this is just the statement in 
Lemma~\ref{lemma_assoc_1}. \eop

In what follows let $X$ be an arbitrary non-empty  set and let $\bullet$ be a binary operation on $X$.
If $I$ is an interval then an \definition{$I$-tuple} with values in $X$ is just a mapping $s \mapsto x_s$ from $I$ to $X$. If the interval is given explicitly 
as $[p,q]$ then a $[p,q]$-tuple will usually be written as $(x_p,\dots,x_q)$.

\begin{proposition}\label{prop_assoc_2} 
Let $[p,q]$ be an interval and let $(x_p,\ldots,x_q)$ be a $[p,q]$-tuple. Then there exists a unique $[p,q]$-tuple $(\pi_p,\ldots, \pi_q)$ with 
$\pi_p = x_p$ and such that $\pi_{f(s)} = \pi_s \bullet x_{f(s)}$ for all $s \in [p,q)$.
\end{proposition}

\proof Let $B$ be the subset of $[p,q]$ consisting of those $s$  for which there exists a unique $[p,s]$-tuple $(\pi_p^s,\ldots,\pi_s^s)$ 
with $\pi_p^s = x_p$ and such that $\pi_{f(y)}^s = \pi_y^s \bullet x_{f(y)}$ for all $y \in [p,x)$. Then $p \in B$ with $(x_p)$ the unique 
$[p,p]$-tuple.  Thus let $s \in B\setminus \{q\}$ and let $(\pi_p^s,\ldots,\pi_s^s)$ be the corresponding unique $[p,s]$-tuple. Then we can extend this 
$[p,s]$-tuple to a $[p,f(s)]$-tuple by putting $\pi_y^{f(s)} = \pi_y^s$ for each $y \in [p,s]$ and  letting $\pi_{f(s)}^{f(s)} = \pi_(s^s \bullet x_{f(s)}$. 
Then $(\pi_p^{f(s)},\ldots,\pi_{f(s)}^{f(s)})$ is the unique $[p,f(s)]$-tuple having the required properties and hence $f(s) \in B$. Therefore by 
Lemma~\ref{lemma_assoc_1} $B = [p,q]$. \eop

We denote the element $\pi_q$ of $X$ defined above by $\bullet(x_p,\ldots,x_q)$.

Let $(x_1,x_2,x_3,x_4)$ be a $[1,4]$-tuple (with $[1,4] \subset \Nat$ and let $(\pi_1,\pi_2,\pi_3,\pi_4)$ be the $[1,4]$-tuple given  by the above proposition. Then 

$\pi_1 = x_1$, $\pi_2 = \pi_1 \bullet x_2 = x_1 \bullet x_2$,
 $\pi_3 = \pi_2 \bullet x_3 = (x_1 \bullet x_2) \bullet x_3$,\\
$\pi_4 = \pi_3 \bullet x_4 = ((x_1 \bullet x_2 )\bullet x_3 )\bullet x_4$.
Hence $\bullet(x_1,x_2,x_3,x_4) = ((x_1 \bullet x_2 )\bullet x_3 )\bullet x_4$.

\begin{theorem}\label{theorem_assoc_1}
Let $p,\,q,\,r \in E$ with $p \le q$ and $f(q) \le r$ and note that $[p,r]$ is the disjoint union of $[p,q]$ and $[f(q),r]$. Let $(x_1,\ldots,x_r)$ be a 
$[p,r]$-tuple and so we have a $[p,q]$-tuple $(x_p,\ldots,x_q)$ and a $[f(q),r]$-tuple $(x_{f(q)}\ldots,x_r)$. Suppose that the operation $\bullet$ is associative. Then
 $\bullet(x_p,,\ldots,x_r) = \alpha \bullet \beta$, where $\alpha = \bullet(x_p,\ldots,x_q)$ and $\beta = \bullet(x_{f(q)},\ldots,x_r)$.
\end{theorem}

\proof For each $s \in [f(q),r]$ put $\beta_s = \bullet (x_{f(q)} \ldots, ,x_s)$ and let
 \[B = \{ s \in [f(q),r] : \bullet(x_p,\ldots,x_s) = \alpha \bullet \beta_s\}.\;\]
Then $f(q) \in B$ since $\bullet(x_p,\ldots,x_{f(q)}) = \bullet(x_p,\ldots,x_q) \bullet x_{f(q)}= \alpha \bullet \beta_{f(q)}$. Thus let $s \in B \cap [f(q),r)$. 
Then $\bullet(x_p,\ldots,s) = \alpha \bullet \beta_s$ and
\[\bullet(x_p,\ldots,x_{f(s)}) =\bullet (x_p,\ldots,x_s) \bullet x_{f(s)} = (\alpha \bullet \beta_s) \bullet x_{f(s)}
= \alpha \bullet(\beta_s \bullet x_{f(s)}) = \alpha \bullet \beta_{f(s)}\]
and hence $f(s) \in B$. Thus by Lemma~\ref{lemma_assoc_1} $B = [f(q),r]$, i.e., $\bullet(x_p,\ldots,x_r) = \alpha \bullet \beta$. \eop

If $I$ is an interval then an \definition{interval partition} of $I$ is a partition $\mathcal{Q}$ such that each element of $\mathcal{Q}$ is also an interval. 
From now on partition always means interval partition. The partition consisting of the singleton sets $[r,r]$, $r \in I$, will be denoted by $\mathcal{Q}_0$ and 
the trivial partition $\{I\}$ by $\mathcal{Q}_T$. 

Let $\mathcal{Q}$ be a partition of $I$ and let $J = [r,s]$ and $K = [u,v]$ be elements of $\mathcal{Q}$. Then $(J,K)$ will be called an \definition{adjacent pair} 
if $f(s) = u$. In this case $J \cup K =  [r,v]$ is the disjoint union of $J$ and $K$. A partition $\mathcal{Q}'$ is said to be an \definition{$I$-reduction} of 
$\mathcal{Q}$ if there exists an adjacent pair $(J,K)$ such that $\mathcal{Q}' = \mathcal{Q} \setminus \{J,K\} \cup \{J \cup K \}$. The partition $\mathcal{Q}'$ 
will be denoted by $\mathcal{Q}(J,K)$. Consider the interval $I = [p,q]$ to be fixed and let $\,\mathbf{Q}$ be a subset of  the set of all partitions of $I$ 
containing  $\mathcal{Q}_0$ and $\mathcal{Q}_T$. Put $\,\mathbf{Q}_S = \mathbf{Q} \setminus \{\mathcal{Q}_T\}$ and let $\Phi : \mathbf{Q} \to \mathbf{Q}$ be a 
mapping with $\Phi(\mathcal{Q}_T) = \mathcal{Q}_T$. Thus there is the iterator $\,\mathbf{R} = (\mathbf{Q},\Phi,\mathcal{Q}_0)$ and $\,\mathbf{R}$ will be  called 
an $I$-\definition{reduction} if it is minimal and $\Phi(\mathcal{Q})$ is a reduction of $\mathcal{Q}$ for each $\mathcal{Q} \in \mathbf{Q}_S$. Note that if 
$p = q$ then there is a unique $I$-reduction with $\,\mathbf{Q} = \mathcal{Q}_T\{I\}$ and with $\Phi$ the identity mapping. 
If $q = f(p)$ then there is also a unique $I$-reduction with $\,\mathbf{Q} = \{\mathcal{Q}_0,\mathcal{Q}_T\}$ and with $\Phi(\mathcal{Q}_0) = \mathcal{Q}_T$ 
and $\Phi(\mathcal{Q}_T) = \mathcal{Q}_T$. In what follows let $\,\mathbf{R} =(\mathbf{Q},\Phi,\mathcal{Q}_0)$ be an $I$-reduction. Thus $\,\mathbf{R}$ is a 
finite minimal iterator with fixed-point $\mathcal{Q}_T$ and so let $\le$ be the total order on $\,\mathbf{Q}$ given in Proposition~\ref{prop_finiteiter_11}. 
Let $<$, $>$ and $\ge$ have their usual meanings. If $A$ and $B$ are finite sets then as before $A \preceq B$ means there is an injective mapping $h : A \to B$ 
and we write $A \prec B$ if both $A \preceq B$ and $A \not\approx B$ hold. Moreover if $B \ne \varnothing$ then we write $A \prec_0 B$ if 
$A \approx B \setminus \{b\}$, where $b$ is any element in $B$. If $\mathcal{Q} \in \mathbf{Q}_S$ then $\Phi(\mathcal{Q}) \prec_0 \mathcal{Q}$.
 
\begin{lemma}\label{lemma_assoc_13}
If $\mathcal{Q},\,\mathcal{Q}' \in \mathbf{Q}$ with $\mathcal{Q} \approx \mathcal{Q}'$ then $\mathcal{Q} = \mathcal{Q}'$. 
\end{lemma}

\proof 
Put $\mathcal{S} = \{\mathcal{Q}_0 \} \cup \{\mathcal{Q} \in \mathbf{Q} \setminus \{\mathcal{Q}_0\} : \mathcal{Q} \prec \mathcal{Q}_0 \}$ and let
$\mathcal{Q} \in \mathcal{S} \setminus \{\mathcal{Q}_0\}$. Then $\mathcal{Q} \prec \mathcal{Q}_0$ and thus either $\mathcal{Q} \in \mathbf{Q}_S$, in
which case $\Phi(\mathcal{Q}) \prec \mathcal{Q}$ or $\mathcal{Q} = \mathcal{Q}_T$, in which case $\Phi(\mathcal{Q}) = \mathcal{Q}$. In both cases
$\Phi(\mathcal{Q})\prec \mathcal{Q}_0$ and so $\Phi(\mathcal{Q}) \in \mathcal{S}$.  Thus $\mathcal{S}$ is $\Phi$-invariant and contains $\mathcal{Q}_0$.
Hence $\mathcal{S} = \mathbf{Q}$, which shows that $\mathcal{Q} \prec \mathcal{Q}_0$ and in particular $\mathcal{Q}\not\approx \mathcal{Q}_0$ for all 
$\mathcal{Q} \in \mathbf{Q} \setminus \{\mathcal{Q}_0\}$.
Now let $\mathcal{Q},\,\mathcal{Q}' \in \mathbf{Q}$ with $\mathcal{Q} \ne \mathcal{Q}'$ and suppose that $\mathcal{Q} \approx \mathcal{Q}'$.
We can assume that $\mathcal{Q} < \mathcal{Q}'$ and that $\mathcal{Q}$ is the least element of $\mathbf{Q}$ such that $\mathcal{Q} \approx \mathcal{Q}'$
for some $\mathcal{Q}'> \mathcal{Q}$. Then $\mathcal{Q} \ne \mathcal{Q}_0$ since $\mathcal{Q}\not\approx \mathcal{Q}_0$ for all 
$\mathcal{Q} \in \mathbf{Q} \setminus \{\mathcal{Q}_0\}$ and so also $\mathcal{Q}' \ne \mathcal{Q}_0$.
Hence by Proposition~\ref{prop_iterators_1} $\mathcal{Q} = \Phi(\mathcal{U})$ and
$\mathcal{Q}' = \Phi(\mathcal{U}'$ for some $\mathcal{U},\, \mathcal{U}' \in \mathbf{Q}$. Now $\mathcal{U} \prec_0 \mathcal{Q}$,
$\mathcal{U}' \prec_0 \mathcal{Q}'$ and $\mathcal{Q} \approx\mathcal{Q}'$ and it follows from Proposition~\ref{prop_finiteiter_11} (4) that 
$\mathcal{U}' > \mathcal{U}$ and $\mathcal{U} \approx \mathcal{U}'$. But $\mathcal{U} < \mathcal{Q}$, which contradicts the minimality of $\mathcal{Q}$.
Therefore if $\mathcal{Q},\,\mathcal{Q}' \in \mathbf{Q}$ with $\mathcal{Q} \ne \mathcal{Q}'$ Then $\mathcal{Q} \not\approx \mathcal{Q}'$ 
and hence if $\mathcal{Q},\,\mathcal{Q}' \in \mathbf{Q}$ with $\mathcal{Q} \approx \mathcal{Q}'$ then $\mathcal{Q} = \mathcal{Q}'$. \eop

\begin{proposition}\label{prop_assoc_14}
For each $J \in \mathcal{P}_0(I)$ there exists a unique $\mathcal{Q} \in \mathbf{Q}$ with $\mathcal{Q} \approx J$.
 \end{proposition}

\proof
Let $\mathcal{S} = \{\varnothing\} \cup \{ J \in \mathcal{P}_0(I) :\mbox{ there exists $\mathcal{Q} \in \mathbf{Q}$ with $\mathcal{Q} \approx J$}\}$.
Let $J \in \proper{\mathcal{S}}$ and $j \in I \setminus J$; put $J' = J \cup \{j\}$. If $J = \varnothing$ and so $J' = \{j\}$ then
$\mathcal{Q}_T \approx J'$.  If $J \ne \varnothing$ then $\mathcal{Q} \approx J$ for some $\mathcal{Q} \in \mathbf{Q}$ and $\mathcal{Q}\ne \mathcal{Q}_0$,
since $\mathcal{Q}_0 \approx I$ and $J$ is a proper subset of $I$. Thus by Proposition~\ref{prop_iterators_1} $\mathcal{Q} = \Phi(\mathcal{Q}')$ for some
$\mathcal{Q}' \in \mathbf{Q}$ and then $\mathcal{Q}' \approx J'$, since $\mathcal{Q} \prec_0 \mathcal{Q}'$. Thus in both cases $J' \in \mathcal{S}$
and therefore $\mathcal{S} = \mathcal{P}(I)$. This shows that for each $J \in \mathcal{P}_0(I)$ there exists $\mathcal{Q} \in \mathbf{Q}$ with 
$\mathcal{Q} \approx J$. The uniqueness follows from Lemma~\ref{lemma_assoc_13}. \eop

\begin{lemma}\label{lemma_assoc_101}
Let $A$ and $B$ be finite sets, let $\mathcal{S} \subset \mathcal{P}(A)$ and suppose that
for each $J \in \mathcal{P}_0(B)$ there exists a unique $C \in \mathcal{S}$ with $C \approx J$.
Then $B \approx \mathcal{S}$. In particular, Proposition~\ref{prop_assoc_14} implies that
$I \approx \mathbf{Q}$.
\end{lemma}
\proof
There is a unique surjective mapping $\alpha : \mathcal{P}_0(B) \to \mathcal{S}$ with $\alpha(J) \approx J$ for each
$J\in \mathcal{P}_0(B)$. Let $\approx'$ be the restriction of the equivalence relation $\approx$ to $\mathcal{P}_0(B)$ and let
$\mathcal{E}(B)$ be the set of equivalence classes. Then the mapping $\alpha$ induces a bijective mapping $\alpha' :\mathcal{E}(B) \to \mathcal{S}$
and hence by Proposition~\ref{prop_fs_888} $B \approx \mathcal{S}$. \eop

\begin{proposition}\label{prop_assoc_101}
Suppose $\,\mathbf{R} = (\mathbf{Q},\Phi,\mathcal{Q}_0)$  satisfies the requirements for being an $I$-reduction except
for the assumption that it be minimal. Then $I \preceq \mathbf{Q}$ and $I \approx \mathbf{Q}$ if and only if $\,\mathbf{R}$ is minimal.
 \end{proposition}

\proof Let $\,\mathbf{Q}_0$ be the least $\Phi$-invariant subset of $\,\mathbf{Q}$ containing $\mathcal{Q}_0$ and let $\Phi_0$ be the restriction of
$\Phi$ to $\,\mathbf{Q}_0$. Then $\,\mathbf{R}_0 = (\mathbf{Q}_0,\Phi_0,\mathcal{Q}_0)$  is minimal and thus is an $I$-reduction and so by
Lemma~\ref{lemma_assoc_101} $I \approx \mathbf{Q}_0$. Hence $I \preceq \mathbf{Q}$ and $I \approx \mathbf{Q}$ if and only if $\,\mathbf{R}$ is minimal.
\eop

In what follows assume that $\,\mathbf{R}$ is minimal.

\begin{proposition}\label{prop_assoc_13}
Let $s \mapsto x_s$ be an $I$-tuple. Then for each $\mathcal{Q} \in \mathbf{Q}$ there is a unique mapping  $\omega_{\mathcal{Q}} : \mathcal{Q} \to X$
with $\omega_{\mathcal{Q}_0}([z,z]) = x_z$ for each $z \in I$ and such that if $\mathcal{Q} \in \mathbf{Q}_S$ with
$\Phi(\mathcal{Q}) = \mathcal{Q}(J,K)$ then $\omega_{\Phi(\mathcal{Q})}(M) = \omega_{\mathcal{Q}}(M)$ if $M \ne J\cup K$ 
and  $\omega_{\Phi(\mathcal{Q})}(J\cup K) = \omega_{\mathcal{Q}}(J)\bullet \omega_{\mathcal{Q}}(K)$.
\end{proposition}

\proof
For each $\mathcal{Q} \in\mathbf{Q}$ let $\prop(\mathcal{Q})$ be the statement that for each $\mathcal{U} \le \mathcal{Q}$ 
there exists a unique mapping $\omega^{\mathcal{Q}}_{\mathcal{U}} : \mathcal{U} \to X$ 
with $\omega^{\mathcal{Q}}_{\mathcal{Q}_0}([z,z]) = x_z$ for each $z \in I$ and such that  
if $\,\mathcal{U} \in \mathbf{Q}_S$ with $\Phi(\mathcal{U}) = \mathcal{U}(J,K)$ then 
$\omega^{\mathcal{Q}}_{\Phi(\mathcal{U})}(M) = \omega^{\mathcal{Q}}_{\mathcal{U}}(M)$ 
if $M \ne J\cup K$ and  $\omega^{\mathcal{Q}}_{\Phi(\mathcal{U})}(J\cup K) = \omega^{\mathcal{Q}}_{\mathcal{U}}(J)\bullet\omega^{\mathcal{Q}}_{\mathcal{U}}(K)$. 
Then clearly $\prop(\mathcal{Q}_0)$ holds and so let $\mathcal{Q} \in \mathbf{Q}$ be such that $\prop(\mathcal{Q})$ holds; put $\mathcal{V} = \Phi(\mathcal{Q})$. 
We must show that $\prop(\mathcal{V})$ holds and if $\mathcal{U} \le \mathcal{V}$ then either $\mathcal{U} \le \mathcal{Q}$ 
or $\mathcal{U} = \mathcal{V}$. If $\mathcal{U} \le \mathcal{Q}$ then put $\omega^{\mathcal{V}}_{\mathcal{U}} =\omega{^{\mathcal{Q}}_{\mathcal{U}}}$. Finally,
if $\mathcal{U} = \mathcal{V} \in \mathbf{Q}_S$ 
with $\Phi(\mathcal{V}) = \mathcal{V}(J,K)$ then put $\omega^{\mathcal{V}}_{\Phi(\mathcal{V})}(M)= \omega^{\mathcal{Q}}_{\mathcal{Q}}(M)$
if $M \ne J \cup K$ and  $\omega^{\mathcal{V}}_{\Phi(\mathcal{V})}(J\cup K) = \omega^{\mathcal{\mathcal{Q}}}_{\mathcal{Q}}(J)\bullet \omega^{\mathcal{Q}}_{\mathcal{Q}}(K)$. 
Then $\{\omega^{\mathcal{V}}_{\mathcal{U}} :\mathcal{U} \le \mathcal{V} \}$ are the unique mappings satisfying the requirements for
$\prop(\mathcal{V})$ and hence $\prop(\mathcal{V})$ holds. Therefore $\prop(\mathcal{Q})$ holds for all $\mathcal{Q} \in \mathbf{Q}$. 
We now define $\omega_{\mathcal{Q}} = \omega_{\mathcal{Q}}^{\mathcal{Q}_T}$ for each $\mathcal{Q} \in \mathbf{Q}$. \eop
 
Note that $\omega_{\mathcal{Q}_T}$ is a mapping from the singleton set $\{I\}$ to $X$ and this element of $X$ will be denoted by
$\bullet_{\mathbf{R}}(s \mapsto x_s)$. If $p=q$ then $\bullet_{\mathbf{R}}(x_p) = x_p$ and if $q = f(p)$ then $\bullet_{\mathbf{R}}(x_p,x_q) = x_p \bullet x_q$, where
in both cases $\,\mathbf{R}$ is the unique $I$-reduction.

The element $\bullet(s \mapsto x_s)$ of $X$ arising from Proposition~\ref{prop_assoc_2} is obtained using the $I$-reduction 
$\,\mathbf{R} = (\mathbf{Q},\Phi,\mathcal{Q}_0)$ 
where $\Phi(\mathcal{Q}) =\mathcal{Q}(L_{\mathcal{Q}}^1,L_{\mathcal{Q}}^2)$ with $(L_{\mathcal{Q}}^1,L_{\mathcal{Q}}^2)$ the adjacent pair consisting of the first and 
second elements in the partition $\mathcal{Q}$.

\begin{theorem}\label{theorem_assoc_2} Suppose that $\bullet$ is associative. Then
\[\bullet_{\mathbf{R}}(s \mapsto x_s) = \bullet(s \mapsto x_s) \]
for each $I$-tuple $(s \mapsto x_s)$ and each reduction $\,\mathbf{R}$. 
\end{theorem}

\proof In what follows we assume that $I$ contains at least two elements and so by Theorem~\ref{theorem_finiteiter_1}~(2) there exists a unique 
$\mathcal{Q}_2 \in \mathbf{Q}$ with $\Phi(\mathcal{Q}_2)  = \mathcal{Q}_T$.
Clearly the partition $\mathcal{Q}_2$ consists of exactly two components. Thus if $I = [p,q]$ then there exists $t \in [p,q]$ with $p \le_E t <_E q$ such that 
$\mathcal{Q}_2 = \{L,R\}$, where $L = [p,t]$ and $R = [f(t),q]$. 

\begin{lemma}\label{lemma_assoc_33} For each $\mathcal{Q} \in \mathbf{Q}_S$ and each $Q \in \mathcal{Q}$ either $Q \subset L$ or 
$Q \subset R$. 
\end{lemma}

\proof Let $\,\mathbf{Q}_0$ denote the set of those $\mathcal{Q} \in \mathbf{Q}$ which contain an element $Q$ which intersects both 
$L$ and $R$. Then $\mathbf{Q}_0$ is $\Phi$-invariant: Let $\mathcal{Q} \in \mathbf{Q}$; then each component of $\Phi(\mathcal{\mathcal{Q}})$ is either equal to a 
component of $\mathcal{Q}$ or is the union of two adjacent components of $\mathcal{Q}$. Thus if $\mathcal{Q} \in \mathbf{Q}_0$ then
$\Phi(\mathcal{Q}) \in \mathbf{Q}_0 $. Note that $\mathcal{Q}_2 \notin \mathbf{Q}_0$ but $\mathcal{Q}_T \in \mathbf{Q}_0$.  
 Suppose $\mathbf{Q}_0 \setminus \{\mathcal{Q}_T\}$ is non-empty ; then by Proposition~\ref{prop_fs_6} it 
contains an element $\mathcal{U}$ which is maximum with respect to the total order $\le$ and $\mathcal{U} < \mathcal{Q}_2$, since
$\mathcal{Q}_2 \notin \mathbf{Q}_0$. But then $\Phi(\mathcal{U}) \in \mathbf{Q}_0\setminus\{\mathcal{Q}_T\}$, which contradicts the maximality of $\mathcal{U}$.
Therefore $\mathbf{Q}_0 \setminus \{\mathcal{Q}_T\} = \varnothing$ and hence for each $\mathcal{Q} \in \mathbf{Q}_S$ and each 
$Q \in \mathcal{Q}$ either $Q \subset L$ or $Q \subset R$. \eop 

For each $\mathcal{Q} \in \mathbf{Q}_S$ let $\mathcal{Q}_L = \{ Q \in \mathcal{Q} : Q \subset L \}$ and 
$\mathcal{Q}_R = \{ Q \in \mathcal{Q} : Q \subset R \}$. Then $\mathcal{Q}_L$ is a partition of $L$ and $\mathcal{Q}_R$  a partition of $R$.
Moreover $\mathcal{Q}_L$ and $\mathcal{Q}_R$ are disjoint subsets of $\mathcal{Q}$ whose union is $\mathcal{Q}$.

Put $\,\mathbf{Q}_Z  = \{ \mathcal{Q} \in \mathbf{Q}_S : \Phi(\mathcal{Q}) \in \mathbf{Q}_S \}$ and so
$\,\mathbf{Q}_Z = \mathbf{Q}_S \setminus \{\mathcal{Q}_2\}$. 
If $\mathcal{Q} \in \mathbf{Q}_Z$ with $\Phi(\mathcal{Q}) = \mathcal{Q}(J,K)$; then either $J\cup K \subset L$ or $J\cup K \subset R$. 
Thus we can define a mapping  $\lambda : \mathbf{Q}_Z  \to \{L,R\}$ by letting $ \lambda(\mathcal{Q}) = L$ if $J\cup K \subset L$ and $ \lambda(\mathcal{Q}) = R$ 
if $J\cup K \subset R$.

Let $\mathcal{Q} \in \mathbf{Q}_Z$ and so $\Phi(\mathcal{Q})$ is a reduction of $\mathcal{Q}$. If $\lambda(\mathcal{Q}) = L$ then   
$\Phi(\mathcal{Q})_L \prec_0\mathcal{Q}_L$ and $\Phi(\mathcal{Q})_R = \mathcal{Q}_R$ 
and if $\lambda(\mathcal{Q}) = R$ then $\Phi(\mathcal{Q})_L = \mathcal{Q}_L$ and $\Phi(\mathcal{Q})_R \prec_0\mathcal{Q}_R$.

Let $\,\mathbf{U}$ be the set of $L$-partitions having the form $\mathcal{Q}_L$ for some $\mathcal{Q} \in \mathbf{Q}_S$. 
Thus $\mathcal{Q}_L \in \mathbf{U}$ for each $\mathcal{Q} \in\mathbf{Q}_S$ but in general the representation as $\mathcal{Q}_L$ cannot be unique.
Note that $\,\mathbf{U}$ contains the elements $\mathcal{U}_0 = (\mathcal{Q}_0)_L$ and $\mathcal{U}_T =(\mathcal{Q}_2)_L = \{L\}$.

For each $\mathcal{U} \in \mathbf{U}$ let $\mathbf{Q}_{\mathcal{U}} = \{ \mathcal{Q} \in \mathbf{Q}_S : \mathcal{Q}_L = \mathcal{U} \}$ and let
$\mathcal{U}^\triangleleft$ be the largest element in $\mathbf{Q}_{\mathcal{U}}$. 
Then $\mathcal{U} = \mathcal{U}^\triangleleft_L$, since 
$\mathcal{U}^\triangleleft \in \mathbf{Q}_{\mathcal{U}}$.

Define a mapping $\Phi_L : \mathbf{U} \to \mathbf{U}$ by letting $\Phi_L(\mathcal{U}) = \Phi(\mathcal{U}^\triangleleft)_L$ if 
$\mathcal{U}^\triangleleft \in \mathbf{Q}_Z$ and letting $\Phi_L(\mathcal{U}) = (\mathcal{U}_T$ if $\mathcal{U}^\triangleleft = \mathcal{Q}_2$. In particular
$\Phi_L(\mathcal{U}_T)= \mathcal{U}_T$.

We thus have the iterator $\,\mathbf{R}_L = (\mathbf{U},\Phi_L,(\mathcal{U}_0)$. 

Let $\,\mathbf{U}_S = \mathbf{U} \setminus \{\mathcal{U}_T\}$, thus $\mathcal{U} \in \mathbf{U}_S$ if and only if $\mathcal{U}^\triangleleft \in \mathbf{Q}_Z$.

\begin{lemma}\label{lemma_assoc_44} (1)\enskip 
Let $\mathcal{U} \in\mathbf{U}_S$ and so  $\mathcal{Q} = \mathcal{U}^\triangleleft \in \mathbf{Q}_Z$. Then 
$\Phi(\mathcal{Q})$
is a reduction of $\mathcal{Q}$ with $\lambda(\mathcal{Q}) = L$ and $\Phi_L(\mathcal{U})$ is an
$L$-reduction of $\mathcal{U}$. 

(2)\enskip Let $\mathcal{Q} \in \mathbf{Q}_Z$ with $\lambda(\mathcal{Q}) = L$. Then $\mathcal{Q} = (\mathcal{Q}_L)^\triangleleft$
and $\Phi_L(\mathcal{Q}_L) = \Phi(\mathcal{Q})_L$.
\end{lemma}

\proof (1) \enskip
We have
$\Phi(\mathcal{U}^\triangleleft) \notin \mathbf{Q}_{\mathcal{U}}$ and so $\Phi(\mathcal{U}^\triangleleft)_L \ne \mathcal{U} = \mathcal{U}^\triangleleft_L$, i.e.,
$\Phi(\mathcal{U}^\triangleleft)_L \ne \mathcal{U}^\triangleleft_L$. Thus
$\Phi(\mathcal{U}^\triangleleft) \ne\mathcal{U}^\triangleleft$ and hence
$\Phi(\mathcal{U}^\triangleleft)$ is a reduction of $\mathcal{U}^\triangleleft$ with $\lambda(\mathcal{U}^\triangleleft) = L$. It follows that 
$\Phi_L(\mathcal{U}) = \Phi(\mathcal{U}^\triangleleft)_L$ is an
$L$-reduction of $\mathcal{U}$. 

(2) \enskip We have $\Phi(\mathcal{Q})$ is a reduction of $\mathcal{Q}$ and $\lambda(\mathcal{Q}) = L$ and so $\Phi(\mathcal{Q})_L \ne \mathcal{Q}_L$. 
Thus $\mathcal{Q}$ is the largest element $\mathcal{V} \in \mathbf{Q}_S$ such that 
$\mathcal{V}_L = \mathcal{Q}_L$ and hence $\mathcal{Q} = (\mathcal{Q}_L)^\triangleleft$.
It follows that $\Phi_L(\mathcal{Q}_L) = \Phi((\mathcal{Q}_L)^\triangleleft) = \Phi(\mathcal{Q})_L$. \eop

Let $\,\mathbf{Q}_L = \{ \mathcal{Q} \in \mathbf{Q}_Z : \mathcal{Q} = \mathcal{U}^\triangleleft \mbox{ for some $\mathcal{U} \in \mathbf{U}_S$} \}$.
If $\mathcal{Q} = \mathcal{U}^\triangleleft$ then $\mathcal{Q}_L = \mathcal{U}$ and so $\mathcal{U}$ is uniquely determined by $\mathcal{Q}$. Thus there is a
mapping $\delta_L : \mathbf{Q}_L \to \mathbf{U}_S$ such that $\delta_L(\mathcal{Q}) = \mathcal{U}$
whenever $\mathcal{Q} = \mathcal{U}^\triangleleft$.
If $\delta_L(\mathcal{Q})= \delta_L(\mathcal{Q}')$ then $\mathcal{Q} = \mathcal{Q}' = \mathcal{U}^\triangleleft$ for some $\mathcal{U} \in \mathbf{Q}_L$
and hence $\delta_L$ is injective. But $\delta_L$ is clearly also surjective and therefore $\delta_L$ is a bijection. The inverse of $\delta_L$ is the
bijective mapping $\gamma_L : \mathbf{U}_S \to \mathbf{Q}_L$ given by $\gamma_L(\mathcal{U}) = \mathcal{U}^\triangleleft$ for all $\mathcal{U} \in \mathbf{U}_S$.

\begin{lemma}\label{lemma_assoc_45} 
$\,\mathbf{Q}_L = \{ \mathcal{Q} \in \mathbf{Q}_Z : \mbox{ $\Phi(\mathcal{Q})$ is a reduction of $\mathcal{Q}$ with $\lambda(\mathcal{Q}) = L$} \}$.
\end{lemma}

\proof
If $\mathcal{Q} \in \mathbf{Q}_L$ then by Lemma~\ref{lemma_assoc_44} (1) 
$\Phi(\mathcal{Q})$ is a reduction of $\mathcal{Q}$ with
$\gamma(\mathcal{Q}) = L$. 
Conversely, if $\Phi(\mathcal{Q})$ is a reduction of $\mathcal{Q}$ with $\gamma(\mathcal{Q}) = L$ then by  Lem's~\ref{lemma_assoc_44} (2)
$\mathcal{Q} = (\mathcal{Q}_L)^\triangleleft$ and so  $\mathcal{Q} \in \mathbf{Q}_L$. \eop

Let $\,\mathbf{V}$ be the set of $R$-partitions having the form $\mathcal{Q}_R$ for some $\mathcal{Q} \in \mathbf{Q}_S$. 
Note that $\,\mathbf{V}$ contains the elements $\mathcal{V}_0 = (\mathcal{Q}_0)_R$ and $\mathcal{V}_T =(\mathcal{Q}_2)_R = \{R\}$.

For each $\mathcal{V} \in \mathbf{V}$ let $\mathbf{Q}_{\mathcal{V}} = \{ \mathcal{Q} \in \mathbf{Q}_S : \mathcal{Q}_R = \mathcal{V} \}$ and let
$\mathcal{V}^\triangleleft$ be the largest element in $\mathbf{Q}_{\mathcal{V}}$. 
Then $\mathcal{V} = \mathcal{V}^\triangleleft_R$, since 
$\mathcal{V}^\triangleleft \in \mathbf{Q}_{\mathcal{V}}$.

Define a mapping $\Phi_R : \mathbf{V} \to \mathbf{V}$ by letting $\Phi_R(\mathcal{V}) = \Phi(\mathcal{V}^\triangleleft)_R$ if 
$\mathcal{V}^\triangleleft \in \mathbf{Q}_Z$ and letting $\Phi_R(\mathcal{V}) = (\mathcal{V}_T$ if $\mathcal{V}^\triangleleft = \mathcal{Q}_2$. In particular
$\Phi_R(\mathcal{V}_T)= \mathcal{V}_T$.

We thus have the iterator $\,\mathbf{R}_R = (\mathbf{V},\Phi_R,(\mathcal{V}_0)$. 

Let $\,\mathbf{V}_S = \mathbf{V} \setminus \{\mathcal{V}_T\}$, thus $\mathcal{V} \in \mathbf{V}_S$ if and only if $\mathcal{V}^\triangleleft \in \mathbf{Q}_Z$.

\begin{lemma}\label{lemma_assoc_54} (1)\enskip 
Let $\mathcal{V} \in\mathbf{V}_S$ and so  $\mathcal{Q} = \mathcal{V}^\triangleleft \in \mathbf{Q}_Z$. Then 
$\Phi(\mathcal{Q})$
is a reduction of $\mathcal{Q}$ with $\lambda(\mathcal{Q}) = R$ and $\Phi_R(\mathcal{V})$ is an
$R$-reduction of $\mathcal{V}$. 

(2)\enskip Let $\mathcal{Q} \in \mathbf{Q}_Z$ with $\lambda(\mathcal{Q}) = R$. Then $\mathcal{Q} = (\mathcal{Q}_R)^\triangleleft$
and $\Phi_R(\mathcal{Q}_R) = \Phi(\mathcal{Q})_R$.
\end{lemma}

\proof This is the same as the proof Lemma~\ref{lemma_assoc_44}. \eop

Let $\,\mathbf{Q}_R = \{ \mathcal{Q} \in \mathbf{Q}_Z : \mathcal{Q} = \mathcal{V}^\triangleleft \mbox{ for some $\mathcal{V} \in \mathbf{V}_S$} \}$.
Then as above thee is a bijective mapping $\delta_R : \mathbf{Q}_R \to \mathbf{V}_S$ such that $\delta_R(\mathcal{Q}) = \mathcal{V}$
whenever $\mathcal{Q} = \mathcal{V}^\triangleleft$. The inverse of $\delta_R$ is the
bijective mapping $\gamma_R : \mathbf{V}_S \to \mathbf{Q}_R$ given by $\gamma_R(\mathcal{V}) = \mathcal{V}^\triangleleft$ for all $\mathcal{V} \in \mathbf{V}_S$.

\begin{lemma}\label{lemma_assoc_55} 
$\,\mathbf{Q}_R = \{ \mathcal{Q} \in \mathbf{Q}_Z : \mbox{ $\Phi(\mathcal{Q})$ is a reduction of $\mathcal{Q}$ with $\lambda(\mathcal{Q}) = R$} \}$.
\end{lemma}

\proof This is the same as the proof of Lemma~\ref{lemma_assoc_45}. \eop

\begin{lemma}\label{lemma_assoc_56} 
The sets $\,\mathbf{U}$ and $\,\mathbf{V}$ are disjoint and $\,\mathbf{U} \cup \mathbf{V} \approx \mathbf{Q}$.
\end{lemma}

\proof
It is clear that $\,\mathbf{U}$ and $\,\mathbf{V}$ are disjoint and by Lemmas \ref{lemma_assoc_45} and \ref{lemma_assoc_55} 
and (since there are bijections $\gamma_L: \mathbf{U}_S \to \mathbf{Q}_L$ and $\gamma_R: \mathbf{V}_S \to \mathbf{Q}_R$) it then follows that 
$\,\mathbf{U}_S \cup \mathbf{V}_S \approx \mathbf{Q}_Z$. Thus $\,\mathbf{U} \cup \mathbf{V} \approx \mathbf{Q}$,
since $\,\mathbf{U} = \mathbf{U}_S \cup \{\mathcal{U}_T\}$, $\,\mathbf{V} = \mathbf{V}_S \cup \{\mathcal{V}_T\}$ and
$\mathbf{Q} = \mathbf{Q}_Z \cup \{\mathcal{Q}_2, \mathcal{Q}_T\}$. \eop

\begin{lemma}\label{lemma_assoc_57} 
The iterator $\,\mathbf{R}_L$ is an $L$-reduction and $\,\mathbf{R}_R$ is an $R$-reduction.
\end{lemma}

\proof By Lemma~\ref{lemma_assoc_44} (1) $\Phi_L(\mathcal{U})$ is an $L$-reduction of $\mathcal{U}$ for each $\mathcal{U} \in\mathbf{U}_S$ 
and by Lemma~\ref{lemma_assoc_54} (1) $\Phi_R(\mathcal{V})$ is an $R$-reduction of $\mathcal{V}$ for each $\mathcal{V} \in\mathbf{V}_S$.
It remains to show that the iterators $\,\mathbf{R}_L$ and $\,\mathbf{R}_R$ are minimal and by  Proposition~\ref{prop_assoc_101} this amounts to
showing that $\,\mathbf{U} \approx L$ and $\,\mathbf{V} \approx R$.
By Proposition~\ref{prop_assoc_101} $\,\mathbf{Q} \approx I$, $L \preceq \mathbf{U}$ and $R \preceq \mathbf{V}$. Also by Lemma~\ref{lemma_assoc_56}
$\,\mathbf{U}$ and $\,\mathbf{V}$ are disjoint with $\,\mathbf{U} \cup \mathbf{V} \approx \mathbf{Q}$, and of course $L$ and $R$ are disjoint with
$L \cup R = I$.
Suppose $L \prec \mathbf{U}$; then $ I= L \cup R \prec \mathbf{U} \cup \mathbf{V} \approx \mathbf{Q}$ and this contradiction implies that
$\,\mathbf{U} \approx L$. In the same way $\,\mathbf{V} \approx R$.
\eop

If $\mathcal{U}$ is an $L$-partition and $\mathcal{V}$  an $R$-partition then the $I$-partition $\mathcal{Q}$ with $\mathcal{Q}_L = \mathcal{U}$ 
and $\mathcal{Q}_R = \mathcal{V}$ will be denoted by $\mathcal{U} + \mathcal{V}$. If $\mathcal{Q} \in \mathbf{Q}_S$  
then $\mathcal{Q} = \mathcal{Q}_L + \mathcal{Q}_R$ and $\mathcal{Q}_L \in \mathbf{U}$ and $\mathcal{Q}_R \in \mathbf{V}$.
Let $\mathcal{Q}$ be an $I$-partition with $\mathcal{Q} = \mathcal{U} + \mathcal{V}$ and let $f_L : \mathcal{U} \to X$ and $f_R : \mathcal{V} \to X$ be mappings. 
Then there is a mapping $f_L \oplus f_R : \mathcal{Q} \to X$ defined by
\[ (f_L \oplus f_R)(J) = \left\{ \begin{array}{cl}
                  f_L(J) &\ \mbox{if}\ J \in \mathcal{U}\;,\\
                  f_R(J)   &\  \mbox{if}\ J \in \mathcal{V}\;.\\
\end{array} \right. \]

Consider the mappings given by Proposition~\ref{prop_assoc_13} for $\,\mathbf{Q},\,\mathbf{U}$ and $\,\mathbf{V}$.
Thus for each $\mathcal{Q} \in \mathbf{Q}$ there is the mapping $\omega_{\mathcal{Q}} : \mathcal{Q} \to X$,
for each $\mathcal{U} \in \mathbf{U}$ there is the mapping $\omega^L_{\mathcal{U}} : \mathcal{U} \to X$
and for each $\mathcal{V} \in \mathbf{V}$ there is the mapping $\omega^R_{\mathcal{V}} : \mathcal{V} \to X$.

\begin{lemma}\label{lemma_assoc_43} 
Let $(x_p,\ldots,x_q)$ be a $[p,q]$-tuple. Then 
$\omega_{\mathcal{Q}} = \omega^L_{\mathcal{Q}_L} \oplus \omega^R_{\mathcal{Q}_R}$ 
for all $\mathcal{Q} \in \mathbf{Q}_S$
and $\omega_{\mathcal{Q}_T} = \omega^L_{\{L\}} \bullet \omega^R_{\{R\}}\:$. Thus
\[\bullet_\mathbf{R}(x_p,\ldots,x_q) = (\bullet_{{\mathbf{R}_L}}(x_p,\ldots,x_t)) \bullet (\bullet_{\mathbf{R}_R}(x_{f(t)},\ldots,x_q)). \] 
\end{lemma} 

\proof Let $\,\mathbf{Q}' = \{ \mathcal{Q} \in \mathbf{Q}_S : \omega_{\mathcal{Q}} = \omega^L_{\mathcal{Q}_L} \oplus \omega^R_{\mathcal{Q}_R}\}$ and so
$\mathcal{Q}_0 \in \mathbf{Q}'$. Let $\mathcal{Q} \in \mathbf\\{Q}' \cap \mathbf{Q}_Z$. If $\lambda(\mathcal{Q}) = L$ then by Lemma~\ref{lemma_assoc_44} (2)
$\Phi(\mathcal{Q}) = \Phi(\mathcal{Q})_L + \Phi(\mathcal{Q})_R =  \Phi_L(\mathcal{Q}_L) + \mathcal{Q}_R$ and thus
$\omega_{\Phi(\mathcal{Q})} = \omega^L_{\Phi_L(\mathcal{Q}_L)} \oplus \omega^R_{\mathcal{Q}_R}
= \omega^L_{{\Phi(\mathcal{Q})_L}} \oplus \omega^R_{{\Phi(\mathcal{Q})_R}}\;$.
In exactly the same way, if $\lambda(\mathcal{Q}) = R$ then by Lemma~\ref{lemma_assoc_54} (2)
$\Phi(\mathcal{Q}) = \Phi(\mathcal{Q})_L + \Phi(\mathcal{Q})_R = \mathcal{Q}_L+ \Phi_R(\mathcal{Q}_R)$ and hence 
$\omega_{\Phi(\mathcal{Q})} = \omega^L_{\mathcal{Q}_L} \oplus \omega^R_{\Phi_R(\mathcal{Q}_R)}  
= \omega^L_{{\Phi(\mathcal{Q})_L}} \oplus \omega^R_{{\Phi(\mathcal{Q})_R}}$.
This shows that $\Phi(\mathcal{Q}) \in \mathbf{Q}'$. Thus $\mathbf{Q}'$ is a $\Phi$-invariant subset of $\,\mathbf{Q}$ containing $\mathcal{Q}_0$
and so $\,\mathbf{Q}' = \mathbf{Q}_S$. Therefore $\omega_{\mathcal{Q}} = \omega^L_{\mathcal{Q}_L} \oplus \omega^R_{\mathcal{Q}_R}$ for all
$\mathcal{Q} \in \mathbf{Q}_S$. 
Finally, $\mathcal{Q}_T = \Phi(\mathcal{Q}_2) = \Phi(\{L\} + \{R\})$
and it follows that $\omega_{\mathcal{Q}_T} = \omega_{\Phi(\mathcal{Q}_2)} = \omega^L_{\{L\}} \bullet \omega^R_{\{R\}}\:$.
\eop

\it{Proof of Theorem~\ref{theorem_assoc_2}}: \enskip
If $B$ is any non-empty finite set then by Lemma~\ref{lemma_fs_7} there exists an $I$-reduction with $I \approx B$.
For each non-empty finite set $B$  let $\prop(B)$ be the statement that $\bullet_\mathbf{R}(x_p,\ldots,x_q) = \bullet(x_p,\ldots,x_q)$
holds for each $I$-tuple $(x_p,\ldots,x_q)$ whenever $\,\mathbf{R}$ is an $I$-reduction with $I \approx B$. Note that if $B_1 \approx B_2$
then $\prop(B_1)$ holds if and only if $\prop(B_2)$ holds.
Suppose there exists a non-empty finite set $B$ for which $\prop(B)$ does not hold and let $\mathcal{S}$ denote the set of non-empty subsets $C$ of $B$
for which $\prop(C)$ does not hold. Then $\mathcal{S}$ is non-empty and so by Proposition~\ref{prop_intro_3} it contains a minimal element
$D$ and $D$ contains at least two elements, since $\prop(\{a\})$ holds for each element $a$, because if $I = \{a\}$ then there is only one $I$-reduction.
If $D'$ is any non-empty finite set with $D' \prec D$ then $\prop(D')$ holds.
 
Now let  $I = [p,q]$ with $I \approx D$. Since $\prop(D)$ does not hold there exists an $I$-reduction $\mathbf{R}$ and a $[p,q]$-tuple $(x_p,\ldots,x_q)$ 
such that $\bullet_\mathbf{R}(x_p,\ldots,x_q) \ne \bullet(x_p,\ldots,x_q)$.
Then by Lemma~\ref{lemma_assoc_43}
$\bullet_\mathbf{R}(x_p,\ldots,x_q) = (\bullet_{{\mathbf{R}_L}}(x_p,\ldots,x_t)) \bullet (\bullet_{\mathbf{R}_R}(x_{t+1},\ldots,x_q))$ 
and $L \prec D$, $R\prec D$. Thus $\prop(L)$ and $\prop(R)$ hold and so
\[\bullet_\mathbf{R}(x_p,\ldots,x_q) = (\bullet(x_p,\ldots,x_t) \bullet (\bullet(x_{t+1},\ldots,x_q)\;. \]
Hence by Theorem~\ref{theorem_assoc_1} 
$\bullet_\mathbf{R}(x_p,\ldots,x_q)= \bullet(x_p,\ldots,x_q)$, and this contradiction shows that $\prop(B)$ holds for each non-empty finite set $B$.
Therefore
 \[\bullet_\mathbf{R}(x_p,\ldots,x_q) = \bullet(x_p,\ldots,x_q)\]
for each $I$-tuple $(x_p,\ldots,x_q)$ and each $I$-reduction $\mathbf{R}$.
\eop


\startsection{Lists}

\label{lists}

In this section we present an approach to dealing with (finite) lists taking their values in some fixed class $E$. We first formulate things using the natural numbers but
then employ a  general Peano iterator. For each $n \in \Nat$ lists of length $n$ are mappings from $L_n =\{0,1.\ldots,n-1\}$. to $E$.
A list $t : L_n \to E$ with $t(k) = e_k$ for all $k \in L_n$ will be represented in the form $[e_0,e_1,\ldots,e_{n-1}]$. In particular, $[]$
represents the empty list, i.e., the list with no elements.

When using the natural numbers the sets $L_n = \{0,1,\ldots,n-1\}$ play an important role and these correspond to the sets in the finite segment iterator associated
with a general Peano iterator. From now on we only use the natural numbers when giving examples. 
Let us fix a Peano iterator
$\,\mathbf{I} = (N,s,0)$ and to reduce the clutter 

we prefer to use $0$ instead of $n_0$ for the third component of $\,\mathbf{I}$. Let $\le$ be the unique total 
order 
on $N$ compatible with $\,\mathbf{I}$ given in Theorem~\ref{theorem_iterators_199}
and let $\,\mathbf{I}_\le = (N_\le,s_\le,\varnothing)$ be the finite segment iterator associated with $\,\mathbf{I}$  given in Theorem~\ref{theorem_iterators_99}.
Therefore $N_\le = \{L_n : n \in N \}$ with $L_n = \{ m \in N : m < n \}$ and $s_\le(L_n) = L_{s(n)}$ for all $n \in N$.

For each $n\in N$ denote by $E^*_n$ the class of all mappings from $L_n$ to $E$. In particular, 
$L_{0} = \varnothing$ and so $E^*_{0}$ consists of the unique mapping from $\varnothing$ to $E$, and this mapping we denote by $\varepsilon$. 

We call $E^*_n$  \definition{the class of lists based on $L_n$ with
values in $E$}. The single element $\varepsilon \in E^*_{0}$ will be referred to as the empty list.
Let $m,\,n \in N$ with $m \ne n$; then $L_m \ne L_n$ and hence $E^*_m$ and $E^*_n$ are disjoint.
Put $E^* = \bigcup_{n\in N} E^*_n$; this is the \definition{class of all lists with values in $E$}. Also let
$E^*_+ = E^* \setminus \{\varepsilon\}$.

We start by defining mappings $\triangleleft : E \times E^* \to E^*_+$ 
and $\triangleright : E^*_+ \to E \times E^*$. These are bijections and  each is the inverse of the other.

In terms of the natural numbers these mappings can be described as follows:
The list $\triangleleft(e,r)$ is obtained by adding the element $e$ to the beginning of the list $r$ and so
$\triangleleft(e,[e_0,e_1,\ldots,e_{n-1}]) = [e,e_0,e_1\ldots,e_{n-1}]$. In particular $\triangleleft(e,[\;]) = [e]$ is a non-empty list.
If $r$ is a non-empty list then the first component of $\triangleright$
is the first element of $r$  (the head of the list) and the second component of $\triangleright(r)$ is the
rest of the list (its tail). Thus $\triangleright([e_0,e_1,\ldots,e_{n-1}]) = (e_0,[e_1,\ldots,e_{n-1}])$.

In the following we make use of Theorem~\ref{theorem_iterators_199} (4), which states 
that $L_{s(n)}$ is the disjoint union of $\{0\}$ and $s(L_n)$.
For each  $n \in N$  we first define a mapping  $\triangleleft_n : E \times E^*_n \to E^*_{s(n)}$.

For $(e,r) \in E \times E^*_n$  let  $\triangleleft \in E^*_{s(n)}$ be given by  $\triangleleft_n(e,r)(0) = e$ and $\triangleleft(e,r)(s(k)) = r(k)$ for all $k \in L_n$. 
Note that $L_{s(0)} =\{0\}$ and so $E^*_{s(0)}$ is the class of all mappings from
$\{0\}$ to $E$. Thus $\triangleleft_{0} : E \times E^*_{0} \to E^*_{s(0)}$ is the mapping with $\triangleleft_{0} (e,\varepsilon)(0) = e$. 

Now define $\triangleleft : E \times E^* \to E^*$ by letting $\triangleleft(e,r) = \triangleleft_n(e,r)$ for all $r \in E^*_n$.

If $r \in E^*_+$ then there exists a unique
$m \in N$ such that $r \in E^*_m$ and $m \ne 0$ since $s \ne \varepsilon$. Thus by Proposition~\ref{prop_iterators_1} there exists a unique
$n \in N$ such that $r \in E^*_{s(n)}$  and so $r$ is a mopping from $\{0\} \cup s(L_n)$ to $E$. For each $n \in N$ we next define a mapping
$\triangleright_n : E^*_{s(n)} \to E \times E^*_n$. For each $v \in E^*_{s(n)}$ let $\triangleright_n(v) = (v(0),r)$, where $r \in E^*_n$ is given by
$r(m) = s(v(m))$ for all $m \in L_n$.

Now define $\triangleright : E^*_+ \to E \times E^*$ by letting $\triangleright(v) = \triangleright_n(v)$ for all $v \in E^*_{s(n)}$.

\begin{proposition}\label{prop_lists_1}
The mappings $\triangleleft : E \times E^* \to E^*_+$ 
and $\triangleright : E^*_+ \to E \times E^*$ are inverse to each other and in particular they are both bijections.
\end{proposition}

\proof
Let $v\in E^*_{s(n)}$; then $\triangleleft(\triangleright(v)) = \triangleleft_n(\triangleright_n(v)) = \triangleright_n(v(0),r)$,where $r \in E^*_n$ is given by
$r(m) = s(v(m))$ for all $m \in L_n$.  Therefore $\triangleleft(\triangleright(v)) = u \in E^*_{s(n)}$, where $u(0) = v(0)$
and $u(s(m)) = r(m) = s(v(m))$ for all $m \in L_n$. Hence $u = v$ which shows that $\triangleleft(\triangleright(v)) = v$.

Now let $(e,r) \in E \times E^*_n$; then $\triangleleft(e,r) = \triangleleft_n(e,r) = v$, where $v \in E^*_{s(n)}$ is given by $v(0) = e$
and $v(r(m)) = r(m)$ for all $m \in L_n$ and hence $\triangleright(\triangleleft((e,r)) = \triangleright(v) = (v(0), t)$, where $t \in E^*_n$ is given by
$t(m) =s(v(m)) = r(m)$ for all $m \in L_n$, i.e., $t = r$. Also $v(0) = e$ and this shows that $\triangleright(\triangleleft(e,r)) = (e,r)$.
\eop

Note that by Proposition~\ref{prop_lists_1} each list in $E^*_{s(n)}$ has a unique representation in the form $\triangleleft(e,r)$ with
$(e,r) \in E \times E^*_n$.

A triple $(X,f,x_0)$ with $X$ a class, $f : E \times X \to X$ a mapping and $x_0$ an element of
$X$ will be called a \definition{list algebra}. (The class $E$ is considered to be fixed here.) Thus $(E^*,\triangleleft,\varepsilon)$ is a list algebra. If
$(X,f,x_0)$ is a list algebra then for each $e \in E$ let $f_e : X \to X$ be the mapping with
$f_e(x) = f(e,x)$ for all $x \in X$. There is then the iterator $(X,f_e,x_0)$.
A subclass $X_0$ of $X$ is said to be \definition{$f$-invariant} if $f_e(X_0) \subset X_0$ for all $e \in E$ and
$(X,f,x_0)$ is said to be \definition{minimal} if $X$ itself is the only $f$-invariant subclass of $X$ containing $x_0$.

\begin{lemma}\label{lemma_lists_1}
The list algebra $(E^*,\triangleleft,\varepsilon)$ is minimal.
\end{lemma} 

\proof
Let $G$ be a $\triangleleft$-invariant subclass of $E^*$ containing $\varepsilon$ and for each $n \in N$ put $G_n = G \cap E^*_n$.
Let $N_0 = \{ n \in N :G_n = E^*_n \}$. In particular, $0 \in N_0$, since $E^*_{0} = \{\varepsilon\}$. Thus consider $n\in N_0$,
and so $G_n = E^*_n$. Let $t \in E^*_{s(n)}$; then by Proposition~\ref{prop_lists_1} there exists $e \in E$ and $r\in E^*_n$ such that
$t = \triangleleft(e,r) = \triangleleft_e(r)$. But $r \in G_n$, since $n \in N_0$ and thus $t \in G_{s(n)}$, since $G$ is $\triangleleft$-invariant. Hence
$G_{s(n)} = E^*_{s(n)}$. Therefore $N_0$ is $f$-invariant, and since $0 \in N_0$, this implies $N_0 = N$, i.e., $G = E^*$.
This shows that $(E^*,\triangleleft,\varepsilon)$ is minimal.\eop

\begin{lemma}\label{lemma_lists_2}
The list algebra $(E^*,\triangleleft,\varepsilon)$ has the following properties:

(1)\enskip The mapping $\triangleleft_e : E^*\to E^*$ is injective for each $e \in E$.

(2)\enskip If $e,\,e' \in E$ with $e \ne e'$ then the classes $\triangleleft_e(E^*)$ and $\triangleleft_{e'}(E^*)$ are disjoint.

(3)\enskip $\varepsilon \notin \triangleleft_e(E^*)$ for all $e \in E$.
\end{lemma} 

\proof (1)\enskip Let $e \in E$ and $u,\,v \in E^*$ with $u \ne v$ and $u \in E^*_m$, $v \in E^*_n$ and  thus $\triangleleft_e(u) \in E^*_{s(m)}$,
$\triangleleft_e(v) \in E^*_{s(n)}$, If $m \ne n$ then $s(m) \ne s(n)$ and in this case $\triangleleft_e(u) \ne \triangleleft_e(v)$ holds
trivially, since $E^*_{s(m)}$ and $E^*_{s(n)}$ are disjoint. We can thus assume that $m = n$ and so $u(k) \ne v(k)$ for some $k \in L_n$.
Now each $r \in E^*_n$ we have  $\triangleleft_e(r) = t$, where  $t \in E^*_{s(n)}$ is given by $t(0) =  e$  and $t(s(j)) = r(j)$ for all $j \in L_n$. 
In particular $\triangleleft_e(u(s(k))) = u(k) \ne v(k) = \triangleleft_e(v(s(k)))$ and hence
$\triangleleft_e(u) \ne \triangleleft_e(v)$. Thus the mapping $\triangleleft_e$ is injective.

(2)\enskip This is clear since $\triangleleft_e(r)(x_0) = e$ for all $r\in \triangleleft_e(E^*)$.

(3)\enskip This is also clear. \eop

Let $(X,f,x_0)$ and $(Y,g,y_0)$ be list algebras.
A  morphism $\pi : (X,f,x_0)) \to (Y,g,y_0)$ is then a mapping $\pi : X \to Y$ with $\pi(x_0) = y_0$ such that
$  g_e \circ\ \pi = \pi \circ f_e$ for all $e \in E$.  Thus $\pi : (X,f,x_0) \to (Y,g,y_0)$ being a morphism means exactly that $\pi (X,f_e,x_0) \to (Y,g_e,y_0)$
is a morphism of iterators for each $e \in E$.
A list algebra $(X,f,x_0)$ is said to be \definition{initial} if for each list algebra $(Y,g,y_0)$ there exists a unique morphism
$\pi : (X,f,x_0) \to (Y,g,y_0)$.

\begin{theorem}\label{theorem_lists_1}
The list algebra $(E^*,\triangleleft,\varepsilon)$ is initial.
\end{theorem}

\proof Let $(X,f,x_0)$ be a list algebra. 
We must show that there is a unique mapping $\pi : E^* \to X$ with $\pi(\varepsilon) = x_0$ and such that 
\[\pi(\triangleleft(e,r)) = f(e,\pi(r))\] 
for all $(e,r) \in E\times E^*$.

For each $n \in N$ put $E^*_{\le n} = \bigcup_{y \le n} E^*_y$ and $E^*_{< n} = \bigcup_{y < n} E^*_y$.
Let $n \in N$; a  mapping $\pi_n : E^*_{\le n} \to X$ with 
$\pi_n(\varepsilon) = x_0$ and
$\pi_n(\triangleleft(e,r)) = f(e,\pi_n(r))$ for all $(e,r) \in E \times E^*_{<n}$ will be called a partial $n$-solution.
Let $N_0$ be the subclass of $N$ consisting of those $n \in N$ for which there exists a unique partial $n$-solution.
There is a unique partial $0$-solution $\pi_{0} : E^*_{\{0\}} \to X$, which must be defined by putting $\pi_{0}\triangleleft(e,\varepsilon) = f(e,x_0)$ for all $e \in E$.
 Hence $0 \in N_0$.

Thus let $n \in  N_0$ with unique partial $n$-solution 
$\pi_n : E^*_{\le n} \to X$.  

We must  define an appropriate mapping
$\pi_{s(n)} : E^*_{\le s(n)} \to X$ and $E^*_{\le s(n)}$ is the disjoint union of $E^*_{\le n}$ and $E^*_{s(n)}$.
Define $\pi_{s(n)}(t)  = \pi_n(t)$ for all $t \in E^*_{\le n}$. 
In particular $\pi_{s(n)}(\varepsilon) = x_0$ and $\pi_{s(n)}(\triangleleft(e,r)) = f(e,\pi_{s(n)}(r))$ for all $(e,r) \in E \times E^*_{<n}$.   
Let $u \in E^*_{s(n)}$. Then $u$ has a unique representation in the form $u = \triangleleft(e,t)$ with
$(e,t) \in E \times E^*_n$ and we put $ \pi_{s(n)}(\triangleleft(e,t)) = f(e,\pi_n(t))$ and so $ \pi_{s(n)}(\triangleleft(e,t)) = f(e,\pi_{s(n)}(t))$,  since
$t \in E^*_{\le n}$. This shows that $\pi_{s(n)}$ is a partial $s(n)$-solution and it is in fact the unique  partial $s(n)$-solution:
The restriction of $\pi_{s(n)}$ to $E^*_{\le n}$ is a partial $n$-solution and so is equal to $\pi_n$, since $n \in N_0$, a requirement which uniquely determines
$\pi_{s(n)}$ on $E^*_{\le n}$. Moreover, it then follows that $\pi_{s(n)}$ is also uniquely determined on $E^*_{s(n)}$. Hence $s(n) \in N_0$ and therefore
$N_0 = N$, since $\,\mathbf{I}$ is minimal. 
We have now shown that for each $n \in N$ there exists a unique partial $n$-solution $\pi_n : E^*_{\le n} \to X$. 

Let $n \le n'$; then the restriction of
$\pi_{n'}$ to $E^*_{\le n}$ is a partial $\pi_n$-solution and so is equal to $\pi_n$.
Define $\pi : E^* \to X$ by  putting $\pi(\varepsilon) = x_0$ and if $n \ne 0$ then letting $\pi(r) = \pi_n(r)$ for all $r \in E^*_n$.

Let $t \in E^*_+$; then there exists a unique $n \in N$ such that $t \in E^*_{s(n)}$ and $t = \triangleleft(e,r)$ with $r \in E^*_n$. Thus 
$\pi(r) = \pi_n(r) = \pi_{s(n)}(r)$ and so
\[\pi(t) = \pi_{s(n)}(t) = \pi_{s(n)}(\triangleleft(e,r)) = f(e,\pi_{s(n)}(r)) = f(e,\pi_n(r))= f(e,\pi(r))\;.\]
Therefore $\pi : E^* \to X$ is such that $\pi(\varepsilon) = x_0$
and $\pi(\triangleleft(e,r)) = f(e,\pi(r))$ for all $(e,r) \in E \times E^*$. 
The uniqueness follows directly from the uniqueness of the partial $n$-solution for each $n \in N$, since the restriction of $\pi$ to $E^*_{\le n}$ is a partial
$n$-solution and is thus equal to $\pi_n$.
\eop

The morphism $\pi : (E^*,\triangleleft,\varepsilon) \to (X,f,x_0)$ is a 'right fold' operation and can be used with the appropriate choices of $(X,f,x_0)$ to obtain several types of mappings defined on lists.

Examples: The list algebra $(X,f,x_0)$ occurring above will be denoted by $\mathbf{L}$.

1.\enskip Here $\,\mathbf{L}  = (N,s',0)$, where $s' : E \times N \to N$ is given by $s'(e,n) = s(n)$. Then
there exists a unique mapping $\ell : E^* \to N$ with $\ell(\varepsilon) = 0$
 and such that $\ell(\triangleleft(e,r)) = s(\ell(r))$ for all $(e,r) \in E \times E^*$. 
 Let $|\cdot| : E^* \to N$ is the mapping with $|r| = n$ for all $r \in E^*_n$ and so $|r|$ 'counts' the number of elements in the list $r$.
Then $|\varnothing| = 0$ and $|\triangleleft(e,r)| = s(|r|)$ for all $(e,r) \in E \times E^*$. 
Thus $|\cdot|$ satisfies the condition which uniquely
determines $\ell$ and therefore $\ell(r) = |r|$ for all $r \in E^*$.

\medskip

2.\enskip Let $F$ be a class, $p : E \to F$ a mapping and let $\,\mathbf{L} = (F^*,\delta,\varepsilon)$, where  the mapping $\delta : E \times F^* \to F^*$ is given 
by $\delta(e,t) = \triangleleft(p(e),t)$. Then there exists a unique mapping $\alpha_p : E^* \to F^*$ with
$\alpha_p(\varepsilon) = \varepsilon$ and such that $\alpha_p(\triangleleft((e,r))) = \triangleleft(p(e),\alpha_p(r)))$ for all $(e,r) \in E\times E^*$. 
The mapping $\alpha_p$ converts each list $r$ with values in $E$ into a list of the same length with values in $F$ by applying the mapping
$p$ to each of the elements of $r$. Thus $\alpha_p([e_0,e_1,\ldots,e_{n-1}]) = [p(e_0),p(e_1),\ldots,p(e_{n-1})]$.

We can give an explicit expression for $\alpha_p$.
Let $\alpha'_p : E^* \to F^*$ be the mapping given by $\alpha'_p(r) = p \circ r$ for all $r \in E^*$. 
Then $\alpha'_p(\varepsilon) = p\circ \varepsilon = \varepsilon$ and if $(e,r) \in E \times E^*$ then
$\alpha'_p(\triangleleft(e,r)) = p \circ \triangleleft(e,r) = \triangleleft(p(e),p \circ r) = \triangleleft(p(e),\alpha'_p(r))$.
Thus $\alpha'_r$ satisfies the condition which uniquely
determines $\alpha_r$ and therefore $\alpha'_r =\alpha_r$ .

\medskip

3.\enskip Let $b : E \to \Bool$ a mapping and let
$q : E \times E^* \to E^*$ be given by $q(e,t) = \triangleleft(e,t)$ if $b(e) = \class{T}$ and $q(e,t) = t$ if $b(e) = \class{F}$. Let $\,\mathbf{L} = (E^*,q,\varepsilon)$.
Then there exists a unique mapping 
$\eta_b : E^* \to E^*$ with
$\eta_b(\varepsilon) = \varepsilon$ and such that $\eta_b(\triangleleft((e,r))) = \triangleleft(e,r)$ if $b(e) = \class{T}$ and $\eta_b((\triangleleft(e,r))) = s$ if $q(e) = \class{F}$.

The mapping $\eta_b$ filters each list $r$ with values in $E$ by removing the elements $e$ in $r$ with $b(e) = \class{F}$.
Thus if $b : \Nat \to \Bool$ with $b(n) = \class{T}$ if and only if $n$ is even then
$\eta_b([0,1,3,6,8,7,6]) = [0,6,8,6]$.

Note that, unlike the first two examples, there does not seem to be an explicit expression for the mapping $\eta_b$.

4.\enskip Let $r \in E^*$ and $\,\mathbf{L} = (E^*,\triangleleft,r)$. Then there is a unique mapping 
$\psi_r : E^* \to E^*$ with
$\psi_r(\varepsilon) = r$ and such that $\psi_r(\triangleleft(e,t)) = \triangleleft(e,\psi_r(t))$ for all $(e,t) \in E \times E^*$.
 
We will see that the mapping $\psi_r$ appends the list $r$ to its argument. Therefore if $r = [b_0,\ldots,b_{n-1}]$ then
$\psi_r([a_0,\ldots,a_{m-1}]) = [a_0,\ldots,b_{m-1},b_0,\ldots,b_{n-1}]$.  In order to see why this is true we give an explicit expression for the mapping $\psi_r$.
First we recall Proposition~\ref{prop_iterators_399}:

For each $n\in N$ put $N_n  = \{ m \in N :n \le m \}$. Then 
$N_n$ is the least $s$-invariant subclass of $N$ 
containing $n$ and $\,\mathbf{N}_n = (N_n,s,n)$ is a Peano iterator.
Moreover, if $\pi_n :\mathbf{N} \to \mathbf{N}_n$ is  the unique isomorphism 
then $L_{\pi_n(m)}$ is the disjoint union of  $ L_n$ and  $\pi_n(L_m)$ for all
$m \in N$.

Let $r \in E^*_m$. Define a mapping $\psi'_r :E^* \to E^*$ as follows:
If $t \in E^*_n$ then $\psi'_r(t)$ is the element of $E^*_{\pi_n(m)}$  given by 
$\psi'_r(t)(k) = t(k)$ if $k \in L_n$ and $\psi'_r(t)(k) = r( k')$ if 
$k = \pi_n(k')$ with $k' \in  L_m$. Thus $\psi'_r(t)$ does append the list $r$ to the list $t$.
\begin{proposition}\label{prop_lists_2}
The mapping  $\psi'_r$ satisfies the condition which uniquely determines
$\psi_r$, i.e., $\psi'_r(\varepsilon) = r$ and  
$\psi'_r(\triangleleft(e,t)) = 
\triangleleft(e,\psi'_r(t))$ 
for all $(e,t) \in E \times E^*$. Therefore $\psi'_r = \psi_r$.
\end{proposition} 

\proof This can be verified with some effort directly from the definition of $\psi'_r$. But it can also be argued as follows: The list $\psi'_r(\triangleleft(e,t)$ 
is obtained by first inserting the element $e$ at the beginning of the list $t$ and then  appending $r$ to the resulting list $\triangleleft(e,t)$, whereas 
the list $\triangleleft(e,\psi'_r(t))$ is obtained by first appending $r$ to the list $t$ and then inserting the element $e$ at the beginning of the resulting  list
$\psi'_r(t)$. The end result in both cases is the same. Therefore $\psi'_r = \psi_r$. \eop

For $r,\,t \in E^*$ we write $t \bowtie r$ instead of $\psi_r(t)$. 
Thus $\psi_r(t) = t \bowtie r$ for all $t, \,r \in E^*$. We consider $\bowtie$ as an infix operation on $E^*$. By the uniqueness of $\psi_r$ it follows that
$\bowtie$ is the unique operation on $E^*$ with   $\varepsilon \bowtie r = r$ for all $r \in E^*$ 
such that $\triangleleft(e,t) \bowtie r = \triangleleft(e, (t \bowtie r))$ for all $(e,t) \in E \times E^*$ and all $r \in E^*$.

Now since  $t \bowtie r$ appends the list $r$ to the list $t$,  it is to be expected that $\bowtie$ is associative, i.e., that
$t\bowtie (u \bowtie v) = (t \bowtie u) \bowtie v$ for all $t,\,u,\,v \in E^*$. This is in fact the case,as we now show.
\begin{proposition}\label{prop_lists_3}
The operation $\bowtie$ is associative and so $(E^*,\bowtie,\varnothing)$ is a monoid.
\end{proposition}

\proof Let 
$E^*_0 = \{ t \in E^* : t\bowtie (u \bowtie v) = (t \bowtie u) \bowtie v\}$ with $u,\,v \in E^*$ considered to be fixed.
In particular $\varepsilon \in E^*_0$, since $\varepsilon \bowtie (u \bowtie v) = (\varepsilon \bowtie u) \bowtie v = u \bowtie v$.
Thus let $t \in E^*_0$ and $e \in E$. Then
\begin{eqnarray*}
\triangleleft(e,t) \bowtie (u \bowtie v) &=& \triangleleft(e,(t \bowtie (u \bowtie v)))\\
 &=& \triangleleft(e,((t \bowtie u) \bowtie v)) = (\triangleleft(e,t) \bowtie u) \bowtie v 
\end{eqnarray*}
where we have twice used the specification which uniquely determines $\bowtie$ and also the assumption that $t \in E^*_0$.
Hence $\triangleleft(e,t) \in E^*_0$. Thus $E^*_0$ is a
$\triangleleft$-invariant subclass of $E^*$ containing $\varepsilon$ and so by Lemma~\ref{lemma_lists_1} $E^*_0 = E^*$. 

Since this holds for all $u,\,v \in E^*$ it follows that $\bowtie$ is associative. It follows that $(E^*,\bowtie,\varnothing)$ is a monoid since also
$\varnothing \bowtie t = t \bowtie \varnothing = t$ for all $t \in E^*$. \eop

\medskip

5.\enskip Here the basic class is not $E^*$ but $E^{**} = (E^*)^*$. The  elements of $E^{**}$ are lists whose elements are lists with values in $E$.
By Theorem~\ref{theorem_lists_1} $(E{^**},\triangleleft^*,\onept)$ is an initial list algebra , where $\triangleleft* : E^* \times E^{^**} \to E^{**}$
appends an element of $E^*$ to each element of $E^{**}$ and $\onept$ is the empty list in $E{^**}$.

Let $\psi : E^* \times E^* \to E^*$ be the mapping with $\psi(r,t) = r \bowtie t$ and put $\,\mathbf{L} = (L^*,\psi,\varepsilon)$.
Then there is a unique mapping $\Psi : E^{**} \to E^*$ with $\Psi(\varepsilon) = \varepsilon$ and such that
$\Psi(\triangleleft^*(r,u)) = r \bowtie \Psi(u)$ for all $(r,u) \in (E^* \times E^{**})$.
The mapping $\Psi$ is a concatenation operator. It takes a list of lists and concatenates them into a single list.

Let $\triangleright_1 : E^*_+ \to E$ and $\triangleright_2 : E^*_+ \to E^*$ denote respectively the first and second components of $\triangleright$.
 Define a mapping $j : E \times E^* \to E^*$ by $j(e,\varepsilon) = \triangleleft(e,\varepsilon)$
and $j(e,r) = \psi_{\triangleright_2(r)}\triangleleft(e,\varepsilon)$ if $r \in E^*_+$, where $\psi$ is as in 4. Then $j(e,r)$
is obtained by moving $e$ to the end of the list and so $j(e,[e_0,\ldots,e_{m-1}] = [e_0,\ldots,e_{m-1},e]$.

6.\enskip Let $\,\mathbf{L} =(E^*,j,\varepsilon)$.  Then there exists a unique mapping $\varrho : E^* \to E^*$ with
$\varrho(\varepsilon) = \varepsilon$ and such that $\varrho(\triangleleft(e,r)) =j(e,\varrho(r))$ for all $(e,r) \in E \times E^*$.
 
The mapping $\varrho$ reverses the order of the elements of a list, so
\[\varrho([e_0,e_1,\ldots,e_{m-1}])= [e_{m-1},\ldots,e_1,e_0]\,\;.\]
To see why this is true,
suppose $\varrho(r) = r'$. Then 

\begin{eqnarray*}
\varrho([e,e_0,e_1,\ldots,e_{m-1}]) &=& \varrho(\triangleleft(e,r))\\
= j(e,\varrho(r)) &=& j(e,r')= j(e,[e_{m-1},\ldots,e_1,e_0]) =[e_{m-1},\ldots,e_1,e_0,e]\;.
\end{eqnarray*}

\medskip

We end the section by giving a characterisation of initial list algebras which corresponds to Theorem~\ref{theorem_iterators_2} for iterators.
We first  state  some simple facts about morphisms which are needed here. These are essentially the same as the corresponding statements for morphism of iterators in
Section~\ref{iterators}.

(1)\enskip
For each list algebra $(X,f,x_0)$ the identity mapping $\id_X$ is a morphism from $(X,f,x_0)$ to 
$(X,f,x_0)$. 

(2)\enskip
If $\pi : (X,f,x_0) \to (Y,g,y_0)$ and $\sigma : (Y,g,y_0) \to (Z,h,z_0)$ are morphisms then $\sigma \circ \pi$ 
is a morphism from $(X,f,x_0)$ to $(Z,h,z_0)$.

If $\pi : (X,f,x_0) \to (Y,g,y_0)$ is a morphism then clearly $\pi \circ \id_X = \pi = \id_Y \circ \pi$, and 
if $\pi,\,\sigma$ and $\tau$ are morphisms for which the compositions are defined then
$(\tau \circ \sigma) \circ \pi = \tau \circ (\sigma \circ \pi)$.

An \definition{isomorphism} is a morphism $\pi : (X,f,x_0) \to (Y,g,y_0)$ for which there exists a morphism 
$\sigma : (Y,g,y_0) \to (X,f,x_0)$ such that $\sigma \circ \pi = \id_X$ and $\pi \circ \sigma = \id_Y$. In this 
case $\sigma$ is uniquely determined by $\pi$.
The morphism $\sigma$ is called the \definition{inverse} of $\pi$.
Moreover,
a morphism $\pi : (X,f,x_0) \to (Y,g,y_0)$ is an isomorphism if and only if the mapping $\pi : X \to Y$ is a 
bijection and the inverse morphism is then the inverse mapping $\pi^{-1} : Y \to X$.

The list algebras $(X,f,x_0)$ and $(Y,g,y_0)$ are said to be \definition{isomorphic} if there exists 
an isomorphism $\pi : (X,f,x_0) \to (Y,g,y_0)$. 

Finally,
if $(X,f,x_0)$ and $(Y,g,y_0)$ are initial list algebras then the unique morphism
$\pi : (X,f,x_0) \to (Y,g,y_0)$ is an isomorphism. In particular, $(X,f,x_0)$ and $(Y,g,y_0)$ are isomorphic.

A list algebra $(X,f,x_0)$ will be called \definition{unambiguous} if the mapping $f_e$ is 
injective for each $e \in E$ and the classes $f_e(X)$, $e \in E$, are disjoint and 
$x_0 \notin \bigcup_{e\in E} f_e(X)$.

By Lemmas~\ref{lemma_lists_1} and \ref{lemma_lists_2} the list algebra $(E^*,\triangleleft,\varepsilon)$ is both minimal and unambiguous
and by Theorem~\ref{theorem_lists_1} $(E^*,\triangleleft,\varepsilon)$ is initial.

\begin{lemma}\label{lemma_lists_3}

Let $(X,f,x_0)$ and $(Y,g,y_0)$ be isomorphic list algebras. Then:

(1) \enskip  $(X,f,x_0)$ is minimal if and only if $(Y,g,y_0)$ is minimal.

(2) \enskip  $(X,f,x_0)$ is unambiguous if and only if $(Y,g,y_0)$ is unambiguous.

\end{lemma}

\proof Let $\pi :(X,f,x_0) \to (Y,g,y_0)$ be an isomorphism.

(1)\enskip  Let $X'$ be an $f$-invariant subclass of $X$ containing $x_0$ Then for all $e \in E$ we have $g_e(\pi(X')) = \pi(f_e((X')) \subset \pi(X'))$
and so $\pi(X')$ is a $g$-invariant subclass of $Y$ containing $y_0$. Thus if $(Y,g,y_0)$ is minimal then $\pi(X') = Y$ which implies that $X' = X$.

This shows that if $(Y,g,y_0)$ is minimal then so is $(X,f,x_0)$. Reversing the roles of $(X,f,x_0)$ and $(Y,g,y_0)$ and using the
isomorphism $\pi^{-1}$ instead of $\pi$ shows that if $(X,f,x_0)$ is minimal then so is $(Y,g,y_0)$.

(2)\enskip Suppose $(Y,g,y_0)$ is unambiguous. Then $f_e$ is injective for all $e \in E$, since $g_e$ is injective and $f_e =   \pi^{-1} \circ g_e \circ\pi$.
Now if $f_e(x) = f_{e'}(x')$ for some $e,\,e' \in E$ and $x,\,x' \in X$
then $g_e(\pi(x)) = \pi(f_e(x)) = \pi(f_{e'}(x')) = g_{e'}(\pi(x'))$. Hence $e = e'$ and $\pi(x) = \pi(x')$, since $(Y,g,y_0)$ is unambiguous and it follows that 
$x = x'$, since $\pi$ is injective. Finally, $\pi(f_e(x)) = g_e(\pi(x)) \ne y_0 = \pi(x_0)$ and so 
$f_e(x) \ne x_0$ for all $x \in X$, $e \in E$.  Thus $(X,f,x_0)$ is unambiguous.
Again reversing the roles of $(X,f,x_0)$ and $(Y,g,y_0)$ and using the
isomorphism $\pi^{-1}$ instead of $\pi$ shows that if $(X,f,x_0)$ is unambiguous then so is $(Y,g,y_0)$.
\eop

\begin{theorem}\label{theorem_lists_2}
A  list algebra is initial if and only if it 
is both minimal and unambiguous.
\end{theorem}

\proof
Let $(X,f,x_0)$ be an initial list algebra. Then by Theorem~\ref{theorem_lists_1} $(E^*,\triangleleft,\varepsilon)$ is initial
and so $(E^*,\triangleleft,\varepsilon)$ and $(X,f,x_0)$ are isomorphic. Thus by Lemmas~\ref{lemma_lists_1}, \ref{lemma_lists_2} and \ref{lemma_lists_3}
$(X,f,x_0)$ is minimal and unambiguous.

The proof of the converse is almost identical to one of the standard proofs of the recursion theorem. Let $(X,f,x_0)$ be an 
unambiguous minimal list algebra and $(Y,g,y_0)$ be any list algebra, and consider the 
list algebra $(X \times Y, f \times_e g, (x_0,y_0))$, where
$f \times_e g : S \times X \times Y \to X \times Y$ is given by
$(f \times_e g)(s,x,y) = (f(s,x),g(s,y))$ for all $e \in E$, $x \in X$, $y \in Y$, and so
$(f \times_e g)_e = f_e \times g_e$ for each $e \in E$. Let $Z$ be the least $(f \times_e g)$-invariant subclass 
of $X \times Y$ containing $(x_0,y_0)$ and let
\[ X_0 = \{\, x \in X : \mbox{there exists exactly one $y \in Y$ such that $(x,y) \in Z$} \,\}\;. \]
It will be shown that $X_0$ is an $f$-invariant subclass of $X$ containing $x_0$, which implies that $X_0 = X$,
since $(X,f,x_0)$ is minimal. We twice need the following fact: If $(x,y) \in Z \setminus \{(x_0,y_0)\}$ then 
there exists $e \in E$ and $(x',y') \in Z$ such that $(f_e(x'),g_e(y')) = (x,y)$. (This follows because
$\{(x_0,y_0)\} \cup \bigcup_{e \in E} (f_e \times g_e)(Z)$ is an $(f \times_e g)$-invariant subclass of 
$X \times Y$ containing $(x_0,y_0)$ and so contains $Z$.)

The element $x_0$ is in $X_0$: Clearly $(x_0,y_0) \in Z$, so suppose also $(x_0,y) \in Z$ for some $y \ne y_0$. 
Then $(x_0,y) \in Z \setminus \{(x_0,y_0)\}$ and hence there exists $(x',y') \in Z$ and $e \in E$ with 
$(f_e(x'),g_e(y')) = (x_0,y)$. In particular $f_e(x') = x_0$, which is not possible, since $(X,f,x_0)$ is 
unambiguous. This shows that $x_0 \in X_0$.

Next let $x \in X_0$ and $e \in E$ and let $y$ be the unique element of $Y$ with $(x,y) \in Z$. Hence 
$(f_e(x),g_e(y)) = (f_e \times g_e)(x,y) \in Z$, since $Z$ is $(f \times_e g)$-invariant. Suppose also  
$(f_e(x),y') \in Z$ for some $y' \in Y$. Then $(f_e(x),y') \in Z \setminus \{(x_0,y_0)\}$, since 
$f_e(x) \ne x_0$, and so $(f_e(x),y') = (f_t(x''),g_t(y''))$ for some $t \in S$ and $(x'',y'') \in Z$. In 
particular $f_t(x'') = f_e(x)$, and this is only possible with $t = s$ and $x'' = x$. since $(X,f,x_0)$ is 
unambiguous. Therefore $y'' = y$, since $x \in X_0$, which implies $y' = g_e(y'') = g_e(y)$. This shows that 
$g_e(y)$ is the unique element $\breve{y} \in Y$ with $(f_e(x),\breve{y}) \in Z$ and in particular that 
$f_e(x) \in X_0$.

We have established that $X_0$ is an $f$-invariant subclass of $X$ containing $x_0$, and so $X_0 = X$. Now define 
a mapping $\pi : X \to Y$ by letting $\pi(x)$ be the unique element of $Y$ such that $(x,\pi(x)) \in Z$ for 
each $x \in X$. Then $\pi(x_0) = y_0$, since $(x_0,y_0) \in Z$ and $\pi(f_e(x)) = g_e(\pi(x))$ for all 
$x \in X$, $e \in E$, since $(f(x),g(y)) \in Z$ whenever $(x,y) \in Z$ and so in particular
$(f_e(x),g_e(\pi(x))) \in Z$ for all $x \in X$, $e \in E$. This gives us a morphism $\pi$ from $(X,f,x_0)$ to 
$(Y,g,y_0)$, and it is easy to see that $f$ being minimal implies that $\pi$ is unique. Therefore the
list algebra $(X,f,x_0)$ is initial.
\eop


\sbox{\ttt}{\textsc{References}}
\thispagestyle{plain}

\bigskip
\bigskip


{\sc Fakult\"at f\"ur Mathematik, Universit\"at Bielefeld}\\
{\sc Postfach 100131, 33501 Bielefeld, Germany}\\
\textit{E-mail address:} \texttt{preston@math.uni-bielefeld.de}\\
\textit{URL:} \texttt{http://www.math.uni-bielefeld.de/\symbol{126}preston}\\


\end{document}